\documentclass[11pt]{amsart} 

\numberwithin{equation}{section} \oddsidemargin=-.0cm
\usepackage{pifont}
\usepackage{empheq} 
\usepackage{bm} 
\usepackage{textcomp}

\usepackage{needspace} 
\usepackage{float}
\usepackage{graphicx}
\usepackage{subfigure}
\usepackage{verbatim}
\usepackage{footmisc}
\usepackage[usenames,dvipsnames]{xcolor}
\usepackage{multirow}
\usepackage{makecell,pict2e}

\usepackage{enumitem} 
\usepackage{framed} 

\usepackage{amsfonts, amssymb, amsmath, amsthm}
\usepackage[colorlinks=true, pdfstartview=FitV, linkcolor=blue,
            citecolor=blue, urlcolor=blue]{hyperref}  

\definecolor{rred}{rgb}{0.7,0,0.1}
\definecolor{greenrb}{rgb}{0.2,0.6,0.2}

\newcommand{\mk}{\color{black}} 
\newcommand{\mkk}{\color{black}}

\newcommand{\sw}{\color{black}}

\newcommand{\att}{\color{black}}
\newcommand{\attn}{\color{black}} 

\newcommand{\hl}{\color{black}} 
\newcommand{\hll}{\color{black}} 

\newcommand{\HL}{\color{black}} 

\newcommand{\HLL}{\color{black}}

\newcommand{\hh}{\color{black}}

\def\bt{\begin{thm}}
\def\et{\end{thm}}
\def\bl{\begin{lem}}
\def\el{\end{lem}}
\def\bd{\begin{defi}}
\def\ed{\end{defi}}
\def\bc{\begin{cor}}
\def\ec{\end{cor}}
\def\bp{\begin{proof}}
\def\ep{\end{proof}}
\def\br{\begin{rem}}
\def\er{\end{rem}}

\def\bprop{\begin{prop}}
\def\eprop{\end{prop}}

\def\y{\mathbf{y}}
\def\z{\mathbf{z}}

\def\s{\mathfrak{s}}
\def\xii{\xi}
\def\c{\mathfrak{c}}
\def\LF{{\rm Lip}(F)}
\def\LFW{{\rm Lip}(F_{\r})}
\def\LH{{\rm Lip}(h_{\lambda})}

\def\Id{{\rm Id}}
\def\Re{{\rm Re \,}}
\def\Im{{\rm Im \,}}

\def\LHW{\mathrm{Lip}(h_{\lambda, \r})}
\def\LHWS{\mathrm{Lip}(h_{\lambda, \r})}

\def\bi{\begin{itemize}}
\def\ei{\end{itemize}}
\def\ben{\begin{enumerate}}
\def\een{\end{enumerate}}

\newcommand{\opnorm}[1]{{\left\vert\kern-0.25ex\left\vert\kern-0.25ex\left\vert #1 
    \right\vert\kern-0.25ex\right\vert\kern-0.25ex\right\vert}}

\def\Forall{\text{ } \forall \:}
\def\d{\mathrm{d}}
\def\Vert{\: \vert \:}
\def\be{\begin{equation}}
\def\ee{\end{equation}}
\def\bes{\begin{equation*}}
\def\ees{\end{equation*}}
\def\bea{\begin{equation} \begin{aligned}}
\def\eea{\end{aligned} \end{equation}}
\def\beas{\begin{equation*} \begin{aligned}}
\def\eeas{\end{aligned} \end{equation*}}

\def\m{q}
\def\p{p}

\def\SS{\S\!}

\def\r{\rho} 

\def\tran{\mathfrak{D}}     
\def\diffusion{\sigma}

\def\M{M}

\newtheorem{thm}{Theorem}[section]
\newtheorem{lem}{Lemma}[section]
\newtheorem{defi}{Definition}[section]

\newtheorem{prop}[thm]{Proposition}

\newtheorem{rem}{Remark}[section]
\newtheorem{cor}{Corollary}[section]

\graphicspath{{./figures/},{../Figures/}}

\title[Stochastic parameterizing manifolds and non-Markovian reduced equations]{On stochastic parameterizing manifolds:\\ Pullback characterization and Non-Markovian reduced equations}

\author[Chekroun]{Micka\"el D. Chekroun}
\address[MC]{Department of Mathematics, University of Hawaii at Manoa, Honolulu, HI 96822, USA, and 
Department of Atmospheric \& Oceanic Sciences, University of California, Los Angeles, CA 90095-1565, USA} 
\email{mdchekroun@math.hawaii.edu}
\email{mchekroun@atmos.ucla.edu}
\author[Liu]{Honghu Liu}
\address[HL]{Department of Atmospheric \& Oceanic Sciences, University of California, Los Angeles, CA 90095-1565, USA}
\email{hliu@atmos.ucla.edu}

\author[Wang]{Shouhong Wang}
\address[SW]{Department of Mathematics, Indiana University, Bloomington, IN 47405, USA}
\email{showang@indiana.edu, \url{http://www.indiana.edu/~fluid}}

\keywords{Random/stochastic invariant manifolds, stochastic partial differential equations, non-Markovian reduced systems, stochastic parameterizing manifolds, pullback limits, random attractors, non-resonance conditions, stochastic inertial manifolds,   stochastic critical manifolds, Taylor approximations of manifolds, memory effects, asymptotic completeness, stochastic Burgers-type equations.}

\subjclass[2010]{34F05, 35B42, 35R60, 37D10, 37L05, 37L10, 37L25, 37L55, 37L65, 60H15}

\date{\today}

\begin{document}

\maketitle

\vspace{-3ex}
\begin{abstract}
Part I of this article is devoted to the leading order approximations of stochastic critical manifolds associated with a broad class of stochastic partial differential equations (SPDEs) which are driven by linear multiplicative white noise.  Stochastic critical manifolds are built naturally as random graphs over a fixed number of critical modes which lose their stability as a control parameter varies. 

Explicit formulas for the leading-order Taylor approximation of such local manifolds about the basic state, are derived.
It is shown that the corresponding approximating manifolds admit  furthermore a pullback characterization, which provides a novel interpretation of such objects in terms of flows. The framework set up in this way allows us, furthermore, to unify the previous approximation approaches from the literature. The existence and attraction properties of one-parameter families of stochastic invariant manifolds are also revisited in this first part. 

In Part II, a point of view more global is adopted. In that respect, a general approach to provide approximate parameterizations of the ``small'' scales by the ``large'' ones, is developed for stochastic 
partial differential equations  driven by linear multiplicative noise. This is accomplished via the concept of {\it parameterizing manifolds} (PMs) that are stochastic manifolds which improve in mean square error  the partial  knowledge of the full SPDE solution $u$ when compared to  the projection of $u$ onto the resolved modes, for a given realization of the noise. 

Backward-forward systems are designed to give access to such PMs in practice.  
The key idea consists of representing  the modes with high wave numbers  (as parameterized by the sought PM) as a pullback limit depending on the time-history of the modes with low wave numbers. 

The resulting manifolds obtained by such a procedure are not subject to a spectral  gap condition such as encountered  in the classical theory.   Instead,   certain PMs  can be determined under weaker {\it non-resonance conditions}: For any given set of resolved modes  for which their self-interactions (through the nonlinear terms) do not vanish when projected against an unresolved mode $e_n$, it is required that  some specific linear combinations of the corresponding  eigenvalues dominate the eigenvalue associated with $e_n$.

Non-Markovian stochastic reduced systems are then derived based on such a PM approach. Such reduced systems take the form of SDEs involving  random coefficients that convey memory effects via the history of the Wiener process, and arise from the nonlinear interactions between the low modes, embedded in the ``noise bath."  These random coefficients follow typically non-Gaussian statistics and exhibit  an exponential decay of correlations whose rate depends explicitly on gaps arising in the non-resonances conditions.

It  is finally  shown on  a stochastic Burgers-type equation, that such PM-based reduced systems can achieve very good performance in reproducing statistical features of the SPDE dynamics projected onto the resolved modes, such as  the autocorrelations and probability functions of the corresponding modes amplitude.  In particular, it is illustrated that the modeling of the large excursions present in the latter can be reproduced with high-accuracy, even when the amount of  the noise is significant. Such success is attributed to the ability of the stochastic PM to capture, for a given realization,  the noise-driven transfer of energy (through the nonlinear terms) to the small scales, as time flows.

\end{abstract}

\tableofcontents

\section{General Introduction} \label{s:intro}

The theory of invariant manifolds for deterministic dynamical systems has been an active research field for a long time, and is now a very well-developed theory; 
see {\it e.g.}, \cite{AW96, AVM96, BJ89, Car81, CFNT89, CH82, CL97, DT96, Far00,FST88, FST89, Hale80, Hale88,  Hen81, HPS77, Kelley67, MR09, MW05, NS11, Pliss64, SY02, Tem97, TW10, TW11,TZ08, Vand89, VI92, VV87}.
Over the past two decades, several important results on random invariant manifolds for stochastically perturbed ordinary as well as partial differential equations (PDEs) {\sw have} been obtained; these results often extend those found in the deterministic setting; see {\it e.g.} ~\cite{Arnold98, AI98, AX95, BF95, Monvel_al98, Boxler91, CCL05, CFD05, CLR01, CG95,  DPD96, DLS03, DLS04, LL10, MS99, MZZ08, NL91, Schmalfuss05}. Even so, the stochastic theory is still much less complete than its deterministic counterpart.  For instance the reduction problem  of a  stochastic partial differential equation (SPDE) to its corresponding stochastic invariant manifolds has been much less studied and only few works in that direction are available \cite{BW10, CLR01, CDZ11, Kan_al12, SDL10, WD07}.

The practical aspects of the reduction problem of a deterministic dynamical system to its corresponding (local) center, center-unstable or unstable manifolds have been well investigated in various finite- and infinite-dimensional settings; see {\it e.g.} \cite{BK98, broer1997algorithms, Car81, DS06, EvP04, Far00, Far01, Hen81, Har08, HI10, johnson1997two,JR05, krauskopf2005survey, Kuznetsov04, MW05, MR09, Pot11, PR06}.  In the stochastic but finite-dimensional context, extensions of center-manifold reduction techniques have been first investigated in  \cite{Boxler89, Boxler91}, and completed by procedures  which consist of deriving  simultaneously both normal forms and center manifold reductions of  stochastic differential equations (SDEs) in \cite{AI98, AX95, NL91,XR96}; see \cite[Sect. ~8.4.5]{Arnold98}.  We mention also \cite{Coullet_al} for prior works on stochastic normal forms.

For stochastic systems of interest for applications, and SPDEs in particular, the  existence problem of stochastic invariant manifolds can be reasonably solved only {\it locally} in practice. This is due to conditions  inherent to the global theory of stochastic invariant manifolds which can be fulfilled only locally via a standard {\it cut-off procedure} of the nonlinear terms  when the latter are not globally Lipschitz; see Section \ref{s:local}.  Such conditions involve typically a gap in the spectrum of the linear part  which has to be large enough {\sw in comparison  with} the variation of the nonlinear terms,  expressed under the form of various {\it spectral gap conditions} in the literature; see Theorems \ref{Lip manfd} and   \ref{attractiveness thm} in Part I.

In situations where the basic state loses its stability as a control parameter $\lambda$ crosses a critical value $\lambda_{\c}$,  such local stochastic invariant manifolds are typically obtained (for $\lambda$ sufficiently close to $\lambda_{\c}$) as random graphs over a neighborhood $\mathcal{V}$ of this basic state, contained in the subspace spanned by the modes which lose their stability, called hereafter the critical modes.  As a consequence, any reduction procedure based on these manifolds can make sense only for $\lambda$ sufficiently close to $\lambda_c$ for which the amplitudes $a_{\lambda}(t)$ of the critical modes remain sufficiently small so that $a_{\lambda}(t) \in \mathcal{V}$.  Such a condition can be satisfied  for all $t$, in the case of deterministic autonomous systems \cite{HI10,MW05}, or even in the case of non-autonomous ones with appropriate time-bounded variations of the vector field  \cite{Aulbach82,DS06,PR10}.   However in a stochastic context, due to large excursions  of the solutions caused by the (white) noise, this condition is expected to be unavoidably violated even when the magnitude of the noise $\sigma$ is small; and is expected to be violated more frequently (as time $t$ flows) when $\sigma$ gets large. We refer to \cite{CM10,Liu_LDP} for results about {\it large deviation principle} for  SPDEs driven by multiplicative noise. 

It is the purpose of this article to propose a new type of stochastic manifolds that are {\it pathwise global} objects\footnote{Defined as graphs of random continuous functions $h(\xi,\omega)$ defined for each realization $\omega$ over the {\it whole} subspace $\mathcal{H}_{\c}$ spanned by the resolved modes, {\it i.e.} for all $\xi\in \mathcal{H}_{\c}$, where  $ \mathcal{H}_{\c}$ is typically the subspace spanned by the first few eigenmodes with low  wave numbers.}  which are able to circumvent this difficulty while allowing  us to derive efficient reduced systems  which are able to achieve good modeling performances of the SPDE dynamics projected onto the resolved (critical) modes, even away from the critical value where the amplitude of these modes gets large.  These manifolds are not invariant in general but provide an approximate parameterization  of the unresolved  variables by the resolved ones, which  improve  in mean square error --- over any (sufficiently large) finite time interval --- the partial  knowledge of the full SPDE solution $u$ when compared with the one obtained from  $u_{\c}:=P_{\c}u $ alone\footnote{Here  $P_{\c}u $ denotes the projection of $u$ onto the resolved modes.}, for a given realization of the noise; see Definition \ref{def:PM}. Such manifolds are naturally named {\it parameterizing manifolds} (PMs) hereafter. A computable criterion coined as the {\it parameterization defect} and inherent to Definition \ref{def:PM}, makes furthermore possible the checking in practice,  that a given stochastic manifold constitutes a stochastic PM; see Section    \ref{ss:Burgers reduction}.

Interestingly, stochastic PMs are not subject to a spectral  gap condition such as encountered  in the classical theory revisited in Part I.   Instead,   certain stochastic PMs  can be determined under weaker {\it non-resonance conditions} in the self-adjoint case:  For any given set of resolved modes  for which their self-interactions (through the nonlinear terms) do not vanish when projected against an unresolved mode $e_n$, it is required that  some specific linear combinations of the corresponding  eigenvalues dominate the eigenvalue associated with $e_n$; see {\it e.g.} ~\eqref{NR} in Section  \ref{ss:PM} and \eqref{NR2} in Section \ref{Sec_Existence_PQ}. No constraints are thus  imposed  on the Lipschitz  constant, explaining why the PM theory  can  overcome the locality issue caused by cut-off arguments.

Stochastic PMs are introduced in Part II of this article where it is shown in particular  that a certain type of  stochastic PMs  coincides actually with the leading order approximation of, for instance,  the local stochastic center or unstable manifolds when restricted to the appropriate neighborhood $\mathcal{V}$ of the basic state; see Section \ref{Sec_app_man+PB}.  The classical theory can thus be reconnected with the important difference that stochastic PMs are global objects and allow thus for large amplitudes.

To further clarify the relationships between stochastic PMs and the classical theory,  Part I is devoted to the derivation of {\sw leading order Taylor approximation} for certain types of local stochastic invariant manifolds
associated with a broad class of SPDEs driven by linear multiplicative white noise. In
particular, approximation formulas --- extending those of \cite[Theorem A.1.1]{MW} --- for local stochastic critical manifolds\footnote{defined as graphs over a neighborhood $\mathcal{V}$  contained in the subspace spanned by the critical modes, for $\lambda$ sufficiently close to $\lambda_{\c}$.} are derived in Section \ref{s:approximation} (see Theorem ~\ref{App CMF} and Corollary ~\ref{Cor:Taylor})
and for local stochastic hyperbolic manifolds, in Section \ref{s:hyperbolic}. Classical theorems about existence and attraction properties of stochastic invariant manifolds are also revisited in Section \ref{s:global} to make the treatment as much  self-contained as possible; see Theorem ~\ref{attractiveness thm} and Corollary ~\ref{attractiveness thm SEE}.   We invite the interested reader to  consult Section ~\ref{Sec_background} to get a more detailed  overview of  the main results 
of Part I.

The notion of (stochastic) PMs such as introduced in this article, raises naturally the question of its relation with the theory of (stochastic) approximate inertial manifolds (AIMs).
The latter has been thoroughly investigated in the deterministic literature, and numerous candidates to  a substitute of inertial manifolds have  been introduced in that purpose; {\mkk see {\it e.g.} \cite{FMT88,FT94,DebTem94,DT95, DT96, NTW01,Tit90}.} In all the cases, the idea was to relax the requirements of the inertial manifold theory so that  the longterm dynamics can be at  least described in some approximate sense, still by some finite-dimensional manifold. Particular efforts have been devoted to developing efficient methods to determine  in practice such manifolds. 
This has led in particular to the so-called nonlinear Galerkin methods; see {\it e.g.} \cite{BJKT90, Chu93,DMT93, FMT88, FST89, JT94, JKT90, JRT01, MT89, NTW01, Tem97, Tit90}.  Approximate inertial manifolds have also been considered in a stochastic context \cite{Chu95}, but contrarily to the deterministic case, only very few algorithms are available to compute stochastic AIMs in practice; see {\it e.g.} \cite{Kan_al12}.

Such a context motivates the important problem examined in Part II concerning the practical  determination of stochastic PMs.  A general approach is introduced in that respect in Section \ref{s:pullback}. 
This approach can be viewed as the cornerstone between the two parts of this article. It consists of obtaining candidates of such PMs as pullback limits of the unresolved variables $\z$
as modeled by auxiliary {\it backward-forward systems}. The key idea consists here of representing  the modes with high wave numbers  as a pullback limit depending on the time-history of the modes with low wave numbers. Such an idea is not new and has been used in the context of 2D-turbulence  \cite{EMS01}, with the essential difference  that   the pullback limits considered in the present article, are associated with backward-forward systems that are {\it partially coupled} in the sense that only the (past values) of the resolved variables $\mathbf{y}$  force   the 
{\sw equations} of the unresolved variables $\mathbf{z}$, without any feedback in the dynamics of the former. 

Due to this partial coupling,  given a realization $\omega$ of the noise, the equations for the resolved variables can be first integrated backward over an interval $[-T,0]$ ($T>0$) from $\xi \in \mathcal{H}_{\c}$, and equations for the unresolved variables are then integrated forward over $[0,T]$. From this operation, the state of the unresolved variables $\z$ at time $s$ ($s>0$)  is thus conditioned on $\xi$ denoted as $\z[\xi]$, and depends  on the state of the resolved variables as well as on the noise at time $s-T$.  The unresolved variables $ \mathbf{z}$ as modeled by such a system is thus a function of   the past history of the resolved variable $\y$, and of the noise.  

The pullback limit of $\z[\xi]$ obtained as $T\rightarrow \infty $, when it exists,   gives generally access to a parameterizing manifold function $h(\xi,\omega)$; but  of course, the design of  appropriate backward-forward systems  is essential for such an operation to be successful.   Different strategies are introduced in Section  \ref{ss:PM} in that respect.  Conditions under which  such systems give access to stochastic PMs  are identified for the stochastic context in Section \ref{Sec_Exist_PMstoch} (Theorems  \ref{thm:hG_PM}  and  \ref{thm:h1_PM}), as well as for the deterministic one in Section \ref{sec_PMdet}. Interestingly,  as for certain AIMs \cite{FMT88},   relations with time analyticity properties  of  the PDE solutions  are involved in order the pullback limit associated with certain backward-forward systems, provides  a PM; see Theorem \ref{thm:h1_PMdet}. The obtention of PMs {\it via} such pullback limits are further analyzed on a stochastic Burgers-type equation in Sections \ref{s:Burgers}  and  \ref{s:Fly}.  

It is worthwhile to mention that such a  pullback characterization of stochastic PMs, constitutes an approach more appealing, from a numerical viewpoint, when compared to methods of approximation of {\it stochastic inertial manifolds} directly rooted in the work of \cite{DPD96} which are also based on backward-forward systems, but this time, {\it fully coupled}; see \cite{Kan_al12}.\footnote{Furthermore, we mention that our particular choice of multiplicative noise  allows us to consider --- {\it via} the cohomology approach (see Section \ref{sec. conjugacy}) ---  transformed versions  ($\omega$ by $\omega$) of our backward-forward systems (such as system \eqref{Eq:bvp-rpde}) so that  we do not have to deal with problems which arise  in trying to solve more general stochastic equations backward in time.}  We mention also that  from a more standard  point of view,  the complementary pullback characterization of  (local) approximating manifolds  presented in Section ~\ref{Sec_app_man+PB} provides a novel interpretation of such objects in terms of  flows. The framework set up in this way allows us, furthermore, to unify the previous approximation approaches from the literature  \cite{BW10, CLR01, CDZ11, WD07}. These features are not limited to the stochastic setting as pointed out in Section \ref{s:pullback},  and are actually dealt with in more details in \cite{CLW13}.

Finally,  two {\it global} stochastic reduction procedures  based on PMs obtained as pullback limits from the auxiliary systems of Section \ref{ss:PM}, are respectively  presented in Section \ref{s:reduction} when analytic expression of such PMs are available, and in Section  \ref{s:Fly} in the general case.  In the case where analytic expressions are not available, a numerical procedure is described in Section  \ref{ss:Burgers-fly} (and in Section \ref{Sec_Existence_PQ}) to determine  ``on the fly''  the  reduced random vector field along a trajectory $\xii(t,\omega)$ generated by the latter as the time is advanced. The cornerstone in this case, is (again) the pullback characterization of  the appropriate PMs, which allows us to update the reduced vector field once $\xii(t,\omega)$  is known at a particular time instance $t$.

In all the cases, a main  feature of the resulting PM-based reduced systems,  comes from the {\it interactions between the stochastic and nonlinear effects} which are shown to contribute to the emergence of {\it memory effects} (see Lemma \ref{lem:Mn} and Section \ref{Sec_Existence_PQ}) conveyed by the stochastic parameterizing manifolds such as built in this article.  These memory effects are shown to play an essential role in the derivation of efficient reduced models able to describe with good accuracy  the main dynamical features of the amplitudes of the resolved modes contained in $\mathcal{H}^{\c}$; see Sections \ref{s:Burgers} and \ref{s:Fly}.  Interestingly, the PM-based stochastic reduction procedure can be seen as an alternative to  the (stochastic) nonlinear Galerkin method where the approximate stochastic  inertial manifolds (AIMs) used therein \cite{Chu95} are substituted here by the parameterizing manifolds introduced above. 

Among the differences with the AIM approach, the PM approach seeks for manifolds which provide modeling error of the evolution of $u_{\c}$ in {\it mean square sense}, over any finite (sufficiently large) time interval; see Proposition  ~\ref{lem:PM error}. This modeling error is controlled by the product of three terms: the energy of the unresolved modes ({\it i.e.} the unknown information), the nonlinear effects (associated with the size of the global random attractor), and the parameterization defect of the stochastic PM employed in the reduction. Interestingly there are cases where parameterization defect can be easily assessed by the theory.  For instance,  when the trivial steady state is unstable, a stochastic inertial manifold (when it exists) is shown to always constitute a stochastic PM; see Theorem \ref{Prop_IMisPM}. The corresponding parameterization defect $Q$ decays  then to zero, 
and  the parameterization of the small scales by the large ones becomes asymptotically  ``exact'' in that case; see again Theorem \ref{Prop_IMisPM}. 
Complementarily, Theorems  \ref{thm:hG_PM}, \ref{thm:h1_PM} and  \ref{thm:h1_PMdet} deal with situations where a stochastic inertial manifold is not known to exist,  and  provide theoretical estimates of  the parameterization defect of various PMs.

The resulting reduced equations are low-dimensional SDEs arising typically with random coefficients which convey  {\it extrinsic memory effects} \cite{HO07, Hai09} expressed in terms of decay of correlations (see Lemma  ~\ref{lem:Mn}),  making the stochastic reduced equations genuinely non-Markovian \cite{Hai09}.  These random coefficients  involve the history of the noise path and  exponentially decaying terms depending in the self-adjoint case on the gap between some linear combinations of the eigenvalues associated with the low modes and the eigenvalues  associated with the high modes. These gaps  correspond exactly to those appearing in the non-resonance condition \eqref{NR} ensuring the pullback limit --- associated with the first backward-forward system introduced in Section  \ref{ss:PM} ---  to exist.  In that case, the memory terms emerge from the nonlinear leading-order interactions between the low modes\footnote{as projected onto the high modes.}, embedded in the ``noise bath.''

Extrinsic memory effects of different type have been encountered in reduction strategies of finite-dimensional SDEs to random center manifolds; see {\it e.g.} ~\cite{Boxler89}. Extrinsic memory effects also arise in procedures  which consist of deriving  simultaneously both normal forms and center manifold reductions of  SDEs; see for instance \cite{AI98, AX95, NL91} and \cite[Sect. ~8.4.5]{Arnold98}. In  such a two-in-one strategy,  anticipating terms may arise  --- as integrals involving the future of the noise path  ---  in  both  the  corresponding  random change of coordinates and the resulting normal form.

In \cite{FS09, LR12, Rob08},  pursuing the works of \cite{AI98,AX95}, reduced stochastic equations involving also extrinsic memory terms have been derived mainly in the context of the stochastic slow manifold; see also \cite{BM13}. By seeking for a random change of variables, which typically involves repeated stochastic convolutions,  reduced equations (different from those  derived in Section \ref{s:reduction}) are obtained to model the dynamics of the slow variables. These reduced equations are also non-Markovian but require a special care in their derivation  to push the anticipative terms (arising in such an approach) to higher order albeit not eliminating them \cite{FS09,Rob08}.

As a comparison, our reduction strategy  is naturally associated with the theory of (stochastic) parameterizing manifolds introduced in this article, and in particular it does not require the existence of  a stochastic slow (or inertial)  manifold. Our approach prevents furthermore the emergence of anticipative terms to any order in the corresponding reduced SDEs.   Memory terms of more elaborated structures than described in Lemma  ~\ref{lem:Mn} (see e.g. \eqref{N_2} or \eqref{N_3} in Section \ref{Sec_Existence_PQ}) can  also arise in  our stochastic reduced equations    built from stochastic PMs  defined as  pullback limits associated with the {\it multilayer backward-forward systems} introduced  in Section  \ref{ss:PM}. As illustrated for the stochastic Burgers-type equation analyzed in Section \ref{s:Fly}, such a multilayer backward-forward system conveys typically  a {\it hierarchy of memory terms} obtained {\it via} repeated compositions of functions involving integrals depending on the past of the noise path driving the SPDE.  Such a hierarchy arises  with a {\it ``matriochka'' of nonlinear self-interactions} between the low modes, as well as with  a sequence of non-resonance conditions,  both of increasing complexity.

 As application, it is shown in Section \ref{ss:num results} that such elaborated memory terms,  along with the corresponding nonlinear cross-interactions brought by the stochastic PM into the reduced systems,  can allow the latter to achieve very good performance in reproducing statistical features of the SPDE dynamics projected onto the resolved modes, such as  the autocorrelations and probability functions of the corresponding modes amplitude.  In particular, it is illustrated that the modeling of the large excursions present in the latter can be reproduced with high-accuracy, and that  such a success can be attributed, for the problem at hand,  to the ability of the corresponding (computed) stochastic PM to capture, for a given realization,  the noise-driven transfer of energy (through the nonlinear term) to the small scales, as time flows.

 \newpage 
\needspace{1\baselineskip}
\vspace{2cm}
\part*{\large \centerline{Part I. Stochastic Invariant Manifolds: Existence, Attraction, and}\\ \centerline{Approximation formulas}\hspace{-2cm}}

\vspace{1cm}

\section{Stochastic Invariant Manifolds: Background and Main Contributions}\label{Sec_background}

The focus of the first part of this article is the derivation of leading order approximations of stochastic invariant manifolds such as stochastic center manifolds by extending, to a stochastic context, the  techniques described in \cite[Chapter 3]{MW05}  and \cite[Appendix A]{MW};  see Section ~\ref{s:approximation}.
New properties of approximating manifolds described in terms of pullback characterization  are {\hh also} identified in Section ~\ref{Sec_app_man+PB}.\footnote{Section ~\ref{Sec_app_man+PB} has been included in part II, but its content concerns the approximating manifolds considered in Part I. This somewhat unconventional presentation has been adopted here in order to articulate, in a unified way,  the {\it pullback characterization} of such manifolds as well as  of  the parameterizing manifolds considered in Part II. }
The framework set up in this way allows us, furthermore, to unify the previous approximation approaches from the literature \cite{BW10, CLR01, CDZ11, Kan_al12}.  These features are not limited to the stochastic setting as pointed out in Section ~\ref{Sec_app_man+PB}.

  {\hh In that respect and} motivated by the study of stochastic bifurcations or more general phase transitions arising in SPDEs\footnote{which is the main purpose of \cite{CGLW}.} \cite{CH93,Mun04}, we first revisit in Sections~\ref{s:global} and \ref{s:local} the  existence and  smoothness properties (Theorems ~\ref{Lip manfd}, \ref{smth manfd}, and ~\ref{thm:local REE})  --- as well as the attraction properties in terms of {\it almost sure asymptotic completeness} (Theorem~\ref{attractiveness thm}) --- of  families of global  stochastic invariant manifolds parameterized by the noise amplitude $\sigma$, and by some control parameter $\lambda$. The latter   is assumed here to vary in some interval $\Lambda$ over which  a uniform decomposition of the spectrum holds; see  \eqref{gap 3} below.   
 The latter condition implies some uniform (partial)-dichotomy estimates that are satisfied by the linearized stochastic flow about the basic state; see ~\eqref{eq:dichotomy}.  

The questions of existence and smoothness are dealt  within a framework rooted in the standard Lyapunov-Perron method \cite{BM61, KB34, Lyapunov47, Perron29}. The techniques follow those used for instance in \cite{CL88, Hen81, Vand89, VI92, VV87}, from which we propose a treatment adapted to the random setting inspired  mainly by the works of \cite{CS01, DLS04}. The related existence and smoothness results are essentially known, but are revisited here in order to set up the precise framework  and to provide the technical tools on which we rely to establish the main results of the first part of this article in Sections ~\ref{s:approximation}, ~\ref{s:hyperbolic}, and ~\ref{Sec_app_man+PB}.

Our treatment of the {\it asymptotic completeness problem}  is inspired by the work  of \cite{CL88} that we adopt in the stochastic framework, and that consists of reformulating this problem as a fixed point problem under constraints. This problem is then recast as an unconstrained fixed point problem associated with a random integral operator that is solved by means of the {\it uniform contraction mapping principle} \cite[Theorems 2.1--2.2]{CH82}.
Various types of attraction properties of random invariant manifolds have been explored in the literature mostly in contexts where the associated SPDEs possess a stable self-adjoint linear part and a bounded and  Lipschitz nonlinearity. For instance, in \cite[Thm. ~3.1]{BF95}, both forward and pullback exponential attractions\footnote{The exponential attraction used therein extends in a random context the classical one encountered in the theory of (deterministic) inertial manifold \cite{FST88}.} of the stochastic inertial manifold are established. Asymptotic completeness in some $n^{\mathrm{th}}$-moment has been established in  \cite[Thm.~2]{CG95} and \cite[Prop.~3.5]{DPD96} (with $n$ being any integer in \cite{CG95} and $n=2$ in \cite{DPD96}) for stochastic inertial manifolds associated with certain types of SPDEs with respectively additive and multiplicative noise.  Almost sure forward asymptotic completeness for deterministic initial data has been established in \cite{CS01} for retarded SPDEs  with additive noise and a stable self-adjoint linear part, adapting also the work of \cite{CL88} to a stochastic context. Almost sure pullback asymptotic completeness of stochastic invariant manifolds has also been investigated in \cite[Thm. ~2.1]{WD07}  for certain type of SPDEs with nonlinearities which do not cause a loss of regularity compared to the ambient space $\mathcal{H}$.\footnote{{\it i.e.}, $F\colon \mathcal{H}\rightarrow \mathcal{H}$ in our notations. Note also that the proof of \cite[Thm. ~2.1]{WD07} provided therein is not complete.}

For SPDEs considered in this article, Theorem ~\ref{attractiveness thm} provides conditions under which the stochastic invariant manifolds ensured by Corollary \ref{Lip manfd SEE} are almost surely forward and pullback asymptotically complete  with  respect to random tempered initial data; see Definition ~\ref{asy. complete}. In particular the existence of  a one-parameter family of global stochastic inertial  manifolds is achieved, and it is shown that the constitutive  manifolds of this family attract exponentially the dynamics at a uniform rate as $\lambda$  varies in $\Lambda$. 
The results obtained in Theorem ~\ref{attractiveness thm} and Corollary ~\ref{attractiveness thm SEE} are not restricted to the case of self-adjoint linear operator and include the cases where unstable modes are present.
 The latter situation is particularly useful to establish in Part II that stochastic inertial manifold always constitute a stochastic PM; see Theorem \ref{Prop_IMisPM} whose proof relies furthermore on some elements contained in the proof of Theorem ~\ref{attractiveness thm}.

 In Section ~\ref{s:local}, we present a local theory of stochastic invariant manifolds associated with the global theory described in Section \ref{s:global}. The ideas 
are standard but the precise framework is detailed here again in view of the main results regarding the approximation formulas of stochastic critical manifolds (Section ~\ref{s:approximation}) and the related pullback characterizations (Section \ref{Sec_app_man+PB}). In particular,  the proof of Theorem ~\ref{thm:local REE} regarding the existence of one-parameter families of local stochastic invariant manifolds,  is provided since some elements  are used in establishing  the main results related to the approximating manifolds considered in Theorem \ref{App CMF}. 

Sections ~\ref{s:approximation}, ~\ref{s:hyperbolic}, and Section ~\ref{Sec_app_man+PB} are  devoted to the main results concerning approximating manifolds of  local stochastic critical manifolds on the one hand (Section \ref{s:approximation}), and local hyperbolic ones, on the other (Section ~\ref{s:hyperbolic}).   They concern the derivation of new approximation formulas of these (local) stochastic  invariant manifolds, and the pullback characterization of the corresponding approximating manifolds such as described in Section ~\ref{Sec_app_man+PB}. More precisely, in  Section~\ref{s:approximation}, we consider the important case for applications where some leading modes lose (once) their stability as $\lambda$ varies in $\Lambda$, which is formulated as the {\it principle of exchange of stabilities (PES)}; see condition \eqref{PES}. It is shown in Lemma \ref{lem:PES} that the  latter implies {\mkk the uniform spectrum decomposition assumed in previous sections}. This allows us in turn to establish in Proposition \ref{prop:cm}, the existence of a family of {\it local stochastic critical manifolds} which are built --- by {\mkk relying on} Section~\ref{s:local} --- as graphs
over some deterministic neighborhood of the origin in the subspace spanned by the critical modes that lose their stability as $\lambda$ varies. By construction, these manifolds {\mkk carry nonlinear dynamical information} associated with the loss of the linear stability of these critical modes; see \cite{CGLW}.

 We then derive  in  Theorem ~\ref{App CMF} and Corollary ~\ref{Cor:Taylor}, {\it explicit random approximation formulas} to the leading order of these local stochastic critical manifolds\footnote{See Definition ~\ref{Def:critical manifold}.} about the origin.  These stochastic critical manifolds are built naturally as graphs over a fixed number of critical modes,  which lose their stability as $\lambda$ varies.  More precisely, the corresponding approximating manifolds are obtained as graphs  --- over some $\lambda$-independent neighborhood $\mathcal{N}$ of zero in the subspace $ \mathcal{H}^{\c}$ spanned by the critical modes --- of the following one-parameter family of random functions:
 \be\label{Eq_app_intro}\tag{AF}
\mathfrak{I}_{\lambda}(\xi, \omega) = \int_{-\infty}^0 e^{ \sigma (k-1)W_t(\omega) \Id} e^{-tL_{\lambda}} P_{\s} F_k(
e^{t L_{\lambda} }\xii) \,\mathrm{d}t, \quad \xi\; \in  \mathcal{N},
\ee
where  $F_k$ denotes the leading-order nonlinear terms of order $k$, $L_{\lambda}$ the corresponding parameterized linear part, $P_{\s}$ the projector upon the non-critical modes, and $W_{t}(\omega)$ the Wiener path associated with the realization $\omega$ of the noise with amplitude $\sigma$.

It is worth mentioning at this stage that the random approximation formulas such as given by \eqref{Eq_app_intro}, contrast with the deterministic ones proposed in ~\cite{BW10} and \cite{CDZ11} for certain types of SPDEs.  In particular, the nonlinearity considered in \cite{BW10} {\mk consists of} a bilinear term, $B(u,u)$, while it consists of power nonlinearity, $u^k$, with $k\ge 2$ in ~\cite{CDZ11}. The error bounds for the approximation of the local random invariant manifold function $h(u,\omega)$ provided in both  \cite{BW10} and \cite{CDZ11} are of the same {\mk order as $\| u\|$ 
and are valid with large probability, and 
for sufficiently small $u$; see \cite[Thm. ~7]{BW10} and \cite[Lemma 4.10]{CDZ11}. 

The class of SPDEs  of type \eqref{SEE} considered below contains the SPDEs {\mk dealt with} in \cite{BW10,CDZ11} as special cases. In contrast with the deterministic approximation formulas obtained in \cite{BW10,CDZ11}, the approximations derived hereafter  are genuine random polynomial functions, which approximate {\mk almost surely} the local {\mk random critical manifolds} and provide (random) Taylor approximations of these manifolds to the leading order; see Corollary ~\ref{Cor:Taylor}.  {\mk More precisely, {\it a priori} error estimates are derived in a general setting} which are of  order {\mk $o(\|u\|_{\alpha}^{k})$} if the {\mk nonlinearity,  $F(u)$,} in {\mk the family of SPDEs}  {\mk is such that $F(u) = O(\|u\|_\alpha^k)$} for some integer $k\ge 2$; see again Theorem ~\ref{App CMF} and {\hll Corollary ~\ref{Cor:Taylor}} for  {\mk precise statements of these} results.

Approximation formulas such as given by  \eqref{Eq_app_intro} are then extended to the case of stochastic hyperbolic manifolds in  Section  ~\ref{s:hyperbolic}, which allows for $\mathcal{H}^{\c}$ to contain a combination of critical  modes, and modes that remain stable as $\lambda$ varies in some interval $\Lambda$.  In that respect, relaxation of the conditions on the spectrum under which the Lyapunov-Perron integral 
$\mathfrak{I}_{\lambda}$ exists, are identified.  In particular, when $L_{\lambda}$ is self-adjoint, it is shown that $\mathfrak{I}_{\lambda}$ exists if a non-resonance condition ~\eqref{NR} is satisfied: For any given set of resolved modes  for which their self-interactions (through the nonlinear term $F_k$) do not vanish when projected against an unresolved mode $e_n$, it is required that  some specific linear combinations of the corresponding  eigenvalues dominate the eigenvalue associated with $e_n$. 

In Section  ~\ref{Sec_app_man+PB}, it is shown that the approximation formulas as given by \eqref{Eq_app_intro} (and derived in Theorem ~\ref{App CMF}, see also {\hll Corollary ~\ref{Cor:Taylor}}) admit an interesting pullback characterization as summarized in Proposition ~\ref{lem:pullback}. More precisely, these formulas can be obtained as the pullback limit of solutions resulting from backward-forward integration of an auxiliary system \eqref{LLL}, which --- unlike the full SPDE \eqref{SEE} --- is only partially coupled in the sense that the critical modes force  the non-critical ones without any feedback into the dynamics of the critical ones. This auxiliary system governs thus the shape of the approximating manifolds proposed in this article.

Remarkably, the approximating manifolds proposed in \cite{BW10, CDZ11} possess also a flow interpretation which allows us  to characterize their shapes as a pullback limit from other auxiliary systems \eqref{eq:BW} and \eqref{Eq:Duan_etal_shape}. In particular it turns out that the approximation formula of \cite{CDZ11} can be interpreted as a special case of ours by setting $\sigma = 0$ in our auxiliary system \eqref{LLL}.

We turn now to the organization of the first part of this article. In Section~\ref{s:preliminary}, we introduce the  class of SPDEs considered throughout this article (including Part II) and describe the main assumptions among which a uniform  decomposition of  the spectrum for the linear constitutes a key ingredient in most of the proofs produced hereafter. We also recall some basic concepts from RDS theory \cite{Arnold98}, and cast such SPDEs into the RDS framework by a classical random change of variables leading to random partial differential equations (RPDEs). The existence, uniqueness, and measurability properties of classical solutions to such RPDEs are recalled in Proposition ~\ref{prop:exist}. For the sake of completeness, the proof and some related results concerning the mild solutions to these RPDEs are presented in Appendix ~\ref{Sect_mild}.

In Section~\ref{s:global}, we revisit the existence and attraction properties of global random/stochastic invariant manifolds within a framework that is suitable for the derivation of certain results regarding the stochastic parameterizing manifolds introduced in Part II; see {\it e.g.} Theorem \ref{Prop_IMisPM}. 
We first derive the existence and smoothness of such manifolds for the transformed RPDEs in Theorems~\ref{Lip manfd}--\ref{smth manfd}. The corresponding results for the original SPDEs are presented in Corollaries ~\ref{Lip manfd SEE}--\ref{smth manfd SEE}. Finally, the almost sure forward-and-pullback asymptotic completeness of these manifolds is examined in Theorem ~\ref{attractiveness thm} and Corollary ~\ref{attractiveness thm SEE}.

In Section~\ref{s:local}, we relax the global Lipschitz condition on the nonlinear term, and derive accordingly the existence of local stochastic invariant manifolds for SPDEs; see Theorem ~\ref{thm:local REE} and Corollary ~\ref{thm:local SEE}.  Section~\ref{s:approximation} is devoted to the main results of the first part of this article regarding the approximation formulas of local stochastic critical manifolds for SPDEs, as summarized in Theorem ~\ref{App CMF} and Corollary ~\ref{Cor:Taylor}. Rigorous error estimates to the leading order are in particular derived.
These results are then extended to the case of (local) stochastic hyperbolic manifolds in  Section  ~\ref{s:hyperbolic}.

\section{Preliminaries} \label{s:preliminary}
In this section, we introduce the functional framework and our  standing hypotheses concerning the abstract stochastic evolution equations of type \eqref{SEE}  below that we will work with. We also recall some basic concepts from the RDS theory \cite{Arnold98,Crauel02}, and introduce a classical random change of variables \cite{DLS04} which will be used to cast  a given parameterized family of SPDEs with abstract formulation as given by Eq. ~\eqref{SEE}  into the RDS framework.

\subsection{Stochastic evolution equations}   \label{ss:SEE}

 We consider the following nonlinear stochastic evolution equation\footnote{Throughout this article, we will often refer to a stochastic evolution equation of type \eqref{SEE} as an SPDE.} driven by linear multiplicative white
noise in the sense of Stratonovich:
\begin{equation} \label{SEE}
\mathrm{d} u = \big( L_\lambda u + F(u) \big) \mathrm{d} t + \diffusion u \circ \mathrm{d}W_t.
\end{equation}
Here, $\{L_\lambda\}$ represents a {\mkk family of linear operators parameterized by a scalar control parameter $\lambda$}, $F(u)$ {\mkk accounts for the} nonlinear {\mkk terms}, $W_t$ is a two-sided one-dimensional Wiener process, and $\diffusion$ is a positive constant which gives a measure of the {\HL ``amplitude"} of the noise.  We make precise below the functional framework that we will adopt {\mkk throughout} this article for the examination {\mkk and approximation} of {\mkk a natural class of} stochastic invariant  manifolds associated with Eq. ~\eqref{SEE}. {\HLL We refer to \cite{CH93,Mun04} for physical contexts {\mkk in which} such equations arise.}

 \vspace{1ex}
{\bf Assumptions  about the operator $L_\lambda$.}  Let $(\mathcal{H}, \|\cdot\|)$ be {\hll an infinite-dimensional} separable {\hll real} Hilbert space. First, let us introduce  a sectorial operator  $A$ on $\mathcal{H}$ \cite[Def. ~1.3.1]{Hen81} with domain 
\bea
\mathcal{H}_1:=D(A)\subset \mathcal{H},
\eea 
and which has compact resolvent. In particular  $\mathcal{H}_1$ is compactly and densely  embedded   in $\mathcal{H}$. We assume  furthermore that $-A$ is stable in the sense that its spectrum satisfies $\Re \sigma(-A) < 0$.

We shall also make use of the fractional powers of $A$ and the associated  interpolated spaces between $\mathcal{H}_1$ and $\mathcal{H}$;  see {\it e.g.} \cite[Sect. ~1.4]{Hen81} and \cite[Sect. ~3.7]{SY02}. 
Let $\mathcal{H}_\gamma:=D(A^{\gamma})$ be such an interpolated space for some $\gamma \in [0,1],$ endowed with the norm $\|\cdot\|_{\gamma}$ induced by the inner product $\langle u, v \rangle_\gamma := \langle A^\gamma u, A^\gamma v \rangle_\mathcal{H}$; in particular $\mathcal{H}_0=\mathcal{H}$, $\mathcal{H}_1$ corresponds to $\gamma=1$, and $\mathcal{H}_1 \subset \mathcal{H}_{\gamma} \subset \mathcal{H}_0$ for $\gamma \in (0,1)$.  Note that in the sequel, $\langle \cdot, \cdot \rangle$, will  be used to denote  the inner-product in the ambient Hilbert space $\mathcal{H}$.

Let us introduce now 
\bea  \label{eq:B-1}
B_\lambda \colon \mathcal{H}_\gamma \rightarrow \mathcal{H}
\eea
a parameterized family of bounded linear operators depending continuously on $\lambda$, with here $\gamma \in [0,1).$
In particular, $-B_\lambda A^{-\gamma}$ is bounded on $\mathcal{H}$ and according to \cite[Cor. ~1.4.5]{Hen81}  the operator $-L_{\lambda}$ is sectorial on   $\mathcal{H}$ with domain $ \mathcal{H}_1 $  where 
\begin{equation} \label{L}
\begin{aligned}
& L_\lambda: = -A + B_\lambda.
\end{aligned}
\end{equation}

Note that $L_\lambda$ has compact resolvent by recalling that $\mathcal{H}_1$ is compactly embedded in $\mathcal{H}$ \cite[Prop. ~II.4.25]{EN00}. 
As a consequence, since  $L_{\lambda} \colon \mathcal{H}_1\rightarrow \mathcal{H}$ is a closed operator\footnote{as a consequence of the sectorial property of $-L_{\lambda}$.}, we have that for each $\lambda$,
the  spectrum of $L_{\lambda}$, $\sigma(L_\lambda)$, consists only of isolated eigenvalues with finite algebraic multiplicities; see \cite[Thm. ~III-6.29]{Kato95} (see also \cite[Corollary IV.1.19]{EN00}).

\vspace{1ex}
{\bf Assumptions about the nonlinearity $F$.} For the nonlinearity,
we assume  that $F \colon \mathcal{H}_\alpha 
\rightarrow \mathcal{H}$ is continuous for some $\alpha\in [0,1)$ which will be fixed throughout this article. 

We assume furthermore that
\be\label{F} 
F(0)=0,\footnotemark
\ee
\footnotetext{Near a nontrivial  steady state $\overline{u}$ of some deterministic  system, one can think $u$ as some deviation from this steady state  (subject to noise fluctuations), and  Eq. \eqref{F} is then satisfied in such situations.}
and in the case where $F$ is at least $C^1$ smooth, the tangent map of $F$ at $0$ is assumed to be the null map, {\it i.e.},
\be \label{DF}
  DF(0) = 0.
\ee
Note that in particular, according to \eqref{F} the noise term in \eqref{SEE} is multiplicative with respect to the basic state; see \cite[p ~473]{Arnold98} for this terminology.

Moreover, for the results on global random invariant manifolds proved in Section ~\ref{s:global}, $F$ will be assumed to be furthermore globally Lipschitz:
\begin{equation} \label{Lip F}
\|F(u) - F(v)\| \le \LF \|u-v\|_\alpha, \quad \forall\, u, v \in
\mathcal{H}_{\alpha},
\end{equation}
where $ \LF$ denotes the smallest positive constant such that \eqref{Lip F} is true. Other assumptions on $F$ will be specified when needed; see {\it e.g.} Theorem ~\ref{smth manfd}, and Sections ~\ref{s:approximation} and \ref{s:reduction}.

\vspace{1ex}
{\bf The spectrum of $L_{\lambda}$ and  the uniform spectrum decomposition.} 
Recall that the spectrum $\sigma(L_\lambda)$ consists only of isolated eigenvalues with finite multiplicities. This combined with the sectorial property of $-L_{\lambda}$ implies that  there are at most finitely many eigenvalues with a given real part. The sectorial property of $-L_{\lambda}$ also implies that $\Re \sigma(L_\lambda)$ is bounded above (see also \cite[Thm. ~II.4.18]{EN00}).  These two properties of $\Re \sigma(L_\lambda)$  allow us in turn to label elements in $\sigma(L_\lambda)$ according to the lexicographical order which we will adopt throughout this article:
\bea  \label{eq:ordering-1}
\sigma(L_\lambda) = \{\beta_n(\lambda) \:|\:  n \in \mathbb{N}^\ast\}, 
\eea
such that for any $1\le n < n'$ we have either 
\bea
\Re \beta_{n}(\lambda) > \Re \beta_{n'}(\lambda), 
\eea
or
\bea  \label{eq:ordering-3}
\Re \beta_{n}(\lambda) = \Re \beta_{n'}(\lambda), \quad \text{ and } \quad \Im \beta_{n}(\lambda) \geq \Im \beta_{n'}(\lambda). 
\eea
Note that we will adopt in this article the convention that each eigenvalue, $\beta_n(\lambda)$, is repeated according to its algebraic multiplicity.

 We assume throughout this article that an open interval $\Lambda$ can be chosen such that 
the following {\it uniform spectrum decomposition of $\sigma(L_\lambda)$ holds over $\Lambda$}:
\begin{equation} \label{gap 3}
\begin{aligned}
\sigma(L_\lambda) = \sigma_{\c}(L_\lambda) \cup \sigma_{\s}(L_\lambda), \; \lambda \in \Lambda, \qquad   \text{ with } \qquad   \eta_{\c} > \eta_{\s},
\end{aligned}
\end{equation}
where 
\bea  \label{eta_cs}
& \eta_{\c} := \inf_{\lambda \in \Lambda} \inf \{ \Re \beta(\lambda) \; \vert \; \beta(\lambda) \in \sigma_{\c}(L_\lambda) \},  \\
& \eta_{\s} := \sup_{\lambda \in \Lambda} \sup\{ \Re \beta(\lambda) \; \vert \; \beta(\lambda) \in \sigma_{\s}(L_\lambda) \},
\eea
and $\sigma_{\c}(L_\lambda)$ consists of the first $m$ eigenvalues (counting multiplicities) in $\sigma(L_\lambda)$: 
\be  \label{m}
\mathrm{card}(\sigma_{\c}(L_\lambda)) = m.
\ee

It is interesting to note that the {\hh uniform spectrum decomposition \eqref{gap 3}} prevents eigenvalues in $\sigma_{\s}(L_\lambda)$ from merging with eigenvalues in $\sigma_{\c}(L_\lambda)$ as $\lambda$ varies in $\Lambda$, while the cardinality of $\sigma_{\c}(L_\lambda)$ remains fixed to be $m$ over $\Lambda$. 
As a consequence, the spaces $\mathcal{H}_\alpha$ and $\mathcal{H}$ can be decomposed into $L_\lambda$-invariant subspaces in such a way that the subspace associated with $\sigma_{\c}(L_\lambda)$,  $\mathcal{H}^{\c}(\lambda)$ defined below in \eqref{Hc}, has fixed dimension $m$ for each $\lambda \in \Lambda$. These subspaces will  be at the basis of the construction of stochastic invariant manifolds considered in later sections. Furthermore, the uniform spectrum decomposition \eqref{gap 3}
will be essential to guarantee the existence of a  neighborhood of the origin  (in $\mathcal{H}^{\c}(\lambda)$) whose diameter is independent of $\lambda \in \Lambda$ over which a family of stochastic local invariant manifolds --- each of dimension $m$ ---  is defined; see Theorem ~\ref{thm:local REE} and Corollary ~\ref{thm:local SEE}. The existence of such a neighborhood will turn out to be particularly useful in the  examination of stochastic bifurcations or more general phase transitions associated with  SPDEs of type \eqref{SEE}; see Remark ~\ref{rmk:critical} and \cite{CGLW}.

\br  \label{rmk:sign sigma}

Note that $\eta_{\c}$ and  $\eta_{\s}$ are allowed to share the same sign in the derivation of all the results of Sections ~\ref{s:global} and ~\ref{s:local}. The results of Section~\ref{s:approximation}  are presented only in the case where $\eta_{\s}<\eta_{\c}<0$, which will be sufficient for the applications dealt with in \cite{CGLW}.

\er

{\bf Related $L_{\lambda}$-invariant subspaces.} We present now decompositions of the spaces $\mathcal{H}$ and $\mathcal{H}_\alpha$  into $L_{\lambda}$-invariant subspaces {\mkk naturally} associated with  the splitting of the spectrum $\sigma(L_\lambda)$ given in \eqref{gap 3}. These decompositions lead naturally to a partial dichotomy of the deterministic {\mkk linear semigroup}\footnote{{\it i.e.} {\mkk the semigroup  associated with}  $\d v =L_{\lambda} v \d t$.} associated with \eqref{SEE}; see \eqref{Proj bounds}.
This partial dichotomy   will turn out to be sufficient in the construction of stochastic invariant manifolds for the full SPDE, as we will see in the forthcoming sections. The reason is that only stochastic invariant manifolds associated with the trivial steady state are considered in this article. For manifolds associated with more general (random) stationary solutions, other types of dichotomy estimates which involve the Lyapunov spectrum are typically required;  see Remark \ref{Rmk_Lyapunov spec} below. To simplify the presentation of our main results, such manifolds will not be considered in this article.

The  decompositions of  $\mathcal{H}$ and $\mathcal{H}_\alpha$ are organized as follows. Let us introduce the complexifications, $\widetilde{\mathcal{H}}_1$ and $\widetilde{\mathcal{H}}$, of the spaces $\mathcal{H}_1$ and $\mathcal{H}$:
\bea  \label{eq: complex H}
\widetilde{\mathcal{H}}_1:= \{u + i v \; \vert \; u, v \in \mathcal{H}_1\}, \quad \widetilde{\mathcal{H}}:= \{u + i v \; \vert \; u, v \in \mathcal{H}\},
\eea 
where $i$ denotes here the imaginary unit. Denote also the complexification of $L_\lambda$ by $\widetilde{L}_\lambda$, {\it i.e.} \bea \label{eq:complex L}
\widetilde{L}_\lambda(u + i v):= L_\lambda u + i L_\lambda v,  \quad \Forall \, u, v \in \mathcal{H}_1.
\eea 
Note that for each $\lambda \in \Lambda$, the set $\sigma_{\c}(L_\lambda)$  in the uniform spectrum decomposition \eqref{gap 3} is bounded. Hence, according to \cite[Thm. ~III-6.17]{Kato95} (see also \cite[Prop. ~IV.1.16]{EN00}),  
there exists a projector $\widetilde{P}_{\c}$ associated with $\sigma_{\c}(L_\lambda)$, such that the space $\widetilde{\mathcal{H}}$ can be decomposed as follows:
\begin{equation} \label{L decomposition cmplx}
\begin{aligned}
& \widetilde{\mathcal{H}} = \widetilde{\mathcal{H}}^{\c}(\lambda) \oplus \widetilde{\mathcal{H}}^{\s}(\lambda),
\end{aligned}
\end{equation} 
where
\bea \label{Def:space cmplx}
 \widetilde{\mathcal{H}}^{\c}(\lambda) &:=  \widetilde{P}_{\c}  \widetilde{\mathcal{H}}, \quad \widetilde{\mathcal{H}}^{\s}(\lambda) := \big(\Id_{\widetilde{\mathcal{H}}} - \widetilde{P}_{\c}\big) \widetilde{\mathcal{H}}. 
\eea
Note that $\widetilde{P}_{\c}$ is given as the Riesz projector defined by
\beas
\widetilde{P}_{\c} := -\frac{1}{2\pi i}\int_{\Gamma_{\c}} (\widetilde{L}_\lambda - \beta\Id)^{-1} \d \beta, 
\eeas
where $\Gamma_{\c}$ is a rectifiable closed curve surrounding the eigenvalues of $\sigma_{\c}(L_\lambda)$ which does not include any elements of $\sigma_{\s}(L_\lambda)$.

Moreover, $\widetilde{\mathcal{H}}^{\c}(\lambda)$ and $\widetilde{\mathcal{H}}^{\s}(\lambda)$ are invariant under $\widetilde{L}_\lambda$ in the following sense \cite[Thm.~5.7-A]{Taylor58}:
\beas
\widetilde{L}_\lambda (\widetilde{\mathcal{H}}^{\c}(\lambda) \cap \widetilde{\mathcal{H}}_1) \subset \widetilde{\mathcal{H}}^{\c}(\lambda), \quad \widetilde{L}_\lambda (\widetilde{\mathcal{H}}^{\s}(\lambda) \cap \widetilde{\mathcal{H}}_1) \subset \widetilde{\mathcal{H}}^{\s}(\lambda);
\eeas
and the restriction of $\widetilde{L}_\lambda$ on $\widetilde{\mathcal{H}}^{\c}(\lambda)$, denoted by $\widetilde{L}_\lambda^{\c}$, is a bounded linear operator on $\widetilde{\mathcal{H}}^{\c}(\lambda)$ according to \cite[Thm.~5.8-A]{Taylor58} or \cite[Thm. ~III-6.17]{Kato95}.

 Now, for $\lambda \in \Lambda$, let us define
\bea  \label{Hc}
\mathcal{H}^{\c}(\lambda) := \{u, v \; \vert \; u, v \in \mathcal{H}, \text{ and } u+ i v \in \widetilde{\mathcal{H}}^{\c}(\lambda)\}.
\eea
Note that $\mathcal{H}^{\c}(\lambda)$ thus defined forms naturally a subspace of $\mathcal{H}$; and it can be checked that\footnote{Eq. ~\eqref{Eq:dimension} can be justified as follows. Since $L_\lambda$ has real coefficients, for any complex eigenvalue $\beta_j(\lambda) \in \sigma_{\c}(L_\lambda)$, its conjugate $\overline{\beta}_j(\lambda)$ is also an eigenvalue which belongs to $\sigma_{\c}(L_\lambda)$. 

After possibly reordering the eigenvalues, the space $\widetilde{\mathcal{H}}^{\c}(\lambda)$ can be further decomposed into $\bigoplus_{j=1}^{l} \big(\widetilde{\mathcal{H}}_{\beta_j(\lambda)} \oplus \widetilde{\mathcal{H}}_{\overline{\beta}_j(\lambda)} \big)  \bigoplus_{j=2l+1}^m \widetilde{\mathcal{H}}_{\beta_j(\lambda)}$, where $\widetilde{\mathcal{H}}_{\beta_j(\lambda)}$ is the eigenspace associated with $\beta_{j}(\lambda)$, and the first $2l$ eigenvalues are (genuinely) complex and the remaining are real. 

For each $\widetilde{\mathcal{H}}_{\beta_j(\lambda)}$, a real vector space $\mathcal{H}_{\beta_j(\lambda)}$ can be defined in the same way as in \eqref{Hc}. 
Note that $\mathcal{H}_{\beta_j(\lambda)} = \mathcal{H}_{\overline{\beta}_j(\lambda)}$, for all $j \in \{1, \cdots, l\}$. 

Then, $\mathcal{H}^{\c}(\lambda)$ admits the following decomposition: $\bigoplus_{j=1}^{l} \mathcal{H}_{\beta_j(\lambda)} \bigoplus_{j=2l+1}^m \mathcal{H}_{\beta_j(\lambda)}$. 
Since $\dim \mathcal{H}_{\beta_j(\lambda)} = \dim_{\mathbb{C}}\widetilde{\mathcal{H}}_{\beta_j(\lambda)} + \dim_{\mathbb{C}}\widetilde{\mathcal{H}}_{\overline{\beta}_j(\lambda)}$ for all $j \in \{1, \cdots, l\}$, and $\dim \mathcal{H}_{\beta_j(\lambda)} = \dim_{\mathbb{C}}\widetilde{\mathcal{H}}_{\beta_j(\lambda)}$ for all $j \in \{2l+1, \cdots, m\}$, the result follows. }
\bea \label{Eq:dimension}
\dim \mathcal{H}^{\c}(\lambda) = \dim_{\mathbb{C}} \widetilde{\mathcal{H}}^{\c}(\lambda) = m, \quad \Forall \lambda \in \Lambda,
\eea 
where $m$ is the cardinality of $\sigma_c(L_\lambda)$ as given in \eqref{m}. 

The space $\mathcal{H}^{\c}(\lambda)$  uniquely determines the closed subspaces $\mathcal{H}^{\s}(\lambda)$ and $\mathcal{H}_\alpha^{\s}(\lambda)$ as the  topological complements in $\mathcal{H}$ and $\mathcal{H}_\alpha$ respectively, {\it i.e.},
\begin{equation} \label{L decomposition}
\begin{aligned}
& \mathcal{H} = \mathcal{H}^{\c}(\lambda) \oplus \mathcal{H}^{\s}(\lambda),  && \mathcal{H}_\alpha = \mathcal{H}^{\c}(\lambda) \oplus \mathcal{H}_\alpha^{\s}(\lambda), \quad \Forall \lambda \in \Lambda.
\end{aligned}
\end{equation}
 Let 
\bea \label{Pc Ps}
P_{\c}(\lambda): \mathcal{H} \rightarrow \mathcal{H}^{\c}(\lambda), \quad P_{\s}(\lambda): \mathcal{H} \rightarrow \mathcal{H}^{\s}(\lambda)
\eea 
be the associated canonical projectors, and we denote
\bea \label{Lc Ls}
L_\lambda^{\c} := L_\lambda P_{\c}(\lambda), \quad L_\lambda^{\s} := L_\lambda P_{\s}(\lambda).
\eea
Note that $L_\lambda$ commutes with $P_{\c}(\lambda)$ and $P_{\s}(\lambda)$,  which follows from the fact that $\widetilde{\mathcal{H}}^{\c}(\lambda)$ and $\widetilde{\mathcal{H}}^{\s}(\lambda)$ are invariant under $\widetilde{L}_\lambda$. {\mkk As a consequence, the subspaces $\mathcal{H}^{\c}(\lambda)$ and $\mathcal{H}^{\s}(\lambda)$ are left invariant by the semigroup $e^{tL_\lambda}$.} {\HL Note also that} similar to $\widetilde{L}_\lambda^{\c}$, the operator $L_\lambda^{\c}$ is a bounded linear operator on $\mathcal{H}^{\c}(\lambda)$, {\HL so that} $e^{tL_\lambda}P_{\c}$ can be extended to $t< 0$, {\HL namely} $e^{tL_\lambda}P_{\c}$ defines a flow on $\mathcal{H}^{\c}(\lambda)$. This fact is used in the partial dichotomy estimate \eqref{Proj-III} below.

As mentioned above, note that by \eqref{Eq:dimension}, the dimension of $\mathcal{H}^{\c}(\lambda)$ is independent  of $\lambda$ as it varies in $\Lambda$, so that  $\mathcal{H}^{\c}(\lambda)$ is unique up to orthogonal transformations. For the sake of concision, this property has led us to suppress the $\lambda$-dependence of the subspaces given in \eqref{L decomposition}, and of the projectors $P_{\c}(\lambda)$ and $P_{\s}(\lambda)$ defined in \eqref{Pc Ps}. The results are derived and presented hereafter according to this convention.

\vspace{1ex}
{\bf Partial-dichotomy estimates.} Thanks to the uniform spectrum decomposition \eqref{gap 3}, for any given numbers $\eta_1$ and $\eta_2$ satisfying
\bea  \label{eta1-2}
\eta_{\c} > \eta_1 > \eta_2 > \eta_{\s},
\eea 
there exists  a constant, $K\ge 1$, such that for all $\lambda \in
\Lambda$ the following  {\it partial-dichotomy}\footnote{The partial aspect of the dichotomy is explained when $\eta_1$ and $\eta_2$ share the same sign which is allowed by \eqref{gap 3}; see also Remark ~\ref{rmk:sign sigma}. In that case, the distinction is made on the magnitude of the rate of contraction (or expansion) {\mkk associated with $\d v =L_{\lambda} v \d t$}; otherwise the concept matches the classical one of exponential dichotomy found in the literature; see {\it e.g.} ~\cite{Hen81, SY02}.} estimates hold for the semigroup generated by $L_{\lambda}$ (see {\it e.g.} ~\cite[Thms.~1.5.3-1.5.4]{Hen81}\footnote{The first two inequalities in \eqref{Proj bounds} follow from \cite[Thm. ~1.5.4]{Hen81} by choosing the operator $A$ thereof to be $-L_\lambda$ and the operator $B$ to be our operator $A$; the third inequality follows from \cite[Thm. ~1.5.3]{Hen81}, the $L_\lambda$-invariance of $\mathcal{H}^{\c}$, and the fact that any norms on $\mathcal{H}^{\c}$ are equivalent since it is finite dimensional.}):
\begin{subequations} \label{Proj bounds}
\begin{align}
& \|e^{t L_\lambda}P_{\s}\|_{L(\mathcal{H}_\alpha, \mathcal{H}_\alpha)} \le K e^{{\eta_2} t}, && t\ge 0, \label{Proj-I} \\
& \|e^{t L_\lambda} P_{\s}\|_{L(\mathcal{H}, \mathcal{H}_\alpha)}  \le K t^{-\alpha}
e^{{\eta_2} t}, &&  t> 0,  \label{Proj-II} \\
& \|e^{t L_\lambda } P_{\c}\|_{L(\mathcal{H}, \mathcal{H}_\alpha)}  \le K
e^{{\eta_1} t}, && t\le 0, \label{Proj-III}
\end{align}
\end{subequations}
where $L(X, Y)$ denotes the space of bounded linear operators from
the Banach space $X$ to the Banach space $Y$. Note that the estimate given in \eqref{Proj-II} accounts for the instantaneous smoothing effects of the semigroup $e^{tL_\lambda}$ for $t>0$ from $\mathcal{H}$ to $\mathcal{H}_\alpha$ where we recall that $\mathcal{H}_\alpha$ has been imposed by the choice of the nonlinearity.

The conditions $\eta_{\c} > \eta_1$ and $\eta_2 > \eta_{\s}$ allow us to absorb the polynomial growth terms in the estimates  \eqref{Proj bounds} that --- because of our assumptions ($L_\lambda$ being not necessarily self-adjoint) --- could be present  in front of the exponential terms with $\eta_{\c}$ (resp. ~$\eta_{\s}$) in place of $\eta_1$ (resp. ~$\eta_2$). As a consequence, $K$ in \eqref{Proj bounds} depends on $\eta_* := \min\{ \eta_{\c} - \eta_1, \eta_2 - \eta_{\s}\}$, and may get larger as $\eta_*$ gets closer to zero in the non-self-adjoint case. Note that however, $K$ is independent of  $\lambda \in \Lambda$ in all the cases.

\br\label{Rmk_Lyapunov spec}
 It is important to note that in {\mkk the first part of this article}, for the sake of simplicity, the stochastic critical manifold theory along with the corresponding approximation formulas are presented {\mkk in Section ~\ref{s:approximation} in the vicinity} of the trivial steady state, as $\lambda$ varies in a neighborhood of the critical parameter $\lambda_c$. {\mkk In Part II, we will consider manifolds which are not subject to such a locality.} 
 
{\mkk In all the cases,  we will} avoid the use of {\mkk a random frame which moves with the cocycle, as determined by} the Lyapunov spectrum and the multiplicative ergodic theory {\mkk in Hilbert or Banach spaces} \cite{Ruelle82, LL10}.  {\mkk The latter are typically} employed when  {\mkk the dynamics is analyzed} in the vicinity of a nontrivial random stationary solution; see {\it e.g.} \cite{Caraballo_al09, LL10, MZZ08}.  {\mkk Certain results of this article could be extended to this more general setting,  but these extensions are not straightforward}. We leave such an enterprise to the interested reader.

\er

\subsection{Random dynamical systems} \label{subsect:RDS}

In this subsection, we recall the definitions of metric dynamical systems (MDSs) and random dynamical systems (RDSs), and specify --- in a measure-theoretic sense --- the canonical MDS associated with the Wiener process  in Eq.~\eqref{SEE} which will be used throughout this article. The interested readers are referred to \cite{Arnold98, Chueshov02, Crauel02} for more details, and to \cite{CSG11} for an intuitive and ``physically-oriented'' {\mkk presentation} of these concepts.

\medskip

{\bf Metric dynamical system}. A family of mappings $\{\theta_t\}_{t\in \mathbb{R}}$ on a probability space $(\Omega, \mathcal{F}, \mathbb{P})$ is called a {\it metric dynamical system} if the following conditions are satisfied:
\begin{itemize}
\item[(i)] $(t,\omega) \mapsto \theta_t \omega$ is $(\mathcal{B}(\mathbb{R})\otimes \mathcal{F}; \mathcal{F})$-measurable, where $\mathcal{B}(\mathbb{R})$ denotes the Borel $\sigma$-algebra on $\mathbb{R}$, and $\mathcal{B}(\mathbb{R})\otimes \mathcal{F}$ denotes the $\sigma$-algebra generated by the direct product of elements of $\mathcal{B}(\mathbb{R})$ and $\mathcal{F}$;

\item[(ii)] $\{\theta_t\}$ satisfies the one-parameter group property, {\it i.e.}, $\theta_0 = \Id_{\Omega}$, and $\theta_{t+s} = \theta_t\circ\theta_s$ for all $t, s \in \mathbb{R}$;

\item[(iii)] $\mathbb{P}$ is invariant with respect to $\theta_t$ for all $t\in \mathbb{R}$,  {\it i.e.}, $(\theta_t)_{\ast} \mathbb{P} = \mathbb{P}$ for all $t\in \mathbb{R}$, where $(\theta_t)_{\ast} \mathbb{P}$ is the {\it push-forward measure} of $ \mathbb{P}$ by $\theta_t$, defined by $(\theta_t)_{\ast} (F) := \mathbb{P}(\theta_{-t}(F))$, for all $F \in \mathcal{F}$.

\end{itemize}

{\bf Continuous random dynamical system}. Given a separable Hilbert space $(H, |\cdot|_H)$ with the associated Borel $\sigma$-algebra denoted by $\mathcal{B}(H)$,  a continuous {\it random dynamical system} acting on $H$ over an MDS, $(\Omega, \mathcal{F}, \mathbb{P}, \{\theta_t\}_{t\in \mathbb{R}})$, is a $(\mathcal{B}(\mathbb{R}^+)\otimes  \mathcal{F} \otimes \mathcal{B}(H); \mathcal{B}(H))$-measurable mapping 
\beas
S\colon \mathbb{R^+}\times \Omega \times H \rightarrow H, \quad (t, \omega, \xi) \mapsto S(t, \omega)\xi,
\eeas 
which satisfies the following properties:
\begin{itemize}

\item[(i)$'$]  $S(0, \omega) = \Id_H$, for all $\omega \in \Omega$,

\item[(ii)$'$] $S$ satisfies the {\it perfect cocycle property}, {\it i.e.}, 
\beas
\hspace*{1.4em}  S(t+s, \omega) = S(t, \theta_s\omega) \circ S(s, \omega), \quad \Forall \, t, s \in \mathbb{R}^+, \text{ and } \omega \in \Omega,
\eeas

\item[(iii)$'$] $S(t, \omega)\colon H \rightarrow H$ is continuous for all $t \in \mathbb{R}^{+}$ and $\omega \in \Omega$.

\end{itemize}

\medskip

We are in position to introduce the MDS that we will work with 
throughout this article. Let us first recall the canonical MDS, $(\Omega, \mathcal{F}, \mathbb{P}, \{\theta_t\}_{t\in \mathbb{R}})$, associated with the Wiener process; see {\it e.g.} ~\cite[Appendices A.2-A.3]{Arnold98} and \cite[Chapter 1]{Chueshov02}. Here the sample space $\Omega$ consists of the sample paths of a two-sided one-dimensional
Wiener process $W_t$ taking zero value at $t=0$, that is,
\beas
\Omega = \{ \omega \in C(\mathbb{R},\, \mathbb{R}) \; \vert \; \omega(0) = 0\}; 
\eeas 
$\mathcal{F}$ is the Borel $\sigma$-algebra associated with the Wiener process; $\mathbb{P}$ is the classical Wiener measure on $\Omega$; and for each $t \in \mathbb{R}$, the map $\theta_t\colon (\Omega, \mathcal{F}, \mathbb{P}) \rightarrow
(\Omega, \mathcal{F}, \mathbb{P})$ is the
measure preserving transformation defined by:
\bea \label{Eq:theta}
\theta_t \omega(\cdot) = \omega(\cdot + t) - \omega(t).
\eea

In order the solution operator associated with Eq. ~\eqref{SEE} to satisfy the perfect cocycle property given in (ii)$'$ above, we will restrict our consideration to some subset in $\Omega$ of full measure which is also $\theta_t$-invariant for all $t \in \mathbb{R}$. In the following, we will identify such a subset, and introduce the restriction of the canonical MDS to this subset.

In that respect, let us consider the following scalar Langevin equation:
\begin{equation} \label{OUP} 
\mathrm{d}z + z \, \mathrm{d}t = \diffusion\,\mathrm{d}W_t.
\end{equation}

It is known that this equation possesses a unique stationary solution $z_\diffusion(\theta_t \omega)$ --- the stationary Ornstein-Uhlenbeck (OU) process --- whose main properties are in particular recalled in the following lemma.

\bl  \label{Lem:OU}

There exists a subset $\Omega^\ast$ of $\Omega$ which is of full measure and is $\theta_t$-invariant for all $t\in \mathbb{R}$, {\it i.e.},
\bea  \label{Eq:z-1}
\mathbb{P}(\Omega^\ast) = 1, \quad \text{and} \quad \theta_t (\Omega^\ast) = \Omega^\ast \quad \Forall \, t \in 
\mathbb{R};
\eea
and the following properties hold on $\Omega^\ast$:

\bi

\item[(i)] For each $\omega \in \Omega^\ast$, $t \mapsto W_t(\omega)$ is $\gamma$-H\"older  for any $\gamma \in (0, 1/2)$.

\item[(ii)] $t\mapsto W_t(\omega)$ has sublinear growth:
\bea  \label{Wt control}
\lim_{t\rightarrow \pm \infty}\frac{W_t(\omega)}{t} = 0, \quad \Forall \, \omega \in \Omega^\ast.
\eea

\item[(iii)] 
For each $\omega \in \Omega^\ast$, $t \mapsto z_\diffusion(\theta_t \omega)$ is $\gamma$-H\"older for any $\gamma \in (0, 1/2)$, and can be written as:
\bea \label{z}
 z_\diffusion(\theta_t
\omega) & = -\diffusion \int_{-\infty}^0 e^\tau W_{\tau}(\theta_t
\omega) \, \mathrm{d} \tau \\
& = -\diffusion \int_{-\infty}^0 e^\tau W_{\tau + t}(\omega) \, \mathrm{d} \tau +\diffusion W_t(\omega), \quad  t\in \mathbb{R},\;\;  \omega \in \Omega^\ast.
\eea

\item[(iv)] The following growth control relations are satisfied:
\begin{equation} \label{Eq:z-4}
\lim_{t\rightarrow \pm \infty}\frac{z_\diffusion(\theta_t \omega)}{t} = 0,
 \text{ and }\lim_{t\rightarrow \pm \infty}
\frac{1}{t}\int_0^t z_\diffusion(\theta_\tau \omega) \, \mathrm{d}\tau = 0, \quad  \Forall \omega \in \Omega^\ast.
\end{equation}

\ei

\el

\bp

It is known that there exists a subset $\Omega' \subset \Omega$ of full measure which is $\theta_t$-invariant for all $t \in \mathbb{R}$, and $(ii)$, $(iv)$ and \eqref{z} hold on $\Omega'$; see {\it e.g. \cite{Caraballo_al04, Chueshov02, DLS03}}.\footnote{Regarding \eqref{Eq:z-4} see also \cite[Lem.~4.1]{Caraballo_al04}, and the Birkhoff-Khinchin ergodic theorem in {\it e.g.} ~\cite[p. 539]{Arnold98}. We will make very often use of these growth control relations in the proofs of Theorems \ref{attractiveness thm} and \ref{App CMF} produced below.}

{\mkk  For the sake of clarity, we explain here how to exhibit a subset $\Omega'' \subset \Omega$ of full measure which is $\theta_t$-invariant and for which $t\mapsto z_\sigma(\theta_t\omega)$  is (locally) $\gamma$-H\"older for any $\gamma \in (0,1/2)$ and any $\omega \in \Omega''.$}
{\mkk Note that by simply integrating Eq. ~\eqref{OUP} we get trivially that:
\bea  \label{OU identity} 
\int_0^t z_\diffusion(\theta_s \omega)\, \mathrm{d}s + z_\diffusion(\theta_t \omega) = z_\diffusion(\omega) + \diffusion W_t(\omega), \quad \Forall t \in \mathbb{R}.
\eea 
From this identity, it is thus sufficient to exhibit such an $\Omega''$ where the $\gamma$-H\"older property holds for $t\mapsto W_t$.
}

It is known  that for each $m\in \mathbb{N}^\ast$ the Wiener process is almost surely $\gamma$-H\"older for any exponent $\gamma < \frac{1}{2} - \frac{1}{2m}$; see {\it e.g.} \cite[Chap. 3C]{Evans01}.  For each $m\ge 1$, let $\Omega_m$ be a corresponding subset of full measure. By possibly redefining $\Omega_m$ to be $\Omega_m \cap \Omega_{m-1}$ so that $\Omega_{m} \subset \Omega_{m-1}$, we deduce naturally that
\bea
\Omega'' := \bigcap_{m=1}^\infty \Omega_m
\eea
is such that $\mathbb{P}(\Omega'') = \lim_{m\rightarrow \infty} \mathbb{P}(\Omega_m) = 1$; see {\it e.g.} ~\cite[Thm ~1.8]{Folland99}. 

By construction of $\Omega''$, any sample path $\omega$ in  $\Omega''$ is $\gamma$-H\"older for any exponent less than $1/2$. Trivially, $\theta_t\omega$ shares  this H\"older property for all $\omega\in \Omega''$  since $\theta_t$ is just a shift on $\Omega$; see \eqref{Eq:theta}. The set  $\Omega''$ is thus $\theta_t$-invariant for all $t \in \mathbb{R}$, and the proof is complete by taking $\Omega^\ast := \Omega' \cap \Omega''$. 

\ep

Now, let $\Omega^\ast$ be the $\theta_t$-invariant subset of $\Omega$ as given in Lemma ~\ref{Lem:OU}, and $\mathcal{F}_{\Omega^\ast}$ be the trace $\sigma$-algebra of $\mathcal{F}$ with respect to $\Omega^\ast$, {\it i.e.},
\beas
\mathcal{F}_{\Omega^\ast} := \{F \cap \Omega^\ast \; \vert \; F  \in \mathcal{F}\}.
\eeas
It can be checked that $(t, \omega) \rightarrow \theta_t\omega$ is $(\mathcal{B}(\mathbb{R})\otimes \mathcal{F}_{\Omega^\ast}; \mathcal{F}_{\Omega^\ast})$-measurable; see for instance \cite[Lemma 3.2]{Caraballo_al09}. It follows that $(\Omega^\ast,\,
\mathcal{F}_{\Omega^\ast}, \mathbb{P}_{\Omega^\ast}, \{\theta_t\}_{t \in \mathbb{R}})$ forms an MDS, where $\mathbb{P}_{\Omega^\ast}$ denotes the restriction of $\mathbb{P}$ on $\Omega^\ast$.

To simplify the notations and the presentation, we will denote hereafter the new sample space $(\Omega^\ast,\,
\mathcal{F}_{\Omega^\ast}, \mathbb{P}_{\Omega^\ast})$ as $(\Omega,\,
\mathcal{F},\, \mathbb{P})$; and we will work with this restricted MDS, $(\Omega, \mathcal{F}, \mathbb{P}, \{\theta_t\}_{t\in \mathbb{R}})$, without confusion with the original MDS. For the sake of concision, we will often omit mentioning such an underlying MDS, thus identifying an RDS with its cocycle part.

\subsection{Cohomologous cocycles and random evolution equations} \label{sec. conjugacy}
Our approach to develop our  theory of stochastic critical manifolds for SPDEs of type  \eqref{SEE} along with the corresponding approximation formulas, will be rooted in the use of a smooth cohomology\footnote{{\mkk The theory of parameterizing manifolds presented in Part II of this article will also make use of such a cohomology.}}, {\it i.e.} a stationary coordinate change, on the phase space $\mathcal{H}_\alpha$. Such a strategy was initially motivated  by the study of certain asymptotic problems arising with finite-dimensional SDEs; see \cite{Arnold98,IL01,IL02,IS01}. The idea is to transform the original SDE into a random differential equation (RDE), {\it i.e.} a randomly parameterized ordinary differential equation, which is {\it de facto} more amenable for the analysis of random attractors or random invariant manifolds since these objects require a  pathwise framework {\mkk from their definitions}.  The cohomology offers then a sort of {\it transfer principle for these asymptotic problems}: if  a result is proved on, say, random invariant manifolds for the transformed RDE, then it holds true for the original SDE up to a stationary conjugacy.  This approach was adopted and extended to the context of SPDEs driven by multiplicative noise in \cite{FL05} for the study of random attractors,  and in \cite{DLS04} for the study of stochastic invariant manifolds.

 We follow here  more specifically the approach of \cite{DLS04} to set up the original SPDE \eqref{SEE} within the RDS framework. In that respect, we transform Eq. ~\eqref{SEE}
into an evolution equation with random coefficients which helps simplify the analysis of the dynamics associated with the SPDE.

{\mk As mentioned above, such a simplification} is rooted into the interpretation --- up to a {\mkk smooth cohomology} --- of the solutions of Eq.~\eqref{SEE} as  $\omega$-wise version of classical solutions of the transformed equation whose existence can be guaranteed from the literature of non-autonomous deterministic evolution equations \cite{Hen81, Robinson1, SY02,Tem97} although the measurability properties of such solutions require a special attention, which is carried out in Appendix ~\ref{Sect_mild}.\footnote{ Note that a direct approach at the level of the SPDE  {\mkk ({\it i.e.} without using random transformations)} could have been used to put Eq.~\eqref{SEE}  within the RDS framework \cite{MZZ08}. We adopted here the cohomology approach \cite{IL01,IL02,IS01} which makes possible the use of standard estimation techniques that do not rely on stochastic analysis  \cite{DPZ08}.  Of course the simple nature of our multiplicative noise makes efficient such a strategy here. For more general multiplicative noise, techniques from \cite{FL05} could be combined with the ones presented hereafter, to extend the main results of {\mkk Sections ~\ref{s:approximation}  and ~\ref{s:hyperbolic}.}} 
}

Let us {\mk now} introduce the following {\mk standard} change of variables:
\bea \label{random transform}
v(t) =e^{-z_\diffusion(\theta_t\omega)}u(t),
\eea
where $z_\diffusion$ is the OU process defined in \eqref{z}.
 
Note that by the It\^o formula (cf. ~\cite[Thm. ~4.1.2]{Oksendal98}), the stochastic process $e^{-z_\diffusion(\theta_t \omega)}$ satisfies
\bea \label{exp z eqn}
\mathrm{d}e^{-z_\diffusion(\theta_t\omega)} &= \Big( z_\diffusion(\theta_t\omega) e^{-z_\diffusion(\theta_t\omega)} + \frac{\diffusion^2}{2} e^{-z_\diffusion(\theta_t\omega)} \Big) \mathrm{d} t - \diffusion e^{-z_\diffusion(\theta_t\omega)}  \mathrm{d} W_t \\
& = z_\diffusion(\theta_t\omega) e^{-z_\diffusion(\theta_t\omega)} \mathrm{d} t - \diffusion e^{-z_\diffusion(\theta_t\omega)} \circ \mathrm{d} W_t,
\eea
where the second equality above follows from the conversion between the It\^o and Stratonovich integrals; cf. \cite[Thm.~2.3.11]{Kunita90}.

Formally, we also have that
\bea
\d v = \d \bigl(e^{-z_\diffusion(\theta_t\omega)}u \bigr) = u \circ \d e^{-z_\diffusion(\theta_t\omega)} + e^{-z_\diffusion(\theta_t\omega)} \circ \d u.
\eea 
Then, by using Eqns. ~\eqref{SEE}  and \eqref{exp z eqn} into the above equation, we find after simplification that  $v$ satisfies formally the following random evolution equation (hereafter referred to as an RPDE):
\begin{equation} \label{REE 1} 
\frac{\mathrm{d} v}{\mathrm{d} t} = L_\lambda v + z_\diffusion(\theta_t\omega)v + G(\theta_t\omega,\, v), 
\end{equation}
where $G(\omega,\, v) := e^{-z_\diffusion(\omega)} F(e^{z_\diffusion(\omega)}v)$. Note that $G(\omega, v)$ is globally Lipschitz in $v$ and has the same Lipschitz constant as $F$; and $G(\omega,\, 0) = 0$ for all $\omega \in \Omega$.

As mentioned above, the existence of pathwise solutions to Eq.~\eqref{REE 1} for each fixed $\omega$ can be guaranteed from the literature of non-autonomous deterministic evolution equations since in particular the OU process $z_\sigma$ is H\"older continuous in time; see Lemma ~\ref{Lem:OU}. Only the measurability properties of such solutions need a particular attention. We summarize the precise results in the following proposition. A full proof is provided in Appendix ~\ref{Sect_mild} for the reader's convenience.

\begin{prop} \label{prop:exist}

Consider the RPDE \eqref{REE 1}.  The assumptions on $L_\lambda$ and $F$ are those of \SS~\ref{ss:SEE} where $F$ is assumed to be globally Lipschitz here; see \eqref{Lip F}.

Then, for each $v_0\in \mathcal{H}_\alpha$, and $\omega\in \Omega$, Eq.~\eqref{REE 1} has a unique classical solution $v_\lambda(t, \omega; v_0)$ with initial datum $v_\lambda(0, \omega; v_0) = v_0$  such that
\be \label{v regularity}
v_\lambda(\cdot, \omega; v_0) \in C((0, T]; \mathcal{H}_1) \cap C([0, T]; \mathcal{H}_\alpha) \cap C^1((0, T]; \mathcal{H}), \quad \Forall T > 0.
\ee

Furthermore, for each $\lambda$,  $v_\lambda$ is $\big(\mathcal{B}(\mathbb{R}^+) \otimes \mathcal{F}\otimes \mathcal{B}(\mathcal{H}_\alpha); \mathcal{B}(\mathcal{H}_\alpha)\big)$-measurable; and for each $\omega$,  $v_\lambda(\cdot, \omega; v_0)$ depends continuously on $\lambda$ and $v_0$.
\end{prop}

It is worthwhile to remark that the conditions under which the conclusions of Proposition \ref{prop:exist} hold are obviously not  optimal but sufficient for our purpose. For instance, the Lipschitz condition on the nonlinearity can be replaced by certain dissipative conditions on the nonlinear terms which hold for a broad class of physical problems. We refer again to the aforementioned works \cite{CH98, Hen81, Robinson1, SY02,Tem97} for more details;  see also \cite{CLR13, DPZ08,Gess13,Fla08}.  

A direct consequence of this proposition is that, for any $\lambda$ and any  $(\mathcal{F}; \mathcal{B}(\mathcal{H}_\alpha))$-measurable random initial datum $v_0(\omega)$, there exists a unique classical solution $v_{\lambda, v_0(\omega)}(t, \omega):= v_\lambda(t, \omega; v_0(\omega))$ of Eq. ~\eqref{REE 1} which is $\big(\mathcal{B}(\mathbb{R}^+) \otimes \mathcal{F}; \mathcal{B}(\mathcal{H}_\alpha)\big)$-measurable, and $v_{\lambda, v_0(\omega)}(\cdot,\omega)$ has the regularity specified in \eqref{v regularity} for each $\omega$.

{\mkk Based} on this proposition, we can define for each $\lambda$, an RDS {\mkk generated by} Eq. ~\eqref{REE 1}, $S_\lambda$, as follows:
\bea \label{eq:S}
S_\lambda \colon  \mathbb{R}^+ \times \Omega \times \mathcal{H}_\alpha \rightarrow \mathcal{H}_\alpha, \quad (t,\omega,  v_0) \mapsto S_\lambda(t,\, \omega)v_0 :=  v_\lambda(t, \omega;  v_0), 
\eea
where $v_\lambda(t, \omega; v_0)$ is the global classical solution to Eq. ~\eqref{REE 1} with initial datum $v_0$, which can be either deterministic or an $\mathcal{H}_\alpha$-valued random variable.

We now introduce the  random (smooth) transformation {acting on the space $\mathcal{H}_\alpha$},
\begin{equation} \label{T transform}
{\tran}(\omega)  \xi :=   \xi e^{-z_\diffusion(\omega)},
\end{equation}
and its inverse random (smooth) transformation,
\begin{equation} \label{T inverse}
{\tran}^{-1}(\omega) \xi :=   \xi e^{z_\diffusion(\omega)},
\end{equation}
where $\xi \in \mathcal{H}_\alpha$, and $\omega \in \Omega$.

Let us now define the mapping $\widehat{S}_\lambda \colon \mathbb{R}^+ \times \Omega \times \mathcal{H}_\alpha \rightarrow \mathcal{H}_\alpha$ via
\bea \label{cohomology relation}
\widehat{S}_\lambda(t,\, \omega)u_0 := {\tran}^{-1}(\theta_t\omega) \circ S_\lambda(t, \omega) \circ {\tran}(\omega)u_0,
\eea
where the symbol, $\circ$, denotes the basic composition operation between self-mappings on $\mathcal{H}_{\alpha}$. The mapping  $\widehat{S}_\lambda$ thus defined is clearly measurable, and defines an RDS acting on $\mathcal{H}_\alpha$.

By a solution to the SPDE \eqref{SEE} with initial datum $u_0 \in \mathcal{H}_\alpha$, we always mean a process $u_\lambda(t, \omega; u_0) := e^{z_\diffusion(\theta_t\omega)}v_\lambda(t, \omega; u_0e^{-z_\diffusion(\omega)})$, where $v_\lambda$ is a classical solution of the RPDE \eqref{REE 1} with initial datum $u_0e^{-z_\diffusion(\omega)}$. In that sense, the RDS, $\widehat{S}_\lambda$, provides solutions to Eq. ~\eqref{SEE} since $u_\lambda(t, \omega; u_0) = \widehat{S}_\lambda(t,\, \omega)u_0$.

 \begin{rem}
{\mkk It is important to note  in \eqref{cohomology relation} that, while the cocycles $S_\lambda(t, \omega)$ and  $\widehat{S}_\lambda(t,\, \omega)$ map the $\omega$-fiber to the $\theta_t \omega$-fiber, the bijective random coordinate transformation $\tran$ maps only each fiber to itself, explaining the shift of fiber appearing in the inverse transformation ${\tran}^{-1}.$ }
\end{rem}

The two random dynamical systems $S_\lambda$ and $\widehat{S}_\lambda$ are {\mkk thus} {\it cohomologous} \cite{Arnold98,IL01} via the random $C^{\infty}$-diffeomorphism ${\tran}$ acting on $\mathcal{H}_{\alpha}$. As {\mkk pointed out above},  the phase portrait of the dynamics associated with Eq.~\eqref{SEE}, in $\mathcal{H}_{\alpha}$,  is {\mkk thus} fiber-wise $C^{\infty}$-diffeomorphic to the phase portrait of the dynamics  associated with Eq.~\eqref{REE 1}. Any regularity as well as topological properties of an invariant manifold for $S_\lambda$ is therefore preserved for $\widehat{S}_\lambda$. Since the random evolution equation  \eqref{REE 1}  is in a more convenient form for the analysis developed hereafter, we will keep proving our results for this equation. The corresponding {\mkk results} associated with Eq.~\eqref{SEE} will be then {\mkk transferred naturally from those proved for \eqref{REE 1} by application of} the cohomology relation \eqref{cohomology relation}; see {\it e.g.} Corollaries  \ref{Lip manfd SEE}--\ref{attractiveness thm SEE}, and \ref{thm:local SEE} below.

\subsection{Linearized stochastic flow and related estimates}  \label{ss:T}

The solution operator $\mathfrak{T}_{\lambda,\sigma}$ associated with the linearized equation (about the trivial steady state),
\bea \label{eq:linear-op}
\frac{\mathrm{d} v}{\mathrm{d} t} =  L_\lambda v + z_\diffusion(\theta_t\omega)v,
\eea
will play a key role in most of the estimates carried out in the forthcoming sections arising in various  {\mkk fixed point problems involving nonlinear} integral equations associated with Eq. ~\eqref{SEE}. We provide here in that respect some related estimates on $\mathfrak{T}_{\lambda,\sigma}$.

First note that since $z_\sigma(\theta_t \omega) \Id$ commutes obviously with $z_\sigma(\theta_{t'} \omega) \Id$ when $t \neq t'$, we obtain thus that $e^{\int_{0}^{t}  z_\sigma(\theta_\tau \omega) \, \mathrm{d}\tau \Id}$ defines the solution operator of the equation $\frac{\mathrm{d} v}{\mathrm{d} t} = z_\diffusion(\theta_t\omega)v$. Now, for any given $ t\ge 0$, by using the basic random change of variables $\widetilde{v} := e^{-\int_{0}^{t}  z_\sigma(\theta_\tau \omega) \, \mathrm{d}\tau \Id} v$, it can be checked that $v(t, \omega; v_0)$ solves Eq. ~\eqref{eq:linear-op} if and only if $\widetilde{v}(t; v_0)$ solves 
\bea \label{eq:linear-op2}
\frac{\mathrm{d} \widetilde{v}}{\mathrm{d} t} =  L_\lambda \widetilde{v},
\eea
where  $\Id \colon \mathcal{H} \rightarrow \mathcal{H}$ is the identity mapping and $v(t, \omega; v_0)$ denotes the solution of Eq. ~\eqref{eq:linear-op} emanating from $v_0$ at $t=0$ in the fiber $\omega$. 

Since $L_\lambda$ generates an analytic semigroup $e^{tL_\lambda}$ on $\mathcal{H}$, Eq. ~\eqref{eq:linear-op2} has a unique solution for each $v_0 \in \mathcal{H}$, so does Eq. ~\eqref{eq:linear-op} and its corresponding solution is given by:
\bea
v(t, \omega; v_0) = e^{\int_{0}^{t}  z_\sigma(\theta_\tau \omega) \, \mathrm{d}\tau \Id}e^{tL_\lambda}v_0, \quad \Forall t \ge 0, \; v_0 \in \mathcal{H}.
\eea
From this observation, we can associate Eq. ~\eqref{eq:linear-op} with a solution operator $\mathfrak{T}_{\lambda,\sigma}({t_f}, {t_s};\omega) \colon \mathcal{H} \rightarrow \mathcal{H}$ given by:
\bea \label{eq:ev-op}
\mathfrak{T}_{\lambda,\sigma}({t_f}, {t_s};\omega) v_0 & := v(t_f - t_s, \theta_{t_s}\omega; v_0) = e^{\int_{0}^{t_f - t_s}  z_\sigma(\theta_\tau \theta_{t_s}\omega) \, \mathrm{d}\tau \Id}e^{(t_f - t_s) L_\lambda}v_0\\
& = e^{\int_{t_s}^{t_f}  z_\sigma(\theta_\tau \omega) \, \mathrm{d}\tau \Id }e^{(t_f - t_s) L_\lambda}v_0, \quad \, \Forall  \, t_s \le  t_f,  \, v_0 \in \mathcal{H}.
\eea

Similar to $e^{tL_\lambda}$, the solution operator $\mathfrak{T}_{\lambda,\sigma}({t_f}, {t_s};\omega)$  leaves invariant the subspaces $\mathcal{H}^{\c}$ and $\mathcal{H}^{\s}$; and $\mathfrak{T}_{\lambda,\sigma}({t_f}, {t_s};\omega) P_{\c}$ can be defined for $t_f < t_s$ since $L_\lambda^{\c}$ is bounded on $\mathcal{H}^{\c}$.

From these simple facts about $\mathfrak{T}_{\lambda,\sigma}({t_f}, {t_s};\omega)$, we can easily derive the partial-dichotomy estimates for $\mathfrak{T}_{\lambda,\sigma}({t_f}, {t_s};\omega)$ from those associated with $e^{tL_\lambda}$ given in \eqref{Proj bounds}, which are summarized in the following

\bl

Let $P_{\c}$ and $P_{\s}$ be the projectors given in \eqref{Pc Ps}, $\mathfrak{T}_{\lambda,\sigma}$ be the solution operator  introduced above, and $\eta_1$ and $\eta_2$ be two constants satisfying \eqref{eta1-2}. Then, for each $\omega$ the following estimates hold:
\begin{subequations} \label{eq:dichotomy}
\begin{align}
& \|\mathfrak{T}_{\lambda,\sigma}({t_f}, {t_s};\omega) P_{\s}\|_{L(\mathcal{H}_\alpha, \mathcal{H}_\alpha)} \le K e^{{\eta_2} (t_f-t_s) + \int_{t_s}^{t_f} z_\sigma(\theta_\tau \omega)\, \mathrm{d}\tau}, &&  {t_s} \le {t_f},   \label{d-I} \\
& \|\mathfrak{T}_{\lambda,\sigma}({t_f},{t_s};\omega)  P_{\s}\|_{L(\mathcal{H}, \mathcal{H}_\alpha)}  \le \frac{K}{({t_f}-{t_s})^{\alpha}}
e^{{\eta_2} ({t_f}-{t_s}) + \int_{t_s}^{t_f} z_\sigma(\theta_\tau \omega)\, \mathrm{d}\tau}, &&  {t_s} < {t_f}, \label{d-II} \\
& \|\mathfrak{T}_{\lambda,\sigma}({t_f},{t_s};\omega) P_{\c}\|_{L(\mathcal{H}, \mathcal{H}_\alpha)}  \le K
e^{{\eta_1} ({t_f}-{t_s})   - \int_{t_f}^{t_s} z_\sigma(\theta_\tau \omega)\, \mathrm{d}\tau}, && {t_f} \le {t_s}.  \label{d-III}
\end{align}
\end{subequations}
\el

\br

Obviously, $e^{-\int_{0}^{t}  z_\sigma(\theta_\tau \omega) \, \mathrm{d}\tau \Id} v = e^{-\int_{0}^{t}  z_\sigma(\theta_\tau \omega) \, \mathrm{d}\tau} v$ for any $v \in \mathcal{H}$. We {\mkk kept} the operator $\Id$ in the above presentation, to emphasize that the linear multiplicative noise that we consider here acts ``diagonally''. More general linear multiplicative noise such as $M u \circ \d W_t$ could be considered where $M$ would be a bounded linear operator on $\mathcal{H}$; see \cite{MZZ08}. For the sake of simplicity, the main results of this article are presented in the case where $M = \Id$, regarding the approximation formulas of stochastic invariant manifolds (Section ~\ref{s:approximation}), the pullback characterization of approximating manifolds (Section ~\ref{s:pullback}), and the theory of parameterizing manifolds introduced in Part II.  The extension of these results to more general multiplicative noise cases will not be treated in this article. 

\er

\section{Existence and Attraction Properties of Global Stochastic Invariant Manifolds} \label{s:global}

In this section, we revisit the existence and attraction properties of families of global stochastic invariant manifolds for the SPDE \eqref{SEE} {\HL in order to set up the precise framework on which we will rely to present the {\mkk main results of the first part of this article} regarding approximation formulas of stochastic critical manifolds (Section ~\ref{s:approximation}) along with the related {\mkk pullback characterizations (Section ~\ref{Sec_app_man+PB}).}

Let us first recall the transformed RPDE associated with the SPDE \eqref{SEE}:
\begin{equation} \label{REE} 
\frac{\mathrm{d}u}{\mathrm{d}t} = L_\lambda u  + z_\diffusion(\theta_t\omega)u + G(\theta_t\omega,\, u),
\end{equation}
where 
\bea \label{eq:G}
G(\omega,\, u) = e^{-z_\diffusion(\omega)} F(e^{z_\diffusion(\omega)}u).
\eea

{\HL Following a classical approach \cite{CL88, CS01, DLS04, Vand89, VI92, VV87}, we will make use --- for the existence theory of random invariant manifolds --- of standard weighted Banach spaces $C_\eta^-$ given by:}
\bea \label{space C-}
C_\eta ^ {-}(\omega) :=\bigl \{ \phi \colon (-\infty,\, 0] \rightarrow \mathcal{H}_\alpha \mid \phi
& \text{ is continuous,  and} \\
& \hspace{0em} \sup_{t\le 0} e^{-\eta t + \int_t^0
z_\diffusion(\theta_\tau \omega)\, \mathrm{d} \tau}\|\phi(t)\|_\alpha < \infty \bigr \}, \quad {\HL \eta \in \mathbb{R}, \, \omega \in \Omega},
\eea
endowed with the norm:
\bea \label{eta norm}
\|\phi\|_{C_\eta ^ {-}(\omega)} := \sup_{t\le 0}e^{-\eta t + \int_t^0
z_\diffusion(\theta_\tau \omega)\, \mathrm{d} \tau}\|\phi(t)\|_\alpha.
\eea
{\HL Similar spaces will also be used in the study of attraction properties of stochastic invariant manifolds in \SS ~\ref{ss:attraction}.}

Throughout this article, we suppress the {\mkk $\omega$-dependence} of the space $C_{\eta}^-$ for simplicity, since for most of the results derived in this article, we work within a fixed fiber $\omega$. Only one exception is made in the proof of Theorem ~\ref{Lip manfd} {\mkk (about existence of random invariant manifolds)} where the $\omega$-dependence is specified explicitly to show the invariance property of the global random invariant manifolds which involves two fibers {\it per se}; see Step ~4 of the proof provided in Appendix ~\ref{appendix section 1}.

Recall that a set valued mapping $\mathfrak{C} \colon \Omega \rightarrow 2^H$ taking values in the closed (resp. ~compact) subsets of a separable Hilbert space $(H, |\cdot|_H)$ is called a {\it random closed (resp. ~compact) set} if for each $\xi \in H$, the map $\omega \rightarrow {\rm dist}(\xi, \mathfrak{C}(\omega))$ is measurable; see {\it e.g.} ~\cite[Def. ~2.1]{Crauel02}. Here and throughout the article, ${\rm dist}(\cdot, \cdot)$ denotes the Hausdorff semi-distance {\mkk associated with} the underlying  space $(H, |\cdot|_H)$, {\it i.e.},
\bea  \label{def:dist}
& {\rm dist}(D, E) := \sup_{a \in D} \inf_{b \in E} |a-b|_{H},  \Forall D, E \subset H.
\eea

We will also make use of the following definition of a  forward random invariant set of an RDS; see \cite[Def. ~6.9]{Crauel02}.
\bd \label{def invariance}

Let $S$ be a continuous RDS acting on some separable Hilbert space $H$ over some MDS, $(\Omega, \mathcal{F}, \mathbb{P}, \{\theta_t\}_{t\in \mathbb{R}})$. A random 
closed set $\mathfrak{B}$ is said to be forward invariant for this RDS if: 
\be\label{Eq_Binv}
S(t, \omega)\mathfrak{B}(\omega) \subset \mathfrak{B}(\theta_t \omega), \quad \Forall \, t >  0, \; \omega \in \Omega. 
\ee
If the inclusion in  \eqref{Eq_Binv} is replaced by a set equality, the random set $\mathfrak{B}$ is said to be forward strictly invariant.

\ed

Let us now introduce the following definition of a random (resp. ~stochastic) invariant manifold:

 \bd \label{def_GM}
Let $S$ be a continuous RDS acting on some separable Hilbert space $H$. A random closed set $\mathfrak{M}$ is called a global random invariant Lipschitz (resp. ~$C^r$ with $r \ge 1$) manifold of $S$ if the following conditions are satisfied:
\bi

\item[(i)] $\mathfrak{M}$ is invariant in the sense of Definition ~\ref{def invariance}.

\item[(ii)] There exist a closed subspace $\mathcal{G} \subset H$ and a measurable mapping $h \colon \mathcal{G} \times \Omega \rightarrow \mathcal{L}$, with $\mathcal{L}$ the topological complement of $\mathcal{G}$ in $H$, such that $\mathfrak{M}(\omega)$ can be represented as the graph of $h(\cdot, \omega)$ for each $\omega \in \Omega$, {\it i.e.}, 
\beas
\mathfrak{M}(\omega) = \{\xii + h(\xii, \omega) \mid \xii \in \mathcal{G}\},  \; \, \omega \in \Omega.
\eeas

\item[(iii)] For each $\omega$, $h(\cdot, \omega)$ is Lipschitz and $h(0, \omega) = 0$. In the case where  $\mathfrak{M}$ is  $C^r$  then  {\hh $h(\cdot, \omega)$ is furthermore $C^r$} and $\mathfrak{M}$ is tangent to $\mathcal{G}$ at the origin, {\it i.e.}, $h(0, \omega) = 0$ and $D_{\xii} h(0, \omega) = 0$, where $D_{\xii} h$ denotes the Fr\'echet derivative of $h$ with respect to $\xii$.

\ei

The function $h$ is called the global random invariant manifold function associated with $\mathfrak{M}$. 

\ed

In the case where the RDS is associated with a stochastic evolution equation of type \eqref{SEE}, we will refer such a manifold as a stochastic invariant Lipschitz (resp. $C^r$) manifold.

\subsection{Existence and smoothness of global stochastic invariant manifolds}  \label{ss:global existence}

In this subsection, we report on results concerning  the existence and smoothness of families of global stochastic (resp. ~random) invariant manifolds for the SPDE \eqref{SEE} (resp. ~the RPDE \eqref{REE}). These results are mainly known, but are framed here with respect to the uniform spectrum decomposition \eqref{gap 3} of the linear part. In particular, the analysis is tailored to make sure the dependence on $\lambda$ remains explicit.

 The approach adopted here is classical  but for the sake of completeness, a proof of the existence theorem below  is provided  in Appendix \ref{appendix section 1}.  
The inclusion of this proof is also motivated by the fact that some of its elements are used in the proofs of Theorem ~\ref{attractiveness thm} and of Theorem ~\ref{App CMF}; the latter constituting  the main result of the first part of this article, regarding the approximation formulas of stochastic critical manifolds.

More specifically, the proofs of the aforementioned  theorems are rooted in the Lyapunov-Perron approach \cite{BM61, KB34, Lyapunov47, Perron29}, and is based on techniques used for instance in \cite{CL88, Hen81, Vand89, VI92, VV87}, from which we propose a treatment adapted to the random setting that follows mainly the works of \cite{CS01, DLS04}.

\bt \label{Lip manfd}

Consider the RPDE \eqref{REE}. The assumptions on $L_\lambda$ and $F$ are those of \SS~\ref{ss:SEE} where $F$ is assumed to be globally Lipschitz with Lipschitz constant $\LF$ as given in \eqref{Lip F}. 
We assume that an open interval $\Lambda$ is chosen such that the uniform spectrum decomposition \eqref{gap 3} holds over $\Lambda$, and  that there exist $\eta_1$ and $\eta_2$ as in \eqref{eta1-2} for which {\hh the following uniform spectral gap condition holds:}
\be \label{thm 3.1 condition}
\boxed{\exists \; \eta \in (\eta_2, \eta_1) \; \; \mathrm{s.t.} \quad \Upsilon_{1}(F):=K \LF \bigl( ({\eta_1} - \eta)^{-1} + \Gamma(1-\alpha) (\eta - {\eta_2})^{\alpha - 1}\bigr) < 1,}
\ee
where $K$ is as given in the partial-dichotomy estimates \eqref{Proj bounds}, and $\Gamma(s) := \int_0^\infty t^{s-1} e^{-t}\, \mathrm{d} t$ is the Gamma function. Let the spaces $\mathcal{H}^{\c}$ (such that $\dim(\mathcal{H}^{\c})=m$) and  $\mathcal{H}_\alpha^{\s}$ be the corresponding subspaces associated with the uniform spectrum decomposition  as defined in \eqref{Hc} and \eqref{L decomposition}, respectively.


Then, for each $\lambda \in \Lambda$, there exists a global random  invariant Lipschitz manifold $\mathfrak{M}_\lambda$ of the RDS, $S_{\lambda}$, associated with Eq. ~\eqref{REE}.

 Each of such manifolds is  m-dimensional and  is given as 
\beas
\mathfrak{M}_\lambda(\omega) := \{\xii + h_\lambda(\xii, \omega) \mid \xii \in \mathcal{H}^{\c} \}, \quad  \omega \in \Omega, \; \lambda\in \Lambda,
\eeas
where, for each $\lambda\in \Lambda$, $h_\lambda$ is the solution of the following integral equation:\footnote{Note that $\Id$ denotes here the identity mapping of $\mathcal{H}$.}
 \begin{equation} \label{h eqn}
\begin{aligned}
 h_\lambda(\xii, \omega) = \int_{-\infty}^0 e^{\int_s^0 z_\diffusion(\theta_\tau \omega)
\,\mathrm{d} \tau \Id}e^{ -s L_\lambda} P_{\s} G\big(\theta_s\omega, u_\lambda(s,  \omega; \xii + h_\lambda(\xii, \omega)) \big)\,\mathrm{d}s, \; \xii \in \mathcal{H}^{\c}, \; \omega \in \Omega,
\end{aligned}
\end{equation}
{\HLL with $u_\lambda(\cdot, \omega; \xii + h_\lambda(\xii, \omega)) \in C_\eta^-$ being the mild solution\footnote{in the sense of Definition ~\ref{Def:mild solution} provided in Appendix ~\ref{Sect_mild}.} of  Eq.~\eqref{REE} on $(-\infty,  0]$  verifying $u_\lambda(0, \omega; \xii + h_\lambda(\xii, \omega)) = \xii + h_\lambda(\xii, \omega)$.}


Furthermore, $h_\lambda(\cdot, \omega) \colon \mathcal{H}^{\c} \rightarrow \mathcal{H}_\alpha^{\s}$ depends continuously on $\lambda$, and its Lipschitz constant satisfies
\be \label{Lip h}
{\HLL \mathrm{Lip}(h_\lambda(\cdot, \omega)) \le \frac{K^2\LF (\eta - {\eta_2})^{\alpha -1}\Gamma(1-\alpha)}{1- \Upsilon_{1}(F)}, \quad \omega \in \Omega.}
\ee

\et

{\HLL Since the Lipschitz constant of $h_\lambda(\cdot, \omega)$ admits a deterministic upper bound for each $\omega$ as seen in \eqref{Lip h}, for the sake of concision, we will denote it by $\LH$ hereafter.}

\needspace{1\baselineskip}
\br \label{rmk: cm dynamical interpretation}
\hspace*{2em}
\bi
\item[1)] The random invariant manifolds, $\mathfrak{M}_\lambda$, as provided by the above theorem, allow us to parametrize the high modes (in $\mathcal{H}_\alpha^{\s}$) of the solutions of Eq.~\eqref{REE} --- evolving on  $\mathfrak{M}_\lambda$ --- by their low modes (in $\mathcal{H}^{\c}$). However in the general case, $\mathfrak{M}_\lambda$ does not provide such a parameterization for all the solutions of Eq.~ \eqref{REE}.
 It turns out indeed that  in order those manifolds  to be genuine random (resp. stochastic) inertial manifolds, extra assumptions on the  spectral gap condition \eqref{thm 3.1 condition} are (as usual) required; see Theorem ~\ref{attractiveness thm}  (resp. Corollary  ~\ref{attractiveness thm SEE}).
 
 \item[2)]  An alternative class of useful manifolds to deal with the parameterization problem of the resolved modes by the unresolved ones, will be introduced in Part II.  These manifolds, called {\it parameterizing manifolds} (PMs) hereafter,  provide an approximate parameterization of the unresolved  variables by the resolved ones in a mean square sense.  In contrast to the classical theory, certain PMs (in the self-adjoint case) can be determined under a  {\it non-resonance conditions} which circumvent the constraints inherent to standard  spectral gap conditions (cf. ~\eqref{NR}-condition in Section  \ref{ss:PM} and \eqref{NR2}-condition in Section \ref{Sec_Existence_PQ});  see also Theorems  \ref{thm:hG_PM},  \ref{thm:h1_PM} (with Corollary \ref{Cor_h1_is_PM}), and  Theorem \ref{thm:h1_PMdet}.
\ei
\er

We provide in the remark below, useful interpretations from a dynamical perspective of $\mathfrak{M}_\lambda$ such as obtained {\it via} Theorem \ref{Lip manfd}.

\needspace{1\baselineskip}
\br\label{Rem_dyn_interp}

\hspace*{2em}

\bi
\item[1)] {\mkk When $\eta_2<\eta_1<0$}, the proof of Theorem ~\ref{Lip manfd} given in Appendix ~\ref{appendix section 1}  shows that $\mathfrak{M}_\lambda(\omega)$ consists of all 
elements in $\mathcal{H}_\alpha$ such that there exists a complete trajectory of Eq. ~\eqref{REE} passing through each such element at $t=0$, which has {\hh a controlled growth rate} as $t\rightarrow -\infty$. {\hh More precisely}, the $\mathcal{H}_\alpha$-norm of  such a solution is controlled by $e^{\eta t - \int_t^0 z_\diffusion(\theta_\tau \omega)\, \mathrm{d} \tau}$ as $t \rightarrow
-\infty$, where $\eta \in (\eta_2, \eta_1)$ is chosen such that the condition \eqref{thm 3.1 condition} holds, and $z_\diffusion(\theta_t \omega)$ is the OU process defined in \eqref{z}.

{\mkk Hyperbolic manifolds can be encountered in this case and are dealt with in details for a stochastic Burgers-type equation in Section ~\ref{s:Burgers}; see also Section ~\ref{s:hyperbolic}. }

\item[2)] {\HL When $\eta_2 < 0 < \eta_1$ in \eqref{Proj bounds}, we are in the case of a classical dichotomy. The dynamical interpretation of $\mathfrak{M}_\lambda$ for this case is more standard.} Indeed, if $\eta$ in \eqref{thm 3.1 condition} is positive, then $\mathfrak{M}_\lambda$ can be interpreted as the random closed set consisting of all {\HLL elements in $\mathcal{H}_\alpha$ such that the corresponding mild solutions decay to $0$ as $t \rightarrow -\infty$} since in this case, $e^{\eta t - \int_t^0 z_\diffusion(\theta_\tau \omega)\, \mathrm{d} \tau} \rightarrow 0$ as  $t \rightarrow -\infty$. The manifold $\mathfrak{M}_\lambda$ is {\mkk just} the random unstable manifold of the trivial steady state. 


\item[3)]  By adapting the proof of Theorem \ref{Lip manfd} accordingly it can be proved that for each $\lambda \in \Lambda$ the RDS associated with Eq. ~\eqref{REE} admits {\HL also a} random invariant manifold given by:
\bea \label{stable manifold}
\mathfrak{M}^{\s}_\lambda(\omega) = \{\zeta + h_\lambda(\zeta, \omega) \mid \zeta \in \mathcal{H}_\alpha^{\s} \}, \quad \omega \in \Omega,
\eea
where the invariant manifold function $h_\lambda(\cdot, \omega)$ is here defined over $\mathcal{H}_\alpha^{\s}$  with, this time, range in  $\mathcal{H}^{\c}$.

If in addition the classical dichotomy is satisfied, and $\eta$ in \eqref{thm 3.1 condition} is negative, then $\mathfrak{M}^{\s}_\lambda$ given in \eqref{stable manifold} corresponds to a random stable manifold which can be interpreted as the random closed set consisting of all initial data through which the mild solutions converges to $0$ as $t\rightarrow +\infty$. {\HL We will not make use of such manifolds in this article.}

\ei

\er

Concerning the smoothness of the random invariant manifold $\mathfrak{M}_\lambda$ given by Theorem ~\ref{Lip manfd}, we have the following

\bt \label{smth manfd}

Consider the RPDE ~\eqref{REE}. Assume that the assumptions of Theorem ~\ref{Lip manfd} hold with $\Lambda$ specified therein. The nonlinearity $F \colon \mathcal{H}_\alpha \rightarrow \mathcal{H}$ is furthermore  assumed to be  $C^p$ smooth for some integer $p\geq 1$ and to satisfy $DF(0) = 0$. Assume also that there exist $\eta \in \mathbb{R}$ and some integer $r \in \{1,...,p\}$, such that
\be
{\eta_2} < j \eta < {\eta_1},  \quad \Forall j \in \{1, \cdots, \,  r\},
\ee 
and {\mkk the following uniform spectral gap conditions hold}
\begin{equation} \label{thm 3.2 condition}
\begin{aligned}
\boxed{\Upsilon_j(F):=K \LF \bigl(({\eta_1} -
 \eta)^{- 1} + \Gamma(1-\alpha) (j \eta - {\eta_2})^{\alpha-1}  \bigr) < 1,   \quad \Forall j \in \{1, \cdots, \, r\}.}
\end{aligned}
\end{equation}

Then, for each $\lambda \in \Lambda$, the random invariant manifold $\mathfrak{M}_\lambda$ obtained in Theorem ~\ref{Lip manfd} is $C^r$ smooth. In particular, {\hh the random invariant manifold function, $(\xi, \omega) \mapsto h_\lambda(\xii, \omega)$, guaranteed therein is $C^r$ smooth with respect to $\xii$ for each $\omega$; and is tangent to the subspace $\mathcal{H}^{\c}$ at the origin}, {\it i.e.}, $h_\lambda(0, \omega)  = 0$, and $D_{\xii} h_\lambda(0, \omega) = 0$.


\et

\bp

The proof follows essentially the same lines of the proof of \cite[Thm.~4.1]{DLS04} and is  omitted here.

\ep

By the cohomology relation \eqref{cohomology relation} between the RDSs associated {\hh respectively} with Eq. ~\eqref{SEE} and Eq. ~\eqref{REE}, {\hh we obtain} the following corollaries regarding the existence and smoothness of families of stochastic invariant manifolds $\widehat{\mathfrak{M}}_\lambda$ for Eq. ~\eqref{SEE}.

 \bc \label{Lip manfd SEE}

Consider the SPDE ~\eqref{SEE}. Assume that the assumptions of Theorem \ref{Lip manfd} hold with $\Lambda$ specified therein. 
Let also $h_\lambda \colon \mathcal{H}^{\c} \times \Omega \rightarrow \mathcal{H}_\alpha^{\s}$ be the random invariant manifold function {\hh associated with Eq. ~\eqref{REE}} provided by Theorem ~\ref{Lip manfd}. Then, for each $\lambda \in \Lambda$, the RDS, $\widehat{S}_\lambda$, as defined in \eqref{cohomology relation} possesses an $m$-dimensional global stochastic invariant Lipschitz manifold $\widehat{\mathfrak{M}}_\lambda$ given by:
\bea \label{M of SEE}
\widehat{\mathfrak{M}}_\lambda(\omega) := \{\xii + e^{z_\diffusion(\omega)}h_\lambda(e^{-z_\diffusion(\omega)}\xii, \omega) \mid \xii \in \mathcal{H}^{\c} \}, \quad  \omega \in \Omega.
\eea

\ec

\bc \label{smth manfd SEE}

Consider the SPDE \eqref{SEE}. Assume that the assumptions of Theorem ~\ref{smth manfd} hold with $r$ specified therein. Then for each $\lambda \in \Lambda$, $\widehat{\mathfrak{M}}_\lambda$ defined in \eqref{M of SEE} is a global
stochastic invariant $C^r$ manifold of $\widehat{S}_\lambda$.

\ec

\subsection{Asymptotic completeness of stochastic invariant manifolds}  \label{ss:attraction}

We {\HL establish} in this subsection the attraction properties of random invariant manifolds as stated in Theorem \ref{attractiveness thm} {\sw below}. As explained in Section ~\ref{Sec_background}, the corresponding results may be viewed as {\mkk complementary to} previous results obtained on the topic \cite{BF95, CG95, CS01, DPD96, Schmalfuss05, SZ12, WD07}. Furthermore, some new insights concerning the asymptotic completeness problem are provided. Our approach is inspired by \cite[Thm.~5.1]{CL88} that we adapt to our framework; see also \cite{CS01}. More precisely, for a given solution $u_\lambda$ to Eq. ~\eqref{REE} we look for a solution $\overline{u}_\lambda$ living on the random invariant manifold $\mathfrak{M}_\lambda$ such that $\|\overline{u}_\lambda(t, \omega) - u_\lambda(t, \omega)\|_\alpha$ decays exponentially as $t\rightarrow \infty$ for almost all $\omega$. The strategy adopted here consists of reformulating this problem as a fixed point problem under the constraint that the sought solution $\overline{u}_\lambda$ emanates from {\HL an initial datum $\overline{u}_0$ which belongs to $\mathfrak{M}_\lambda$ and is} well-prepared with respect to the given initial datum $u_0$. This problem is then recast as an unconstrained fixed point problem associated with a random integral operator that is solved by means of the uniform contraction mapping principle \cite[Theorems 2.1--2.2]{CH82}. The price to pay is an additional condition \eqref{thm 3.5 condition 2} to the existence theory which involves the spectral gap and the Lipschitz constant associated with $\mathfrak{M}_\lambda$. The results so obtained show that the manifolds obtained in \SS ~\ref{ss:global existence} under this additional condition are, for $\eta < 0$, almost surely forward and pullback  asymptotically complete at a {\HL {\it uniform rate}} with respect to the parameter $\lambda$.  By doing so, we establish naturally the existence of stochastic inertial manifolds, which are almost surely exponentially attracting at a uniform rate in  both a forward and a pullback sense.
Before presenting the results, let us introduce the following definition of forward (resp. ~pullback) asymptotic completeness of a stochastic invariant manifold, which is motivated by its deterministic analogue;  see {\it e.g.} ~\cite{CFNT89, Rob96}.

\bd  \label{asy. complete}

Let $S$ be a continuous RDS acting on some separable Hilbert space $H$ over some MDS, $(\Omega, \mathcal{F}, \mathbb{P}, \{\theta_t\}_{t\in \mathbb{R}})$, and $\mathfrak{M}$ be a random invariant manifold of $S$. The manifold $\mathfrak{M}$ is said to be forward  asymptotically complete if there exist $\kappa > 0$ and a positive random variable $C_{\kappa}$ such that for any $H$-valued tempered random initial datum\footnote{See {\it e.g.} ~\cite[Def.~1.3.3]{Chueshov02} for the definition of a tempered random variable.} $u_0$ there exists a tempered random variable $v_0$ on the manifold $\mathfrak{M}$ satisfying
\bea\label{Eq_attract_fwd}
& \| S(t, \omega)u_0(\omega) - S(t, \omega)v_0(\omega)\| \le {\hll C_{\kappa}(\omega)} \|u_0(\omega) - v_0(\omega)\| e^{- \kappa t}, \quad  t\ge 0, \: \omega \in \Omega.\\
\eea

It is said to be pullback asymptotically complete if under the same conditions the following attraction property holds:
\bea\label{Eq_attract_pba}
&   \| S(t, \theta_{-t}\omega)u_0(\theta_{-t} \omega) - S(t, \theta_{-t} \omega)v_0(\theta_{-t} \omega)\| \\
& \hspace{6em} \le {\hll C_{\kappa}(\omega)} \|u_0(\theta_{-t} \omega) - v_0(\theta_{-t}\omega)\|  e^{-\kappa t} , \quad  t\ge 0, \: \omega \in \Omega.
\eea

The rate $\kappa$ is called an attraction rate. The supremum of such attraction rates is called the critical attraction rate.\footnote{Note that  such supremum is not an attraction rate in general.}

If furthermore $\mathfrak{M}(\omega)$ is finite-dimensional of fixed dimension for all $\omega$, then $\mathfrak{M}$ is called a pullback (resp. forward) random inertial manifold  in case where  \eqref{Eq_attract_pba} (resp. \eqref{Eq_attract_fwd}) holds.
\ed

\bt \label{attractiveness thm}

Consider the RPDE ~\eqref{REE}. Assume that the assumptions of Theorem ~\ref{Lip manfd} hold with $\Lambda$ and $\Upsilon_{1}(F)$ specified therein. Assume also that $\eta$ in condition \eqref{thm 3.1 condition} is negative. If the invariant manifold function $h_\lambda$ guaranteed by Theorem ~\ref{Lip manfd} satisfies furthermore 
\bea \label{thm 3.5 condition 2}
 & \frac{K{\rm Lip}(h_{\lambda})}{1- \Upsilon_{1}(F)} < 1,
\eea
then for each $\lambda \in \Lambda$ the random invariant manifold $\mathfrak{M}_\lambda$ obtained by Theorem \ref{Lip manfd} is both a pullback and a forward random inertial manifold with critical attraction rate $|\eta|$. 

\et
Once more, by using the cohomology relation between the two RDSs associated respectively with Eqns. ~\eqref{SEE} and \eqref{REE}, we obtain the following results on the forward and pullback asymptotic completeness of the stochastic invariant manifolds $\widehat{\mathfrak{M}}_\lambda$ for  Eq.~\eqref{SEE} {\HL as provided by} Corollary ~\ref{Lip manfd SEE}.

\bc \label{attractiveness thm SEE}

Consider the SPDE ~\eqref{SEE}. Assume that the assumptions of Theorem ~\ref{attractiveness thm} hold with $\Lambda$ and $\eta$ specified therein. Then for each $\lambda \in \Lambda$ the stochastic invariant manifold $\widehat{\mathfrak{M}}_\lambda$ guaranteed by Corollary ~\ref{Lip manfd SEE} is both a pullback and a forward stochastic inertial manifold with critical attraction rate $|\eta|$. 


\ec

\br 
We warn the reader on the fact that for a local invariant manifold, $\mathfrak{M}$,  defined as a random graph over a (deterministic) bounded neighborhood $\mathcal{V}$ of this basic state, contained in $\mathcal{H}^{\c}$ (see Corollary ~\ref{thm:local SEE}),  a pathwise asymptotic completeness property such as \eqref{Eq_attract_fwd} cannot be expected to hold in the the stochastic context. 
The reasons are due  to the theory of existence which imply that $\mathfrak{M}(\omega)$  (over $\mathcal{V}$) is uniformly bounded in $\omega$\footnote{This property results from  \eqref{Lip h} in Theorem ~\ref{Lip manfd} since in that case, $\|h_{\lambda}(\xi,\omega)\|_{\alpha} \leq \gamma  \mathrm{diam}(\mathcal{V})$ with $\gamma$ denoting the RHS of \eqref{Lip h}.}, on the one hand, and to the large deviations caused by the white noise \cite{CM10,Liu_LDP}, on the other.  The latter  imply that for any solution $u \neq 0$ of \eqref{SEE}, $\mathrm{dist} (u(t,\omega),\mathfrak{M}(\omega))$ 
cannot remain bounded as $t$ flows (for all $\omega$), in contradiction with \eqref{Eq_attract_fwd}. In particular \cite[Lemma 9]{CRD12} should be reconsidered in a probabilistic sense. 
\er

We first present the main ideas of the proof. 

\medskip
{\bf Ideas of the proof of Theorem ~\ref{attractiveness thm}.}  As mentioned above, we describe here what fixed point problem under constraint is {\HL naturally} related to the problem of asymptotic completeness  and how this problem can be solved by recasting it as an unconstrained fixed point problem. 

In that respect, we first point out in Step ~1 a natural integral equation to be satisfied by the difference, $v_\lambda:= \overline{u}_\lambda - u_\lambda$, between any solution $\overline{u}_\lambda$ of Eq. ~\eqref{REE} and a solution $u_\lambda$  emanating from a given $\mathcal{H}_\alpha$-valued tempered random initial datum $u_0$ 
\beas 
v_\lambda(t,\omega) & = \mathfrak{T}_{\lambda, \sigma}(t, t_0; \omega) v_\lambda(t_0, \omega)  +  \int_{t_0}^t \mathfrak{T}_{\lambda, \sigma}(t, s; \omega) \delta G(\theta_s\omega, v_\lambda(s, \omega)) \,\mathrm{d} s,   \quad 0 \le t_0 < t, \, {\HLL \omega \in \Omega,}
\eeas
where $\delta G(\theta_s\omega, v_\lambda(s, \omega)) = G(\theta_s\omega, u_\lambda(s, \omega) + v_\lambda(s, \omega))  - G(\theta_s\omega, u_\lambda(s, \omega))$, {\HL and $\mathfrak{T}_{\lambda, \sigma}$ is the linearized stochastic flow introduced in \SS ~\ref{ss:T}.

The above integral equation is then rewritten into the following fixed point problem:
\bea \label{int eqn 2-intro}
v_\lambda(t,\omega) = \mathcal{L}^{\omega, \lambda}_{\m} [v_\lambda](t),  \quad  t \ge 0, 
\eea
where
\bea \label{operator L-intro}
\mathcal{L}^{\omega, \lambda}_{\m} [v_\lambda](t) & := \mathfrak{T}_{\lambda, \sigma}(t, 0; \omega)  \m +  \int_{0}^t \mathfrak{T}_{\lambda, \sigma}(t, s; \omega)  P_{\s} \delta G(\theta_s\omega, v_\lambda(s,\omega))\,\mathrm{d} s \\
& \hspace{2em}  -  \int_t^{+\infty} \mathfrak{T}_{\lambda, \sigma}(t, s; \omega)  P_{\c} \delta G(\theta_s\omega, v_\lambda(s,\omega))\,\mathrm{d} s,  \quad  t \ge 0, 
\eea
with $q = P_{\s}v_\lambda(0, \omega)$. {\HL The advantage of doing so is to dispose of an operator $\mathcal{L}^{\omega, \lambda}_{\m}$ which is written under a form naturally adapted for establishing its contraction property based on the partial-dichotomy estimates \eqref{eq:dichotomy}; see Steps ~2 and 3.}

Now, {\HL given $v_\lambda[q]$ which solves \eqref{int eqn 2-intro}}, since $\overline{u}_\lambda = v_\lambda[q] + u_\lambda$, in order to show that $\overline{u}_\lambda$ lives on $\mathfrak{M}_\lambda$, it is sufficient to show that there exists $q \in \mathcal{H}_\alpha^{\s}$ such that $v_\lambda[q](0,\omega) + u_0(\omega) \in \mathfrak{M}_\lambda(\omega)$ for each $\omega \in \Omega$, due to the invariance property of $\mathfrak{M}_\lambda$. In  that respect, we introduce a weighted Banach space $C_\eta^+$,
\be \label{space C+}
C_\eta ^{+} :=\{ \phi \colon [0,\, \infty) \rightarrow \mathcal{H}_\alpha \mid \phi
\text{ is continuous and } \sup_{t\ge 0} e^{-\eta t - \int_0^t
z_\diffusion(\theta_\tau \omega)\, \mathrm{d} \tau}\|\phi(t)\|_\alpha < \infty \},\footnote{We suppress the {\mkk $\omega$-dependence of $C_\eta ^ +$  for the same reasons that were invoked for $C_\eta ^ -$  in \eqref{space C-}.}}
\ee
which is endowed with the {\HL following norm}:
$$|\phi|_{C_\eta ^{+}} := \sup_{t\ge 0}e^{-\eta t - \int_0^t z_\diffusion(\theta_\tau \omega)\, \mathrm{d} \tau}\|\phi(t)\|_\alpha.$$

If furthermore we show that $v_\lambda[q](\cdot, \omega)$ can be obtained for each $\omega$ as a fixed point of $\mathcal{L}^{\omega, \lambda}_{\m}$ in the space $C_\eta^+$ {\HL for some} $\eta < 0$  (independent of $\omega$) chosen according to \eqref{thm 3.1 condition}, then we get naturally that $v_\lambda[q](t, \omega) = \overline{u}_\lambda(t, \omega) - u_\lambda(t,\omega)$ approaches $0$ exponentially as $t \rightarrow \infty$, so that the {\HL asymptotic completeness problem is solved}.

The problem of searching for the desired solution $\overline{u}_\lambda$ which solves the asymptotic completeness problem is then recast into the problem of searching for a fixed point $v_\lambda[q](\cdot, \omega)$ of the integral operator $\mathcal{L}^{\omega, \lambda}_{\m}$ given by \eqref{operator L-intro} in the weighted Banach space $C_\eta^+$ under the constraint that $\overline{u}_0(\omega):=v_\lambda[q](0,\omega) + u_0(\omega) \in \mathfrak{M}_\lambda(\omega)$ for each $\omega \in \Omega$, namely, the following problem
\bea \label{FPC} 
\begin{cases}
\mathcal{L}^{\omega, \lambda}_{\m}[v] = v, \quad   v \in C_\eta^+, \quad {\HL \eta < 0,}\\
v(0,\omega) + u_0(\omega) \in \mathfrak{M}_\lambda(\omega), \quad {\HL \omega \in \Omega}. 
\end{cases}
\eea
The problem \eqref{FPC} is solved in two stages.

In the first stage, we deal with the fixed point problem $\mathcal{L}^{\omega, \lambda}_{\m}[v] = v$ (in $C_\eta^+$) disregarding first  the constraint, which is done in Steps ~2 and 3. It is shown in Step ~2 that $\mathcal{L}^{\omega, \lambda}_{\m}$ leaves $C_\eta^+$ stable, and in Step ~3 that $\mathcal{L}^{\omega, \lambda}_{\m}$ has a unique fixed point $v_\lambda[q](\cdot, \omega) \in C_\eta^+$ parameterized by $q = P_{\s}v_\lambda[q](0, \omega) \in \mathcal{H}_\alpha^{\s}$ using the uniform contraction mapping principle \cite[Theorems 2.1--2.2]{CH82}. Also derived in Step ~3 is a control of the $C_\eta^{+}$-norm of $v_\lambda[q](\cdot, \omega)$ in terms of $\|q\|_\alpha$ which is uniform with respect to $\omega$; see \eqref{v bound}. This will turn out to be important to ensure the temperedness of the initial datum $\overline{u}_0$ examined in Step ~5.

In the second stage, with the preparation carried out in the first stage we solve the constrained fixed point problem \eqref{FPC} by transforming it into a fixed point problem without constraint. This is done in Step ~4. 
First note that {\HL given a $q \in \mathcal{H}_\alpha^{\s}$ and the corresponding fixed point $v_\lambda[q](\cdot, \omega)$ of $\mathcal{L}^{\omega, \lambda}_{\m}$}, {\HL an initial datum  $\overline{u}_0(\omega)$ given by  $v_\lambda[q](0, \omega) + u_0(\omega)$ belongs to $\mathfrak{M}_\lambda(\omega)$ if $\overline{u}_0(\omega) = \p + h_\lambda(\p, \omega)$ for some $\p \in \mathcal{H}^{\c}$}, where $h_\lambda$ is the random invariant manifold function of Eq. ~ \eqref{REE} as {\HL provided by} Theorem ~\ref{Lip manfd}. {\HL Since $q = P_{\s} v_\lambda[\m](0, \omega)$}, it is thus natural to seek for $q$ of the following form:
\beas
q = P_{\s} (\overline{u}_0(\omega) - u_0(\omega)) = h_\lambda(\p, \omega) - P_{\s}u_0(\omega).
\eeas
The constraint $\overline{u}_0(\omega) \in \mathfrak{M}_\lambda(\omega)$ is then written as a functional equation in terms of $\p$:
\beas 
\p  = P_{\c}\overline{u}_0(\omega) & = P_{\c}v_\lambda[q](0, \omega) + P_{\c}u_0(\omega) \\
& = P_{\c}v_\lambda[h_\lambda(\p, \omega) - P_{\s}u_0(\omega)](0, \omega) + P_{\c}u_0(\omega).
\eeas
{\HL As a consequence}, the problem \eqref{FPC} can be rewritten as
\bea \label{FPC-2} 
\begin{cases}
\mathcal{L}^{\omega, \lambda}_{q(p)}[v] = v, \quad   v \in C_\eta^+, \\
\p = P_{\c}v(0, \omega) + P_{\c}u_0(\omega),
\end{cases}
\eea
where $q(p)= h_\lambda(\p, \omega) - P_{\s}u_0(\omega)$. Note also that for the fixed point $v_\lambda[q(p)](\cdot, \omega)$ of $\mathcal{L}^{\omega, \lambda}_{q(p)}$, $P_{\c}v_\lambda[q(p)](0, \omega)$ admits a natural {\HL integral representation} given by:
\bea \label{eq:op_I-intro}
P_{\c} v_\lambda[q(p)](0, \omega)  = - \int_{0}^{+\infty} \mathfrak{T}_{\lambda, \sigma}(0, s; \omega) P_{\c} \delta G (\theta_s\omega, v_\lambda[q(p)](s, \omega))\,\mathrm{d} s, \quad \Forall \; p \in \mathcal{H}^{\c},
\eea
which is obtained by applying the projector $P_{\c}$ to $\mathcal{L}^{\omega, \lambda}_{q(p)}$ given in \eqref{operator L-intro} and setting $t$ to $0$. In virtue of {\HL the integral representation} \eqref{eq:op_I-intro}, the problem \eqref{FPC-2} can be further transformed into the following fixed point problem in terms of $p$ but without any constraint:
\bea \label{eq:p}
\p = - \int_{0}^{+\infty} \mathfrak{T}_{\lambda, \sigma}(0, s; \omega) P_{\c} \delta G (\theta_s\omega, v_\lambda[q(p)](s, \omega))\,\mathrm{d} s + P_{\c}u_0(\omega)=: \mathcal{J}_{\omega, \lambda}(p),
\eea
{\HL where $v_\lambda[q(p)](\cdot, \omega)$ is the fixed point of $\mathcal{L}^{\omega, \lambda}_{q(p)}$ with $q(p) = h_\lambda(\p, \omega) - P_{\s}u_0(\omega)$.}

The unknown $p$ is then found using the contraction mapping principle applied to the operator $\mathcal{J}_{\omega, \lambda}$ defined on $\mathcal{H}^{\c}$.

In Step ~5, it is shown that the random initial datum $\overline{u}_0$ obtained from the previous step is tempered by using an estimate on $v_\lambda[q](0,\omega)$ derived in Step ~3 and the assumption that $u_0$ is tempered.

\bp [{\bf Proof of Theorem \ref{attractiveness thm}}]

{\HL We proceed in five steps following the ideas outlined above to prove both the forward and pullback asymptotic completeness properties of $\mathfrak{M}_\lambda$.}

\medskip

{\bf Step 1. Fixed point problem satisfied by $v_\lambda$.}  
%
Let us first point out a natural integral equation to be satisfied by the difference, {\HL $v_\lambda:= \overline{u}_\lambda - u_\lambda$}, between any solution {\HL $\overline{u}_\lambda$} of Eq. ~\eqref{REE} and a given solution $u_\lambda$  emanating from a given $\mathcal{H}_\alpha$-valued tempered random initial datum $u_0$: 
\bea \label{diff eqn}
v_\lambda(t,\omega)  = \mathfrak{T}_{\lambda, \sigma}(t, t_0; \omega) v_\lambda(t_0, \omega)  +  \int_{t_0}^t \mathfrak{T}_{\lambda, \sigma}(t, s; \omega) \delta G(\theta_s\omega, v_\lambda(s, \omega)) \,\mathrm{d} s, \quad 0 \le t_0 < t,
\eea
where 
\bea  \label{delta G}
\delta G(\theta_s\omega, v_\lambda(s, \omega)) := G(\theta_s\omega, u_\lambda(s, \omega) + v_\lambda(s, \omega))  - G(\theta_s\omega, u_\lambda(s, \omega)),
\eea
and $\mathfrak{T}_{\lambda, \sigma}$ is the solution operator defined in \SS~\ref{ss:T}.

For any given $\m \in \mathcal{H}_\alpha^{\s}$ and $\omega \in \Omega$, let us define then a mapping $\mathcal{L}^{\omega, \lambda}_{\m}$ acting on $C_\eta^{+}$ as follows:
\bea \label{int eqn 2}
\mathcal{L}^{\omega, \lambda}_{\m} [v](t) & := \mathfrak{T}_{\lambda, \sigma}(t, 0; \omega)  \m +  \int_{0}^t \mathfrak{T}_{\lambda, \sigma}(t, s; \omega)  P_{\s} \delta G(\theta_s\omega, v(s,\omega))\,\mathrm{d} s \\
& \hspace{2em}  -  \int_t^{+\infty} \mathfrak{T}_{\lambda, \sigma}(t, s; \omega)  P_{\c} \delta G(\theta_s\omega, v(s,\omega))\,\mathrm{d} s,  \quad  t \ge 0, \, v(\cdot, \omega) \in C_\eta^{+},
\eea 
where $\delta G$ is defined in \eqref{delta G}, {\HL and $C_\eta^+$ is the weighted Banach space defined in \eqref{space C+}.}

Similar to Step ~1 of the Proof of Theorem ~\ref{Lip manfd} given in Appendix ~\ref{appendix section 1}, one can show that  $v_\lambda(\cdot,\omega)$ is a solution to \eqref{diff eqn} such that  $v_\lambda(\cdot,\omega) \in C_\eta ^ +$ if and only if $v_\lambda(\cdot,\omega)$ is a fixed point of the operator $\mathcal{L}^{\omega, \lambda}_{\m}$  in $ C_\eta ^ +$, with $\m = P_{\s} v_\lambda(0, \omega)$.

{\HL The rest of the proof is devoted to solving the following fixed point problem:
\bea \label{FPC-recall} 
\begin{cases}
\mathcal{L}^{\omega, \lambda}_{\m}[v] = v, \quad   v \in C_\eta^+, \quad \eta < 0,\\
v(0,\omega) + u_0(\omega) \in \mathfrak{M}_\lambda(\omega), \quad  \omega \in \Omega,
\end{cases}
\eea
which as explained above and in Step ~5 below is sufficient to solve the asymptotic completeness problem.
}

\medskip
{\bf Step 2. $\mathcal{L}^{\omega, \lambda}_{\m}$ leaves $C_\eta^+$ stable.} For each $\m\in \mathcal{H}_\alpha^{\s}$ and $\omega \in \Omega$, we show in this step that 
\bea \label{eq:L stable}
\mathcal{L}^{\omega, \lambda}_{\m} C_\eta ^ + \subset C_\eta ^ +.
\eea 

Since $\eta > \eta_2$ according to condition \eqref{thm 3.1 condition} and $\m \in \mathcal{H}_\alpha^{\s}$, it follows from the partial-dichotomy estimate \eqref{d-I} that
\bea \label{T est}
 \sup_{t\ge 0}e^{-\eta t - \int_{0}^t z_\diffusion(\theta_s \omega)\,\mathrm{d}s} \| \mathfrak{T}_{\lambda, \sigma}(t, 0; \omega) \m \|_\alpha 
 \le  \sup_{t\ge 0} Ke^{({\eta_2} - \eta)t }\|\m\|_\alpha = K \|\m\|_\alpha.
\eea

Recalling that $G$ and $F$ have the same Lipschitz constant $\LF$, we then derive from \eqref{delta G} that $\| \delta G(\theta_s\omega, v(s, \omega))\| \le \LF \|v(s, \omega)\|_\alpha$ for any $s\ge 0$ and $v(\cdot, \omega) \in C_\eta^+$. This together with \eqref{d-II} implies that
\bea  \label{T map-1}
&e^{-\eta t - \int_{0}^t z_\diffusion(\theta_s \omega)\,\mathrm{d}s} \Big\|\int_{0}^t \mathfrak{T}_{\lambda, \sigma}(t, s; \omega)  P_{\s} \delta G(\theta_s\omega, v(s, \omega))\,\mathrm{d} s\Big\|_\alpha \\
& \le K e^{-\eta t - \int_{0}^t z_\diffusion(\theta_s \omega)\,\mathrm{d}s}\int_{0}^t \frac{e^{{\eta_2} (t-s) + \int_s^t z_\diffusion(\theta_\tau \omega) \,\mathrm{d} \tau}}{(t-s)^\alpha} \| \delta G(\theta_s\omega, v(s, \omega))\| \,\mathrm{d} s \\
& \le K \LF e^{-\eta t - \int_{0}^t z_\diffusion(\theta_s \omega)\,\mathrm{d}s} \int_{0}^t \frac{e^{{\eta_2} (t-s) + \int_s^t z_\diffusion(\theta_\tau \omega) \,\mathrm{d} \tau}}{(t-s)^\alpha}\|v(s, \omega)\|_\alpha \,\mathrm{d} s \\
& \le K \LF e^{-\eta t - \int_{0}^t z_\diffusion(\theta_s \omega)\,\mathrm{d}s} \int_{0}^t \frac{e^{{\eta_2} (t-s) + \eta s + \int_0^t z_\diffusion(\theta_\tau \omega) \,\mathrm{d} \tau}}{(t-s)^\alpha} \,\mathrm{d} s \|v(\cdot, \omega)\|_{C^+_\eta} \\
& = K \LF \int_{0}^t \frac{e^{({\eta_2} - \eta) (t-s)}}{(t-s)^\alpha}\,\mathrm{d} s \|v(\cdot, \omega)\|_{C^+_\eta} \\
& \le \frac{K \LF\Gamma(1-\alpha)}{(\eta - {\eta_2})^{1-\alpha}} \|v(\cdot, \omega)\|_{C^+_\eta}, \quad \Forall t \ge 0, \, {\HLL v(\cdot, \omega) \in C_\eta^+}.
\eea
Similarly, we have for all $t\ge 0$ {\HLL and $v(\cdot, \omega) \in C_\eta^+$} that
\bea  \label{T map-2}
e^{-\eta t - \int_{0}^t z_\diffusion(\theta_s \omega)\,\mathrm{d}s} \Big\| \int_t^{+\infty} \mathfrak{T}_{\lambda, \sigma}(t, s; \omega)  P_{\c} \delta G(\theta_s\omega, v(s, \omega))\,\mathrm{d} s \Big \|_\alpha \le \frac{K \LF}{{\eta_1} - \eta} \|v(\cdot, \omega)\|_{C^+_\eta}.
\eea

Now, \eqref{eq:L stable} follows from \eqref{T est}--\eqref{T map-2}.

\medskip
{\bf Step 3. Unique fixed point of $\mathcal{L}^{\omega, \lambda}_{\m}$.} In this step, we show that $\mathcal{L}^{\omega, \lambda}_{\m}$ has a unique fixed point $v_\lambda[\m](\cdot, \omega) \in C_\eta ^ +$ for each $\m \in \mathcal{H}_\alpha^{\s}$ and $\omega \in \Omega$; and that $v_\lambda[\m] \colon [0, \infty) \times \Omega \rightarrow \mathcal{H}_\alpha$ is $(\mathcal{B}([0, \infty)) \otimes \mathcal{F}; \mathcal{B}(\mathcal{H}_\alpha))$-measurable. The Lipschitz continuity of $v_\lambda[\m](t, \omega)$ with respect to $\m$ is also examined.

For any $v_1(\cdot, \omega)$, $ v_2(\cdot, \omega) \in C_\eta^+$, following the same type of estimates as given in \eqref{T map-1}--\eqref{T map-2}, we obtain
\bea \label{Lip L}
&\|\mathcal{L}^{\omega, \lambda}_{\m}[v_1]  - \mathcal{L}^{\omega, \lambda}_{\m}[v_2] \|_{C^+_\eta}  \\
& \hspace{0em} \le \sup_{t\ge 0} e^{-\eta t - \int_{0}^t z_\diffusion(\theta_s \omega)\,\mathrm{d}s} \Big\| \int_{0}^t \mathfrak{T}_{\lambda, \sigma}(t, s; \omega) P_{\s} \big [\delta G(\theta_s\omega, v_1) - \delta G(\theta_s\omega, v_2) \big] \,\mathrm{d} s \Big \|_\alpha \\
& \hspace{2em} + \sup_{t\ge 0} e^{-\eta t - \int_{0}^t z_\diffusion(\theta_s \omega)\,\mathrm{d}s} \Big\| \int_t^{+\infty} \mathfrak{T}_{\lambda, \sigma}(t, s; \omega) P_{\c} \big[\delta G(\theta_s\omega, v_1) - \delta G(\theta_s\omega,  v_2) \big ] \,\mathrm{d} s \Big \|_\alpha \\
& \le  K \LF  \big(({\eta_1} - \eta)^{-1}+ \Gamma(1-\alpha) (\eta- {\eta_2})^{\alpha-1} \big) \|v_1(\cdot, \omega) - v_2(\cdot, \omega)\|_{C^+_\eta} \\
& = \Upsilon_{1}(F) \|v_1(\cdot, \omega) - v_2(\cdot, \omega)\|_{C^+_\eta},
\eea
where $\Upsilon_{1}(F) \in (0, 1)$ is given in \eqref{thm 3.1 condition}. This implies that $\mathcal{L}^{\omega, \lambda}_{\m}$ is a uniform contraction in the space $C_\eta^{+}$ with respect to $\m$. By the uniform contraction mapping principle, see {\it e.g.} ~\cite[Theorems 2.1--2.2]{CH82}, the operator $\mathcal{L}^{\omega, \lambda}_{\m}$ has a unique fixed point $v_\lambda[\m](\cdot, \omega) \in C_\eta^{+}$ for any $\m \in \mathcal{H}_\alpha^{\s}$.  

The measurability of $v_\lambda[\m]\colon [0, \infty) \times \Omega \rightarrow \mathcal{H}_\alpha$ can be obtained in the same fashion as in Step ~2 of the proof in Appendix ~\ref{appendix section 1}, {\HL by relying on the Picard scheme naturally associated with the fixed point problem $\mathcal{L}^{\omega, \lambda}_{\m}[v] = v$.}

We turn now to the Lipschitz continuity of $v_\lambda[\m]$ with respect to $\m$. According to the definition of $\mathcal{L}^{\omega, \lambda}_{\m}$ given in \eqref{int eqn 2}, we get
\bea \label{Lip L 2}
\|\mathcal{L}^{\omega, \lambda}_{\m_1}[v]  - \mathcal{L}^{\omega, \lambda}_{\m_2}[v]\|_{C^+_\eta} & =  \|\mathfrak{T}_{\lambda, \sigma}(\cdot, 0; \omega) \m_1 - \mathfrak{T}_{\lambda, \sigma}(\cdot, 0; \omega) \m_2  \|_{C^+_\eta} \\
&   =\sup_{t\ge 0} e^{-\eta t - \int_0^t
z_\diffusion(\theta_\tau \omega)\, \mathrm{d} \tau}
\|\mathfrak{T}_{\lambda, \sigma}(t, 0; \omega)(\m_1 - \m_2)  \|_\alpha \\
& \le (\text{ by \eqref{T est}}) \\
& \le K \|\m_1 - \m_2\|_\alpha, \; \Forall \m_1, \m_2 \in \mathcal{H}_\alpha^{\s}, \, v(\cdot, \omega) \in C_\eta^+.
\eea
Hence, 
\beas
 \|v_\lambda[\m_1](\cdot, \omega)  & - v_\lambda[\m_2](\cdot, \omega)\|_{C^+_\eta} \\
& = \|\mathcal{L}^{\omega, \lambda}_{\m_1}[v_\lambda[\m_1]]  - \mathcal{L}^{\omega, \lambda}_{\m_2} [v_\lambda[\m_2]] \|_{C^+_\eta} \\
&\le\|\mathcal{L}^{\omega, \lambda}_{\m_1}[v_\lambda[\m_1]] - \mathcal{L}^{\omega, \lambda}_{\m_1}[ v_\lambda[\m_2]]\|_{C^+_\eta} +  \|\mathcal{L}^{\omega, \lambda}_{\m_1}[v_\lambda[\m_2]]   - \mathcal{L}^{\omega, \lambda}_{\m_2}[v_\lambda[\m_2]] \|_{C^+_\eta}  \\
&\le (\text{ by \eqref{Lip L} and \eqref{Lip L 2}}) \\
& \le \Upsilon_{1}(F) \|v_\lambda[\m_1](\cdot, \omega) - v_\lambda[\m_2](\cdot, \omega)\|_{C^+_\eta} + K\|\m_1 - \m_2\|_\alpha.
\eeas
We then obtain
\bea \label{Lip v}
 \|v_\lambda[\m_1](\cdot, \omega)  - v_\lambda[\m_2](\cdot, \omega)\|_{C^+_\eta} \le \frac{K}{1- \Upsilon_{1}(F)}\|\m_1 - \m_2\|_\alpha, \quad \Forall \m_1, \m_2 \in \mathcal{H}_\alpha^{\s}.
\eea

Note that when $\m = 0$, $v\equiv 0$ satisfies \eqref{int eqn 2}. By the uniqueness of the fixed point, we have $v_\lambda[0](t, \omega) \equiv 0$. This together with \eqref{Lip v} implies that
\be \label{v bound}
\|v_\lambda[\m](\cdot, \omega)\|_{C^+_\eta} \le \frac{K}{1- \Upsilon_{1}(F)}\|\m\|_\alpha, \; \Forall \m \in \mathcal{H}_\alpha^{\s}.
\ee

\medskip
{\bf Step 4. Unconstrained formulation of the problem \eqref{FPC-recall}.}  We show in this step that for each $\omega \in \Omega$ there exists $q\in \mathcal{H}_\alpha^{\s}$ such that the constraint $\overline{u}_0(\omega):=v_\lambda[\m](0,\omega) + u_0(\omega) \in \mathfrak{M}_\lambda(\omega)$ is met {\HL for any given $u_0$.}

Since $v_\lambda[\m](\cdot, \omega)$ is the fixed point of $\mathcal{L}^{\omega, \lambda}_{\m}$ given in \eqref{int eqn 2}, it holds that
\bea \label{eq:q}
q = P_{\s} v_\lambda[\m](0, \omega). 
\eea
{\HL Note also that an initial datum  $\overline{u}_0(\omega)$ given by  $v_\lambda[q](0, \omega) + u_0(\omega)$ belongs to $\mathfrak{M}_\lambda(\omega)$ if $\overline{u}_0(\omega) = \p + h_\lambda(\p, \omega)$ for some $\p \in \mathcal{H}^{\c}$, where $h_\lambda$ is the random invariant manifold function of Eq. ~ \eqref{REE} as {\HL provided by} Theorem ~\ref{Lip manfd}. {\HL Since $q = P_{\s} v_\lambda[\m](0, \omega)$}, it is thus natural to seek for $q$ of the following form}:
\bea
q = P_{\s} v_\lambda[\m](0, \omega) = P_{\s} (\overline{u}_0(\omega) - u_0(\omega)) = h_\lambda(\p, \omega) - P_{\s}u_0(\omega).
\eea
{\HL The constraint $\overline{u}_0(\omega) \in \mathfrak{M}_\lambda(\omega)$ is then written as a functional equation in terms of $\p$}:
\bea \label{functional eqn}
\p  = P_{\c}\overline{u}_0(\omega) & = P_{\c}v_\lambda[q](0, \omega) + P_{\c}u_0(\omega) \\
& = P_{\c}v_\lambda[h_\lambda(\p, \omega) - P_{\s}u_0(\omega)](0, \omega) + P_{\c}u_0(\omega).
\eea
Note also that $P_{\c} v_\lambda[h_\lambda(\p, \omega) - P_{\s}u_0(\omega)](0, \omega)$ has a natural integral representation. 
Indeed, since $v_\lambda[\m](\cdot, \omega)$ is the fixed point of $\mathcal{L}^{\omega, \lambda}_{\m}$ for each $q \in \mathcal{H}_\alpha^{\s}$, we obtain by applying  $P_{\c}$ to  \eqref{int eqn 2} and setting $t$ to $0$ that
\bea \label{eq:op_I}
P_{\c} v_\lambda[\m](0, \omega)  = - \int_{0}^{+\infty} \mathfrak{T}_{\lambda, \sigma}(0, s; \omega) P_{\c} \delta G (\theta_s\omega, v_\lambda[\m](s, \omega))\,\mathrm{d} s, \quad \Forall \; q \in \mathcal{H}_\alpha^{\s},
\eea
{\HL where we used the fact that the solution operator $\mathfrak{T}_{\lambda, \sigma}$ leaves invariant the subspaces $\mathcal{H}^{\c}$ and $\mathcal{H}^{\s}$ as pointed out in \SS~\ref{ss:T}.}  

It follows from \eqref{functional eqn} and \eqref{eq:op_I} that {\hh the original fixed point problem \eqref{FPC-recall}, {\it i.e.} $\mathcal{L}^{\omega, \lambda}_{\m}[v] = v$ with the constraint $v(0, \omega) + u_0(\omega) \in \mathfrak{M}_\lambda(\omega)$}, can be transformed into the following unconstrained fixed point problem for $p$:
\bea  \label{problem p}
p = \mathcal{J}_{\omega, \lambda}(p),
\eea
where the operator $\mathcal{J}_{\omega, \lambda}$ is defined by:
\bea \label{eq:op_J}
\mathcal{J}_{\omega, \lambda}(\p) &:=  P_{\c} v_\lambda[q(p)](0, \omega)  + P_{\c} u_0(\omega) \\
& = - \int_{0}^{+\infty} \mathfrak{T}_{\lambda, \sigma}(0, s; \omega) P_{\c} \delta G (\theta_s\omega, v_\lambda[q(p)](s, \omega))\,\mathrm{d} s + P_{\c}u_0(\omega), \qquad \p \in \mathcal{H}^{\c},
\eea
{\HL with $v_\lambda[q(p)](\cdot, \omega)$ being the fixed point of $\mathcal{L}^{\omega, \lambda}_{q(p)}$ and $q(p) = h_\lambda(\p, \omega) - P_{\s}u_0(\omega)$.}

It is then sufficient to solve \eqref{problem p} in $\mathcal{H}^{\c}$ in order to solve the constrained fixed point problem \eqref{FPC-recall} and therefore the asymptotic completeness problem; the temperedness of $\overline{u}_0$ being checked in Step ~5.


We show in the following that $\mathcal{J}_{\omega, \lambda}$ is a contraction mapping on $\mathcal{H}^{\c}$, leading then to a unique fixed point {\HL $\p_{\omega, \lambda}(u_0)$} of $\mathcal{J}_{\omega, \lambda}$ {\HL for each given $u_0$}.

For any given $p_1$ and $p_2 \in \mathcal{H}^{\c}$, we get
%
\bea  \label{J-est-1}
& \|\mathcal{J}_{\omega, \lambda}( \p_1) - \mathcal{J}_{\omega, \lambda}(\p_2)\|_\alpha \\
& = \| P_{\c} v_\lambda[q(p_1)](0, \omega) - P_{\c} v_\lambda[q(p_2)](0, \omega)\|_\alpha \\
&  = {\HLL \Bigl \|\int_{0}^{+\infty} \mathfrak{T}_{\lambda, \sigma}(0, s; \omega) P_{\c} \Bigl ( \delta G \bigl( \theta_s\omega, v_\lambda[q(p_1)](s, \omega) \bigr) - \delta G \bigl (\theta_s\omega, v_\lambda[q(p_2)](s, \omega) \bigr ) \Bigr) \,\mathrm{d} s \Bigr \|_\alpha} \\
& \le K\LF\int_0^{\infty} e^{-{\eta_1} s - \int_0^s z_\diffusion(\theta_\tau \omega)
\,\mathrm{d} \tau} \|v_\lambda[q(p_1)](s, \omega) - v_\lambda[q(p_2)](s,\omega) \|_\alpha \,\mathrm{d} s,
\eea
where we used the partial-dichotomy estimate \eqref{d-III} and the fact that $G$ and $F$ have the same Lipschitz constant $\LF$
to obtain the last inequality above.

Recalling that $q(p_i) = h_\lambda(\p_i, \omega) - P_{\s}u_0(\omega)$, we obtain the following estimates by applying \eqref{Lip v} to the RHS of the last inequality above:
\bea \label{J-est-2}
& K\LF\int_0^{\infty} e^{-{\eta_1} s - \int_0^s z_\diffusion(\theta_\tau \omega)
\,\mathrm{d} \tau} \|v_\lambda[q(p_1)](s, \omega) - v_\lambda[q(p_2)](s, \omega) \|_\alpha \,\mathrm{d} s \\
& \le \frac{K^2\LF ({\eta_1} - \eta)^{-1}}{1- \Upsilon_{1}(F)}\|h_\lambda(\p_1, \omega) - h_\lambda(\p_2, \omega)\|_\alpha \\
& \le \frac{K^2\LF \LH({\eta_1} - \eta)^{- 1}}{1- \Upsilon_{1}(F)}\|\p_1 - \p_2\|_\alpha,
\eea
{\HL which --- by noting that $K \LF ({\eta_1} - \eta)^{-1} < 1$ from \eqref{thm 3.1 condition}  --- leads to}
\beas
\|\mathcal{J}_{\omega, \lambda}( \p_1) - \mathcal{J}_{\omega, \lambda}(\p_2)\|_\alpha \le \frac{K \LH}{1- \Upsilon_{1}(F)}\|\p_1 - \p_2\|_\alpha, \qquad  \Forall \, p_1, \, p_2 \in \mathcal{H}^{\c}.
\eeas
Hence, $\mathcal{J}_{\omega, \lambda}$ is a contraction mapping on $\mathcal{H}^{\c}$ by assumption \eqref{thm 3.5 condition 2}. 

{\HL From what precedes, given $u_0$ and the corresponding fixed point $p_{\omega, \lambda}(u_0)$ of the operator $\mathcal{J}_{\omega, \lambda}$, we get therefore the existence of $\overline{u}_0(\omega) = v_\lambda \bigl[ q(\p_{\omega, \lambda}(u_0)) \bigr](0, \omega) + u_0(\omega)$ that belongs to $\mathfrak{M}_\lambda(\omega)$ so that  $\overline{u}_\lambda = v_\lambda[q(\p_{\omega, \lambda}(u_0))] + u_\lambda$ is a mild solution --- hence also a classical solution according to Proposition ~\ref{prop_mild} --- of Eq. ~\eqref{REE} on the manifold, and $\overline{u}_\lambda(\cdot, \omega) - u_\lambda(\cdot, \omega) \in C_\eta^+$ for each $\omega$. We are just left with the verification that $\overline{u}_0$ is tempered.}

\medskip
{\bf Step 5. Temperedness of $\overline{u}_0$ and the asymptotic completeness.} 
We saw in the previous step that given an $\mathcal{H}_\alpha$-valued tempered random initial datum $u_0$, the sought solution $\overline{u}_{\lambda}$, which solves the problem \eqref{FPC-recall}, emanates from $\overline{u}_0$ given by:
\be\label{Eq_good_ubar_0}
\overline{u}_0(\omega)= \p_{\omega, \lambda}(u_0) + h_\lambda(\p_{\omega, \lambda}(u_0), \omega),
\ee
where $p_{\omega, \lambda}(u_0)$ is the corresponding fixed point of the operator $\mathcal{J}_{\omega, \lambda}$ given in \eqref{eq:op_J}.
In what follows we check that such a $\overline{u}_0$ is an  $\mathcal{H}_{\alpha}$-valued tempered random variable.

First, we check that $\overline{u}_0$ is a random variable. Since
$\overline{u}_0(\omega) = \p_{\omega, \lambda}(u_0) + h_\lambda(\p_{\omega, \lambda}(u_0), \omega)$ for each $\omega$ and $h_\lambda$ is measurable, we only need to check that $\omega \mapsto \p_{\omega, \lambda}(u_0)$ is $(\mathcal{F}; \mathcal{B}(\mathcal{H}^{\c}))$-measurable. In that respect, note that the function $v_\lambda[\m]$ involved in the definition of the operator $\mathcal{J}_{\omega, \lambda}$ is known to be $(\mathcal{B}([0, \infty)) \otimes \mathcal{F}; \mathcal{B}(\mathcal{H}_\alpha))$-measurable according to Step ~3. 
The desired measurability of {\HLL $\omega \mapsto \p_{\omega, \lambda}(u_0)$} can be thus obtained by relying on the Picard scheme mentioned before associated with here the fixed point problem $\mathcal{J}_{\omega, \lambda}(p) = p$.

We  check  now that $\overline{u}_0$ is tempered. This is simply obtained from an appropriate control of  $\|\overline{u}_0(\omega)\|_\alpha$ by $\| u_0(\omega)\|_\alpha$ which is 
made possible from the construction of $\overline{u}_0$ given in \eqref{Eq_good_ubar_0}.

{\HLL In that respect, first note that by the construction of $\mathcal{J}_{\omega, \lambda}$ given in  \eqref{eq:op_J}, it holds that
\beas 
 \|\p_{\omega, \lambda}(u_0)\|_\alpha = \|\mathcal{J}_{\omega, \lambda}(\p_{\omega, \lambda}(u_0))\|_\alpha  & \le  \|P_{\c} v_\lambda\bigl[q(\p_{\omega, \lambda}(u_0))\bigr](0, \omega)\|_\alpha + \|P_{\c} u_0(\omega)\|_\alpha \\
& \le \|P_{\c} v_\lambda\bigl[ q(\p_{\omega, \lambda}(u_0)) \bigr](\cdot, \omega)\|_{C_\eta^+} + \|P_{\c} u_0(\omega)\|_\alpha.
\eeas
By applying \eqref{v bound} we obtain then:}
\bea
\|\p_{\omega, \lambda}(u_0)\|_\alpha & \le   \frac{K}{1- \Upsilon_{1}(F)}\| q(\p_{\omega, \lambda}(u_0)) \|_\alpha   + \|P_{\c} u_0(\omega)\|_\alpha \\
& = \frac{K}{1- \Upsilon_{1}(F)}\|h_\lambda(\p_{\omega, \lambda}(u_0), \omega) - P_{\s}u_0(\omega)\|_\alpha   + \|P_{\c} u_0(\omega)\|_\alpha \\
& \le \frac{K}{1- \Upsilon_{1}(F)}\big( \LH\|\p_{\omega, \lambda}(u_0)\|_\alpha + \|P_{\s}u_0(\omega)\|_\alpha \big) +  \|P_{\c} u_0(\omega)\|_\alpha, \; \Forall \, \omega \in \Omega.
\eea
Using assumption \eqref{thm 3.5 condition 2} in the last inequality, we obtain after simplification that
\bea \label{p est}
\|\p_{\omega, \lambda}(u_0)\|_\alpha & \le \frac{1 - \Upsilon_{1}(F)}{1 - \Upsilon_{1}(F) - K \LH} \Bigl(\frac{K}{1 - \Upsilon_{1}(F)}\|P_{\s} u_0(\omega)\|_\alpha + \|P_{\c} u_0(\omega)\|_\alpha \Bigr).
\eea
By noting that $\|P_{\c} u_0(\omega)\|_\alpha < \frac{K}{1 - \Upsilon_{1}(F)} \|P_{\c} u_0(\omega)\|_\alpha$ (recalling that $K\ge 1$ and $0<\Upsilon_{1}(F)<1$), we deduce:
\beas
\|\p_{\omega, \lambda}(u_0)\|_\alpha & \le \frac{2K}{1 - \Upsilon_{1}(F) - K \LH} \|u_0(\omega)\|_\alpha,  \quad \Forall \, \omega \in \Omega.
\eeas
{\mkk Hence, from the definition of $\overline{u}_0$ in \eqref{Eq_good_ubar_0}:}
\bea \label{temperedness}
\|\overline{u}_0(\omega)\|_\alpha & \le \|\p_{\omega, \lambda}(u_0)\|_\alpha + \|h_\lambda(\p_{\omega, \lambda}(u_0), \omega)\|_\alpha \\
& \le (1+ \LH)\|\p_{\omega, \lambda}(u_0)\|_\alpha \\
& \le \frac{2 K (1+ \LH)}{1 - \Upsilon_{1}(F) - K \LH} \|u_0(\omega)\|_\alpha, \quad \Forall \, \omega \in \Omega.
\eea
The temperedness property of $\overline{u}_0$ follows then from the assumption that the initial datum $u_0$ is  tempered. 

Finally, we conclude about the asymptotic completeness property. Let $\overline{u}_{\lambda}$ be the solution to Eq. ~\eqref{REE} which solves the problem \eqref{FPC-recall}, emanating from the initial datum as given in \eqref{Eq_good_ubar_0}.

{\HLL
From the definition of $C_\eta^{+}$-norm, we obtain for each $t\ge 0$ that
\beas
\|\overline{u}_\lambda(t,\omega)  - u_\lambda(t, \omega)\|_\alpha &= \|v_\lambda[q(\p_{\omega, \lambda}(u_0))](t, \omega)\|_\alpha \le e^{\eta t + \int_0^t z_\diffusion(\theta_s \omega)\,\mathrm{d} s}\|v_\lambda[q(\p_{\omega, \lambda}(u_0))](\cdot, \omega)\|_{C_\eta^+}.
\eeas
{\mkk Using \eqref{v bound} and $q(\p_{\omega, \lambda}(u_0)) = h_\lambda(\p_{\omega, \lambda}(u_0), \omega) - P_{\s}u_0(\omega)$ in the above inequality, we obtain}
\be 
\|\overline{u}_\lambda(t,\omega)  - u_\lambda(t, \omega)\|_\alpha  \le \frac{Ke^{\eta t + \int_0^t z_\diffusion(\theta_s \omega)\,\mathrm{d} s}}{1 - \Upsilon_{1}(F)}\|h_\lambda(\p_{\omega, \lambda}(u_0), \omega) - P_{\s}u_0(\omega)\|_\alpha,
\ee
which leads to
\be\label{attractive 1}
\|\overline{u}_\lambda(t,\omega)  - u_\lambda(t, \omega)\|_\alpha \leq  \frac{K e^{\eta t + \int_0^t z_\diffusion(\theta_s \omega)\,\mathrm{d} s}}{1 - \Upsilon_{1}(F)}  \|P_{\s}\overline{u}_0(\omega) - P_{\s}u_0(\omega)\|_\alpha, 
\ee
by simply applying $P_{\s}$ to \eqref{Eq_good_ubar_0}.}

For any $\epsilon \in (0, |\eta|)$, let us now introduce the following positive random variables
\bea \label{att coef}
& {\mkk C_{\epsilon,\sigma}}(\omega) := \sup_{t\ge 0}  \frac{K e^{-\epsilon t + \int_0^t z_\diffusion(\theta_s \omega)\,\mathrm{d} s}}{1 - \Upsilon_{1}(F)}, \quad  C'_{\epsilon, \sigma}(\omega) := \sup_{t\ge 0}  \frac{K e^{-\epsilon t + \int_{-t}^0 z_\diffusion(\theta_{s} \omega)\,\mathrm{d} s}}{1 - \Upsilon_{1}(F)},
\eea
where $\Upsilon_{1}(F)$ is defined in \eqref{thm 3.1 condition}.

Note that $C_{\epsilon, \sigma}(\omega)$ and $C'_{\epsilon,\sigma}(\omega)$ are finite for all $\omega$ because of the growth control relation of $z_\diffusion(\theta_t \omega)$ given in \eqref{Eq:z-4}. We obtain then from \eqref{attractive 1} that 
\bea \label{attractive 2}
 & \|\overline{u}_\lambda(t,\omega) - u_\lambda(t, \omega)\|_\alpha \le C_{\epsilon,\sigma}(\omega)  e^{(\eta + \epsilon) t}  \|P_{\s}\overline{u}_0(\omega) - P_{\s}u_0(\omega)\|_\alpha, \; \Forall t \ge 0,\, \omega \in \Omega.
\eea
By a simple change of fiber, we also obtain from \eqref{attractive 1} that:
\bea  \label{attractive 2-pb}
& \|\overline{u}_\lambda(t,\theta_{-t}\omega) - u_\lambda(t, \theta_{-t}\omega)\|_\alpha  \\
& \hspace{3em}  \le  C'_{\epsilon,\sigma}(\omega)  e^{(\eta + \epsilon) t}  \|P_{\s}\overline{u}_0(\theta_{-t}\omega) - P_{\s}u_0(\theta_{-t}\omega)\|_\alpha, \; \Forall t \ge 0,\, \omega \in \Omega.
\eea
Recalling that $\eta<0$, the forward and pullback asymptotic completeness of $\mathfrak{M}_\lambda$ follow respectively from \eqref{attractive 2} and \eqref{attractive 2-pb} with $\overline{u}_0(\omega)$ defined in \eqref{Eq_good_ubar_0}. In both \eqref{attractive 2} and \eqref{attractive 2-pb}, the rate of attraction is given by $|\eta| - \epsilon$ for any $\epsilon>0$ so that the critical attraction rate is $|\eta|$. The proof is now complete.

\ep

\needspace{1\baselineskip}
\br

\hspace*{2em}
\bi

\item[1)] The proof of Corollary ~\ref{attractiveness thm SEE} follows from the cohomology relation \eqref{cohomology relation} and the estimate \eqref{attractive 1}. Indeed, the same estimates as given in \eqref{attractive 2} and \eqref{attractive 2-pb} for the SPDE case can be derived from \eqref{attractive 1} using the cohomology relation, but with the following random coefficients $\widehat{C}_{\epsilon,\sigma}(\omega)$ and $\widehat{C}'_{\epsilon,\sigma}(\omega)$: 
\bea\label{eq:widehat_C}
\widehat{C}_{\epsilon,\sigma}(\omega) := \sup_{t\ge 0}  \frac{K e^{-\epsilon t  + \diffusion W_t(\omega)}}{1 - \Upsilon_{1}(F)}, \quad \widehat{C}'_{\epsilon,\sigma}(\omega) := \sup_{t\ge 0}  \frac{K e^{-\epsilon t - \diffusion W_{-t}(\omega)}}{1 - \Upsilon_{1}(F)}.
\eea

\item[2)] Note that the random positive variable $\widehat{C}_{\epsilon,\sigma}$ exhibits fluctuations which get larger as $\sigma$ gets bigger. In particular, for any $\gamma>0$, it can be shown that the time $t^\ast_{\sigma}(\omega)$, after which 
$\|\overline{u}_\lambda(t,\omega) - u_\lambda(t, \omega)\|_{\alpha}$ gets smaller than $\gamma$, is such that its expected value $m_{\sigma}$ increases with $\sigma$. However, when $\sigma$ tends to zero, $m_{\sigma}$ converges to the corresponding attraction time associated with $\sigma=0$. Hence, although the critical rate of attraction is independent of $\sigma$, the expected attraction time to the inertial manifold (for a given precision $\gamma$) is not.  


\item[3)] {\HL Note that from the proof of Theorem ~\ref{attractiveness thm}, we see that in the case where $\eta>0$} the difference between the given solution $u_\lambda$ of Eq. ~\eqref{REE} and the constructed solution $\overline{u}_\lambda$ living on the random invariant manifold $\mathfrak{M}_\lambda$ has a growth controlled by $C_{\epsilon, \sigma}(\omega)e^{\eta + \epsilon t}$ in the forward case (resp. ~$C'_{\epsilon, \sigma}(\omega)e^{\eta + \epsilon t}$ in the pullback case) for any $\epsilon > 0$.

\ei

\er

\section{Local Stochastic Invariant Manifolds: Preparation to Critical Manifolds} \label{s:local}
 In this section we present a local theory of stochastic invariant manifolds associated with the global theory described in Section \ref{s:global}. The ideas 
are standard but the precise framework is detailed here again in view of the main results regarding the approximation formulas of stochastic critical manifolds (Section ~\ref{s:approximation}) and the related pullback characterizations (Section \ref{Sec_app_man+PB}).  
In particular,  the proof of Theorem ~\ref{thm:local REE} is provided below where some elements will be used in establishing Theorem \ref{App CMF}.

We introduce now the following notion of a local random (resp. ~stochastic) invariant manifold on which we will rely to build our theory of stochastic critical manifolds in Section ~\ref{s:approximation}.

\bd  \label{Def. local inv. manfld}
Let $S$ be a continuous RDS acting on some separable Hilbert space $H$, {\hl and $r$ be a given positive integer.} A random closed set $\mathfrak{M} ^{\mathrm{loc}}$ is called a local random invariant $C^r$ manifold of $S$ if the following conditions hold:
\bi
\item[(i)] $\mathfrak{M} ^{\mathrm{loc}}$ is locally invariant in the sense that for each $\omega$ and each $u_0 \in \mathring{\mathfrak{M}}^{\mathrm{loc}}(\omega)$ ($\neq \emptyset$),\footnote{Here, $\mathring{\mathfrak{M}}^{\mathrm{loc}}(\omega)$ denotes the interior of $\mathfrak{M}^{\mathrm{loc}}(\omega)$.} there exists $t_{u_0,\omega}>0$, such that $S(t, \omega)u_0 \in \mathring{\mathfrak{M}}^{\mathrm{loc}}(\theta_t\omega)$ for all $t\in [0, t_{u_0, \omega})$.

\item[(ii)] There exist a closed subspace $\mathcal{G} \subset H$ and a measurable function $h\colon \mathcal{G}\times \Omega \rightarrow \mathcal{L}$, with $\mathcal{L}$ the topological complement of $\mathcal{G}$ in $H$, such that $\mathfrak{M}^{\mathrm{loc}}(\omega)$ can be represented as the graph of $h(\cdot,\omega)\vert_{\mathfrak{B}(\omega)}$ for each $\omega$, {\it i.e.},
\beas
\mathfrak{M}^{\mathrm{loc}}(\omega) = \{\xii + h(\xii, \omega) \mid \xii\in \mathfrak{B}(\omega)\},
\eeas
where  $\mathfrak{B} \colon \Omega \rightarrow 2^{\mathcal{G}}$ is some random closed ball, such that $\mathfrak{B}(\omega)$
 is centered at zero for each $\omega$.
\item[(iii)] $h(\cdot, \omega)\vert_{\mathfrak{B}(\omega)}$ is $C^r$ for all $\omega$.

\item[(iv)] $\mathfrak{M}^{\mathrm{loc}}(\omega)$ is tangent to the subspace $\mathcal{G}$ at the origin:
\beas
h(0, \omega) = 0, \quad D_{\xii} h(0, \omega) = 0, \, \quad \Forall \, \omega \in \Omega.
\eeas 

\ei
The function $h$ is called the local random invariant manifold function associated with $\mathfrak{M}^{\mathrm{loc}}$.
\ed 

As in the global version, we will refer the associated local manifold as a local stochastic invariant $C^r$ manifold when we will deal with Eq. ~\eqref{SEE}. 

\br \label{rmk def local cm}
\hspace*{2em}
\bi
{\mkk
\item[1)] The random closed ball $\mathfrak{B}$ in part (ii) of the above definition can  be deterministic; see {\it e.g.} ~Corollary ~\ref{thm:local SEE} for such a situation.
\item[2)] Note that this notion of local random invariant manifold is a natural generalization of the classical one already encountered for  non-autonomous dynamical systems; see {\it e.g.} ~\cite[Def. 6.1.1]{Hen81}. Here, ``local''  has to be understood in terms of both the time variable $t$ {\hh and} the phase-space variable $\xii$ as specified in (i) and (ii) in the definition above. As we will see in Step ~3 of the proof of Theorem ~\ref{thm:local REE}, this ``local'' property arises as a direct consequence of the cut-off procedure.
\item[3)] Note also that Definition ~\ref{Def. local inv. manfld} does not exclude that 
$$S(t, \omega)u_0 \in \mathring{\mathfrak{M}}^{\mathrm{loc}}(\theta_t\omega),$$
for $t \in [t^1_{u_0, \omega}, t^2_{u_0, \omega}]$ with $t^2_{u_0, \omega} > t^1_{u_0, \omega}> t_{u_0, \omega}$. Actually, the noise effects can make this to happen on {\mkk infinitely-many} such time intervals; see \cite{CGLW} for more details.}
\ei
\er

Let us introduce now a cut-off version of the nonlinearity $F$. Let $\zeta \colon
\mathbb{R}^+ \rightarrow [0,1]$ be a $C^\infty$ decreasing function, such that
\begin{equation} \label{cut-off}
\zeta(s) = \begin{cases} 1 \quad \text{if } s \le 1, \\
0 \quad \text{if } s \ge 2.
\end{cases}
\end{equation}
For any given positive constant $\r$, let us define a mapping $F_{\r} \colon \mathcal{H}_\alpha \rightarrow  \mathcal{H}$ {\it via}
\be \label{F local version 1}
F_{\r}(u):=\zeta\Big(\frac{\|u\|_\alpha}{\r}\Big)F(u).
\ee


{\mkk We are now in position to formulate the results regarding  the existence and smoothness of families of local random invariant manifolds for Eq.~\eqref{REE}, when the latter generates a continuous (global) RDS acting on  $\mathcal{H}_\alpha$.\footnote{For applications we have in mind, it is often met that for any given initial datum $u_0\in \mathcal{H}_\alpha$ there exists a unique classical solution to Eq. ~\eqref{REE} for all $t \ge 0$ in the sense given in Proposition ~\ref{prop:exist}, and hence Eq. ~\eqref{REE} generates an RDS in the sense given in \SS~\ref{subsect:RDS}. Such property holds  for a broad class of random dissipative evolution equations; see {\it e.g.} ~\cite{CF94, Gess13, Kloeden_al07}.}
}

\bt \label{thm:local REE}

Consider the RPDE \eqref{REE}. The linear operator $L_\lambda$ is assumed to satisfy the corresponding assumptions in \SS ~\ref{ss:SEE}. The nonlinearity $F\colon \mathcal{H}_\alpha \rightarrow \mathcal{H}$ is assumed to be $C^p$ smooth for some integer $p\geq 2$, and to satisfy \eqref{F}--\eqref{DF}.\footnote{but not necessarily required to be globally Lipschitz here.} Assume also that Eq. ~\eqref{REE} generates a continuous RDS acting on the space $\mathcal{H}_\alpha$.

{\mkk We assume furthermore that {\hh an open interval $\Lambda$ is chosen} such that the uniform spectrum decomposition \eqref{gap 3} holds over $\Lambda$, and such that there exist $\eta_1$ and $ \eta_2$ satisfying 
\bea  \label{eta1-2b}
\eta_{\c} > \eta_1 > \eta_2 > \eta_{\s},
\eea 
for which
\be \label{thm 4.1 condition}
\exists \; \eta \in (\eta_2, \eta_1) \; \; \mathrm{s.t.} \quad \eta_2 < j\eta < \eta_1, \quad \Forall j \in \{1, \cdots, r\},
\ee
for some integer $r \in \{1,...,p\}$. As above,  $\eta_{\c}$ and $\eta_{\s}$ are defined in \eqref{eta_cs}, and  $\mathcal{H}^{\c}$ and  $\mathcal{H}_\alpha^{\s}$ denote  the subspaces associated with the uniform spectrum decomposition, defined in \eqref{Hc} and \eqref{L decomposition} with $\dim(\mathcal{H}^{\c})=m$.}

{\hh {\mkk Under these conditions}, there exists $\r^\ast > 0$ such that for each $\r \in (0, \r^\ast)$ the following uniform spectral gap conditions hold:
\be \label{Upsilon_1 local}
\boxed{\Upsilon_1(F_\r) = K \text{Lip}(F_{\r}) \bigl(({\eta_1} -  \eta)^{- 1} + \Gamma(1-\alpha) (\eta - {\eta_2})^{\alpha-1}  \bigr)  < \frac{1}{2}, }
\ee
and
\be  \label{Upsilon_j local}
\boxed{\Upsilon_j(F_\r) = K \text{Lip}(F_{\r}) \bigl(({\eta_1} - j \eta)^{- 1} + \Gamma(1-\alpha) (j \eta - {\eta_2})^{\alpha-1}  \bigr) < 1, \quad \Forall j \in \{2, \cdots, r\}, }
\ee
where  $F_{\r}$ is defined by \eqref{F local version 1} and $K$ is as given in the partial-dichotomy estimates \eqref{Proj bounds}.}

Moreover, for each $\r \in (0, \r^\ast)$, Eq. ~\eqref{REE} possesses a family of local  random invariant $C^r$ manifolds $\{\mathfrak{M}^{\mathrm{loc}}_\lambda\}_{\lambda \in \Lambda}$, where each of such manifolds is  m-dimensional and  is given by
\bea  \label{eq:def M local}
\mathfrak{M}^{\mathrm{loc}}_\lambda(\omega) := \{\xii + h_\lambda(\xii, \omega) \mid \xii \in \mathfrak{B}_{\r}(\omega) \}, \quad   \omega \in \Omega, \; \lambda\in \Lambda,
\eea
with $h_\lambda \colon \mathcal{H}^{\c}\times \Omega \rightarrow \mathcal{H}_\alpha^{\s}$ denoting the corresponding random invariant manifold function. The associated random closed ball  $\mathfrak{B}_{\r}$ is
given by:
\be  \label{eq:set B}
\mathfrak{B}_{\r}(\omega) := \overline{B_\alpha(0, \delta_{\r}(\omega))} \cap \mathcal{H}^{\c}, \; \omega\in \Omega, \footnote{Throughout this article, for any given positive random variable $r(\omega)$, $B_\alpha(0, r)$ denotes the random open ball in $\mathcal{H}_\alpha$ centered at the origin with random radius $r$; {\HL see {\it e.g.} \cite[Def. ~2.1]{Crauel02}}.}
\ee
with 
\be  \label{delta in local cm}
\delta_{\r}(\omega) :=  \frac{e^{-z_\diffusion(\omega)}\r}{1+ K},  \; \omega\in \Omega,
\ee
where $z_\diffusion(\omega)$ is the OU process defined in \eqref{z}.

 \et

Thanks to the cohomology relation between the two RDSs associated with Eqns. ~\eqref{SEE} and \eqref{REE} as recalled in \SS~\ref{sec. conjugacy}, we deduce the following corollary.

\bc \label{thm:local SEE}
Consider the SPDE \eqref{SEE}. Assume that the same assumptions as in Theorem ~\ref{thm:local REE} hold with $\Lambda$ specified therein. Let also $\r^\ast$ be the positive constant as in Theorem ~\ref{thm:local REE}. 

Then for each $\r \in (0, \r^\ast)$, Eq. ~\eqref{SEE} possesses a family of local  stochastic invariant $C^r$ manifolds $\{\widehat{\mathfrak{M}}^{\mathrm{loc}}_\lambda\}_{\lambda \in \Lambda}$, where each of such manifolds is  $m$-dimensional and  is given by
\bea \label{SEE cm local version}
\widehat{\mathfrak{M}}^{\mathrm{loc}}_\lambda(\omega) := \{\xii + \widehat{h}_\lambda(\xii, \omega) \mid \xii \in\widehat{\mathfrak{B}}_{\r} \}, \quad  \omega \in \Omega, \; \lambda\in \Lambda,
\eea
with $\widehat{h}_\lambda \colon \mathcal{H}^{\c}\times \Omega \rightarrow \mathcal{H}_\alpha^{\s}$ denoting the corresponding stochastic invariant manifold function for each $\lambda\in \Lambda$.

In this case,  the associated ``random closed ball" is simply a deterministic ball, $\widehat{\mathfrak{B}}_{\r}$, which is given by 
\bea  \label{eq:set B hat}
\widehat{\mathfrak{B}}_{\r}:=\overline{{B_\alpha(0, \widehat{\delta}_{\r})}} \cap \mathcal{H}^{\c}, 
\eea
with 
\be  \label{delta in SEE local cm}
\widehat{\delta}_{\r}:=  \frac{\r}{(1+ K)}.
\ee
Moreover, the following relation holds:
\bea \label{h and h hat}
\widehat{h}_\lambda (\xii, \omega)= e^{z_\diffusion (\omega)}h_{\lambda}(e^{-z_\diffusion (\omega)}\xii, \omega),   \quad  \xi \in \mathcal{H}^{\c}, \,  \omega \in \Omega,
\eea 
where $h_\lambda(\xii, \omega)$ is specified in Theorem ~\ref{thm:local REE}.

\ec

\needspace{1\baselineskip}

\br \label{rmk:radius}

\hspace*{2em}

\bi
 
\item[1)] Note that from the proof of Theorem ~\ref{thm:local REE} given below ({\HL see {\it e.g.} \eqref{eq:xi range1}}), for each $\r\in(0, \r^\ast)$, the local random invariant manifold $\mathfrak{M}^{\mathrm{loc}}_\lambda$ can actually be obtained as the random graph over a larger random closed ball with random radius given by 
\bea
\widetilde{\delta}_{\r}(\omega):= \frac{e^{-z_\sigma(\omega)} \r}{1+ \LHWS},
\eea
{\HL since} indeed $\delta_{\r}(\omega) < \widetilde{\delta}_{\r}(\omega)$ for all $\omega$ {\HL according to \eqref{Lip h local cm 2}.} 

Similarly, the local stochastic invariant manifold $\widehat{\mathfrak{M}}^{\mathrm{loc}}_\lambda$ given in Corollary ~\ref{thm:local SEE} can also be defined over a larger deterministic ball with radius given by $\frac{\r}{1+ \LHWS}$. {\mkk In all the cases, it is important to note that such deterministic  radii  are  
artificial and result from the techniques employed in the proofs of Theorem ~\ref{thm:local REE} and Corollary ~\ref{thm:local SEE}, such as {\it e.g.} ~the cohomology approach; see also Remark ~\ref{Rmk_artificial_bounds2}. }

\item[2)] The radius $\delta_\r$ provided in Theorem ~\ref{thm:local REE} results from a control of $\LHW$ by $K$ made possible by {\HLL the choice of $\r^\ast$}; {\HL see Step ~3 below}.  The interest is that the resulting manifolds $\mathfrak{M}^{\mathrm{loc}}_\lambda$ live then above some $\lambda$-independent neighborhoods, which help simplify certain algebraic manipulation in the proof of Theorem ~\ref{App CMF} that relies on Theorem ~\ref{thm:local REE}.

\ei

\er

\bp[{\bf Proof of Theorem ~\ref{thm:local REE}}]

As mentioned before, the proof is accomplished via a cut-off procedure. We proceed in three steps.

\medskip
{\bf Step 1. Modification of the nonlinearity $F$ through a cut-off function.} We show in this first step that the modified nonlinearity $F_{\r}$ defined in \eqref{F local version 1} is globally Lipschitz while still being $C^p$ smooth.


We first examine the smoothness of $F_\r$. Note that the $\|\cdot\|_\alpha$-norm is $C^\infty$ smooth at any nonzero element in $\mathcal{H}_\alpha$ as norm induced by the inner product $\langle \cdot, \cdot \rangle_\alpha$ given in Section~\ref{s:preliminary}. Then, a {\HL basic} composition argument leads to $C^p$ smoothness of $F_\r$ at any nonzero element. It is also clear that $F_\r$ has this smoothness at the origin by the construction of $\zeta$; {\HL and $F_\r$ is thus} $C^p$ smooth on $\mathcal{H}_\alpha$.

Now, we show that $F_\r$ is globally Lipschitz when $\r$ is sufficiently small.  Note that $DF$ is continuous since $F$ is $C^p$ smooth with $p\ge 2$ by our assumption. Note also that $DF(0) = 0$ by \eqref{DF}. Then, there exists $\r_1 > 0$ and a positive function $C(\r)$ defined on $[0, \r_1]$ such that 
\bea \label{Cr est}
C(\r)\rightarrow 0  \quad \text{ as }  \quad {\HLL \r \rightarrow 0}, 
\eea
and
\bea \label{Eq:op bdd-1}
\|DF(u)\|_{\alpha,0} \le C(\r), \quad \Forall \,  \r \in (0, \r_1], \, u \in B_\alpha(0, \r) \subset \mathcal{H}_\alpha,
\eea
where $\|DF(u)\|_{\alpha,0}$ stands for the operator norm of $DF(u)$ as a  linear operator from $\mathcal{H}_\alpha$ to $\mathcal{H}$.

For any $u_1, u_2 \in \mathcal{H}_\alpha$, by the mean value theorem, there exists $s \in [0, 1]$ such that
\bea \label{Eq:MVT-1}
F(u_2) = F(u_1) + DF(s u_1 + (1-s)u_2)(u_2 - u_1). 
\eea
Now, it follows from \eqref{Eq:op bdd-1} and \eqref{Eq:MVT-1} that
\be \label{eq:F est}
\|F(u_1) - F(u_2)\| \le C(\r) \|u_1 - u_2 \|_\alpha, \quad \Forall \, \r \in (0, \r_1], \, u_1, u_2 \in B_\alpha(0, \r).
\ee

Recall that $D^2F$ is also continuous. In particular, it is continuous at the origin. There exists then $\r_2>0$ and a bounded positive function $\widetilde{C}(\r)$ defined on $[0, \r_2]$ such that
\bea \label{Eq:bdd-2}
\|D^2F(u)\|_{\mathrm{op}} \le \widetilde{C}(\r), \quad \Forall \, \r \in (0, \r_2], \, u \in B_\alpha(0,\r),
\eea
where $\|D^2F(u)\|_{\mathrm{op}}$ is the operator norm of $D^2F(u)$ as a bilinear operator from $\mathcal{H}_\alpha \times \mathcal{H}_\alpha$ to $\mathcal{H}$.

Again, since $F$ is $C^p$ smooth with $p\ge 2$, $F(0) = 0$, and $DF(0) = 0$, it follows then from the Taylor-Lagrange theorem that
\bea \label{Eq:MVT-2}
F(u) = D^2F(s u)(u,u), \quad \, \Forall u \in \mathcal{H}_\alpha,
\eea
for some $s \in [0,1]$ which depends on $u$.

We obtain from \eqref{Eq:bdd-2} and \eqref{Eq:MVT-2} that 
\be \label{local bound F}
\|F(u)\| \le \widetilde{C}(\r) \|u\|_\alpha^2, \quad \Forall \, \r \in (0, \r_2], \,  u \in B_\alpha(0,\r).
\ee

Now, let 
\bea
\rho_3 := \min \left\{\frac{\r_1}{2}, \frac{\r_2}{2}\right\}.
\eea
We show that for each fixed $\r\in (0, \r_3)$ it holds that
\bea \label{Lip F rho}
 \|F_{\r}(u_1) - F_{\r}(u_2)\| & \le \big( C(2\r) + 4 \r \widetilde{C}(2 \r)  \text{Lip}(\zeta) \big) \|u_1-u_2\|_\alpha, \;  \Forall u_1, u_2 \in \mathcal{H}_\alpha,
\eea
where $\mathrm{Lip}(\zeta)$ denotes the global Lipschitz constant of the cut-off function $\zeta$: 
\bea \label{eq:Lip zeta}
|\zeta(s_1) - \zeta(s_2)| \le \mathrm{Lip}(\zeta) |s_1 - s_2|, \quad \Forall s_1, s_2 \in \mathbb{R}^{+}.
\eea

There are three cases to be considered. If $\|u_1\|_\alpha, \|u_2\|_\alpha \ge 2\r$, then \eqref{Lip F rho} holds trivially by the definition of $F_\r$ and the construction of $\zeta$.

If $\|u_1\|_\alpha < 2\r$ and $\|u_2\|_\alpha \ge 2\r$, we get
\bea
\|F_{\r}(u_1) - F_{\r}(u_2)\| & = \|F_{\r}(u_1) \| = \Big\| \zeta\Big(\frac{\|u_1\|_\alpha}{\r} \Big)F(u_1) \Big \| \le \widetilde{C}(2\r) \|u_1\|_\alpha^2 \Big \| \zeta\Big (\frac{\|u_1\|_\alpha}{\r}\Big) \Big \| \\
& \le 4 \r^2 \widetilde{C}(2\r) \Big \| \zeta\Big (\frac{\|u_1\|_\alpha}{\r}\Big) \Big \| =  4 \r^2 \widetilde{C}(2\r)  \Big \| \zeta\Big (\frac{\|u_1\|_\alpha}{\r}\Big )  - \zeta\Big(\frac{\|u_2\|_\alpha}{\r}\Big)\Big \| \\
& \le 4 \r \widetilde{C}(2\r) \text{Lip}(\zeta) (\|u_2\|_\alpha - \| u_1\|_\alpha) \\
& \le 4 \r \widetilde{C}(2\r) \text{Lip}(\zeta) \|u_1 -  u_2\|_\alpha,
\eea
where we used \eqref{local bound F} and $2\r < 2\r_3 \le \r_2$ to derive the first inequality above. So \eqref{Lip F rho} holds in this case.  By exchanging the role of $u_1$ and $u_2$, the result also holds for the case when {\HL $\|u_2\|_\alpha < 2\r$ and $\|u_1\|_\alpha \ge 2\r$}.

Finally, if $\|u_1\|_\alpha< 2 \r$ and $\|u_2\|_\alpha < 2 \r$, we get
\bea
& \|F_{\r}(u_1) - F_{\r}(u_2)\| = \Big\| \zeta\Big(\frac{\|u_1\|_\alpha}{\r} \Big)F(u_1) - \zeta\Big(\frac{\|u_2\|_\alpha}{\r}\Big)F(u_2)\Big \| \\
& \le \Big \| \zeta\Big (\frac{\|u_1\|_\alpha}{\r}\Big)F(u_1) -  \zeta\Big(\frac{\|u_1\|_\alpha}{\r}\Big)F(u_2)\Big \| + \Big \| \zeta\Big (\frac{\|u_1\|_\alpha}{\r}\Big )F(u_2) -  \zeta\Big(\frac{\|u_2\|_\alpha}{\r}\Big)F(u_2)\Big \| \\
& \le C(2\r)\|u_1 - u_2\|_\alpha + {\HL \text{Lip}(\zeta) \Bigl | \frac{\|u_1\|_\alpha}{\r} - \frac{\|u_2\|_\alpha}{\r} \Bigr|\cdot\|F(u_2)\|} \\
& \le C(2\r)\|u_1 - u_2\|_\alpha + 4 \r \widetilde{C}(2\r) \text{Lip}(\zeta) \|u_1 - u_2\|_\alpha,
\eea
where we used \eqref{eq:Lip zeta} and \eqref{eq:F est} to derive the second last inequality above, and used \eqref{local bound F} and $\|u_2\|_\alpha < 2\r$ to obtain the last inequality above. Again, \eqref{Lip F rho} follows. Thus, $F_{\r}$ is indeed globally Lipschitz for each fixed $\r \in (0, \r_3)$.
 
Now, let
\bea  \label{local Lip def}
\LFW :=  C(2\r) +  4 \r \widetilde{C}(2\r) \text{Lip}(\zeta), \quad \r \in (0, \r_3).
\eea
Then, it follows from \eqref{Lip F rho} that
\bea \label{local Lip}
 \|F_{\r}(u_1) - F_{\r}(u_2)\| & \le  \LFW \|u_1 - u_2\|_\alpha, \Forall u_1, u_2 \in \mathcal{H}_\alpha \text{ and } \r \in (0, \r_3).
\eea

\medskip
{\bf Step 2. Existence of global random invariant manifolds for the modified equation.} Now, for each fixed $\r>0$, let us consider the following cut-off version of Eq.~\eqref{REE}: 
\begin{equation} \label{REE local version}
\frac{\mathrm{d}u}{\mathrm{d}t} = L_\lambda u + z_\diffusion(\theta_t\omega)u + G_{\r}(\theta_t\omega,\, u),
\end{equation}
where
\bea \label{G modified}
& G_{\r}(\omega,\, u) := e^{-z_\diffusion(\omega)} F_{\r}(e^{z_\diffusion(\omega)}u) = e^{-z_\diffusion(\omega)} \zeta\Big( \frac{e^{z_\diffusion(\omega)} \|u\|_\alpha}{\r}\Big)F(e^{z_\diffusion(\omega)}u).
\eea

We first check condition \eqref{thm 3.2 condition} for Eq. ~\eqref{REE local version} when $\r$ is sufficiently small. Note that the Lipschitz constant of $F_\r$ given in \eqref{local Lip def} satisfies that
\bea \label{Lip F asym}
\LFW \rightarrow 0  \quad \text{ as } \quad \r \rightarrow 0,
\eea 
which follows from \eqref{Cr est} and the boundedness of $\widetilde{C}(\r)$ for sufficiently small $\r$. 

For $\eta$ given in \eqref{thm 4.1 condition}, we infer from \eqref{Lip F asym} that there exists a positive constant $\r_4$ such that for all $\r \in (0, \r_4)$, the following conditions are satisfied:
\bea \label{eq: F local condition}
\Upsilon_j(F_\r) = K \text{Lip}(F_{\r}) \bigl(({\eta_1} - j \eta)^{- 1} + \Gamma(1-\alpha) (j \eta - {\eta_2})^{\alpha-1}  \bigr) < 1, \quad \Forall j \in \{1, \cdots, r\}, 
\eea
where $r \in\{1, \cdots, p\}$ is as given in \eqref{thm 4.1 condition}. Condition \eqref{thm 3.2 condition} is thus verified with $F_\r$ in place of $F$.

It is clear that Eq. ~\eqref{REE local version} also satisfies other assumptions required in Theorems ~\ref{Lip manfd}--\ref{smth manfd}. Hence, for each $\r \in (0, \r_4)$, there exists a family of global random invariant $C^r$ manifolds $\{\mathfrak{M}_{\lambda,\r}\}_{\lambda \in \Lambda}$ for Eq. ~\eqref{REE local version}, which is given by:
\bea \label{M cutoff}
\mathfrak{M}_{\lambda,\r}(\omega) = \{\xii + h_{\lambda, \r}(\xii, \omega) \mid \xii \in \mathcal{H}^{\c} \}, \, \quad \Forall \, \omega \in \Omega, \, \lambda \in \Lambda, \, \r \in (0, \r_4), 
\eea
where the random invariant manifold function $h_{\lambda, \r} \colon \mathcal{H}^{\c} \times \Omega \rightarrow \mathcal{H}^{\s}_\alpha$ is $C^r$ smooth and globally Lipschitz with respect to $\xii$ for each $\omega$, and it is tangent to $\mathcal{H}^{\c}$ at the origin:
\be \label{tangency local cm}
h_{\lambda, \r}(0, \omega)  = 0, \quad D_{\xii} h_{\lambda, \r}(0, \omega) = 0, \, \quad \Forall \, \omega \in \Omega.
\ee
Moreover, by the Lipschitz estimate given in \eqref{Lip h}, we also have that
\be \label{Lip h local cm 1}
\LHW \le \frac{K^2\LFW (\eta - {\eta_2})^{\alpha -1}\Gamma(1-\alpha)}{1- \Upsilon_{1}(F_\r)},  \quad \r \in (0, \r_4),
\ee
where $\Upsilon_{1}(F_\r)$ is given in \eqref{Upsilon_1 local}.

\medskip
{\bf Step 3. Existence of local random invariant manifolds for Eq. ~\eqref{REE}.}  Now, we show that there exists $\r^\ast > 0$ sufficiently small such that for any given $\r \in (0, \r^\ast)$, we can find  a random closed ball centered at the origin in $\mathcal{H}^{\c}$ such that for each $\lambda \in \Lambda$ the random graph of $h_{\lambda, \r}$ over this random closed ball is a local random invariant $C^r$ manifold for Eq. ~\eqref{REE}. The construction of such manifolds relies on the obvious fact that a solution of the original equation \eqref{REE} coincides with a solution of the cut-off version Eq. ~\eqref{REE local version} over some time interval $[0, t]$ provided that they both emanate from the same sufficiently small initial datum.

We first remark that for each $\r> 0$ by comparing $G_\r$ given in \eqref{G modified} with $G$ given in \eqref{eq:G}, we have that 
\bea \label{eq:G equiv}
\Forall \, u \in \mathcal{H}_\alpha, \quad \Bigl( \|u\|_\alpha \le e^{-z_\sigma(\omega)} \r \Bigr)  \Longrightarrow \Bigl( G_\r(\omega, u) = G(\omega, u) \Bigr).
\eea
Note also that by the continuous dependence of the solutions of Eq. ~\eqref{REE local version} with respect to the initial data and the continuity of the OU process $t\mapsto z_\sigma(\theta_t\omega)$, we have:
\bea \label{eq: u est}
 \Forall \, u_0 \in B_\alpha(0,  e^{-z_\sigma(\omega)} \r),  & \ \exists \ t^{*}_{u_0, \omega} > 0, \quad \mathrm{ s.t.} \\
&  \|u_{\lambda, \r}(t, \omega; u_0)\|_\alpha < e^{-z_\sigma(\theta_t\omega)} \r, \quad \Forall \; t \in [0, t^{*}_{u_0, \omega}).
\eea

For each $\omega \in \Omega$ and $u_0 \in \mathcal{H}_\alpha$, let us denote the solution to Eq. ~\eqref{REE} emanating from $u_0$ (in fiber $\omega$) by $u_\lambda(t, \omega; u_0)$, which is guaranteed to exist  since Eq. ~\eqref{REE} is assumed to generate an RDS. {\HL From \eqref{eq:G equiv} and \eqref{eq: u est} we then obtain that}
\bea \label{eq:u and u local}
\Forall \, u_0 \in B_\alpha(0,  e^{-z_\sigma(\omega)} \r), \  \exists \; t^{*}_{u_0, \omega} > 0, \quad \mathrm{s.t.} \quad   u_\lambda(t, \omega; u_0) = u_{\lambda, \r}(t, \omega; u_0), \quad  t \in [0, t^{*}_{u_0, \omega}).
\eea

{\HL The idea is then to build the sought local random invariant manifold for Eq. ~\eqref{REE} as a subset $\mathfrak{M}_{\lambda}^{\mathrm{loc}}$ that lives in $\overline{B_\alpha(0,  e^{-z_\sigma(\omega)} \r)}$ of the (global) invariant manifold $\mathfrak{M}_{\lambda, \r}$ for the modified equation \eqref{REE local version}. By taking indeed any initial datum $u_0$ in the interior of such a subset, it holds $\|u_0\|_\alpha < e^{-z_\sigma(\omega)} \r$, so that $u_\lambda = u_{\lambda, \r}$ on $[0, t_{u_0, \omega}^*)$ by \eqref{eq:u and u local}. As we will see, this will imply the local invariance property of $\mathfrak{M}_{\lambda}^{\mathrm{loc}}$ from the invariance property of $\mathfrak{M}_{\lambda, \r}$ under $u_{\lambda, \r}$.

We focus then first on how to ensure $\|u_0\|_\alpha \le e^{-z_\sigma(\omega)} \r$ when $u_0 \in \mathfrak{M}_{\lambda, \r}$. For each such $u_0$, and each $\r \in (0, \r_4)$, we have that
\beas
u_0 = P_{\c}u_0 + h_{\lambda, \r}(P_{\c}u_0, \omega),
\eeas
which implies
\bea  \label{eq:xi range1}
\|u_0\|_\alpha & \le \|P_{\c} u_0\|_\alpha  + \|h_{\lambda, \r}(P_{\c} u_0, \omega)\|_\alpha  \le (1 + \LHW) \|P_{\c} u_0\|_\alpha.
\eea

Recalling that $\LFW \rightarrow 0$ as $\r \rightarrow 0$, we have
\bea  \label{Fountain est}
\exists \; \r_5 > 0, \quad \mathrm{s.t.} \quad  \Upsilon_{1}(F_\r) \le \frac{1}{2}, \quad \r \in (0, \r_5].
\eea
Therefore, by noting that $K \LFW (\eta - {\eta_2})^{\alpha -1}\Gamma(1-\alpha) \le \Upsilon_{1}(F_\r)$ from \eqref{Upsilon_1 local}, we get from \eqref{Lip h local cm 1} that
\be \label{Lip h local cm 2}
\LHWS \le K, \qquad \Forall \, \r \in (0, \r^\ast),
\ee
where \bea  \label{eq:rho}
\r^\ast := \min\{\r_4, \r_5\},
\eea
leading in turn to 
\bea \label{eq:Pc u0}
\|u_0\|_\alpha & \le (1+K) \|P_{\c} u_0\|_\alpha, \quad \Forall \,  \r \in (0, \r^\ast). 
\eea
It is then sufficient that $\|P_{\c} u_0\|_\alpha < \frac{e^{-z_\sigma(\omega)} \r}{1+K}$ to ensure $\|u_0\|_\alpha \le e^{-z_\sigma(\omega)} \r$. Let us now introduce 
\bea
\mathfrak{B}_{\r}(\omega) := \overline{{B_\alpha(0, \delta_{\r}(\omega))}} \cap \mathcal{H}^{\c},  \quad  \r \in (0, \r^\ast), \; \omega\in \Omega,
\eea
with} 
\bea \label{delta-recall}
\delta_{\r}(\omega) := \frac{e^{-z_\sigma(\omega)}\r}{1+K},  \quad  \omega \in \Omega.
\eea

We claim that for each $\r \in (0, \r^\ast)$, $\mathfrak{M}_\lambda^{\mathrm{loc}}$ defined by
\bea  \label{M local-recall}
\mathfrak{M}^{\mathrm{loc}}_\lambda(\omega) := \{\xii + h_\lambda(\xii, \omega) \mid \xii \in \mathfrak{B}_{\r}(\omega) \}, \quad   \omega \in \Omega, \; \lambda\in \Lambda,
\eea
{\HL corresponds to} a local random invariant $C^r$ manifold associated with Eq. ~\eqref{REE}, where the corresponding local random invariant manifold function $h_\lambda$ is taken to be $h_{\lambda, \r}$ as obtained in Step ~2.

The fact that $\mathfrak{M}_{\lambda}^{\mathrm{loc}}$ is a random closed set can be obtained by considering the sequence $\{\gamma_n\}_{n \in \mathbb{N}}$ of measurable mappings $\gamma_n \colon \Omega \rightarrow \mathcal{H}_\alpha$ defined by $\gamma_n(\omega) := e^{-z_\sigma(\omega)}u_n + h_\lambda(e^{-z_\sigma(\omega)}u_n, \omega)$, where $\{u_n\} \in (\mathcal{H}^{\c})^{\mathbb{N}}$ is a sequence dense in $B_\alpha(0, \frac{\r}{1+K}) \cap \mathcal{H}^{\c}$. Then similar to Step ~3 of Appendix ~\ref{appendix section 1}, we get that $\mathfrak{M}_{\lambda}^{\mathrm{loc}}(\omega)=\overline{\{\gamma_n(\omega)\vert n\in \mathbb{N}\}}^{\mathcal{H}_{\alpha}}$, which implies that $\mathfrak{M}_{\lambda}^{\mathrm{loc}}$ is a random closed set by application of the Selection Theorem \cite[Thm.~III.9]{CV77}.

Now, we show that $\mathfrak{M}_{\lambda}^{\mathrm{loc}}$ is locally invariant under the RDS associated with Eq. ~\eqref{REE}. Note that $\|P_{\c} u_0\|_\alpha < \delta_{\r}(\omega)$ for all $u_0$ in the interior, $\mathring{\mathfrak{M}}_{\lambda}^{\mathrm{loc}}(\omega)$,  of $\mathfrak{M}_{\lambda}^{\mathrm{loc}}(\omega)$, which implies that $\|u_0\|_\alpha < e^{-z_\sigma(\omega)} \r$ thanks to \eqref{eq:Pc u0} and \eqref{delta-recall}. Then, it follows from \eqref{eq:u and u local} that there exists $t^{*}_{u_0, \omega} > 0$ such that 
\bea \label{eq:u and u local-2}
u_\lambda(t, \omega; u_0) = u_{\lambda, \r}(t, \omega; u_0), \quad  t \in [0, t^{*}_{u_0, \omega}).
\eea
Since $\mathfrak{M}_{\lambda, \r}$ is invariant under $u_{\lambda, \r}$ {\HL for all $t >0$} and $h_\lambda = h_{\lambda, \r}$, we also have for each $ \omega \in \Omega$ that
\bea
u_{\lambda, \r}(t, \omega; u_0) = P_{\c} u_{\lambda, \r}(t, \omega; u_0) + h_{\lambda}(P_{\c} u_{\lambda, \r}(t, \omega; u_0), \theta_t\omega), \quad \Forall \, t \ge 0,  \, u_0 \in \mathfrak{M}^{\mathrm{loc}}_{\lambda}(\omega).
\eea 
Therefore, 
\bea \label{eq:u representation}
u_{\lambda}(t, \omega; u_0) = P_{\c} u_{\lambda}(t, \omega; u_0) + h_{\lambda}(P_{\c} u_{\lambda}(t, \omega; u_0), \theta_t\omega), \quad  \Forall \, u_0 \in \mathring{\mathfrak{M}}_{\lambda}^{\mathrm{loc}}(\omega), \, t \in [0, t^{*}_{u_0, \omega}).
\eea
Using again $\|P_{\c} u_0\|_\alpha < \delta_{\r}(\omega)$ for all $u_0 \in \mathring{\mathfrak{M}}_{\lambda}^{\mathrm{loc}}(\omega)$, it follows from the continuous dependence of the solutions with respect to the initial data and the continuity of the OU process $t \mapsto z_\sigma(\theta_t\omega)$ that there exists $t^{**}_{u_0, \omega}> 0$ such that 
\bea \label{eq:uc bound}
\|P_{\c} u_{\lambda}(t, \omega; u_0)\|_\alpha < \delta_{\r}(\theta_t\omega), \quad \Forall \, t \in [0, t^{**}_{u_0, \omega}).
\eea
{\HL This together with \eqref{eq:u representation} implies that}
\bea
u_{\lambda}(t, \omega; u_0) \in \mathring{\mathfrak{M}}_{\lambda}^{\mathrm{loc}}(\theta_t\omega), \quad \Forall \, t \in [0, t_{u_0, \omega}),
\eea
where
\beas
t_{u_0, \omega} := \min \{t^{*}_{u_0, \omega}, t^{**}_{u_0, \omega}\}.
\eeas
The local invariance property of $\mathfrak{M}_{\lambda}^{\mathrm{loc}}$ follows.

Note also that $h_\lambda$ clearly satisfies the required conditions in Definition ~\ref{Def. local inv. manfld} since $h_{\lambda, \r}$ does from Step ~2. The proof is now complete.

\ep

\section{Local Stochastic Critical Manifolds: Existence and Approximation Formulas} \label{s:approximation}

In this section, we focus on a typical situation where the control parameter $\lambda$ varies in an interval that contains the critical value, $\lambda_c$, at which the trivial steady state of Eq. ~\eqref{SEE} changes its linear stability stated as  a {\it principle of exchange of stabilities (PES)} given in \eqref{PES} below. We show in Lemma \ref{lem:PES} that this PES implies that  the uniform spectrum decomposition \eqref{gap 3} introduced in Section ~\ref{s:preliminary} is naturally satisfied. It allows us in turn to establish in Proposition \ref{prop:cm} the existence and smoothness of a family of local stochastic critical manifolds in the sense of Definition ~\ref{Def:critical manifold}, which are built --- by relying on Section~\ref{s:local} --- as random graphs over some deterministic neighborhoods of the origin in the subspace spanned by the critical modes that lose their stability as $\lambda$ crosses $\lambda_c$. By construction, these manifolds carry nonlinear dynamical information associated with this loss of linear stability.


We then derive  in  Theorem ~\ref{App CMF} and Corollary ~\ref{Cor:Taylor} our main results concerning explicit random approximation formulas to the leading order of these local stochastic critical manifolds about the origin. 
%

A pullback characterization of the approximation formulas of these stochastic critical manifolds will be derived in Section~\ref{Sec_app_man+PB}. As we will see such a characterization provides geometric insights of approximating manifolds associated with these approximation formulas that will be furthermore useful to build up more general manifolds for the parameterization problem of the ``small'' scales by the ``large'' ones; see Section \ref{ss:PM}.

\subsection{Standing hypotheses}  \label{subsec:hyp}

Consider the SPDE ~\eqref{SEE}. We recall and reformulate some of our assumptions from Section ~\ref{s:preliminary}. The assumptions about the linear operator $L_\lambda$ are those of Section ~\ref{s:preliminary}. {The nonlinearity $F \colon \mathcal{H}_\alpha \rightarrow \mathcal{H}$ is assumed here to be $C^p$ smooth and to take the following form:
\begin{equation} \label{F Taylor}
F(u)= F_k(\underbrace{u, \cdots, u}_{k \text{ times}}) + O(\|u\|^{k+1}_\alpha), 
\end{equation}
where $\mathcal{H}_\alpha$ is as before the space associated with the fractional power $A^\alpha$ for some $\alpha \in [0, 1)$,  $k$ and $p$ are integers such that 
\bea \label{k and p}
p > k \ge 2,
\eea 
and
\bea \label{k-linear def}
F_k \colon \underbrace{\mathcal{H}_\alpha \times \cdots \times \mathcal{H}_\alpha}_{k
\text{ times}} \rightarrow \mathcal{H}
\eea
is a continuous $k$-linear operator. Without any confusion we will often write  $F_k(u)$ instead of $F_k(u,\,\cdots \, , u)$ to simplify the presentation. Note that \eqref{F Taylor} encompasses a large class of nonlinearities satisfying \eqref{F} and \eqref{DF}.

Throughout this section, we also assume that the RPDE \eqref{REE} associated with the SPDE \eqref{SEE} has a unique classical solution for any given initial datum in $\mathcal{H}_\alpha$ in the sense specified in Proposition ~\ref{prop:exist}. We denote as before the corresponding RDS acting on $\mathcal{H}_\alpha$  by $S_\lambda$. The resulting RDS associated with Eq. ~\eqref{SEE} via the cohomology relation \eqref{cohomology relation} is denoted by $\widehat{S}_\lambda$.

We assume furthermore that there exists a critical value $\lambda_c$ and an integer $m >0$ such that the following principle of exchange of stabilities (PES) holds at $\lambda_c$ {\hll for the spectrum $\sigma(L_\lambda)$}:
\begin{equation} \label{PES}
\begin{aligned}
& {\hll \Re \beta_j(\lambda)}
\begin{cases} <0 & \mbox{if } \lambda < \lambda_c, \\ =0 & \mbox{if } \lambda = \lambda_c,\\ >0 & \mbox{if } \lambda > \lambda_c,
\end{cases} && \forall \; j \in \{1, \cdots, m\}, \\
&{\hll \Re \beta_j(\lambda_c)} < 0, && \forall \, j\ge m+1.
\end{aligned}
\end{equation}
 Note that the critical value $\lambda_c$ in \eqref{PES} is the value of $\lambda$ at which the trivial steady state of Eq. ~\eqref{SEE} (or Eq. ~\eqref{REE}) changes its linear stability. 

In what follows we will make use of the following decomposition of $\sigma(L_\lambda)$ associated with the PES condition \eqref{PES}:
\bea \label{splitting}
 & \sigma_{\c}(L_\lambda) := \{\beta_j(\lambda) \:|\: j = 1, \ 2,\ \cdots, m\},  \quad \sigma_{\s}(L_\lambda) := \{\beta_j(\lambda) \:|\: j = m+1, \ m+2,\ \cdots\},
\eea
where $m$ is as given in \eqref{PES}.

Recall that according to the convention adopted in this article, eigenvalues are counted with multiplicity;  so repetitions are allowed  in \eqref{splitting}. In particular  $\sigma_{\c}(L_\lambda)$ may, for instance, be constituted of only one eigenvalue of multiplicity $m$; such situations will be analyzed in details in \cite{CGLW}.

\subsection{Existence of local stochastic critical manifolds} \label{subsec:critical}

Under the assumptions given in the previous subsection, we establish in Proposition ~\ref{prop:cm} below the existence and smoothness of local stochastic (resp. ~random) invariant manifolds for the RDSs associated with Eq. ~\eqref{SEE} (resp. ~Eq. ~\eqref{REE}) for $\lambda$ {\mkk in some neighborhood of}  the critical value $\lambda_c$. The resulting manifolds will be called local stochastic (resp. ~random) critical manifolds; see Definition ~\ref{Def:critical manifold}.

First, let us introduce the following lemma. 

\bl \label{lem:PES}
Assume that the standing hypotheses of \SS~ \ref{subsec:hyp} about the linear operator $L_{\lambda}$ and its spectrum  hold.

Then, for any $r \in \mathbb{N}^\ast$, there exists an open interval $\Lambda_r$ containing the critical value $\lambda_c$ specified in \eqref{PES}, such that
\be  \label{general gap}
0 > r \eta_{\c}(r) > \eta_{\s}(r),
\ee
where
\bea \label{general etac etas}
& \eta_{\c}(r) := \inf_{\lambda \in \Lambda_r} \inf_{j = 1, \cdots, m} \{\Re \beta_j(\lambda)\} , & \eta_{\s}(r):= \sup_{\lambda \in \Lambda_r} \sup_{j \ge m+1} \{\Re \beta_j(\lambda)\}.
\eea

\el

\needspace{1\baselineskip}
\br \label{rmk:PES}

{\hspace*{2em}}

\begin{itemize}

\item[1)] It follows from \eqref{general gap} that $\eta_{\s}(r) < \eta_{\c}(r)$. This shows in particular that the PES condition implies that $L_\lambda$ satisfies the uniform spectrum decomposition \eqref{gap 3}  over $\Lambda_r$ for each $r \in \mathbb{N}^\ast$. The spectrum of $L_{\lambda}$ is thus separated into two disjoint parts, $\sigma_\c(L_{\lambda})$ and $\sigma_\s(L_{\lambda})$ defined in \eqref{splitting}, with $\mathrm{card}(\sigma_c(L_{\lambda}))=m$, as $\lambda$ varies in $\Lambda_r$. This allows us to introduce subspaces $\mathcal{H}^{\c}(\lambda)$ and $\mathcal{H}^{\s}_\alpha(\lambda)$ associated with the decomposition \eqref{splitting} for each $\lambda \in \Lambda_r$ as done in \eqref{Hc} and \eqref{L decomposition} with therefore $\mathcal{H}^{\c}(\lambda)$ of fixed dimension $m$.  For the same reasons as given in Section ~\ref{s:preliminary}, we will {\mkk omit hereafter} to point out the $\lambda$-dependence of these subspaces.

In the sequel, $\mathcal{H}^\c$ (resp.~$\mathcal{H}_\alpha^{\s}$) will be called the critical subspace (resp. ~non-critical subspace) of $\mathcal{H}_\alpha$. The corresponding eigenvectors that span $\mathcal{H}^{\c}$ (resp. ~$\mathcal{H}_\alpha^{\s}$) will be called the critical modes (resp.~non-critical modes).

\item[2)] It is important to note also that condition \eqref{general gap}  --- and thus the PES condition in virtue of Lemma \ref{lem:PES} --- prevents eigenvalues from $\sigma_{\s}(L_\lambda)$ given in \eqref{splitting} to cross the imaginary axis as $\lambda$ varies  in $\Lambda_r$ for any $r \in \mathbb{N}^\ast$. Hence, no eigenvalues other than those of $\sigma_\c(L_\lambda)$ change sign in each such $\Lambda_r$. The approximation formulas provided in Theorem ~\ref{App CMF} and Corollary ~\ref{Cor:Taylor} are subject {\mkk to condition \eqref{general gap}} with $r = 2k$, where $k$ is the leading order of $F$; see \eqref{F Taylor}.

\item[3)] Note that conditions similar to \eqref{general gap} has been often met in the literature when dealing with regularity properties of (deterministic) invariant manifolds; see for instance \cite{Chow_al91, Gallay93, Koc97}. Lemma \ref{lem:PES} shows that such regularity properties may be in fact derived from the PES condition which is easier to check than  \eqref{general gap} in practice.

\end{itemize}

\er

\bp

It follows from the PES condition \eqref{PES} that
 \be\label{eq:tau}
\tau := \Re \beta_{m+1}(\lambda_c) < 0.
\ee
Since $L_\lambda$ is a family of closed operators on $\mathcal{H}$ which depends continuously on $\lambda$, it follows that each eigenvalue $\beta_n(\lambda)$ depends continuously on $\lambda$; see {\it e.g.} \cite[Thm. ~IV 3.16]{Kato95}. This together with \eqref{eq:tau} implies that there exists an open bounded interval $\Lambda$ containing $\lambda_c$ such that
\bea
\sup_{\lambda \in \Lambda} \Re \beta_{m+1}(\lambda) \le \frac{3}{4}\tau.
\eea
Then, according to the ordering of the eigenvalues as specified in \eqref{eq:ordering-1}--\eqref{eq:ordering-3}, we have
\bea \label{gap infinite part}
\sup_{\lambda \in \Lambda} \sup_{j \ge m+1} \Re \beta_j(\lambda)  =  \sup_{\lambda \in \Lambda} \Re \beta_{m+1}(\lambda) \le \frac{3}{4}\tau.
\eea

Since $\Re \beta_j(\lambda_c) = 0$ for $j \in \{ 1,  \cdots,  m\}$ (from the PES condition), by the continuous dependence of the corresponding eigenvalues on $\lambda$ again, we have  for each given $r \in \mathbb{N}^\ast$ that there exists a bounded open interval $\widetilde{\Lambda}_r$ containing $\lambda_c$ such that 
\be \label{gap finite part}
\inf_{\lambda \in \widetilde{\Lambda}_r} \inf_{j = 1, \cdots, m} \Re \beta_j(\lambda) \ge \frac{\tau}{4r}.
\ee

Now, by introducing
\beas \label{Lambda}
\Lambda_r := \Lambda  \cap \widetilde{\Lambda}_r,
\eeas
it follows from \eqref{gap infinite part} and \eqref{gap finite part} that 
\bea \label{gap 3-2}
& \eta_{\c}(r) = \inf_{\lambda \in \Lambda_r} \inf_{j = 1, \cdots, m} \{\Re \beta_j(\lambda)\} \ge \frac{\tau}{4r}, \\
& \eta_{\s}(r) = \sup_{\lambda \in \Lambda_r} \sup_{j \ge m+1} \{\Re \beta_j(\lambda)\} \le \frac{3\tau}{4}.
\eea
Recalling that $\tau < 0$ from \eqref{eq:tau},  we have that $ \eta_{\s}(r) \le \frac{3\tau}{4} < \frac{\tau}{4} \le r \eta_{\c}(r)$. Note also from ~\eqref{PES} that $\Re \beta_j(\lambda) < 0$ for $j \in \{ 1,  \cdots,  m\}$ if $\lambda < \lambda_c$,  so that  $\eta_{\c}(r)$ as given in \eqref{gap 3-2} is negative.  The proof is complete.

\ep

\begin{prop} \label{prop:cm}

Assume that the standing hypotheses of \SS~ \ref{subsec:hyp} hold. For each $r \in \mathbb{N}^\ast$, let $\Lambda_{r}$ be the open interval provided by Lemma ~\ref{lem:PES} with  $\eta_{\s}(r)$ and  $\eta_{\c}(r)$ as specified therein which in particular satisfy
\be
0 > r \eta_{\c}(r) > \eta_{\s}(r).
\ee

{\hh For each $\r > 0$, let $F_{\r}$ be the modified nonlinearity associated with $F$ as given by \eqref{F local version 1}.}

Let us also introduce:
\bea
r^\ast := \min\{r, p\},
\eea
where $p$ indicates the regularity of $F$  given in \SS~\ref{subsec:hyp} as a $C^p$-smooth function.


Then for each 
\bea  \label{eq:eta range}
\eta \in \Bigl( \frac{\eta_{\s}(r)}{r}, \eta_{\c}(r) \Bigr),
\eea 
there exists $\r^\ast> 0$ 
{\hh such that for each $\r \in (0, \r^\ast)$ the following uniform spectral gap conditions hold:
\be \label{Upsilon_1 localb}
\Upsilon_1(F_\r) = K \text{Lip}(F_{\r}) \bigl(({\eta_1} -  \eta)^{- 1} + \Gamma(1-\alpha) (\eta - {\eta_2})^{\alpha-1}  \bigr)  < \frac{1}{2}, 
\ee
and
\be  \label{Upsilon_j localb}
\Upsilon_j(F_\r) = K \text{Lip}(F_{\r}) \bigl(({\eta_1} - j \eta)^{- 1} + \Gamma(1-\alpha) (j \eta - {\eta_2})^{\alpha-1}  \bigr) < 1, \quad \Forall j \in \{2, \cdots, r^\ast\}, 
\ee
where $K$ is as given in the partial-dichotomy estimates \eqref{Proj bounds}.}

Moreover, for each $\r \in (0, \r^\ast)$ and each $\lambda \in \Lambda_{r}$ the following assertions hold:
\bi

\item[(i)] The RDS, $S_{\lambda}$, associated with Eq. ~\eqref{REE} admits  an $m$-dimensional local  random invariant {\mkk $C^{r^\ast}\hspace{-1ex}-$manifold}, $\mathfrak{M}_{\lambda}^{\mathrm{loc}}$, provided by Theorem \ref{thm:local REE}, which is given by:
\bea
\mathfrak{M}^{\mathrm{loc}}_\lambda(\omega) := \{\xii + h_\lambda(\xii, \omega) \mid \xii \in \mathfrak{B}_{\r}(\omega) \}, \quad   \omega \in \Omega,
\eea
where $h_\lambda$ and $\mathfrak{B}_{\r}$ are as specified therein.

\item[(ii)]  

The RDS, $\widehat{S}_\lambda$, associated with Eq. ~\eqref{SEE} via the cohomology relation \eqref{cohomology relation} admits an $m$-dimensional local stochastic invariant {\mkk $C^{r^\ast}\hspace{-1ex}-$manifold}, $\widehat{\mathfrak{M}}_{\lambda}^{\mathrm{loc}}$,  provided by Corollary \ref{thm:local SEE}, which is given by:
\bea  \label{critical manifold SPDE}
\widehat{\mathfrak{M}}^{\mathrm{loc}}_\lambda(\omega) := \{\xii + \widehat{h}_\lambda(\xii, \omega) \mid \xii \in\widehat{\mathfrak{B}}_{\r} \}, \quad  \omega \in \Omega.
\eea  
\ei

\end{prop}

\br
In case (i), recall from Remark \ref{Rem_dyn_interp} (1) and the proof of Theorem  ~\ref{thm:local REE} (Step 3), that the local random invariant manifold, $\mathfrak{M}_\lambda^{\mathrm{loc}}$, is characterized as the random set consisting of all elements $u_0$ in $\mathcal{H}_\alpha$ such that there exists a complete trajectory of the truncated version of Eq. ~\eqref{REE} given by Eq. ~\eqref{REE local version}, passing through $u_0$ at $t=0$ and which has a controlled growth as $t\rightarrow -\infty$.  This control growth is given by $e^{\eta t - \int_t^0 z_\diffusion(\theta_\tau \omega)\, \mathrm{d} \tau}$ which bounds the $\mathcal{H}_\alpha$-norm of  such a trajectory  as $t \rightarrow -\infty$, for $\eta \in ( \frac{\eta_{\s}(r)}{r}, \eta_{\c}(r))$ chosen such that the conditions \eqref{Upsilon_1 localb} and \eqref{Upsilon_j localb} hold.  Similar statement holds in case (ii).
\er

\bp
Lemma \ref{lem:PES} shows that the standing hypotheses of \SS~\ref{subsec:hyp} about the linear operator $L_{\lambda}$  imply that the uniform spectrum decomposition \eqref{gap 3} holds {\HL over $\Lambda_{r}$ for each $r \in \mathbb{N}^\ast$}. This proposition is thus a direct consequence of Theorem ~\ref{thm:local REE} and Corollary ~\ref{thm:local SEE}.  The {\HL $C^{r^\ast}\hspace{-1ex}-$smoothness} of the manifolds is achieved by noting that for  each $\eta$ as given by \eqref{eq:eta range}, we can choose $\eta_1$ and $\eta_2$ such that 
\bea \label{gap 4-2}
0 > \eta_{\c}(r) > \eta_1 > \eta_2 > \eta_{\s}(r),  \quad 0 > r \eta_{1} > \eta_{2},  
\eea
and
\bea \label{eq:eta-111}
\eta \in \Bigl( \frac{\eta_2}{r}, \eta_1 \Bigr),
\eea
which implies that $\eta_2 < r \eta < \eta < \eta_1$ (since $\eta<0$) leading in turn to 
\bea \label{eta ineq 1}
 \eta_2 < j\eta < \eta_1, \quad \Forall j \in \{1, \cdots, r^\ast \},
\eea
{\mkk so that condition \eqref{thm 4.1 condition} of Theorem ~\ref{thm:local REE}  is satisfied.}
\ep

This proposition allows us to introduce the following notion of local stochastic (resp. ~random) critical manifold associated with the  SPDE ~\eqref{SEE} (resp. ~the RPDE ~\eqref{REE}).

\bd \label{Def:critical manifold}

For each $\Lambda_r$ provided by Lemma ~\ref{lem:PES} and each $\lambda \in \Lambda_r$,  a local random invariant {\mkk $C^{r^\ast}\hspace{-1ex}-$manifold} associated with Eq. ~\eqref{REE} guaranteed by Proposition ~\ref{prop:cm} is called a local random critical $C^{r^\ast}\hspace{-1ex}-$manifold, and is denoted by  $\mathfrak{M}_{\lambda}^{\mathrm{crit}}$. The corresponding local invariant manifold function is called a local critical manifold function.

Similarly, for each $\lambda \in \Lambda_r$, a local stochastic invariant $C^{r^\ast}\hspace{-1ex}-$manifold associated with Eq.~\eqref{SEE} guaranteed by Proposition ~\ref{prop:cm} is called a local stochastic critical $C^{r^\ast}\hspace{-1ex}-$manifold, and is denoted by $\widehat{\mathfrak{M}}_\lambda^{\mathrm{crit}}$.

\ed

Hereafter, a local stochastic (resp. ~random) critical manifold will be simply referred to as a local critical manifold for the sake of concision.

\needspace{1\baselineskip}
\br  \label{rmk:critical}


 Note that in the deterministic case, the notion of a critical manifold has been used in \cite[Sect. ~6.3]{Hen81} to refer to a classical center manifold, {\it i.e.}, only when $\lambda = \lambda_c$. Our notion of critical manifold includes but is not limited to center manifold. For instance, when $\lambda > \lambda_c$ a critical manifold corresponds to an unstable manifold.  When $\lambda < \lambda_c$,  it corresponds to the manifold constituted by initial data in $\mathcal{H}_\alpha$  from which emanate solutions that are pullback attracted by the origin, whose rate of attraction is entirely  determined by the decay rate of $\|P_{\c} u\|_{\alpha}$ to zero. These manifolds  constitute thus,   in a certain sense, the natural objects that reflect  locally the PES at the level of the nonlinear dynamics. In \cite{CGLW}, it will be shown that the parameterizing manifolds introduced in Part II of this article can serve to provide global analogues of critical  manifolds  that will turn out to be particularly useful in the study of stochastic bifurcations arising in SPDEs driven by (linear) multiplicative noise.  

\er

\subsection{Approximation of local stochastic critical manifolds} \label{subsec:approx}

We present in this subsection approximation formulas of the local critical manifolds provided by Proposition ~\ref{prop:cm}.

First, let us introduce the following Landau notations for random functions.

\medskip
{\bf  Landau notations.} Let $X$, $Y_1$, and $Y_2$ be real Banach spaces, and $\{\Omega, \mathcal{F}, \mathbb{P}\}$ be a probability space. Consider two measurable mappings $f_i \colon X \times \Omega \rightarrow Y_i$, $i=1,2$. We write 
\be \label{big-O}
f_1(\xi, \omega) = O( f_2(\xi, \omega))
\ee
to mean that there exist a constant $M>0$ and a random open ball $B(0, r(\omega)) \subset X$ centered at $0$, such that the following holds
\bes
\|f_1(\xi, \omega)\|_{Y_1} \le  M \|f_2(\xi, \omega)\|_{Y_2}, \quad \Forall \xi \in B(0,\, r(\omega)), \, \omega \in \Omega.
\ees

We denote by
\be \label{little-o}
f_1(\xi, \omega) = o( f_2(\xi, \omega))
\ee
to mean that for each constant $\epsilon>0$, there exists a random open ball $B(0, r_\epsilon(\omega)) \subset X$ centered at $0$, such that 
\bes 
\|f_1(\xi, \omega)\|_{Y_1} \le  \epsilon \|f_2(\xi, \omega)\|_{Y_2}, \quad \Forall \xi \in B(0,\, r_\epsilon(\omega)), \, \omega \in \Omega.
\ees

Similar notations apply to the case when $f_2$ is a deterministic mapping and $f_2$ is still random. Note that when $f_1$ is also deterministic, these notations agree with the classical Landau notations. It should be clear from the context which one is meant.

\medskip

We turn now to the main {\mkk abstract} results of this article.  For this purpose, let us first introduce the following {\it Lyapunov-Perron integral}:
\bea \label{LP integral}
\mathfrak{I}_{\lambda}(\xi, \omega)=\int_{-\infty}^0 &  e^{\diffusion (k-1)W_s(\omega) \Id}e^{- s L_\lambda } P_{\s} F_k(
e^{s L_\lambda}\xii) \, \mathrm{d}s, \quad \Forall \xii \in  \mathcal{H}^{\c}, \, \omega \in \Omega.
\eea

{\hh Note that this integral is well defined for all $\lambda \in \Lambda$ if $\Lambda \subset \mathbb{R}$ is chosen such that the following condition holds: 
\bea \label{gap for LP integral-1}
\eta_{\s}<  \eta_{\c} \quad \text{ and } \quad \eta_{\s}< k \eta_{\c},
\eea
where $\eta_{\c}$ and $\eta_{\s}$ are as given in \eqref{eta_cs}, {\mkk and $k$ denotes the order of the nonlinear terms $F_k$}. }

Indeed, for such a $\Lambda$, we can find $\eta_1$ and $\eta_2$ such that
\bea \label{gap for LP integral}
\eta_{\s}< \eta_2 < \eta_1 < \eta_{\c}, \quad \text{ and } \quad \eta_2 < k \eta_1.
\eea
Then, by using the dichotomy estimate \eqref{Proj-II}, we get
\beas
\Bigl \|\int_{-\infty}^0 &  e^{\diffusion (k-1)W_s(\omega) \Id}e^{- s L_\lambda } P_{\s} F_k(
e^{s L_\lambda}\xii) \, \mathrm{d}s \Bigr \|_\alpha &\le 
K \int_{-\infty}^0  \frac{e^{\diffusion (k-1)W_s(\omega) - \eta_2 s}}{|s|^\alpha} \Bigl \|F_k(
e^{s L_\lambda}\xii) \Bigr \| \, \mathrm{d}s.
\eeas
Note also that by continuous $k$-linearity of $F_k$ and the dichotomy estimate \eqref{Proj-III}, there exists $C > 0$ such that
\beas
\|F_k(e^{s L_\lambda}\xii) \| \le C  \|e^{s L_\lambda}\xii\|_\alpha^{k} \le CK^k e^{k \eta_1 s}\|\xi\|_\alpha^k,\quad  \Forall \xi \in \mathcal{H}^{\c}.
\eeas
We obtain then
\beas
\Bigl \|\int_{-\infty}^0 &  e^{\diffusion (k-1)W_s(\omega) \Id}e^{- s L_\lambda } P_{\s} F_k(
e^{s L_\lambda}\xii) \, \mathrm{d}s \Bigr \|_\alpha &\le CK^{k+1} \int_{-\infty}^0  \frac{e^{\diffusion (k-1)W_s(\omega)  + (k\eta_1 - \eta_2) s}}{|s|^\alpha} \mathrm{d}s \|\xi\|_\alpha^k,
\eeas
which is finite thanks to the condition $k\eta_1 > \eta_2$ and the fact that $s\mapsto W_s(\omega)$ has sublinear growth as recalled in Lemma ~\ref{Lem:OU}.

The theorem below  shows  that  the Lyapunov-Perron integral, $\mathfrak{I}_\lambda$, {\mkk provides in fact} the leading-order approximation to the stochastic critical manifold function $\widehat{h}_\lambda$  for $\|\xi\|_{\alpha}$ sufficiently small\footnote{Note however that the random radius $\widehat{r}_\epsilon(\omega)$   where the approximation ~\eqref{SEE App CMF formula est-1} holds such as deduced from the proof of Theorem ~\ref{App CMF}, is not optimal.}, and for $\Lambda$ chosen to be $\Lambda_{2k}$ {\mkk such as provided by} Lemma  ~\ref{lem:PES}.

\bt \label{App CMF}


Assume the standing hypotheses of \SS~\ref{subsec:hyp} hold. Let $k$ be the leading order of the nonlinearity specified in \eqref{F Taylor}, and $\Lambda_{2k}$ be the open interval provided by Lemma ~\ref{lem:PES} with $r=2k$.

Let us define, for each $\lambda \in \Lambda_{2k}$,  the mapping $h^\mathrm{app}_\lambda \colon \mathcal{H}^{\c} \times \Omega \rightarrow  \mathcal{H}_\alpha^{\s}$ as the following Lyapunov-Perron integral:
\bea \label{App CMF formula 1}
h^\mathrm{app}_\lambda(\xii, \omega) := e^{(k-1)z_\diffusion(\omega)}\int_{-\infty}^0 &  e^{\diffusion (k-1)W_s(\omega) \Id}e^{- s L_\lambda } P_{\s} F_k(
e^{s L_\lambda}\xii) \, \mathrm{d}s, \quad \Forall \xii \in  \mathcal{H}^{\c}, \, \omega \in \Omega.
\eea

Then 
the local critical manifold function $h_\lambda$ associated with $\mathfrak{M}_\lambda^{\mathrm{crit}}$ 
provided by Proposition ~\ref{prop:cm} (i) (with $r=2k$) can be approximated to the leading order, $k$, by  $h^\mathrm{app}_\lambda$ in  the following sense:

\begin{itemize}
\item For any $\epsilon > 0$, there exists a random open ball $B_\alpha(0, r_\epsilon(\omega))$ such that
\bea \label{App CMF formula est}
 \|h_\lambda(\xii, \omega) - h^\mathrm{app}_\lambda(\xii, \omega)\|_\alpha \le \epsilon \|\xii\|_\alpha^k, \quad \xi \in B_\alpha(0, r_\epsilon(\omega)) \cap \mathcal{H}^{\c}, \, \lambda \in \Lambda_{2k},  \, \omega \in \Omega.
\eea
Furthermore  $r_\epsilon(\omega) \le \delta_{\r}(\omega)$ for all $\omega$, where $\delta_{\r}$ is the random variable given by 
\be  \label{delta in local cmb}
\delta_{\r}(\omega) :=  \frac{e^{-z_\diffusion(\omega)}\r}{1+ K},  \; \omega\in \Omega,
\ee
for $\rho$ sufficiently small\footnote{More precisely, for  $\rho \in (0,\overline{\rho})$ where $\overline{\rho}$ is defined in \eqref{eqn: rho limited} in the proof below.}.


\end{itemize}
 
Similarly, for each $\lambda \in \Lambda_{2k}$, let us define 
$\widehat{h}^\mathrm{app}_\lambda(\xii, \omega) := e^{z_\diffusion(\omega)} h^\mathrm{app}_\lambda(e^{-z_\diffusion(\omega)}\xii, \omega)$, which leads to
\be \label{AF_man}\tag{AF}
\boxed{\widehat{h}^\mathrm{app}_\lambda(\xii, \omega) = \int_{-\infty}^0  e^{ \diffusion (k-1)W_s(\omega)\Id}e^{-s L_\lambda} P_{\s} F_k(
e^{s L_\lambda}\xii) \,\mathrm{d}s, \quad \Forall \xii \in  \mathcal{H}^{\c}, \omega \in \Omega.}
\ee

Then, the local critical manifold function $\widehat{h}_\lambda$ associated with $\widehat{\mathfrak{M}}_\lambda^{\mathrm{crit}}$ provided by Proposition ~\ref{prop:cm} (ii)  (with $r=2k$) can be approximated to the leading order, $k$, by  $\widehat{h}^\mathrm{app}_\lambda$ in the following sense:


\bi
\item For any $\epsilon > 0$, there exists a random open ball $B_\alpha(0, \widehat{r}_\epsilon(\omega))$ such that
\bea  \label{SEE App CMF formula est-1}
 \|\widehat{h}_\lambda(\xii, \omega) - \widehat{h}^\mathrm{app}_\lambda(\xii, \omega)\|_\alpha \le \epsilon \|\xii\|_\alpha^k, \quad \xi \in B_\alpha(0, \widehat{r}_\epsilon(\omega)) \cap \mathcal{H}^{\c}, \, \lambda \in \Lambda_{2k},  \, \omega \in \Omega.
\eea
Furthermore, $\widehat{r}_\epsilon(\omega) \le \widehat{\delta}_{\r}$ for all $\omega$, where $\widehat{\delta}_{\r}$ is the deterministic constant given by
\be  \label{delta in SEE local cmb}
\widehat{\delta}_{\r}:=  \frac{\r}{(1+ K)},
\ee
for $\rho$ sufficiently small.\footnote{Same remark given above applies here.}
\ei
\et

\br\label{Rmk_artificial_bounds2}

Note that the above theorem provides therefore conditions under which the  local critical manifold function $\widehat{h}_\lambda$ associated with $\widehat{\mathfrak{M}}_\lambda^{\mathrm{crit}}$ is approximated by $\widehat{h}^\mathrm{app}_\lambda$ given by \eqref{AF_man},  in the sense that 
\bea \label{SEE App CMF formula est}
 \|\widehat{h}_\lambda(\xii, \omega) - \widehat{h}^\mathrm{app}_\lambda(\xii, \omega)\|_\alpha = o(\|\xii\|_\alpha^k), \quad \Forall \lambda \in \Lambda_{2k},
\eea
where $o(\|\xii\|_\alpha^k)$ is the Landau notation defined in \eqref{little-o}. 

It is important to note that neither these conditions nor the radius over which such an approximation is valid are optimal here. In particular, the deterministic bound provided by \eqref{delta in SEE local cmb} is 
artificial and results from the techniques adopted in the proof of Theorem  \ref{App CMF}; see also Remark ~\ref{rmk:radius}.

\er

\medskip

We assume for the rest of this subsection that $L_\lambda$ is self-adjoint. In this case, {\HL as we will see} the approximation formulas $h^\mathrm{app}_\lambda$  and $\widehat{h}^\mathrm{app}_\lambda$ given respectively in \eqref{App CMF formula 1} and \eqref{AF_man} can be written as random homogeneous polynomials of order $k$ in the critical state variables, and thus constitute genuine leading-order Taylor approximations of the {\hh corresponding} families of local critical manifolds. 

{\mkk For that purpose, first note that $L_\lambda$ being} self-adjoint, all its eigenvalues are real. Since $\mathcal{H}_1$, the domain of $L_\lambda$, is compactly and densely embedded into $\mathcal{H}$, it follows that the set of normalized eigenvectors $\{ e_i \Vert i \in \mathbb{N}^\ast\}$\footnote{As done for the subspaces $\mathcal{H^{\c}}$ and $\mathcal{H}_\alpha^{\s}$,  we suppress the $\lambda$-dependence of the eigenvectors $e_i$.} forms a Hilbert basis of $\mathcal{H}$.\footnote{This is a direct consequence of the spectral properties of symmetric compact operators (see {\it e.g.} \cite[Appendix E, Thm. ~7]{Evans10}), and the fact that there exists a positive constant $a$ such that $(-L_\lambda + a \Id)^{-1}$ exists and is a symmetric compact operator on $\mathcal{H}$.} Then, the critical subspace $\mathcal{H}^{\c}$ and the non-critical subspace $\mathcal{H}_\alpha^{\s}$ are given {\mkk respectively} by:
\bea  \label{eq: spaces adjoint}
& \mathcal{H}^{\c}:= \mathrm{span} \{e_i \; \vert \; i = 1, \cdots, m\}, \; \mathrm{ and } \\
& \mathcal{H}_\alpha^{\s}:= \overline{ \mathrm{span}\{e_j \; \vert \; j = m+1, m+2, \cdots, \}}^{\mathcal{H}_\alpha}.
\eea

Let us also introduce the following notion of random homogeneous  polynomials adapted to our framework, which allows us to make precise what we mean by Taylor approximations of critical manifolds.

\bd \label{Def:random polynomial}

Let $\mathcal{H}^{\c}$ be the $m$-dimensional critical subspace and $\mathcal{H}_\alpha^{\s}$ be the non-critical subspace as given in \eqref{eq: spaces adjoint}.
%
%
%
%
%
A random  {\mkk homogeneous} polynomial function of order $p$ in the critical variable $\xii$ with range in $\mathcal{H}_\alpha^{\s}$ is a function of the form
\bea
g \colon \mathcal{H}^{\c} \times \Omega \rightarrow \mathcal{H}_\alpha^{\s}, \quad (\xii, \omega) \mapsto \sum_{n=m+1}^{\infty} g_n(\xii, \omega) e_n,
\eea
with $g_n(\xii,\cdot):=\langle g(\xi,\cdot),e_n\rangle$,\footnote{where $\langle \cdot, \cdot \rangle$ denoting  the inner-product in the ambient Hilbert space $\mathcal{H}$.}  which satisfies furthermore the following conditions:
\bi
\item[(i)] For each fixed $\xii \in \mathcal{H}^{\c}$, {\hh $g_n(\xii, \cdot)\colon \Omega \rightarrow \mathbb{R}$ is $\bigl(\mathcal{F}; \mathcal{B}(\mathbb{R})\bigr)$-measurable,} where $\mathcal{F}$ is the $\sigma$-algebra associated with the Wiener process;
\item[(ii)] $g_n(\cdot, \omega) \colon \mathcal{H}^{\c} \rightarrow \mathbb{R}$ is a {\mkk homogeneous} polynomial in $\xii_1, \cdots, \xii_m$ of order $p$ for  all $\omega$, where $\xii_i =  \langle \xii, e_i\rangle$ for $i=1, \cdots, m$. 
\ei

\ed

The approximation formula $\widehat{h}^\mathrm{app}_\lambda$ (resp. ~$h^\mathrm{app}_\lambda$) given in \eqref{AF_man} (resp. ~\eqref{App CMF formula 1}) {\HL will be} called the leading-order Taylor approximation of the corresponding local critical manifold function $\widehat{h}_\lambda$ (resp.~$h_\lambda$) if {\HL the estimate \eqref{SEE App CMF formula est-1} (resp. ~\eqref{App CMF formula est}) holds} and if {\mkk  $\widehat{h}^\mathrm{app}_\lambda$ (resp. ~$h^\mathrm{app}_\lambda$) is} furthermore a random {\mkk homogeneous} polynomial of order $k$ in the critical state variable {\mkk $\xi \in \mathcal{H}^{\c}$.}

\bc \label{Cor:Taylor}

Assume the standing hypotheses of \SS~\ref{subsec:hyp} hold. Assume furthermore that $L_\lambda$ is self-adjoint. Then, the approximation $h^\mathrm{app}_\lambda$ given in \eqref{App CMF formula 1} constitutes the leading-order Taylor approximation of the corresponding $h_\lambda$. {\mkk More precisely, this approximation is given by:} 
\begin{equation} \label{diagonal case}
h^\mathrm{app}_\lambda(\xii,\, \omega) = \sum_{n=m+1}^\infty h^{\mathrm{app},n}_\lambda(\xii,\,
\omega) e_n,  \quad \Forall \, \xi \in \mathcal{H}^{\c}, \, \omega \in \Omega,
\end{equation}
where $m$ is the dimension of the critical subspace $\mathcal{H}^{\c}$, and $h^{\mathrm{app},n}_\lambda(\xii,\, \omega)$ is given for each $n \ge m+1$ by
\bea  \label{App CMF formula 2}
h^{\mathrm{app},n}_\lambda(\xii, \omega) &:= e^{(k-1)z_\diffusion(\omega)}  \sum_{(i_1, \cdots, i_k )\in \mathcal{I}^k }\xii_{i_1} \cdots \xii_{i_k} \langle F_k(e_{i_1}, \cdots, e_{i_k}), e_n \rangle \M_n^{i_1,\cdots,i_k}(\omega, \lambda).
\eea
Here, $\mathcal{I}=\{1,\cdots,m\}$,   $\xii_i =  \langle \xii, e_i \rangle$ for each $i \in \mathcal{I}$ with  $\langle \cdot, \cdot \rangle$ denoting  the inner-product in the ambient Hilbert space $\mathcal{H}$, and for each $(i_1, \cdots, i_k ) \in \mathcal{I}^k$, the $\M_n^{i_1,\cdots,i_k}(\omega, \lambda)$-term is given by:
\be \label{Mn_gen}
 \boxed{\M_n^{i_1,\cdots,i_k}(\omega, \lambda) := \int_{-\infty}^0 e^{ \bigl( \sum_{j=1}^{k}\beta_{i_j}(\lambda) -
\beta_n(\lambda) \bigr)s + \diffusion (k-1)W_s(\omega)}\, \mathrm{d}s,}
\ee
{\mkk where  each $\beta_{i_j}(\lambda)$ denotes the eigenvalue associated with the corresponding mode $e_{i_j}$ in $\mathcal{H}^{\c}$, and $\beta_n(\lambda)$  denotes the eigenvalue associated with the high mode $e_n$ in $\mathcal{H}^{\s}$.}


Under the same assumptions, the approximation $\widehat{h}^\mathrm{app}_\lambda$ given in \eqref{AF_man} constitutes the leading-order Taylor approximation of the corresponding {\mkk stochastic manifold function,} $\widehat{h}_\lambda$.  {\mkk More precisely, this approximation is given by:} 
\bea \label{diagonal case SEE}
\widehat{h}^\mathrm{app}_\lambda(\xii,\, \omega) & = \sum_{n=m+1}^\infty \widehat{h}^{\mathrm{app},n}_\lambda(\xii,\, \omega) e_n,  \quad \Forall \, \xi \in \mathcal{H}^{\c}, \, \omega \in \Omega,
\eea
where $\widehat{h}^{\mathrm{app},n}_\lambda(\xii, \omega)$ is given for each $n \ge m+1$ by
\bea  \label{SEE App CMF formula 2}
\widehat{h}^{\mathrm{app},n}_\lambda(\xii, \omega) &:= \sum_{(i_1, \cdots, i_k )\in \mathcal{I}^k} \xii_{i_1} \cdots \xii_{i_k}\langle F_k(e_{i_1}, \cdots, e_{i_k}), e_n \rangle \M_n^{i_1,\cdots,i_k}(\omega, \lambda).
\eea

\ec

\needspace{1\baselineskip}

\br \label{rmk finiteness of xi}

\hspace*{2em}

\bi
\item[1)]  Note that $\M_n^{i_1,\cdots,i_k}(\omega, \lambda)$ is finite for each $(i_1, \cdots, i_k ) \in \mathcal{I}^k$, $\lambda \in \Lambda_{2k}$, and $\omega \in \Omega$. This follows from the fact that $\lim_{t \rightarrow \pm \infty} \frac{W_t(\omega)}{t} = 0$ for each $\omega$ (see Lemma ~\ref{Lem:OU} {\mkk again}), and that $\sum_{j=1}^{k}\beta_{i_j}(\lambda)  > \beta_n(\lambda)$. 
This latter inequality {\mkk holds since} according to \eqref{general gap} and \eqref{general etac etas}, we have {\mkk that} $\sum_{j=1}^{k}\beta_{i_j}(\lambda) \ge k \eta_{\c}(2k) > 2k \eta_{\c}(2k) > \eta_{\s}(2k) \ge \beta_n(\lambda)$ for all $\lambda \in \Lambda_{2k}$. Besides, it can be  checked that $\M_n^{i_1,\cdots,i_k}(\omega, \lambda)$ is a tempered random variable. In Section~\ref{s:reduction}, it will be shown that such random variables convey memory effects that will be expressed in terms of decay of correlations provided that $\sigma$ {\mkk lives} in some admissible range; see Lemma ~\ref{lem:Mn}. 

\item[2)]  {\mkk It is worthwhile noting that the $\M_n^{i_1,\cdots,i_k}(\omega, \lambda)$-terms come with $\langle F_k(e_{i_1}, \cdots, e_{i_k}), e_n \rangle$,  {\it i.e.} ~the nonlinear, leading-order self-interactions of the (corresponding) critical modes, as projected against the non-critical mode $e_n$.} 
Note that by definition,  any permutation  of a $k$-tuple $(i_1,\cdots,i_k)$ corresponds to a same  $\M_n^{i_1,\cdots,i_k}(\omega, \lambda)$-term but may correspond to different {\mkk interactions projected against the $e_n$-mode}.

\item[3)] Note that if $\sigma_{\c}(L_\lambda)$ as defined in \eqref{splitting} consists of one eigenvalue with multiplicity $m$, {\it i.e.,}
\bea  \label{beta-c}
\beta_1(\lambda) = \cdots = \beta_m(\lambda) =: \beta_\ast (\lambda), \Forall \lambda \in \Lambda_{2k},
\eea
then each $\M_n^{i_1,\cdots,i_k}(\omega, \lambda)$-term  reduces to the following simplified form:
\be \label{Mn}
\M_n(\omega, \lambda) := \int_{-\infty}^0 e^{[k\beta_{\ast}(\lambda) -
\beta_n(\lambda)]s + \diffusion (k-1)W_s(\omega)}\, \mathrm{d}s.
\ee

\item[4)]  When $\sigma = 0$, {\sw in the general case for $L_\lambda$}, the approximation result stated in Theorem ~\ref{App CMF} recovers those obtained in \cite[Lemma 6.2.4]{Hen81}, \cite[Chap. ~3]{MW05}, and \cite[Appendix A]{MW}. The proof presented in \SS~\ref{ss:proofs} below works literally for this case by simply setting $\sigma$ and hence $z_\sigma$ to zero; see also Lemma ~\ref{lem:det}. 

{\sw In the self-adjoint case for $L_\lambda$}, the approximation formula ~\eqref{diagonal case SEE} becomes deterministic {\sw when $\sigma = 0$}, where the $\M_n$-terms given in \eqref{Mn_gen} are then reduced to 
\be \label{Mn_gen_det}
\M_n^{i_1,\cdots,i_k}(\lambda) = \int_{-\infty}^0 e^{\bigl( \sum_{j=1}^{k}\beta_{i_j}(\lambda) -
\beta_n(\lambda) \bigr)s}\, \mathrm{d}s = \frac{1}{\sum_{j=1}^{k}\beta_{i_j}(\lambda) -
\beta_n(\lambda)},
\ee
for each $n \ge m+1$ and $\lambda \in \Lambda_{2k}$. {\mkk In this case also, the $\M_n^{i_1,\cdots,i_k}$-terms  come  with the nonlinear, 
 leading-order self-interactions of the critical modes, as projected against the non-critical mode $e_n$.}
\ei 

\er

\subsection{Proofs of Theorem ~\ref{App CMF} and Corollary ~\ref{Cor:Taylor}} \label{ss:proofs}

We turn first to the description of the main ideas of the proof of Theorem ~\ref{App CMF}.

\medskip
{\bf Skeleton of the Proof of Theorem ~\ref{App CMF}.} We focus on the approximation of $h_\lambda$, the results for $\widehat{h}_\lambda$ following then from the cohomology relation given in \eqref{h and h hat}.

The estimate given in \eqref{App CMF formula est} will be {\HL achieved} by using an appropriate pivot quantity $h^\mathrm{app}_{\lambda,\r}$ in the basic triangle inequality:
\bea
\|h_\lambda(\xii, \omega) - h^\mathrm{app}_{\lambda}(\xii, \omega)\|_\alpha \le \|h_\lambda(\xii, \omega) - h^\mathrm{app}_{\lambda,\r}(\xii, \omega)\|_\alpha + 
\|h^\mathrm{app}_{\lambda,\r}(\xii, \omega) - h^\mathrm{app}_{\lambda}(\xii, \omega)\|_\alpha.
\eea

To describe this pivot quantity we first recall {\HL from Section ~\ref{s:local}} that $h_\lambda = h_{\lambda, \r}$, where $h_{\lambda, \r}$ is obtained as the fixed point of the following integral equation
\bea \label{int h1-intro}
 h_{\lambda,\r}(\xii, \omega) = \int_{-\infty}^0 \mathfrak{T}_{\lambda, \sigma}(0, s; \omega) P_{\s} G_{\r}\bigl(\theta_s\omega, u_{\lambda, \r}(s,  \omega; \xii + h_{\lambda,\r}(\xii, \omega))\bigr) \mathrm{d}s
\eea
associated with a modified equation of Eq. ~\eqref{REE} given by
\begin{equation} \label{REE local version-recall}
\frac{\mathrm{d}u}{\mathrm{d}t} = L_\lambda u + z_\diffusion(\theta_t\omega)u + G_{\r}(\theta_t\omega,\, u),
\end{equation}
where $\mathfrak{T}_{\lambda, \sigma}$ is the solution operator introduced in \SS ~\ref{ss:T}, $G_{\r}$ is an appropriate cut-off version of $G$, and {\HLL $u_{\lambda, \r}$ is a mild solution to \eqref{REE local version-recall} on $(-\infty, 0]$ taking value $\xi +  h_{\lambda,\r}(\xii, \omega)$ at $t= 0$}; see {\HL also} \eqref{REE local version}--\eqref{G modified} and \eqref{h eqn}.

The pivot $h^\mathrm{app}_{\lambda,\r}$ is then obtained by replacing $G_\r$ in the integral equation \eqref{int h1-intro} with $G_{\r,k}$ (the leading order term of $G_\r$ as defined in \eqref{F local version 2})  and by replacing the mild solution $u_{\lambda, \r}$ of Eq. ~\eqref{REE local version-recall} involved in \eqref{int h1-intro} with the backward solution $v(s):=\mathfrak{T}_{\lambda, \sigma}(s, 0; \omega) \xii$, $s \le 0$, associated with the random linear equation:
\bea \label{linear RPDE}
\frac{\mathrm{d} v}{\mathrm{d} t} =  L_\lambda v + z_\diffusion(\theta_t\omega)v.
\eea 
{\HLL Note that $v(s)$ is well-defined since $\xi \in \mathcal{H}^{\c}$; see \SS ~\ref{ss:T}.}

In other words, the pivot quantity $h^\mathrm{app}_{\lambda, \r}$ reads:
\bea \label{pivot-intro}
h^\mathrm{app}_{\lambda,\r}(\xii, \omega) := \int_{-\infty}^0 &  \mathfrak{T}_{\lambda, \sigma}(0, s; \omega) P_{\s}G_{\r, k}(\theta_s \omega, \mathfrak{T}_{\lambda, \sigma}(s, 0; \omega)\xii) \mathrm{d}s, \quad  \xii \in \mathcal{H}^{\c}, \, \omega\in \Omega.
\eea

Let us first describe how we control $\|h^\mathrm{app}_{\lambda,\r}(\xii, \omega) - h^\mathrm{app}_{\lambda}(\xii, \omega)\|_\alpha$. Actually, $h^\mathrm{app}_{\lambda, \r}$ differs from the approximation formula $h^\mathrm{app}_{\lambda}$ as rewritten in \eqref{App CMF formula 1 rewritten} by replacing $G_{k}$ with $G_{\r, k}$, so that $h^\mathrm{app}_{\lambda,\r} - h^\mathrm{app}_{\lambda}$ is given by:
\beas 
h^\mathrm{app}_\lambda(\xii, \omega) -h^\mathrm{app}_{\lambda,{\r}}(\xii, \omega)  = \int_{-\infty}^0 \mathfrak{T}_{\lambda, \sigma}(0, s; \omega)  \Big( 1 - \zeta\Big( \frac{e^{z_\diffusion(\theta_s \omega)} \|v(s)\|_\alpha}{{\r}}\Big) \Big) P_{\s} G_k(\theta_s \omega,
v(s)) \, \mathrm{d}s;
\eeas
see \eqref{diff h 22}. 

By noting that  $\zeta\Big( \frac{e^{z_\diffusion(\theta_s \omega)} \|v(s)\|_\alpha}{{\r}}\Big) = 1$ if $e^{z_\diffusion(\theta_s \omega)} \|v(s)\|_\alpha \le {\r}$, the last identity above becomes:
\beas
h^\mathrm{app}_\lambda(\xii, \omega) -h^\mathrm{app}_{\lambda,{\r}}(\xii, \omega)  = \int_{-\infty}^{s_0(\xi)} \mathfrak{T}_{\lambda, \sigma}(0, s; \omega)  \Big( 1 - \zeta\Big( \frac{e^{z_\diffusion(\theta_s \omega)} \|v(s)\|_\alpha}{{\r}}\Big) \Big) P_{\s} G_k(\theta_s \omega, v(s)) \, \mathrm{d}s,
\eeas
where $s_0(\xi)<0$ denotes the first time at which the backward solution $v(s)$ of Eq. ~\eqref{linear RPDE} emanating from some $\xi \in \mathcal{H}^{\c}$ leaves $B_\alpha(0, e^{-z_\diffusion(\theta_s \omega)} \r)$.  This remark allows us in Step ~5 to control  $\| h^\mathrm{app}_\lambda(\xii, \omega) -h^\mathrm{app}_{\lambda,{\r}}(\xii, \omega)\|_\alpha$ as follows by usage of the partial-dichotomy estimates \eqref{eq:dichotomy} and continuous $k$-linear properties of $G_k$:
\beas
\| h^\mathrm{app}_\lambda(\xii, \omega) -h^\mathrm{app}_{\lambda,{\r}}(\xii, \omega)\|_\alpha \le \int_{-\infty}^{s_0(\xi)} \kappa_0(s, \omega)\d s \|\xi\|^k, 
\eeas
where $\kappa_0(s, \omega)$ is some exponentially decaying factor; see \eqref{diff h 2}.\footnote{The OU process $z_\sigma(\theta_s\omega)$ is involved in the factor $\kappa_0(s, \omega)$, and its control is made possible thanks to the growth control estimates given in \eqref{Eq:z-4}.} 

By taking $\xi$ in a sufficiently small random ball $B_\alpha(0, r_{\epsilon}(\omega))\cap \mathcal{H}^{\c}$, it is then shown that 
\beas
\int_{-\infty}^{s_0(\xi)} \kappa_0(s, \omega)\d s  \le \frac{\epsilon}{2},
\eeas
leading to a good control of $\| h^\mathrm{app}_\lambda(\xii, \omega) -h^\mathrm{app}_{\lambda,{\r}}(\xii, \omega)\|_\alpha$. The last inequality above follows from the observation that $s_0(\xi)$ gets more negative as $\|\xi\|_\alpha$ gets smaller.

\medskip
The control of $\|h_\lambda(\xi, \omega) - h^\mathrm{app}_{\lambda,\r}(\xi, \omega)\|_\alpha$ is more challenging and is carried out in Steps ~2, 3, and 4. Note that according to \eqref{int h1-intro}, \eqref{pivot-intro}, and the fact $h_\lambda = h_{\lambda, \r}$, it holds that
\bea  \label{main task intro}
h_\lambda(\xii, \omega) - h^\mathrm{app}_{\lambda,\r}(\xii, \omega) = \int_{-\infty}^0 \mathfrak{T}_{\lambda, \sigma}(0, s; \omega) P_{\s} \Bigl( G_{\r}(\theta_s \omega, u(s))  - G_{\r, k}(\theta_s \omega, v(s)) \Bigr)
\mathrm{d}s,
\eea
where we have introduced $u(s):= u_{\lambda, \r}(s,  \omega; \xii + h_{\lambda,\r}(\xii, \omega)$ to simplify the presentation. The main task is then to obtain a good control of the quantity $\|G_{\r}(\theta_s \omega, u(s)) -  G_{\r,k}(\theta_s \omega, v(s))\|$, which can be split further as follows by introducing another pivot $G_{\r,k}(\theta_s \omega, u(s))$:
\bea \label{G diff-intro}
 \|G_{\r}(\theta_s \omega, u(s)) -  G_{\r,k}(\theta_s \omega, v(s))\| & \le \|G_{\r}(\theta_s \omega, u(s)) - G_{\r,k}(\theta_s \omega, u(s))\| \\
&  \hspace{1em} + \|G_{\r,k}(\theta_s \omega, u(s)) - G_{\r,k}(\theta_s \omega, v(s))\|.
\eea

We then point out in Step ~2 at the level of ``vector fields'' key estimates of $\|G_{\r}(\omega, u) -  G_{\r,k}(\omega, u)\|$ and $\|G_{\r,k}(\omega, u_1) -  G_{\r,k}(\omega, u_2)\|$ that result from continuous $k$-linear properties of the leading order term $F_k$ (see again \eqref{G modified} and \eqref{F local version 2} for the definitions of the notations):
\beas 
 \|G_{\r}(\omega, u) - G_{\r,k}(\omega, u)\| & \le C_1(\omega) \|u\|^{k+1}_\alpha,  && \Forall u \in \mathcal{H}_\alpha, \, \omega \in \Omega, \\
\|G_{\r,k}  (\omega, u_1) - G_{\r,k} (\omega, u_2)\| & \le C_2(\omega) \big( \|u_1\|_\alpha^{k-1} + \|u_2\|_\alpha^{k-1}\big) \|u_1 - u_2\|_\alpha \\
& \hspace{4em} + C_3(\omega) \|u_2\|_\alpha^k\|u_1 - u_2\|_\alpha,  && \Forall u_1, u_2 \in \mathcal{H}_\alpha, \, \omega \in \Omega,
\eeas
where the $C_i(\omega)$'s are some positive random constants; see \eqref{G est. 2-1}--\eqref{G est. 2-2}. Using these estimates and the partial-dichotomy estimates ~\eqref{eq:dichotomy} in \eqref{main task intro}, we then obtain within the same step that
\beas
\|h_\lambda(\xi, \omega) - h^\mathrm{app}_{\lambda,\r}(\xi, \omega)\|_\alpha  \le I_1(\xi, \omega) + I_2(\xi, \omega) + I_3(\xi, \omega),
\eeas
where
\bea \label{I-intro}
I_1(\xi, \omega) &= \int_{-\infty}^0 \kappa_1(s, \omega) \|u(s)\|_\alpha^{k+1}\d s, \\
I_2(\xi, \omega) &= \int_{-\infty}^0 \kappa_2(s, \omega) (\|u(s)\|_\alpha^{k-1} + \|v(s)\|_\alpha^{k-1}) \|u(s) - v(s)\|_\alpha \d s, \\
I_3(\xi, \omega) &= \int_{-\infty}^0 \kappa_3(s, \omega) \|v(s)\|_\alpha^{k} \|u(s)- v(s)\|_\alpha \d s,
\eea
with the $\kappa_i(s, \omega)$'s denoting some exponentially decaying factors,  and $u(s)$, $v(s)$ the solutions  associated with respectively Eq. ~\eqref{REE local version-recall} and Eq. ~\eqref{linear RPDE} as before; see \eqref{I} for more details. Here, $I_1(\xi, \omega)$ provides a control of
\beas
\Bigl \|\int_{-\infty}^0  \mathfrak{T}_{\lambda, \sigma}(0, s; \omega) P_{\s} \Big (G_{\r}(\theta_s \omega, u(s)) - G_{\r, k}(\theta_s \omega, u(s)) \Big) \d s \Bigr \|_\alpha,
\eeas
 and $I_1(\xi, \omega) + I_3(\xi, \omega)$ provides a control of 
\beas
\Bigl\|\int_{-\infty}^0  \mathfrak{T}_{\lambda, \sigma}(0, s; \omega) P_{\s} \Big (G_{\r, k}(\theta_s \omega, u(s)) - G_{\r, k}(\theta_s \omega, v(s)) \Big) \d s \Bigr\|_\alpha.
\eeas

In Step ~3, the following estimates of $\|u(s)\|_\alpha$, $\|v(s)\|_\alpha$, and $\|u(s) - v(s)\|_\alpha$ are then carried out:
\begin{align}  \label{flow estimate intro}
& \|u(s)\|_\alpha \le \kappa_4(s, \omega) \|\xii\|_\alpha, && \|v(s)\|_\alpha  \le \kappa_5(s,\omega) \|\xii\|_\alpha,  && \|u(s) - v(s)\|_\alpha \le  \kappa_6(s, \omega) \|\xii\|_\alpha^k,  
\end{align}
where the $\kappa_i(s,\omega)$'s are positive but not necessarily decaying factors; see \eqref{u est.}--\eqref{u-v est.}. 

The estimate for $\|v(s)\|_\alpha$ follows directly from the partial-dichotomy estimate \eqref{d-III}. The controls of $\|u(s)\|_\alpha$ and $\|u(s) - v(s)\|_\alpha$ are subject to a control of $\|u(\cdot)\|_{C_\eta^{-}}$ by $2K\|\xi\|_\alpha$ as pointed out at the beginning of Step ~3 resulting from application of the estimate \eqref{Lip in x} derived in the proof of Theorem ~\ref{Lip manfd}. The rest of Step ~3 is devoted to showing that 
\beas
\|u(s) - v(s)\|_\alpha \le C_4(\omega)\kappa_7(s,\omega)\|u(\cdot)\|_{C_\eta^-}^k,
\eeas 
where $\kappa_7(s,\omega)$ is a positive factor and $C_4(\omega)$ is a positive random constant. This quantity $C_4(\omega)$ is obtained by appropriately controlling the integrands arising from application of the dichotomy estimates to the integral equation \eqref{VCF-1} given below that $u(t)$ satisfies; see \eqref{u-v}--\eqref{u-v 2} and \eqref{u-v 5}. A basic algebraic lemma (Lemma ~\ref{lem:algebra}) is then used to achieve this appropriate control.

In Step ~4, by using \eqref{flow estimate intro} in \eqref{I-intro}, we obtain with the help of Lemma ~\ref{lem:algebra} that
\beas
\|h_\lambda(\xi, \omega) - h^\mathrm{app}_{\lambda,\r}(\xi, \omega)\|_\alpha & \le I_1(\xi, \omega) + I_2(\xi, \omega) + I_3(\xi, \omega) \\
& \le   C_5(\omega)\|\xi\|_\alpha^{k+1} +  C_6(\omega)\|\xi\|_\alpha^{2k-1} + C_7(\omega)\|\xi\|_\alpha^{2k}, \quad \Forall \xi \in \mathcal{H}^{\c},
\eeas
where the $C_i(\omega)$'s are some positive random constants; see \eqref{diff h 1 est.}. It then follows that $\|h_\lambda(\xi, \omega) - h^\mathrm{app}_{\lambda,\r}(\xi, \omega)\|_\alpha \le \frac{\epsilon}{2}\|\xi\|_\alpha^k$ when $\xi$ is in some sufficiently small random ball $B_\alpha(0, r_\epsilon(\omega)) \cap \mathcal{H}^{\c}$.


Finally, the estimate \eqref{App CMF formula est} follows by combining estimates for $\|h_{\lambda}(\xii, \omega) - h^\mathrm{app}_{\lambda,\r}(\xii, \omega)\|_\alpha$ and $\|h^\mathrm{app}_{\lambda,\r}(\xii, \omega) - h^\mathrm{app}_{\lambda}(\xii, \omega)\|_\alpha$, which is done in Step ~6.


\bp[{\bf Proof of Theorem ~\ref{App CMF}}] We proceed in six steps following the ideas outlined above.

\medskip
{\bf Step 1. Integral equation for $h_\lambda$ and the pivot $h^\mathrm{app}_{\lambda,\r}$.} In this step, we recall an integral equation that $h_\lambda$ satisfies, and introduce an intermediate approximation formula $h^\mathrm{app}_{\lambda,\r}$. Corresponding error estimates between $h_\lambda$ and $h^\mathrm{app}_{\lambda,\r}$ will be {\HLL obtained} in the next three steps. 

First, let us recall that from the proof of Theorem ~\ref{thm:local REE}, we have
\beas
h_\lambda = h_{\lambda, \r},
\eeas
where $h_{\lambda, \r}$ is the global random invariant manifold function associated with the modified equation \eqref{REE local version}. By application of Theorem ~\ref{Lip manfd} to Eq. ~\eqref{REE local version}, we know that $h_{\lambda, \r}$ satisfies the following integral equation:
 \begin{equation} \label{int h1}
\begin{aligned}
 h_{\lambda,\r}(\xii, \omega) = \int_{-\infty}^0 \mathfrak{T}_{\lambda, \sigma}(0, s; \omega) P_{\s} G_{\r}\big(\theta_s\omega, u_{\lambda, \r}(s,  \omega; \xii + h_{\lambda,\r}(\xii, \omega))\big)\,\mathrm{d}s,
\end{aligned}
\end{equation}
where $\mathfrak{T}_{\lambda, \sigma}$ is as given in \SS~\ref{ss:T} and $u_{\lambda, \r}(\cdot,  \omega; \xii + h_{\lambda,\r}(\xii, \omega))$ is a mild solution of Eq. ~\eqref{REE local version} defined on $(-\infty, 0]$ taking value $\xii + h_{\lambda,\r}(\xii, \omega)$ at $t=0$. Recall from the construction of such an $h_{\lambda,\r}$ (see Step ~1 in the proof of Theorem ~\ref{Lip manfd}) that the corresponding mild solution $u_{\lambda, \r}(\cdot,  \omega; \xii + h_{\lambda,\r}(\xii, \omega))$ is obtained as the unique fixed point in $C_\eta^{-}$ of the integral operator $\mathcal{N}_{\xi, \r}^{\omega, \lambda}$ defined by:
\bea \label{VCF-1} 
\mathcal{N}_{\xii,\r}^{\omega, \lambda}[u](s) & :=  \mathfrak{T}_{\lambda, \sigma}(s, 0; \omega) \xii  - \int_s^0 \mathfrak{T}_{\lambda, \sigma}(s, s'; \omega) P_{\c} G_{\r}(\theta_{s'}\omega, u(s'))\,\mathrm{d} s'
\\
& \hspace{2em} + \int_{-\infty}^s \mathfrak{T}_{\lambda, \sigma}(s, s'; \omega) P_{\s} G_{\r}(\theta_{s'}\omega, u(s'))\,\mathrm{d} s', \qquad s \le 0.
\eea
{\HLL This fact will be used in Step ~3 below.}
 
Now, for each $\r > 0$ and $\lambda \in \Lambda_{2k}$, let us introduce a mapping $h^\mathrm{app}_{\lambda,\r}\colon \mathcal{H}^{\c} \times \Omega \rightarrow \mathcal{H}_\alpha^{\s}$ defined by:
\bea \label{int h2}
h^\mathrm{app}_{\lambda,\r}(\xii, \omega) := \int_{-\infty}^0 &  \mathfrak{T}_{\lambda, \sigma}(0, s; \omega) P_{\s}G_{\r, k}(\theta_s \omega, \mathfrak{T}_{\lambda, \sigma}(s, 0; \omega)\xii) \,
\mathrm{d}s, \quad  \xii \in \mathcal{H}^{\c}, \, \omega\in \Omega,
\eea
where
\bea \label{F local version 2}
& G_{\r, k}(\omega,\, u) := e^{-z_\diffusion(\omega)} F_{\r,k}(e^{z_\diffusion(\omega)}u) \quad \text{with }  \quad F_{\r,k}(u):=\zeta\Big(\frac{\|u\|_\alpha}{\r}\Big)F_k(u),
\eea
and $\zeta$ is the cut-off function defined in \eqref{cut-off}.

Compared with $h_{\lambda,\r}(\xii, \omega)$ defined in \eqref{int h1}, the nonlinearity $G_{\r}$ is replaced here by $G_{\r, k}$, and the pathwise solution $u_{\lambda, \r}$ is replaced by the solution $\mathfrak{T}_{\lambda, \sigma}(s, 0; \omega) \xii$ of the linear problem $\frac{\mathrm{d} v}{\mathrm{d} t} =  L_\lambda v + z_\diffusion(\theta_t\omega)v$ with initial datum $\xii \in \mathcal{H}^{\c}$.

We aim to show in the next three steps that there exists $\overline{\r}>0$, such that for each fixed $\r \in (0, \overline{\r})$ and any $\epsilon >0$, there exists a random open ball $B_\alpha(0, r^{\ast}_{\epsilon}(\omega))\subset \mathcal{H}_\alpha$ verifying 
\bea  \label{diff h aim 1}
\| h_{\lambda,\r}(\xii, \omega) - h^\mathrm{app}_{\lambda,\r}(\xii, \omega)\|_\alpha \le \frac{\epsilon}{2} \|\xii\|_\alpha^k, \quad \forall \, \xii \in B_\alpha(0, r^{\ast}_{\epsilon}(\omega)) \cap \mathcal{H}^{\c},\, \lambda \in \Lambda_{2k}, \, \omega \in \Omega.
\eea


Since in what follows all estimates are uniform in $\lambda$ for all $\lambda \in \Lambda_{2k}$, and are produced for each fixed $\omega$, to simplify the presentation, we introduce for all $s \le 0$ {\HLL and $\xi \in \mathcal{H}^{\c}$} the following two notations keeping only the time dependence explicit:\footnote{It is also safe to suppress the $\r$-dependence of the mild solution $u_{\lambda, \r}$ because all the estimates on $u_{\lambda, \r}$ are done for a fixed $\r \in (0, \overline{\r})$, with $\overline{\r}$ specified later in \eqref{eqn: rho limited}. {\HLL The dependence on $\sigma$ is also removed for the sake of readability.}}
\bea \label{def. u & v}
&u(s) := u_{\lambda, \r}(s,  \omega; \xii + h_{\lambda,\r}(\xii, \omega)), \\
& v(s) := \mathfrak{T}_{\lambda, \sigma}(s, 0; \omega) \xii.
\eea
Then from \eqref{int h1} and \eqref{int h2} we obtain 
\bea \label{h diff}
h_{\lambda,\r}(\xii, \omega) - h^\mathrm{app}_{\lambda,\r}(\xii, \omega) = \int_{-\infty}^0  \mathfrak{T}_{\lambda, \sigma}(0, s; \omega) \Big (P_{\s}G_{\r}(\theta_s \omega, u(s)) - P_{\s}G_{\r, k}(\theta_s \omega, v(s)) \Big) \d s.
\eea

{\HLL The integral representation of $h_{\lambda,\r} - h^\mathrm{app}_{\lambda,\r}$ involved in \eqref{h diff} motivates} the estimates presented in the next three steps, which are organized as follows. In Step ~2, we first point out some key estimates related to $G_{\r}$ and $G_{\r, k}$ at the level of ``vector fields,'' which are then used to control \eqref{h diff} via three integral terms $I_1(\xi, \omega)$, $I_2(\xi, \omega)$, and $I_3(\xi, \omega)$. These integral terms involve the mild solution $u(s)$ and the linear flow $v(s)$; see \eqref{I}. In Step ~3, we derive {\HLL some related estimates of $\|u(s)\|_\alpha$, $\|v(s)\|_\alpha$, and $\|u(s) - v(s)\|_\alpha$}, which will be used in Step ~4 to derive the desired error estimate announced in \eqref{diff h aim 1}.

\medskip
{\bf Step 2. Estimates of the nonlinear terms at the level of {\HL ``vector fields.''}} In this step, we first establish some estimates about $G_{\r}$ and $G_{\r, k}$ {\HLL as ``vector fields'' from $(\mathcal{H}_\alpha, \|\cdot\|_\alpha)$ to $(\mathcal{H}, \|\cdot\|)$}, and then derive a control of $\| h_{\lambda,\r}(\xii, \omega) - h^\mathrm{app}_{\lambda,\r}(\xii, \omega) \|_\alpha$ as given in \eqref{I}.

{\HLL We first note that:}
\bea \label{G diff}
 \|P_{\s}G_{\r}(\theta_s \omega, u(s)) -  P_{\s}G_{\r,k}(\theta_s \omega, v(s))\| & \le \|P_{\s}G_{\r}(\theta_s \omega, u(s)) - P_{\s}G_{\r,k}(\theta_s \omega, u(s))\| \\
&  \hspace{1em} + \|P_{\s}G_{\r,k}(\theta_s \omega, u(s)) - P_{\s}G_{\r,k}(\theta_s \omega, v(s))\|.
\eea

To proceed further from \eqref{G diff}, we establish some estimates about $F_{\r}$ and $F_{\r, k}$, where 
\bea \label{F local version}
& F_\r(u):=\zeta\Big(\frac{\|u\|_\alpha}{\r}\Big)F(u), \quad F_{\r,k}(u):=\zeta\Big(\frac{\|u\|_\alpha}{\r}\Big)F_k(u),
\eea
and $\zeta$ is still the cut-off function defined in \eqref{cut-off}.

In the sequel, $C>0$ will denote a generic positive constant which may or may not depend on $\r$. We claim that there exist $\widetilde{\r}> 0$ and a positive constant $C$ such that:
\begin{subequations} \label{F high order est.} 
\begin{align}
 \|F_{\r}(u)\|  & \le  C \|u\|^{k}_\alpha,  \qquad \Forall \, \r \in (0, \widetilde{\r}),  \, u \in \mathcal{H}_\alpha, \label{eq:F-1}\\
 \|F_{\r}(u) - F_{\r,k}(u)\|  & \le C \|u\|^{k+1}_\alpha, \quad \Forall \, \r \in (0, \widetilde{\r}), \, u \in \mathcal{H}_\alpha. \label{eq:F-2}
\end{align}
\end{subequations}

Let us first check \eqref{eq:F-1}. It follows from the assumption \eqref{F Taylor} that $F(u) = O(\|u\|_\alpha^k)$. Hence, there exists $\widetilde{\r}> 0$ and a positive constant $C$ such that 
\beas
\|F(u)\|  & \le C \|u\|^{k}_\alpha,  \qquad \Forall \, u  \in B_\alpha(0, 2 \widetilde{\r}).
\eeas
Then, for each $\r \in (0, \widetilde{\r})$, \eqref{eq:F-1} holds for all $u \in B_\alpha(0, 2\r)$. Note also that \eqref{eq:F-1} holds obviously if $\|u\|_\alpha \ge 2\r$ by the definition of the cut-off function $\zeta$ in \eqref{cut-off}. 

From \eqref{F Taylor}, we also have that $F(u) - F_k(u) = O(\|u\|_\alpha^{k+1})$. Following the same arguments as above, we obtain \eqref{eq:F-2} by choosing possibly a smaller $\widetilde{\r}$ and a larger constant $C$. 

Moreover, since $F_k$ is a continuous $k$-linear operator {\HLL and $\zeta$ is Lipschitz, we have that for each $\r > 0$, there exists $C> 0$} such that
\bea \label{F_k est.}
&  \|F_{\r,k}(u_1) - F_{\r,k}(u_2)\|  =  \Big\|\zeta\Big(\frac{\|u_1\|_\alpha}{\r}\Big) F_{k}(\underbrace{u_1, \cdots, u_1}_{k
\text{ times}}) - \zeta\Big(\frac{\|u_2\|_\alpha}{\r}\Big) F_{k}(\underbrace{u_2, \cdots, u_2}_{k
\text{ times}}) \Big\| \\
& \le \Big \| \zeta\Big(\frac{\|u_1\|_\alpha}{\r}\Big) F_{k}(u_1, \cdots, u_1) - \zeta\Big(\frac{\|u_1\|_\alpha}{\r}\Big) F_{k}(\underbrace{u_1, \cdots, u_1}_{k-1
\text{ times}}, u_2) \Big\| \\
& \hspace{1em} + \Big\| \zeta\Big(\frac{\|u_1\|_\alpha}{\r}\Big) F_{k}(\underbrace{u_1, \cdots, u_1}_{k-1
\text{ times}}, u_2) - \zeta\Big(\frac{\|u_1\|_\alpha}{\r}\Big) F_{k}(\underbrace{u_1, \cdots, u_1}_{k-2
\text{ times}}, u_2, u_2) \Big \| + \cdots \\
& \hspace{1em} + \Big \| \zeta\Big(\frac{\|u_1\|_\alpha}{\r}\Big) F_{k}(u_1, u_2, \cdots, u_2) - \zeta\Big(\frac{\|u_1\|_\alpha}{\r}\Big) F_{k}(u_2, \cdots, u_2) \Big \| \\
& \hspace{1em} + \Big \| \zeta\Big(\frac{\|u_1\|_\alpha}{\r}\Big) F_{k}(u_2, \cdots, u_2) - \zeta\Big(\frac{\|u_2\|_\alpha}{\r}\Big) F_{k}(u_2, \cdots, u_2) \Big \| \\
& \le C \sum_{i = 0}^{k-1} \|u_1\|_\alpha^{k-1-i}\|u_2\|_\alpha^{i} \|u_1 - u_2\|_\alpha +  C \|u_2\|_\alpha^k\|u_1 - u_2\|_\alpha \\
& \le C \sum_{i = 0}^{k-1} \big( \|u_1\|_\alpha^{k-1} + \|u_2\|_\alpha^{k-1}\big) \|u_1 - u_2\|_\alpha + C \|u_2\|_\alpha^k\|u_1 - u_2\|_\alpha \\
& =  k C \big( \|u_1\|_\alpha^{k-1} + \|u_2\|_\alpha^{k-1}\big) \|u_1 - u_2\|_\alpha + C \|u_2\|_\alpha^k\|u_1 - u_2\|_\alpha, \Forall u_1, u_2 \in \mathcal{H}_\alpha.
\eea

Now, let us define $\overline{\r}$ by:
\bea \label{eqn: rho limited}
\overline{\r} := \min \{\widetilde{\r}, \r^\ast\},
\eea
where $\widetilde{\r}$ is as specified in \eqref{F high order est.}, and $\r^\ast$ is as given in Proposition ~\ref{prop:cm} to ensure the existence of a critical manifold for $\lambda \in \Lambda_{2k}$; see also Remark ~\ref{rmk:radius}. From now on, we fix an arbitrary $\r \in (0, \overline{\r})$.

Recalling from \eqref{G modified} and \eqref{F local version 2} that
\beas
G_{\r}(\omega,\, u) = e^{-z_\diffusion(\omega)} F_{\r}(e^{z_\diffusion(\omega)}u), \quad G_{\r, k}(\omega,\, u) = e^{-z_\diffusion(\omega)} F_{\r,k}(e^{z_\diffusion(\omega)}u), 
\eeas
it follows then from \eqref{F high order est.} that
\bea \label{G est. 2-1}
 \|G_{\r}(\omega, u)\|  & \le C e^{(k-1) z_\diffusion(\omega)} \|u\|^{k}_\alpha, && \Forall u \in \mathcal{H}_\alpha, \, \omega \in \Omega, \\
 \|G_{\r}(\omega, u) - G_{\r,k}(\omega, u)\| & \le C e^{k z_\diffusion(\omega)} \|u\|^{k+1}_\alpha, && \Forall u \in \mathcal{H}_\alpha, \, \omega \in \Omega.
\eea
Similarly, it follows from \eqref{F_k est.} that
\bea  \label{G est. 2-2}
\|G_{\r,k} (\omega, u_1) - G_{\r,k} (\omega, u_2)\|  & \le k C e^{(k-1)z_\diffusion(\omega)}\big( \|u_1\|_\alpha^{k-1} + \|u_2\|_\alpha^{k-1}\big) \|u_1 - u_2\|_\alpha \\
& \hspace{1.5em} + C e^{k z_\diffusion(\omega)} \|u_2\|_\alpha^k\|u_1 - u_2\|_\alpha, \quad \Forall u_1, u_2 \in \mathcal{H}_\alpha, \, \omega \in \Omega. 
\eea

It follows directly from \eqref{G est. 2-1} and \eqref{G est. 2-2} that
\bea \label{G est. 3-1}
& \|P_{\s}G_{\r}(\theta_s \omega, u(s)) - P_{\s}G_{\r,k}(\theta_s \omega, u(s))\| \le C e^{k z_\diffusion(\theta_s \omega)} \|u(s)\|^{k+1}_\alpha, \Forall s \le 0, 
\eea
and
\bea  \label{G est. 3-2}
\|P_{\s}G_{\r,k}(\theta_s \omega, u(s)) & - P_{\s}G_{\r,k}(\theta_s \omega, v(s))\| \\
& \le k C e^{(k-1)z_\diffusion(\theta_s \omega)}\big( \|u(s)\|_\alpha^{k-1} + \|v(s)\|_\alpha^{k-1}\big) \|u(s) - v(s)\|_\alpha \\
& \hspace{3em} + C e^{k z_\diffusion(\omega)} \|v(s)\|_\alpha^k\|u(s) - v(s)\|_\alpha, \quad \Forall s \le 0.
\eea

Now, we apply \eqref{G diff} and \eqref{G est. 3-1}--\eqref{G est. 3-2} together with the partial-dichotomy estimates given in \eqref{eq:dichotomy} to estimate \eqref{h diff}. For this purpose, let us choose $\eta_1$ and $\eta_2$ {\HLL as in} the proof of Proposition ~\ref{prop:cm} with $2k$ in place of $r$. In particular, 
\bea  \label{gap 4-3}
0 > \eta_{\c}^{\ast} > \eta_1 > \eta_2 > \eta_{\s}^{\ast},
\eea
{\HL recalling that $\eta_{\c}^\ast =\eta_{\c}(2k)$ and $\eta_{\s}^{\ast}=\eta_{\s}(2k)$.}

{\HLL The bounds in \eqref{G est. 3-1}--\eqref{G est. 3-2} allow us to control the terms in \eqref{G diff} which in turn allow us to control \eqref{h diff} after application of} the partial-dichotomy estimate \eqref{d-II}:
\bea \label{I}
 \| h_{\lambda,\r}(\xii, \omega)  & - h^\mathrm{app}_{\lambda,\r}(\xii, \omega) \|_\alpha \le K C \int_{-\infty}^0  \frac{e^{-\eta_2 s  +
\int_s^0 z_\diffusion(\theta_\tau \omega)\, \mathrm{d}\tau + k z_\diffusion(\theta_s \omega)}}{|s|^\alpha} \|u(s)\|^{k+1}_\alpha \d s \\
& + k KC \int_{-\infty}^0  \Big ( \frac{e^{ -\eta_2 s+ \int_s^0 z_\diffusion(\theta_\tau \omega)\, \mathrm{d}\tau+ (k-1)z_\diffusion(\theta_s \omega)}}
{|s|^\alpha} \\
& \hspace{6em} \times \big( \|u(s)\|_\alpha^{k-1} + \|v(s)\|_\alpha^{k-1}\big) \|u(s) - v(s)\|_\alpha \Big)  \d s \\
& + KC \int_{-\infty}^0  \frac{e^{ -\eta_2 s+ \int_s^0 z_\diffusion(\theta_\tau \omega)\, \mathrm{d}\tau+ k z_\diffusion(\theta_s \omega)}}
{|s|^\alpha} \|v(s)\|_\alpha^{k} \|u(s) - v(s)\|_\alpha  \d s \\
& =: I_1(\xii,\omega) + I_2(\xii,\omega) + I_3(\xii,\omega).
\eea

To estimate the integral terms above, we first derive in the next step some estimates of $\|u(s)\|_\alpha$, $\|v(s)\|_\alpha$, and $\|u(s) - v(s)\|_\alpha$ for all $s \le 0$.

\medskip
{\bf Step 3. Estimates {\HL about} the stochastic flows $u(s)$ and $v(s)$.} We show that for each $s \le 0$, $\xii \in \mathcal{H}^{\c}$, and $\omega\in \Omega$, the following estimates hold:
\begin{align} 
& \|u(s)\|_\alpha \le 2K e^{\eta s - \int_s^0 z_{\diffusion}(\theta_{\tau} \omega)\,\mathrm{d} \tau} \|\xii\|_\alpha,  \label{u est.} \\ 
& \|v(s)\|_\alpha  \le K e^{\eta_1 s - \int_s^0 z_\diffusion(\theta_\tau \omega)\, \mathrm{d}\tau} \|\xii\|_\alpha, \label{v est.} \\
& \|u(s) - v(s)\|_\alpha \le  (2K)^k\mathfrak{C}_1(\omega) e^{(k \eta - (k-1)\epsilon_1) s - \int_s^0 z_\diffusion(\theta_\tau
\omega) \,\mathrm{d} \tau} \|\xii\|_\alpha^k,  \label{u-v est.}
\end{align}
where $\mathfrak{C}_1(\omega)$ and $\epsilon_1$ will be defined below in \eqref{C1} and  \eqref{epsilon1}, respectively.

Recall that $u(\cdot) = u_{\lambda, \r}(\cdot,  \omega; \xii + h_{\lambda,\r}(\xii, \omega))$ {\HLL is the mild solution of Eq. ~\eqref{REE local version} defined on $(-\infty, 0]$ taking value $\xi + h_{\lambda,\r}(\xii, \omega)$ at $t=0$}, which is furthermore the unique fixed point in $C_{\eta}^{-}$ of the operator $\mathcal{N}_{\xii, \r}^{\omega, \lambda}$ defined in \eqref{VCF-1}. In particular, $P_{\c} u(0) = \xii$. It is also clear that $u(s)\equiv 0$ if $\xii=0$. {\HLL As a consequence}, by using the estimate \eqref{Lip in x} applied to the mild solution $u(s)$ (with $\xi_1 = \xi$ and $\xi_2 = 0$), we obtain
\bea \label{eta norm of u est.}
\|u(\cdot)\|_{C_\eta^-} & \le \frac{K}{1- \Upsilon_{1}(F_\r)}\|\xii\|_\alpha.
\eea

Recall from Theorem ~\ref{thm:local REE} that $\Upsilon_{1}(F_\r)$ can be made less than $\frac{1}{2}$ by taking $\r < \r^\ast$.
%
%
%
Using this in \eqref{eta norm of u est.}, we obtain for $\r < \overline{\r}$ ({\hh with $\overline{\r}$ given in \eqref{eqn: rho limited}}):
\bea  \label{eta norm of u est.-2}
\|u(\cdot)\|_{C_\eta^-} & \le 2K \|\xii\|_\alpha.
\eea
Now, \eqref{u est.} follows directly from \eqref{eta norm of u est.-2} and the definition of the $\|\cdot\|_{C_{\eta}^-}$-norm given in  \eqref{eta norm}.

The estimate about $v(s)$ in \eqref{v est.} follows directly from the definition of $v$ in \eqref{def. u & v} and the partial-dichotomy estimate \eqref{d-III}.

The remaining part of this step is devoted to deriving \eqref{u-v est.}. Since $u(\cdot)$ is a fixed point in $C_{\eta}^{-}$ of the integral operator $\mathcal{N}_{\xii,\r}^{\omega, \lambda}$ defined in \eqref{VCF-1} and $v(\cdot) =\mathfrak{T}_{\lambda, \sigma}(\cdot, 0; \omega) \xii$, we get naturally:
\bea \label{u-v}
u(s)-v(s) & =  - \int_s^ 0 \mathfrak{T}_{\lambda, \sigma}(s, s'; \omega)  P_{\c} G_{\r}(\theta_{s'}\omega, u(s'))\,\mathrm{d} s'
\\
& \hspace{2em} + \int_{-\infty}^s \mathfrak{T}_{\lambda, \sigma}(s, s'; \omega)   P_{\s} G_{\r}(\theta_{s'}\omega, u(s'))\,\mathrm{d} s', \quad s \le 0, \, \omega \in \Omega.
\eea
Applying the partial-dichotomy estimates \eqref{eq:dichotomy} to \eqref{u-v}, we obtain that
\beas 
\|u(s)-v(s)\|_\alpha & \le K  \int_s^0 e^{\eta_1 (s-s') - \int_s^{s'} z_\diffusion(\theta_\tau
\omega) \,\mathrm{d} \tau} \|P_{\c} G_{\r}(\theta_{s'}\omega, u(s'))\|\,\mathrm{d} s'
\\
& \hspace{1em} + K \int_{-\infty}^s \frac{e^{\eta_2 (s-s') + \int_{s'}^s z_\diffusion(\theta_\tau \omega)
\,\mathrm{d} \tau}}{|s-s'|^\alpha} \|P_{\s} G_{\r}(\theta_{s'}\omega, u(s'))\|\,\mathrm{d} s', \quad s \le 0, \, \omega \in \Omega.
\eeas
Using then \eqref{G est. 2-1} in the above estimate, we obtain
\bea  \label{u-v 2}
\|u(s)-v(s)\|_\alpha & \le   K C \int_s^0 e^{\eta_1 (s-s') - \int_s^{s'} z_\diffusion(\theta_\tau
\omega) \,\mathrm{d} \tau + (k-1)z_\diffusion(\theta_{s'}
\omega)} \|u(s')\|_\alpha^k \,\mathrm{d} s'
\\
& \hspace{1em} + K C \int_{-\infty}^s \frac{e^{\eta_2 (s-s') + \int_{s'}^s z_\diffusion(\theta_\tau \omega)
\,\mathrm{d} \tau + (k-1)z_\diffusion(\theta_{s'}
\omega)}}{|s-s'|^\alpha} \|u(s')\|_\alpha^k\,\mathrm{d} s'.
\eea

{\HLL Before pursuing the analysis, we summarize in the following lemma some basic algebraic inequalities, which will be used to ensure mainly the existence of certain  integrals arising in various places from the error estimates; see for instance \eqref{u-v 3}, \eqref{I_2 est.}, and \eqref{I_3 est.} below. }

{\mkk In that respect, let us recall that $\eta_1$ and $\eta_2$  chosen from the proof of Proposition ~\ref{prop:cm} with $2k$ in place of $r$
  satisfy furthermore 
\bea  \label{gap 4-3-2}
 0 > 2k \eta_{1} > \eta_{2}, \quad \text{and} \quad \eta \in \Bigl( \frac{\eta_2}{2k}, \eta_1 \Bigr).
\eea
}

\bl  \label{lem:algebra}

Let $\epsilon_1$ be the positive constant defined by:
\be \label{epsilon1}
\epsilon_1 := \frac{2k \eta - \eta_2}{2(2k-1)}.
\ee
Then, the following set of inequalities hold:
\bea  \label{eta relation}
 \eta_2 & < 2k \eta - (2k-1)\epsilon_1 < (2k-1) \eta - 2(k-1) \epsilon_1 \\
&  \le (k+1) \eta - k \epsilon_1 < k\eta - (k-1)\epsilon_1< \eta_1.
\eea

\el

\bp

These inequalities can be checked directly using \eqref{gap 4-3}, \eqref{gap 4-3-2}, \eqref{epsilon1}, and the fact that $k\ge 2$. For instance, the first inequality in \eqref{eta relation} can be established as follows by using \eqref{epsilon1} and $2k\eta > \eta_2$: 
\bes
2k \eta - (2k-1)\epsilon_1  = 2k\eta - \frac{2k\eta - \eta_2}{2} = \frac{2k\eta + \eta_2}{2} > \eta_2.
\ees 

To obtain the last inequality in \eqref{eta relation}, we note that 
\beas
\eta_1 - (k\eta  - (k-1)\epsilon_1) = (\eta_1 - \eta) + (k-1)(\epsilon_1  - \eta).
\eeas
The desired result follows since $0 > \eta_1 > \eta$, $k-1 > 0$, and $\epsilon_1 > 0$.

\ep

For $\epsilon_1>0$ given in \eqref{epsilon1}, let us introduce the following random variable:
\be \label{random const}
C_{\epsilon_1}(\omega) := \sup_{s'\le 0}  e^{\epsilon_1 s' + |\int_{s'}^0 z_{\diffusion}(\theta_{\tau} \omega)\,\mathrm{d} \tau| +|z_{\diffusion}(\theta_{s'} \omega)|}.
\ee
Note that $C_{\epsilon_1}(\omega)$ is finite for each $\omega \in \Omega$ thanks to the growth control estimates \eqref{Eq:z-4} of  $z_\diffusion$.

With the help of Lemma ~\ref{lem:algebra}, the first integral on the RHS of \eqref{u-v 2} can be then controlled as follows:
\bea \label{u-v 3}
\int_s^0 & e^{\eta_1 (s-s') - \int_s^{s'} z_\diffusion(\theta_\tau
\omega) \,\mathrm{d} \tau + (k-1)z_\diffusion(\theta_{s'}
\omega)} \|u(s')\|_\alpha^k \,\mathrm{d} s' \\
& \le \int_s^0 e^{\eta_1 (s-s') + k \eta s' - \int_s^0 z_\diffusion(\theta_\tau
\omega) \,\mathrm{d} \tau  + (k-1)[z_\diffusion(\theta_{s'}
\omega) + \int_0^{s'} z_\diffusion(\theta_\tau
\omega) \,\mathrm{d} \tau]} \,\mathrm{d} s' \|u(\cdot)\|_{C_\eta^-}^k \\
& \le C_{\epsilon_1}^{k-1}(\omega) e^{-\int_s^0 z_\diffusion(\theta_\tau
\omega) \,\mathrm{d} \tau} \int_s^0 e^{\eta_1 (s-s') + (k \eta - (k-1)\epsilon_1) s'
} \,\mathrm{d} s' \|u(\cdot)\|_{C_\eta^-}^k \\
& \le \frac{e^{(k \eta - (k-1)\epsilon_1) s - \int_s^0 z_\diffusion(\theta_\tau
\omega) \,\mathrm{d} \tau} C_{\epsilon_1}^{k-1}(\omega)}{\eta_1 - k \eta + (k-1)\epsilon_1} \|u(\cdot)\|_{C_\eta^-}^k, \quad \Forall s \le 0, 
\eea
{\HLL where the last integral above is well-defined since} $\eta_1 - k \eta + (k-1)\epsilon_1>0$ thanks to \eqref{eta relation}.
Similarly,
\bea \label{u-v 4}
& \int_{-\infty}^s \frac{e^{\eta_2 (s-s') + \int_{s'}^s z_\diffusion(\theta_\tau \omega)
\,\mathrm{d} \tau + (k-1)z_\diffusion(\theta_{s'}
\omega)}}{|s-s'|^\alpha} \|u(s')\|_\alpha^k\,\mathrm{d} s'  \\
& \hspace{3em} \le \frac{e^{(k \eta - (k-1)\epsilon_1) s - \int_s^0 z_\diffusion(\theta_\tau
\omega) \,\mathrm{d} \tau} C_{\epsilon_1}^{k-1}(\omega) \Gamma(1-\alpha)}{(k \eta - \eta_2 - (k-1)\epsilon_1)^{1-\alpha}} \|u(\cdot)\|_{C_\eta^-}^k, \quad \Forall s \le 0. 
\eea

By using the controls obtained in \eqref{u-v 3} and \eqref{u-v 4}, we deduce from \eqref{u-v 2} that
\bea \label{u-v 5}
\|u(s)-v(s)\|_\alpha & \le \mathfrak{C}_1(\omega) e^{(k \eta - (k-1)\epsilon_1)s - \int_s^0 z_\diffusion(\theta_\tau
\omega) \,\mathrm{d} \tau}  \|u(\cdot)\|_{C_\eta^-}^k, \; \Forall s \le 0, \; \omega \in \Omega,
\eea
where
\be \label{C1}
\mathfrak{C}_1(\omega) := K C C_{\epsilon_1}^{k-1}(\omega) \Big( \frac{1}{\eta_1 - k \eta +  (k-1)\epsilon_1} + \frac{ \Gamma(1-\alpha)}{(k\eta - \eta_2 - (k-1)\epsilon_1)^{1 - \alpha}}\Big).
\ee
Then, \eqref{u-v est.} follows from \eqref{u-v 5} by controlling $\|u(\cdot)\|_{C_\eta^-}$ using \eqref{eta norm of u est.-2}.

\medskip
{\bf Step 4. Estimates of {\HL $\|h_{\lambda,\r}(\xii, \omega) - h^\mathrm{app}_{\lambda,{\r}}(\xii, \omega)\|_\alpha$}.} Now we are ready to estimate $I_1(\xii, \omega)$, $I_2(\xii, \omega)$, and $I_3(\xii, \omega)$ as given in \eqref{I} in order to derive the estimate \eqref{diff h aim 1} at the end of this step. 

Let us begin {\HLL by estimating} $I_1(\xi, \omega)$. Using \eqref{u est.} in $I_1(\xii,\omega)$, we obtain
\bea \label{I_1 est. 1}
I_1(\xii, \omega) & =  K C \int_{-\infty}^0  \frac{e^{-\eta_2 s  +
\int_s^0 z_\diffusion(\theta_\tau \omega)\, \mathrm{d}\tau + k z_\diffusion(\theta_s\omega)}}{|s|^\alpha} \|u(s)\|^{k+1}_\alpha \d s \\
& \le 2^{k+1}K^{k+2} C \int_{-\infty}^0  \frac{e^{[(k+1)\eta -\eta_2] s - k\int_s^0 z_\diffusion(\theta_{\tau}\omega) \d \tau + k z_\diffusion(\theta_s\omega)}}{|s|^\alpha} \d s \|\xii\|_\alpha^{k+1} \\
& \le 2^{k+1} K^{k+2} C C_{\epsilon_1}^k(\omega) \int_{-\infty}^0  \frac{e^{[(k+1)\eta -\eta_2 - k \epsilon_1] s}}{|s|^\alpha} \d s \|\xii\|_\alpha^{k+1} \\
& \le 2^{k+1} K^{k+2} C C_{\epsilon_1}^k(\omega) \Gamma(1-\alpha) \big( (k+1) \eta - \eta_2 -k \epsilon_1\big)^{\alpha - 1} \|\xii\|_\alpha^{k+1} \\
& =:\mathfrak{C}_2(\omega)\|\xii\|_\alpha^{k+1}, \qquad \Forall \xii \in \mathcal{H}^{\c}, \, \omega \in \Omega,
\eea
{\HLL where the last integral above is well-defined since 
$(k+1) \eta - \eta_2 -k \epsilon_1>0$ thanks to \eqref{eta relation}.}

For $I_2(\xii, \omega)$, first note that since $\eta < \eta_1 < 0$, by {\hh using the controls of $\|u(s)\|_\alpha$ and $\|v(s)\|_\alpha$ given respectively in} \eqref{u est.} and \eqref{v est.} we obtain
\bea \label{u+v est.}
\|u(s)\|_\alpha^{k-1} + \|v(s)\|_\alpha^{k-1} \le (2^{k-1} + 1)K^{k-1}e^{(k-1)\eta s - (k-1)\int_s^0 z_\diffusion(\theta_{\tau}\omega) \d \tau}\|\xii\|_\alpha^{k-1}, \quad \Forall \, s \le 0.
\eea

Now, $I_2(\xii,\omega)$ as defined in \eqref{I} can be controlled as follows by using \eqref{u+v est.} and \eqref{u-v est.}:
\bea \label{I_2 est.}
 I_2(\xii, \omega) & \le 2^k (2^{k-1} + 1)k K^{2k} C \mathfrak{C}_1(\omega) \\
& \hspace{2em} \times \int_{-\infty}^0  \frac{e^{[(2k-1)\eta - \eta_2 -(k-1)\epsilon_1] s  + (k-1) 
[z_\diffusion(\theta_s\omega) -\int_s^0 z_\diffusion(\theta_\tau \omega)\, \mathrm{d}\tau]}}{|s|^\alpha} \d s \|\xii\|_\alpha^{2k-1}\\
& \le 2^k (2^{k-1} + 1)k K^{2k} C \mathfrak{C}_1(\omega) C_{\epsilon_1}^{k-1}(\omega)\int_{-\infty}^0  \frac{e^{[(2k-1)\eta - \eta_2 -2(k-1)\epsilon_1] s }}{|s|^\alpha} \d s \|\xii\|_\alpha^{2k-1}\\
&= \frac{2^k (2^{k-1} + 1)k K^{2k} C \mathfrak{C}_1(\omega) C_{\epsilon_1}^{k-1}(\omega) \Gamma(1-\alpha)}{\big( (2k-1) \eta - \eta_2 -2(k-1) \epsilon_1\big)^{1-\alpha}}  \|\xii\|_\alpha^{2k-1}\\
& =: \mathfrak{C}_3(\omega)\|\xii\|_\alpha^{2k-1}, \qquad \Forall \xii \in \mathcal{H}^{\c}, \, \omega \in \Omega,
\eea
{\HLL where the last integral above is well-defined since $(2k-1) \eta - \eta_2 -2(k-1) \epsilon_1 > 0$ thanks to \eqref{eta relation}.}

By the same type of estimates, we obtain
\bea \label{I_3 est.}
I_3(\xii,\omega) & = KC \int_{-\infty}^0  \frac{e^{ -\eta_2 s+ \int_s^0 z_\diffusion(\theta_\tau \omega)\, \mathrm{d}\tau+ k z_\diffusion(\theta_s \omega)}}
{|s|^\alpha} \|v(s)\|_\alpha^{k} \|u(s) - v(s)\|_\alpha  \d s \\
& \le \frac{2^k K^{2k+1} C \mathfrak{C}_1(\omega) C_{\epsilon_1}^{k}(\omega) \Gamma(1-\alpha)}{\big( k \eta_1 + k \eta - \eta_2 - (2k-1) \epsilon_1\big)^{1-\alpha}}  \|\xii\|_\alpha^{2k}\\
& =: \mathfrak{C}_4(\omega)\|\xii\|_\alpha^{2k}, \qquad \Forall \xii \in \mathcal{H}^{\c}, \, \omega \in \Omega.
\eea
Note that $\mathfrak{C}_4(\omega)$ is positive since $k \eta_1 + k \eta - \eta_2 - (2k-1) \epsilon_1 > 2 k \eta - \eta_2 - (2k-1) \epsilon_1 > 0$ from \eqref{eta relation}.

Now, using the controls of $I_1(\xi,\omega)$ in \eqref{I_1 est. 1}, $I_2(\xi,\omega)$ in \eqref{I_2 est.}, and $I_3(\xi,\omega)$ in \eqref{I_3 est.},  we obtain then from \eqref{I} that
\be \label{diff h 1 est.}
\| h_{\lambda,{\r}}(\xii, \omega) - h^\mathrm{app}_{\lambda,{\r}}(\xii, \omega) \|_\alpha \le (\mathfrak{C}_2(\omega) + \mathfrak{C}_3(\omega)\|\xii\|_\alpha^{k-2} + \mathfrak{C}_4(\omega)\|\xii\|_\alpha^{k-1} )\|\xii\|_\alpha^{k+1},
\ee
for all $\xii \in \mathcal{H}^{\c}$, $\lambda \in \Lambda_{2k}$, and all $\omega$.

Let $r^{\ast}_\epsilon$ be now the positive random variable defined by
\bea \label{radius for diff h 1}
r^{\ast}_\epsilon(\omega)  := \min \Big\{ \frac{\epsilon}{2\big( \mathfrak{C}_2(\omega) + \mathfrak{C}_3(\omega) + \mathfrak{C}_4(\omega) \big)}, \ 1 \Big \}.
\eea
Then, it follows that 
\beas
\Big( \mathfrak{C}_{2}(\omega) + \mathfrak{C}_3(\omega)\|\xii\|_\alpha^{k-2} + \mathfrak{C}_4(\omega)\|\xii\|_\alpha^{k-1}\Big) \|\xii\|_\alpha \le \frac{\epsilon}{2}, \quad \Forall \xii \in B_\alpha(0, r^{\ast}_{\epsilon}(\omega))\cap \mathcal{H}^{\c},\,  \omega \in \Omega.
\eeas
This together with \eqref{diff h 1 est.} implies that for each $\r \in (0, \overline{\r})$ with $\overline{\r}$ given in \eqref{eqn: rho limited}, the following error estimate holds:
\be 
\| h_{\lambda,{\r}}(\xii, \omega) - h^\mathrm{app}_{\lambda,{\r}}(\xii, \omega) \|_\alpha \le \frac{\epsilon}{2}\|\xii\|_\alpha^k,\, \Forall \xii \in B_\alpha(0, r^{\ast}_{\epsilon}(\omega))\cap \mathcal{H}^{\c}, \, \lambda \in \Lambda_{2k},\,  \omega \in \Omega,
\ee
and \eqref{diff h aim 1} is proved. 

\medskip

{\bf Step 5. Estimates of {\HL $\|h^\mathrm{app}_\lambda(\xii, \omega) -h^\mathrm{app}_{\lambda,{\r}}(
\xi, \omega)\|_\alpha$}.}  In this step, we get rid of the cut-off function from the expression of $h^\mathrm{app}_{\lambda,{\r}}(\xii, \omega)$ given in \eqref{int h2}, and provide a corresponding error estimate. The purpose {\HLL here} is to show that, for any $\r \in (0, \overline{\r})$ and {\HLL any $\epsilon>0$}, there exists a random open ball $B_\alpha(0, r^{\ast \ast}_{\epsilon}(\omega)) \subset \mathcal{H}_\alpha$ such that for all $\lambda \in \Lambda_{2k}$ the following estimate holds:
\be \label{diff h aim 2}
\|h^\mathrm{app}_\lambda(\xii, \omega) -h^\mathrm{app}_{\lambda,{\r}}(\xii, \omega)\|_\alpha \le \frac{\epsilon}{2}\|\xii\|_\alpha^k,\, \Forall \xii \in B_\alpha(0, r^{\ast \ast}_{\epsilon}(\omega))\cap \mathcal{H}^{\c}, \, \omega \in \Omega,
\ee
where $h^\mathrm{app}_{\lambda}$ and $h^\mathrm{app}_{\lambda,{\r}}$ are as defined in \eqref{App CMF formula 1} and \eqref{int h2}, respectively.

To do so, we first evaluate the difference $h^\mathrm{app}_\lambda -h^\mathrm{app}_{\lambda,{\r}}$. In that respect, let us  introduce 
\bea \label{G_k}
& G_{k}(\omega,\, u) := e^{-z_\diffusion(\omega)} F_k(e^{z_\diffusion(\omega)}u) = e^{(k-1)z_\diffusion(\omega)} F_k(u).
\eea
Then, by using \eqref{OU identity}, \eqref{G_k} and $k$-linear properties of $F_k$, the formula \eqref{App CMF formula 1} can be rewritten as 
\bea  \label{App CMF formula 1 rewritten}
h^\mathrm{app}_\lambda(\xii, \omega) = \int_{-\infty}^0 & \mathfrak{T}_{\lambda, \sigma}(0, s; \omega)  P_{\s} G_k(\theta_s \omega, \mathfrak{T}_{\lambda, \sigma}(s, 0; \omega) \xii) \, \mathrm{d}s.
\eea

Note that $G_{{\r}, k}(\omega, u) =  \zeta\Bigl( \frac{e^{z_\diffusion(\omega)} \|u\|_\alpha}{\r} \Bigr) G_{k}(\omega,\, u)$. Recall also that $v(s)  = \mathfrak{T}_{\lambda, \sigma}(s, 0; \omega) \xi$ {\HLL is well-defined for $s<0$ since $\xi \in \mathcal{H}^{\c}$; see \SS ~\ref{ss:T}.} It then follows from \eqref{int h2} and \eqref{App CMF formula 1 rewritten} that 
\bea \label{diff h 22}
 h^\mathrm{app}_\lambda(\xii, \omega) -h^\mathrm{app}_{\lambda,{\r}}(\xii, \omega)  = \int_{-\infty}^0 \! \mathfrak{T}_{\lambda, \sigma}(0, s; \omega)  \Bigl( 1 - \zeta\Bigl( \frac{e^{z_\diffusion(\theta_s \omega)} \|v(s)\|_\alpha}{{\r}}\Bigr) \Bigr) P_{\s} G_k(\theta_s \omega,v(s)) \, \mathrm{d}s. 
\eea
The term $1 - \zeta\Big( \frac{e^{z_\diffusion(\theta_s \omega)} \|v(s)\|_\alpha}{{\r}}\Big)$ requires a special attention in order to control $\|h^\mathrm{app}_\lambda(\xii, \omega) -h^\mathrm{app}_{\lambda,{\r}}(\xii, \omega)\|_\alpha$. From the definition of the cut-off function $\zeta$ given in \eqref{cut-off}, we have that  
\bea  \label{eq:zeta one}
\Bigl( e^{z_\diffusion(\theta_s \omega)} \|v(s)\|_\alpha \le {\r} \Bigr) \Longrightarrow \Bigl( \zeta\Big( \frac{e^{z_\diffusion(\theta_s \omega)} \|v(s)\|_\alpha}{{\r}}\Big) = 1 \Bigr). 
\eea
Note furthermore that
\bes
e^{z_\diffusion(\theta_s \omega)} \|v(s)\|_\alpha \le K e^{\eta_1 s - \int_s^0 z_\diffusion(\theta_\tau \omega)\, \mathrm{d}\tau + z_\diffusion(\theta_s \omega)} \|\xii\|_\alpha, \Forall s \le 0,
\ees
which follows directly from \eqref{v est.}. Hence,
\be \label{v est. 22}
e^{z_\diffusion(\theta_s \omega)} \|v(s)\|_\alpha \le K C_{\epsilon_1}(\omega) e^{(\eta_1 - \epsilon_1) s} \|\xii\|_\alpha,  \Forall s \le 0, 
\ee
where $C_{\epsilon_1}(\omega)$ is as defined in \eqref{random const}.
 
Recalling that $\eta_1 < 0$ thanks to \eqref{gap 4-3}, we then have $\eta_1 - {\epsilon_1} < 0$. This together with \eqref{v est. 22} implies that 
\be \label{s0}
\Bigl( e^{z_\diffusion(\theta_s \omega)} \|v(s)\|_\alpha > {\r}  \Bigr)   \Longrightarrow  \Bigl( s <  \min \left\{ 0, \; s_0(\xi) \right \} \Bigr ),
\ee
where
\be \label{s0-2}
s_0(\xi) := \frac{\log(\frac{{\r}}{K C_{\epsilon_1}(\omega)  \|\xii\|_\alpha})}{\eta_1 - {\epsilon_1}}.
\ee
 From now on, we consider $\xii$ such that  
\beas
\|\xii\|_\alpha < \frac{{\r}}{K C_{\epsilon_1}(\omega)}.
\eeas
For such $\xii$, we have $s_0(\xi) < 0$. Then, it follows from \eqref{s0} that $e^{z_\diffusion(\theta_s \omega)} \|v(s)\|_\alpha \le {\r}$ for all $s \in [s_0(\xi), 0]$ by contraposition. Hence, $\zeta\Big( \frac{e^{z_\diffusion(\theta_s \omega)} \|v(s)\|_\alpha}{{\r}}\Big) = 1$ for all $s \in [s_0(\xi), 0]$ according to \eqref{eq:zeta one}.

We obtain thus from \eqref{diff h 22} the following control of $\|h^\mathrm{app}_\lambda(\xii, \omega) -h^\mathrm{app}_{\lambda,{\r}}(\xii, \omega)\|_\alpha$:
\bea \label{diff h}
\|h^\mathrm{app}_\lambda(\xii, \omega) -h^\mathrm{app}_{\lambda,{\r}}(\xii, \omega)\|_\alpha & \le \Big \| \int_{-\infty}^{s_0(\xi)}  \mathfrak{T}_{\lambda, \sigma}(0, s; \omega)  P_{\s} G_k(\theta_s \omega,
v(s)) \, \mathrm{d}s \Big \|_\alpha \\
& \le K \int_{-\infty}^{s_0(\xi)}  \frac{ e^{-\eta_2 s  +
\int_s^0  z_\diffusion(\theta_\tau \omega)\, \mathrm{d}\tau}}{|s|^\alpha} \|P_{\s} G_k(\theta_s \omega,
v(s))\| \, \mathrm{d}s.
\eea

The continuous $k$-linear property of $F_k$ implies trivially the existence of $C>0$ such that:
\bes
\|F_k(u)\| \le C \|u\|_\alpha^k, \qquad \Forall u \in \mathcal{H}_\alpha,
\ees
leading thus to
\beas
\|G_k(\theta_s \omega, v(s))\| & = e^{-z_\diffusion(\theta_s \omega)}\|F_k( e^{z_\diffusion(\theta_s \omega)}v(s))\| \le C e^{(k-1)z_\diffusion(\theta_s \omega)}\|v(s)\|_\alpha^{k} \\
& \le C K^k e^{k \eta_1 s + (k-1)z_\diffusion(\theta_s \omega)  -
k \int_s^0 z_\diffusion(\theta_\tau \omega)\, \mathrm{d}\tau} \|\xii\|_\alpha^k, \quad \Forall s \le 0,
\eeas
where the control of $\|v(s)\|_\alpha$ given in \eqref{v est.} is used to derive the last inequality above.
 
We obtain then from \eqref{diff h}:
\bea \label{diff h 2}
\|h^\mathrm{app}_\lambda(\xii, \omega) -h^\mathrm{app}_{\lambda,{\r}}(\xii, \omega)\|_\alpha & \le C K^{k+1} \int_{-\infty}^{s_0(\xi)}  \frac{ e^{(k\eta_1-\eta_2) s  +
(k-1)[z_\diffusion(\theta_s \omega) - \int_s^0  z_\diffusion(\theta_\tau \omega)\, \mathrm{d}\tau ] } }{|s|^\alpha} \, \mathrm{d}s \|\xii\|_\alpha^k \\
& =: J(\omega) \|\xii\|_\alpha^k.
\eea
We {\HLL conclude now about} the control of $\|h^\mathrm{app}_\lambda(\xii, \omega) -h^\mathrm{app}_{\lambda,{\r}}(\xii, \omega)\|_\alpha$ by showing that the random constant $J(\omega)$ given in \eqref{diff h 2} can be actually dominated by $\epsilon/2$ provided that $\|\xi\|_\alpha$ is sufficiently small. This control relies on the fact that $s_0(\xi)$ {\HLL as defined in \eqref{s0-2}} becomes more negative as $\|\xi\|_\alpha$ gets smaller.


By noting that 
\beas
(k\eta_1-\eta_2) s  +
(k-1) \Bigl [z_\diffusion(\theta_s \omega) - \int_s^0  z_\diffusion(\theta_\tau \omega)\, \mathrm{d}\tau 
\Bigr ] & \le (k\eta_1-\eta_2) s - (k-1)\epsilon_1 s \\
& + (k-1) \Bigl (\epsilon_1 s + |z_\diffusion(\theta_s \omega)| + \Bigl | \int_s^0  z_\diffusion(\theta_\tau \omega)\, \mathrm{d}\tau\Bigr| \Bigr),
\eeas 
the constant $ C_{\epsilon_1}(\omega)$ given in \eqref{random const} appears in the control of $J(\omega)$ as follows:
\beas  
 J(\omega) & \le \frac{C K^{k+1}C^{k-1}_{\epsilon_1}(\omega)}{|s_0(\xi)|^\alpha} \int_{-\infty}^{s_0(\xi)} e^{(k\eta_1-\eta_2 - (k-1){\epsilon_1}) s } \, \mathrm{d}s  \\
& = \frac{C K^{k+1} C^{k-1}_{\epsilon_1}(\omega)}{(k\eta_1-\eta_2 - (k-1){\epsilon_1}) |s_0(\xi)|^\alpha}  e^{(k\eta_1-\eta_2 - (k-1){\epsilon_1}) s_0(\xi) },
\eeas
where $k\eta_1-\eta_2 - (k-1){\epsilon_1}> k\eta -\eta_2 - (k-1){\epsilon_1} >0$ thanks to \eqref{eta relation}.

{\HLL Now by recalling the definition of $s_0(\xi)$ in \eqref{s0-2}, we obtain then:}
\beas 
J(\omega) & \le \frac{C K^{k+1}C_{\epsilon_1}^{k-1}(\omega) |\eta_1 - {\epsilon_1}|^\alpha}{(k\eta_1-\eta_2 - (k-1){\epsilon_1})}\frac{1}{|\log(\frac{{\r}}{K C_{\epsilon_1}(\omega)\|\xii\|_\alpha})|^\alpha} \Big(\frac{{\r}}{K C_{\epsilon_1}(\omega)\|\xii\|_\alpha}\Big)^{\eta_\epsilon},
\eeas
with
\beas 
\eta_\epsilon :=  \frac{k\eta_1-\eta_2 - (k-1){\epsilon_1}}{\eta_1 - {\epsilon_1}}.
\eeas
Let us introduce furthermore  the following positive random variable
\beas 
\widetilde{r}_{\epsilon}(\omega) & := \frac{{\r}}{KC_{\epsilon_1}(\omega)} \Big(\frac{\epsilon [k\eta_1-\eta_2 - (k-1){\epsilon_1}]}{2 CK^{k+1}C_{\epsilon_1}^{k-1}(\omega) |\eta_1 - {\epsilon_1}|^\alpha}\Big)^{-\frac{1}{\eta_\epsilon}}.
\eeas
{\HLL We obtain then for any $\xi \in \mathcal{H}^{\c}$ such that $\|\xi\|_\alpha < \widetilde{r}_{\epsilon}(\omega)$:}
\beas
\Big(\frac{{\r}}{K C_{\epsilon_1}(\omega)\|\xii\|_\alpha}\Big)^{\eta_\epsilon} < \Big(\frac{{\r}}{K C_{\epsilon_1}(\omega) \widetilde{r}_{\epsilon}(\omega)} \Big)^{\eta_\epsilon} = \frac{\epsilon [k\eta_1-\eta_2 - (k-1){\epsilon_1}]}{2 CK^{k+1}C_{\epsilon_1}^{k-1}(\omega) |\eta_1 - {\epsilon_1}|^\alpha},
\eeas
{\HLL since $\eta_{\epsilon} < 0$ due to the fact that $\eta_1 < 0$ from \eqref{gap 4-3}, and $k\eta_1-\eta_2 - (k-1){\epsilon_1}>0$ from what precedes.}

As a consequence,
\beas 
& \frac{CK^{k+1}C_{\epsilon_1}^{k-1}(\omega) |\eta_1 - {\epsilon_1}|^\alpha}{(k\eta_1-\eta_2 - (k-1){\epsilon_1})} \Big(\frac{{\r}}{K C_{\epsilon_1}(\omega)\|\xii\|_\alpha}\Big)^{\eta_\epsilon} < \frac{\epsilon}{2},  \quad \Forall \xii \in B_\alpha(0, \widetilde{r}_{\epsilon}(\omega)) \cap \mathcal{H}^{\c}.
\eeas
Note also that 
\beas  
\log \Bigl(\frac{{\r}}{K C_{\epsilon_1}(\omega)\|\xii\|_\alpha}\Bigr) > 1, \quad \Forall \xii \in B_\alpha\left(0, \ \frac{{\r}}{eKC_{\epsilon_1}(\omega)}\right) \cap \mathcal{H}^{\c},
\eeas
where $e$ denotes the Euler constant.

Now, {\HLL by introducing}
\bea  \label{r double star}
r^{\ast \ast}_{\epsilon}(\omega) & := \min \biggl\{  \frac{{\r}}{eKC_{\epsilon_1}(\omega)},\ \ \widetilde{r}_{\epsilon}(\omega) \biggl\},
\eea
{\HLL we have that}
\bea  \label{final J est}
J(\omega) < \frac{\epsilon}{2}, \quad \Forall \xii \in B(0, r^{\ast \ast}_{\epsilon}(\omega)), \, \omega \in \Omega.
\eea
{\HLL The desired control \eqref{diff h aim 2} on $\|h^\mathrm{app}_\lambda(\xii, \omega) -h^\mathrm{app}_{\lambda,{\r}}(\xii, \omega)\|_\alpha$ is thus achieved from \eqref{diff h 2}.}

\medskip

{\bf Step 6. Estimates of {\HL $\|h_{\lambda}(\xii, \omega) -  h^\mathrm{app}_\lambda(\xii, \omega)\|_\alpha$}.} Let 
\bea
r_{\epsilon} := \min \Big\{
r^{\ast}_{\epsilon}, r^{\ast \ast}_{\epsilon} \Big\},
\eea
where $r^{\ast}_{\epsilon}$ and $r^{\ast \ast}_{\epsilon}$ are given by \eqref{radius for diff h 1} and \eqref{r double star}, respectively.

It follows then from \eqref{diff h aim 1} and \eqref{diff h aim 2} that for each $\r \in (0, \overline{\r})$: 
\bea  \label{diff h 3}
\| h_{\lambda,{\r}}(\xii, \omega) -  h^\mathrm{app}_\lambda(\xii, \omega)\|_\alpha \le \epsilon \|\xii\|_\alpha^k, \quad \forall \, \xii \in B_\alpha(0, r_{\epsilon}(\omega)) \cap \mathcal{H}^{\c}, \, \lambda \in \Lambda_{2k}, \, \omega \in \Omega.
\eea
The estimate in \eqref{App CMF formula est} follows then from \eqref{diff h 3} by recalling from Step ~1 that $h_\lambda = h_{\lambda, \r}$. 

{\HLL By recalling that $\widehat{h}_\lambda (\xii, \omega)= e^{z_\diffusion (\omega)}h_{\lambda}(e^{-z_\diffusion (\omega)}\xii, \omega)$ from \eqref{h and h hat}, the corresponding error estimate for $\widehat{h}^\mathrm{app}_\lambda$ given in \eqref{SEE App CMF formula est-1} follows then directly from \eqref{App CMF formula est}.} The proof is complete.

\ep

\bp[{\bf Proof of Corollary ~\ref{Cor:Taylor}}]

Since the eigenvectors of $L_\lambda$ form a Hilbert basis of $\mathcal{H}$, the approximation $h^{\mathrm{app}}_\lambda$ defined in \eqref{App CMF formula 1} can be expanded as {\hh given} in \eqref{diagonal case} with the random coefficients $h^{\mathrm{app},n}_\lambda$ determined  for all $n \ge m+1$ as follows:
\bea \label{h_n terms}
& h^{\mathrm{app},n}_\lambda(\xii, \omega) =  \langle h^\mathrm{app}_\lambda(\xii, \omega), e_n \rangle \\
&=e^{(k-1)z_\diffusion(\omega)} \int_{-\infty}^0 \Big \langle e^{\diffusion (k-1)W_s(\omega)\Id}e^{-s L_\lambda}  P_{\s} F_k\Big(e^{sL_\lambda}\xi \Big) , e_n \Big \rangle  \,\mathrm{d}s, \quad \xi \in \mathcal{H}^{\c}, \, \omega \in \Omega.
\eea
We check now that $ h^{\mathrm{app},n}_\lambda$ can be written in the form as given by \eqref{App CMF formula 2}.

First note that, for any $\xi \in \mathcal{H}^{\c}$, it can be written as $\xi = \sum_{i=1}^m \xi_i e_i$,
with $\xii_i=\langle \xii, e_i \rangle$, $i = 1, \cdots, m$. Then, 
\bea \label{linear approx.}
e^{s L_\lambda}\xii  = \sum_{i=1}^m \xii_i e^{s L_\lambda} e_i  = \sum_{i=1}^m e^{\beta_{i}(\lambda) s}\xii_i e_i, \quad \Forall\, \xii \in \mathcal{H}^{\c}, \; s \le 0.
\eea

Note also that for any $s \le 0$, $u \in \mathcal{H}$, and $n> m$, we have
\bea \label{linear approx. 2}
\bigl \langle e^{-s L_\lambda} P_{\s} u, e_n \bigr \rangle =  \bigl \langle P_{\s} u, e^{-s L_\lambda} e_n \bigr \rangle = \bigl \langle P_{\s} u, e^{-s \beta_n(\lambda)} e_n \bigr \rangle  = e^{-s \beta_n(\lambda)} \bigl \langle u, e_n \bigr \rangle,
\eea
where in the last equality above we used $\langle P_{\c} u, e_n \rangle = 0$ since $n > m$.

Now, using the identities \eqref{linear approx.} and \eqref{linear approx. 2} in \eqref{h_n terms} we obtain
\bea  \label{eq:cor taylor-1}
h^{\mathrm{app},n}_\lambda(\xii, \omega) &=e^{(k-1)z_\diffusion(\omega)} \int_{-\infty}^0 \Big \langle e^{\diffusion (k-1)W_s(\omega)\Id}e^{-s L_\lambda}  P_{\s} F_k\Big(\sum_{i=1}^m e^{\beta_{i}(\lambda) s} \xii_i e_i\Big) , e_n \Big \rangle  \,\mathrm{d}s \\
& = e^{(k-1)z_\diffusion(\omega)} \int_{-\infty}^0 \Big \langle e^{\diffusion (k-1)W_s(\omega) - s \beta_n(\lambda)} F_k\Big(\sum_{i=1}^m e^{\beta_{i}(\lambda) s} \xii_i e_i\Big) , e_n \Big \rangle  \,\mathrm{d}s.
\eea
Since $F_k$ is $k$-linear, we have
\beas
F_k\Big(\sum_{i=1}^m e^{\beta_{i}(\lambda) s} \xii_i e_i\Big) = \sum_{(i_1, \cdots, i_k )\in \mathcal{I}^k} \xii_{i_1} \cdots \xii_{i_k}F_k(e_{i_1}, \cdots, e_{i_k}) e^{ \sum_{j=1}^{k}\beta_{i_j}(\lambda)s},
\eeas
where $\mathcal{I} = \{1, \cdots, m\}$. We get then
\bea  \label{eq:cor taylor}
h^{\mathrm{app},n}_\lambda(\xii, \omega) &=  e^{(k-1)z_\diffusion(\omega)} \sum_{(i_1, \cdots, i_k )\in \mathcal{I}^k } \xii_{i_1} \cdots \xii_{i_k}\langle F_k(e_{i_1}, \cdots, e_{i_k}), e_n \rangle \\
& \hspace{6em} \times    \int_{-\infty}^0 e^{\bigl ( \sum_{j=1}^{k}\beta_{i_j}(\lambda) -
\beta_n(\lambda) \bigr) s +  \diffusion (k-1)W_s(\omega)}\, \mathrm{d}s \\
&=  e^{(k-1)z_\diffusion(\omega)} \sum_{(i_1, \cdots, i_k )\in \mathcal{I}^k } \xii_{i_1} \cdots \xii_{i_k}\langle F_k(e_{i_1}, \cdots, e_{i_k}), e_n \rangle \M_n^{i_1,\cdots,i_k}(\omega, \lambda), 
\eea
where
$$\M_n^{i_1,\cdots,i_k}(\omega, \lambda) := \int_{-\infty}^0 e^{\bigl ( \sum_{j=1}^{k}\beta_{i_j}(\lambda) -
\beta_n(\lambda) \bigr) s + \diffusion (k-1)W_s(\omega)}\, \mathrm{d}s.$$
The formula \eqref{App CMF formula 2} is now obtained. 

It is clear from \eqref{eq:cor taylor} that $h^{\mathrm{app},n}_\lambda(\xi,\cdot)$ is $\bigl( \mathcal{F}; \mathcal{B}(\mathbb{R})\bigr)$-measurable for each fixed $\xii \in \mathcal{H}^{\c}$, and $h^{\mathrm{app},n}_\lambda(\cdot,\omega)$ is a homogeneous polynomial in $\xi_1, \cdots, \xi_m$ of order $k$ for each $\omega$. Hence, $h^{\mathrm{app}}_\lambda$ is a random homogeneous polynomial of order $k$ in the sense given in Definition ~\ref{Def:random polynomial}. This together with the estimate \eqref{App CMF formula est} implies that $h^{\mathrm{app}}_\lambda$ constitutes the leading-order Taylor approximation of $h_\lambda$.

The corresponding results for $\widehat{h}^\mathrm{app}_\lambda$ given in \eqref{diagonal case SEE}--\eqref{SEE App CMF formula 2} can be derived using the relation $\widehat{h}_\lambda (\xii, \omega)= e^{z_\diffusion (\omega)}h_{\lambda}(e^{-z_\diffusion (\omega)}\xii, \omega)$; {\HLL see \eqref{h and h hat}}. The proof is now complete.

\ep

\section{Approximation of Stochastic Hyperbolic Invariant Manifolds}  \label{s:hyperbolic} 

The approximation formulas provided by Theorem ~\ref{App CMF}  and Corollary ~\ref{Cor:Taylor}  were presented in the case where the subspace $\mathcal{H}^{\c}$ contains only critical modes which lose their stability,  formulated in terms of the PES condition \eqref{PES}, as the control parameter $\lambda$ varies; see also Remark ~\ref{rmk:PES}. {\mkk In practice, it can be also of interest to} consider other situations where the subspace $\mathcal{H}^{\c}$ contain a combination of critical  modes and modes that remain stable as $\lambda$ varies in some interval $\Lambda$. It is then  natural to ask whether {\mkk the formulas provided by Theorem ~\ref{App CMF}  and Corollary ~\ref{Cor:Taylor}} still provide approximation to the leading order of the corresponding {\mkk local} stochastic invariant manifolds.\footnote{{\mkk According to Corollary ~\ref{thm:local SEE}, the latter always exist in a sufficiently small neighborhood of the origin}. }

 First, note that by adding stable modes into $\mathcal{H}^{\c}$ which remain stable as $\lambda$ varies, 
then when $\lambda>\lambda_c$ the subspace $\mathcal{H}^{\c}$ is spanned by a mixture of stable and unstable modes, making thus {\it hyperbolic} the corresponding {\mkk local} stochastic invariant manifold. This property results from the fact that $\mathcal{H}^{\c}$ --- which is tangent to the local stochastic manifold at the origin --- is then decomposed as the direct sum of the unstable subspace and the subspace spanned by the stable modes contained in $\mathcal{H}^{\c}$. We are thus concerned here with the approximation of such hyperbolic stochastic invariant manifolds, which is the content of the following corollary. Here, the PES condition required in Theorem ~\ref{App CMF} and Corollary ~\ref{Cor:Taylor} is  replaced by the condition \eqref{gap condition hyperbolic case} below.

\bc \label{cor:hyperbolic}

Consider the SPDE ~\eqref{SEE}. Assume {\mkk that all the assumptions given in \SS ~\ref{subsec:hyp} are fulfilled} except the PES condition \eqref{PES}. 
We assume furthermore that an open interval $\Lambda$ is chosen such that the uniform spectrum decomposition as given in \eqref{gap 3} holds over $\Lambda$, and  that there exist $\eta_1$ and $\eta_2$ such that
\bea \label{gap condition hyperbolic case}
\eta_{\s} < \eta_2 < \eta_1 < \eta_{\c}, \qquad  \eta_2 < 2 k \eta_1 < 0,
\eea
where $k\ge 2$ denotes the leading order of the nonlinear term $F$; see \eqref{F Taylor}. Let $\mathcal{H}^{\c}$ and  $\mathcal{H}_\alpha^{\s}$ be the corresponding subspaces associated with the uniform spectrum decomposition as defined in \eqref{Hc} and \eqref{L decomposition} with $\dim(\mathcal{H}^{\c})=m$.

Then, for each $\eta \in (\frac{\eta_2}{2k}, \eta_1)$, there exists $\r^\ast > 0$, such that for each $\r \in (0, \r^\ast)$, Eq. ~\eqref{SEE} possesses a family of local  stochastic invariant $C^1$-manifolds $\{\widehat{\mathfrak{M}}^{\mathrm{loc}}_\lambda\}_{\lambda \in \Lambda}$, where each of such manifolds is  $m$-dimensional and takes the {\mkk abstract} form given by \eqref{SEE cm local version}. 

Moreover, the corresponding local stochastic manifold function $\widehat{h}_\lambda$ is approximated\footnote{in the sense given by \eqref{SEE App CMF formula est}.} to the leading order, $k$, by $\widehat{h}^{\mathrm{app}}_\lambda$ as defined in \eqref{AF_man}, that is
\be \label{AF_manbis}
\boxed{\widehat{h}^\mathrm{app}_\lambda(\xii, \omega) = \int_{-\infty}^0  e^{ \diffusion (k-1)W_s(\omega)\Id}e^{-s L_\lambda} P_{\s} F_k(
e^{s L_\lambda}\xii) \,\mathrm{d}s, \quad \Forall \xii \in  \mathcal{H}^{\c}, \omega \in \Omega.}
\ee

If the linear operator $L_\lambda$ is furthermore assumed to be self-adjoint, then the approximation $\widehat{h}^{\mathrm{app}}_\lambda$  constitutes the leading-order Taylor approximation of $\widehat{h}_\lambda$, which can be rewritten into the form given by \eqref{diagonal case SEE}--\eqref{SEE App CMF formula 2}.

\ec 

{\mkk
\br
Note that the existence of $\rho^*$ in the above theorem  is related to the uniform spectral gap conditions \eqref{Upsilon_1 local} and \eqref{Upsilon_j local}, which are subordinated to the choice of  $\eta$ in $ (\frac{\eta_2}{2k}, \eta_1)$.
\er
}
\bp The existence of local  stochastic invariant $C^1$-manifolds follows directly from Corollary ~\ref{thm:local SEE}. The proof consists then of noting that in {\mkk the derivation of Theorem ~\ref{App CMF}, the PES condition \eqref{PES} is only used to ensure --- via Proposition ~\ref{prop:cm} ---  that  the condition $\eta_2 < 2 k \eta_1 < 0$  holds in order to apply the technical Lemma \ref{lem:algebra} for the control of certain integrals as pointed out in the description of the skeleton of the proof of Theorem ~\ref{App CMF}; see also  \eqref{gap 4-3-2}.

 The assumption \eqref{gap condition hyperbolic case} allows us  to get rid of the PES condition, and the approximation results can thus be derived by following the same lines of proofs of Theorem ~\ref{App CMF}  and Corollary ~\ref{Cor:Taylor}. }
\ep

 As just explained in the proof above, the condition $\eta_2 < 2 k \eta_1 < 0$ of Corollary ~\ref{cor:hyperbolic} is required here for purely technical reasons in order to guarantee certain integrals emerging from the estimates to converge; see {\it e.g.}~ Lemma ~\ref{lem:algebra} and Step ~4 of the proof of Theorem ~\ref{App CMF}. 

When this condition is not satisfied, it is reasonable to conjecture that the Lyapunov-Perron integral {\mkk $\mathfrak{I}_{\lambda}$  given in \eqref{LP integral}}
still gives, when it exists, the leading order approximation of $\mathfrak{M}^{\mathrm{loc}}_\lambda$. {\mkk In the general case we saw in \SS ~\ref{subsec:critical} that the condition $\eta_2 <  k \eta_1 < 0$ was sufficient to ensure the existence of $\mathfrak{I}_{\lambda}$ when  $\lambda$ varies in some interval $\Lambda$, so that the conclusions of Corollary ~\ref{cor:hyperbolic} should still hold under this weaker condition}. 
{\mkk In the case where $L_\lambda$ is self-adjoint, a necessary  and sufficient condition for $\mathfrak{I}_{\lambda}$ to exist can be formulated as a {\it non-resonance condition} given in \eqref{NR} below.}

{\mkk As an illustration,}  we will see in Sections ~\ref{s:Burgers} and \ref{s:Fly} in the case of a stochastic Burgers-type equation, that {\mkk when the  \eqref{NR}-condition is satisfied,} the formula \eqref{diagonal case SEE} {\mkk given in Corollary ~\ref{Cor:Taylor}} {\mkk provides still an efficient tool to derive reduced equations for  the amplitudes  of the modes in $\mathcal{H}^{\c}$}, {\mkk in the case where stable modes are included in $\mathcal{H}^{\c}$ while the condition $\eta_2 < 2 k \eta_1 < 0$ is violated.}

The aforementioned non-resonance condition for the self-adjoint case {\mkk can be described as follows}. If the leading order nonlinear interactions between the low modes in $\mathcal{H}^{\c}$ when projected against a given high mode is not zero, then the corresponding eigenvalues associated with these low modes and the given high mode should satisfy {\mkk the following  weak form of non-resonance}:\footnote{See \cite[Sect. ~22 A]{Arnold83} for {\mkk a more standard definition of the notion of non-resonance}.}

\begin{equation}  \label{NR} \tag{NR}
\boxed{
\begin{aligned} 
& \Forall \, (i_1, \cdots, i_k ) \in \mathcal{I}^k,  \ n > m, \, \lambda \in \Lambda,  \text{ it holds that} \\
 & \Bigl (\langle F_k(e_{i_1}, \cdots, e_{i_k}), e_n \rangle \neq 0 \Bigr) \Longrightarrow \biggl ( \sum_{j=1}^{k}\beta_{i_j}(\lambda) - \beta_n(\lambda) > 0 \biggr),  
\end{aligned}
}
\end{equation}
where $m$ is the dimension of  $\mathcal{H}^{\c}$ spanned by the first $m$ eigenvectors of $L_\lambda$, {\it i.e.}:
\bes
\mathcal{H}^{\c} :=\mathrm{span}\{e_1, \cdots, e_m\},
\ees
{\mkk and where $\mathcal{I}=\{1,\cdots,m\}$, and $\langle \cdot, \cdot \rangle$ denotes  the inner-product in the ambient Hilbert space $\mathcal{H}$.} 
{\mkk As for Corollary  ~\ref{Cor:Taylor}, each $\beta_{i_j}(\lambda)$ denotes the eigenvalue associated with the corresponding mode $e_{i_j}$ in $\mathcal{H}^{\c}$, and $\beta_n(\lambda)$ denotes the eigenvalue associated with the mode $e_n$ in $\mathcal{H}^{\s}$.} In fact, the \eqref{NR}-condition is needed for the corresponding $\M_n^{i_1,\cdots,i_k}$-terms,  {\mkk given in \eqref{Mn_gen} of Corollary ~\ref{Cor:Taylor}}, to be finite so that each $\widehat{h}^{\mathrm{app},n}_\lambda$ defined by \eqref{SEE App CMF formula 2} exists.

\br
We mention that cross-non-resonances similar to \eqref{NR} arise in the power series expansion of invariant manifolds and normal forms of deterministic (and finite-dimensional) vector fields near an equilibrium; see {\it e.g.} ~\cite{Har08}. 
\er

We will see  in the applications that when the  \eqref{NR}-condition is met  the formula \eqref{diagonal case SEE}  {\mkk characterizes  in certain cases a (hyperbolic) manifold which does not necessarily approximate an invariant manifold, but still conveys very useful dynamical  information}; see Sections ~\ref{s:Burgers} and \ref{s:Fly}. This is formulated below via the concept of {\it parameterizing manifold} introduced in Section \ref{ss:PM}. As we will see, the construction of such manifolds  can be achieved in terms of {\it pullback characterization} of these objects that we describe in the next part of this article.

\newpage
\part*{\large \centerline{Part II. Stochastic Reduction Based on Parameterizing Manifolds:} \\ \centerline{Theory and Applications}\hspace{-1cm}}}
\vspace{1cm}

\section{Pullback Characterization of Approximating, and Parameterizing Manifolds} \label{s:pullback}


In this section 
we show that the {\it stochastic approximating manifold}\footnote{Such a manifold in our terminology is aimed to approximate some targeted stochastic invariant manifold (that can be local). It should not be confused with the notion of {\it stochastic approximate manifold} \cite{Chu95} which makes sense even when no stochastic invariant manifold is guaranteed to exist.} of a critical manifold as obtained in Theorem ~\ref{App CMF} can be interpreted as  the pullback limit associated with an auxiliary backward-forward system introduced in \eqref{LLL} below;  see Lemma \ref{lem:pullback}.  In terms of a Fourier series representation of the solution to this system, the key idea consists of representing  the modes with high wave numbers  (living in the range of $\widehat{h}^{\mathrm{app}}_\lambda$) as a pullback limit depending on the time-history of the modes with low wave numbers; see also \cite{EMS01}. Such an idea has also been used to construct stochastic inertial (or simply invariant) manifolds in \cite{DPD96,CLR01}, but to the best of the authors knowledge {\sw it} has not been explored for the construction of stochastic approximating manifolds.

  In  \SS ~\ref{ss:PM},  this idea is extended to construct  more general  manifolds: the parameterizing manifolds such as introduced in \SS ~\ref{ss:PM_def}. As shown in Theorems \ref{thm:hG_PM} and \ref{thm:h1_PM}, and as illustrated in Sections \ref{s:reduction}, \ref{s:Burgers}  and \ref{s:Fly},  several  auxiliary backward-forward systems can be designed to get efficiently access 
 {\mkk --- still via functionals of time-history of their low modes ---} 
  to useful approximate parameterizations $h^{\mathrm{pm}}$ of the unresolved variables  $u_{\s}(t,\omega)$ in terms of the  resolved ones $u_{\c}(t,\omega)$, for a given realization $\omega$.\footnote{See also Corollary \ref{Cor_h1_is_PM} for more precise rigorous results in the deterministic context.}  Such parameterizations lead to a mean squared error, $\int_0^T \bigl \|u_{\s}(t,\omega) -  h^{\mathrm{pm}}(u_{\c}(t,\omega),\theta_t \omega) \bigr \|_\alpha^2 \, \d t $, smaller than the variance of $u_{\s}$,  $\int_0^T  \|u_{\s}(t,\omega)\|_\alpha^2 \; \d t$, whenever $T$ is sufficiently large; see Definition ~\ref{def:PM}.

Such parameterizing manifolds will turn out to be useful in the design of efficient reduced models for the amplitudes of the low (resolved) modes of an SPDE solution, even when these amplitudes are large; see Sections \ref{s:reduction}, \ref{s:Burgers}  and \ref{s:Fly}.}

\subsection{Pullback characterization of approximating manifolds}\label{Sec_app_man+PB}
For each $T>0$, let us introduce the following {\it backward-forward system}:
\begin{subequations}\label{LLL}
\begin{empheq}[box=\fbox]{align}
& \mathrm{d} \widehat{u}^{(1)}_{\c} =  L_\lambda^{\c} \widehat{u}^{(1)}_{\c}\, \mathrm{d} s + \sigma \widehat{u}^{(1)}_{\c} \circ \mathrm{d} W_s,  
 && s \in[ -T, 0], \label{Eq:spde-I} \\
& \mathrm{d} \widehat{u}^{(1)}_{\s} = \bigl( L_\lambda^{\s} \widehat{u}^{(1)}_{\s}  +  P_{\s} F_k(\widehat{u}^{(1)}_{\c}(s-T, \omega)) \bigr) \mathrm{d} s + \sigma \widehat{u}_{\s}^{(1)} \circ \mathrm{d} W_{s-T} ,
\label{Eq:spde-II}
&&  
s \in [0, T],  \\
& \textrm{with } \widehat{u}^{(1)}_{\c}(s, \omega)\vert_{s=0} = \xii, \textrm{ and } \widehat{u}^{(1)}_{\s}(s, \theta_{-T}\omega)\vert_{s=0}= 0.
\end{empheq}
\end{subequations}
where $\xii \in \mathcal{H}^{\c}$, $L_\lambda^{\c}$ and $L_\lambda^{\s}$ are respectively the projections of $L_\lambda$ to the critical and non-critical subspaces  as defined in \eqref{Lc Ls}, and $F_k$ is the $k$-linear operator given in \eqref{F Taylor}.

In the system above, the  initial value of $\widehat{u}_{\c}^{(1)}$ is prescribed in fiber $\omega$, and the {\HL initial} value of $\widehat{u}_{\s}^{(1)}$ in fiber $\theta_{-T}\omega$. The solution of this system is obtained by using a {\it backward-forward integration procedure} {\mkk made possible due to the partial coupling  present in \eqref{LLL} where $\widehat{u}^{(1)}_{\c}$ forces the evolution equation of $\widehat{u}^{(1)}_{\s}$ but not reciprocally}. 

More precisely, we first integrate \eqref{Eq:spde-I} from fiber $\omega$ backward up to fiber $\theta_{-T}\omega$ to obtain $\widehat{u}^{(1)}_{\c}(s, \omega; \xii)$ for $s \in [-T, 0]$. The second equation is then integrated forward for $s\geq0$ as follows. Since the initial datum for this equation is specified at $s=0$ in the fiber $\theta_{-T}\omega$, the random forcing $\widehat{u}_{\c}^{(1)}$ as well as the noise terms  have to start --- as $s$ evolves --- from the same fiber, explaining the $(s-T)$-dependence of $\widehat{u}_{\c}^{(1)}$ and of $\d W$ which appear in Eq.  \eqref{Eq:spde-II} by noting that $\d W_{s-T} (\omega)=\d W_s (\theta_{-T} \omega) $ for $s\in[0,T].$ It is then important to remark that  the process $\widehat{u}_{\s}^{(1)}$ thus obtained depends on $\xi$ via $\widehat{u}_{\c}^{(1)}$ which emanates (backward) from $\xi$. For that reason, we will emphasize this dependence as $\widehat{u}_{\s}^{(1)}[\xi]$ in what follows.

\bprop \label{lem:pullback}
Let $\widehat{h}^{\mathrm{app}}_\lambda$ given by  Theorem ~\ref{App CMF} be the  approximation of the stochastic critical manifold $\widehat{\mathfrak{M}}_{\lambda}^{\mathrm{crit}}$. 
Then $\widehat{h}^{\mathrm{app}}_\lambda$  can be characterized as the pullback limit of the process $\widehat{u}_{\s}^{(1)}[\xii]$ associated with \eqref{Eq:spde-II}, {\it i.e.}:
\be  \label{Eq_PBA}\tag{CL}
\boxed{\widehat{h}^{\mathrm{app}}_\lambda(\xii, \omega) = \lim_{T \rightarrow +\infty} \widehat{u}_{\s}^{(1)}[\xii](T, \theta_{-T}\omega; 0), \quad \Forall \, \xii \in \mathcal{H}^{\c}, \; \omega \in \Omega,}
\ee
where the limit above is taken in the topology of $\mathcal{H}_{\alpha}^{\s}$.

\eprop


\br



Note that this  proposition can be related to the idea --- introduced in the context of the closure problem associated with some stochastic Navier-Stokes equations in \cite{EMS01} --- of representing the high modes as a pullback limit depending on the whole history of the low modes.

\er

\bp

Let us first introduce two processes $u_{\c}^{(1)}$ and $u_{\s}^{(1)}[\xii]$ built from $\widehat{u}_{\c}^{(1)}$ and $\widehat{u}_{\s}^{(1)}[\xi]$ as follows:
\begin{subequations} \label{Eq:transform}
\begin{align}
& u_{\c}^{(1)}(s, \omega; \xii) :=e^{-z_\diffusion(\theta_s\omega)}\widehat{u}_{\c}(s, \omega; e^{z_\diffusion(\omega)} \xii), && s \in [-T, 0], &&  \xi \in \mathcal{H}^{\c}, \label{Eq:transform-a}\\
& u_{\s}^{(1)}[\xii](s, \theta_{-T}\omega; 0) :=e^{-z_\diffusion(\theta_{s-T}\omega)}\widehat{u}_{\s}[e^{z_\diffusion(\omega)}\xii](s,\theta_{-T}\omega; 0), &&  s \in [0, T]. \label{Eq:transform-b}
\end{align}
\end{subequations}
Here, the process $u_{\c}^{(1)}$ is obtained by application of the cohomology relation \eqref{cohomology relation} to $\widehat{u}_{\c}$ solving \eqref{Eq:spde-I}. The process $ u_{\s}^{(1)}[\xii]$ is obtained from the process $\widehat{u}_{\s}[e^{z_\diffusion(\omega)} \xi]$ still by application of the cohomology relation where the origin of time has been shifted to the fiber $\theta_{-T} \omega$.

{\HL Then by using the change of variables \eqref{Eq:transform} and following similar derivations  to \SS ~\ref{sec. conjugacy},  the system \eqref{LLL} is transformed into the following backward-forward system  of abstract random differential }equations: 
\bea \label{Eq:bvp-rpde}
& \frac{\mathrm{d} u_{\c}^{(1)}}{\mathrm{d} s} = L_\lambda^{\c} u_{\c}^{(1)} + z_\sigma(\theta_s \omega) u_{\c}^{(1)}, &&  s \in [-T, 0],\\
& \frac{\mathrm{d} u_{\s}^{(1)}}{\mathrm{d} s} = L_\lambda^{\s} u_{\s}^{(1)} + z_\sigma(\theta_{s-t}\omega) u_{\s}^{(1)} + e^{(k-1)z_\diffusion(\theta_{s-T}\omega)} P_{\s} F_k(u_{\c}^{(1)}(s-T, \omega)), && s \in [0, T],
\eea
supplemented by the initial conditions:
\be
u_{\c}^{(1)}(s, \omega)\vert_{s=0} = \xii \textrm{ and } u_{\s}^{(1)}(s, \theta_{-T}\omega)\vert_{s = 0} = 0.
\ee

This system is integrated according to the same  backward-forward integration procedure as for the system \eqref{LLL}, {\it i.e.} the first equation is integrated backward from fiber $\omega$ (corresponding to $s=0$) to fiber $\theta_{-T}\omega$ (corresponding to $s=-T<0$), and the second equation is then integrated forward from fiber $\theta_{-T}\omega$ up to fiber $\omega$.
Similar to \SS ~\ref{ss:T}, we obtain that
\bea \label{Eq:uc} 
u_{\c}^{(1)}(s, \omega; \xii) & = e^{- \int_s^0 z_\sigma(\theta_\tau \omega)\, \mathrm{d}\tau \Id} e^{s L_\lambda^{\c}}  \xii,   \quad \Forall \, s \in [-T, 0] \quad \text{ for }  \quad T > 0;
\eea
{\HL and by using the variation-of-constants formula we obtain that}
\bea \label{Eq:us}
u_{\s}^{(1)}[\xii](T, \theta_{-T}\omega; 0)  =  \int_{0}^{T}  e^{(k-1)z_\diffusion(\theta_{s-T}\omega)} e^{ \int_{s}^{T} z_\sigma(\theta_{\tau-T} \omega)\, \mathrm{d}\tau  \Id} e^{(T-s) L_\lambda^{\s}} P_{\s} F_k(u_{\c}^{(1)}(s-T, \omega; \xii))\, \mathrm{d}s.    
\eea
Since $F_k$ is $k$-linear, we infer from Eq. ~\eqref{Eq:uc} that for $s \in [0, T]$:
\bea \label{Eq:Fk}
F_k(u_{\c}^{(1)}(s-T, \omega; \xii)) & = F_k\Bigl(e^{-  \int_{s-T}^0 z_\sigma(\theta_{\tau} \omega)\, \mathrm{d}\tau \Id} e^{(s-T) L_\lambda^{\c}}  \xii\Bigr) \\
&=  F_k\Bigl( e^{-\int_{s-T}^0 z_\sigma(\theta_{\tau} \omega)\, \mathrm{d}\tau } e^{(s-T) L_\lambda^{\c}} \xii\Bigr) \\
& =  e^{-k \int_{s-T}^0 z_\sigma(\theta_{\tau} \omega)\, \mathrm{d}\tau} F_k\bigl( e^{(s-T) L_\lambda^{\c}} \xii\bigr).
\eea
Note also that
\bea \label{Eq:L}
e^{ \int_{s}^{T} z_\sigma(\theta_{\tau-T} \omega)\, \mathrm{d}\tau  \Id} e^{(T-s) L_\lambda^{\s}} & = e^{\int_{s}^{T} z_\sigma(\theta_{\tau-T} \omega)\, \mathrm{d}\tau} e^{(T-s) L_\lambda^{\s}}  = e^{\int_{s-T}^{0} z_\sigma(\theta_{\tau'} \omega)\, \mathrm{d}\tau'} e^{(T-s) L_\lambda^{\s}}.
\eea

Using \eqref{Eq:Fk} and \eqref{Eq:L} in Eq. ~\eqref{Eq:us}, we obtain after simplifications that
\bea \label{Eq:us-2}
{\HLL u_{\s}^{(1)}[\xii](T, \theta_{-T}\omega;0) = \int_{0}^{T} e^{(k-1) [ z_\diffusion(\theta_{s-T}\omega) - \int_{s-T}^0 z_\sigma(\theta_{\tau} \omega)\, \mathrm{d}\tau ] } e^{(T-s) L_\lambda^{\s}} P_{\s} F_k(e^{(s-T) L_\lambda^{\c}} \xii)\, \mathrm{d}s. }
\eea

Now, by using the identity \eqref{OU identity} and introducing the variable $s':= s - T$, we can rewrite \eqref{Eq:us-2} as 
\bea \label{Eq:us-3}
u_{\s}^{(1)}[\xii](T, \theta_{-T}\omega; 0) & = \int_{-T}^{0} e^{(k-1) [ z_\diffusion(\omega) + \sigma W_{s'}(\omega) ] } e^{-s' L_\lambda^{\s}} P_{\s} F_k(e^{s' L_\lambda^{\c}} \xii)\, \mathrm{d}s'.  
\eea

By comparing \eqref{Eq:us-3} with the expression of $h^{\mathrm{app}}_\lambda(\xii, \omega)$ given in \eqref{App CMF formula 1}, we get 
\be \label{Eq:rpde-limit}
h^{\mathrm{app}}_\lambda(\xii, \omega) = \lim_{T \rightarrow +\infty} u_{\s}^{(1)}[\xii](T, \theta_{-T}\omega; 0), \quad \Forall \, \xii \in \mathcal{H}^{\c}, \, \omega \in \Omega.
\ee

By {\HL recalling} the definition of $\widehat{h}^\mathrm{app}_\lambda$ given in \eqref{AF_man}, we have  $\widehat{h}^\mathrm{app}_\lambda (\xii, \omega)= e^{z_\diffusion (\omega)}h^{\mathrm{app}}_{\lambda}(e^{-z_\diffusion (\omega)}\xii, \omega)$. This together with \eqref{Eq:rpde-limit} implies that

\beas 
\widehat{h}^{\mathrm{app}}_\lambda(\xii, \omega)  &= e^{z_\diffusion (\omega)} h^{\mathrm{app}}_{\lambda}(e^{-z_\diffusion (\omega)}\xii, \omega) = \lim_{T \rightarrow +\infty} e^{z_\diffusion (\omega)} u_{\s}^{(1)}[e^{-z_\diffusion (\omega)} \xii](T, \theta_{-T}\omega; 0), \quad \Forall \, \xii \in \mathcal{H}^{\c},\, \omega \in \Omega.
\eeas
It then follows from \eqref{Eq:transform-b} that
\bea 
\widehat{h}^{\mathrm{app}}_\lambda(\xii, \omega) & = \lim_{T \rightarrow +\infty} \widehat{u}_{\s}^{(1)}[\xii](T, \theta_{-T}\omega; 0), \quad \Forall \, \xii \in \mathcal{H}^{\c},\, \omega \in \Omega.
\eea 

The proof is now complete.

\ep

It is worthwhile to note at this stage that the deterministic approximating manifolds proposed in \cite{BW10} and \cite{CDZ11} also possess a flow interpretation which helps pointing out the similarities and differences   with the random approximating manifolds proposed in this article.

By  working  within a functional framework slightly different than the one of the present article,  it was derived in \cite{BW10, CDZ11} approximation of stochastic  manifolds associated with SPDEs of type \eqref{SEE} with a bilinear term in \cite{BW10} and nonlinearity of power type in \cite{CDZ11}. We recall below first the corresponding approximation formula of \cite{BW10}. By adapting the framework of \cite{BW10} to fit our notations, given a linear operator $L$  and {\HL the corresponding} $L$-invariant subspaces $\mathcal{H}^{\s}$ and $\mathcal{H}^{\c}$, the approximating manifold of \cite{BW10} is obtained as the graph of the following deterministic function
\bea\label{Eq_BWapp}
h_B(\xii) := (-L_{\s})^{-1}P_{\s}B(\xii, \xii), \qquad \xi\in \mathcal{H}^{\c},
\eea
where $L_{\s}$ and $P_{\s}B$ denote respectively the projections associated with $\mathcal{H}^{\s}$ of the operator $L$ and the bilinear operator $B$; see \cite[Thm. ~3]{BW10} for more details.  

By reworking the framework  of \cite{CDZ11} and following the main steps in the proof of \cite[Theorem 4.1]{CDZ11}, it can be shown that the resulting deterministic approximation formula of the stochastic unstable manifold as considered in \cite{CDZ11} is given by:
\bea\label{Eq_Duan_alapp}
h_N(\xii) :=\int_{-\infty}^0 e^{-\tau L_{\s}} P_{\s}N_p(e^{\tau L_{\c}}\xi)\d \tau, \qquad \xi\in \mathcal{H}^{\c},
\eea
where the nonlinear operator $N_p$ is defined as $[N_p(\xi)](x):=(\xi(x))^p,$ with $x$ denoting the spatial variables in the original SPDE. Here $\mathcal{H}^{\c}$ denotes the unstable subspace, $P_{\s}$ the projector onto the (center)-stable modes, and $L_{\s}=P_{\s} L$. 

Using either \eqref{Eq_BWapp} or \eqref{Eq_Duan_alapp} it can be shown that the error between the original stochastic manifold and the corresponding deterministic approximation is  controlled in the appropriate norms by $\mathcal{O}(\| \xi\|)$ terms\footnote{Using standard Landau notations. Note that there is a typo in the error bound reported in \cite[Thm. ~7]{BW10}. The correct error bound, when using their notations, should read 
$\|h(\omega, \xi) - e^{z(0)} L_s^{-1}B_s(\xi, \xi)\| \le C\big[R^2 + \big(\sigma K_2(\omega) + R^2e^{2\sigma K(\omega)} + \sigma K_2(\omega)e^{2\sigma K(\omega)}\big)\|\xi\| + e^{2\sigma K(\omega)}\|\xi\|^2\big]\|\xi\|$, after combining the estimates derived by the authors in \cite[Thm.~6]{BW10}, Eqns. ~(39), (42), and (43) thereof. }, {{\sw and holds} typically over a probability set that is never of full measure but can be large.

In comparison, for SPDEs  which would fit our functional setting with same nonlinear terms as in \cite{BW10, CDZ11},  the error between the approximating manifold and the stochastic critical manifold\footnote{not necessarily restricted to the case of stochastic unstable manifold considered in \cite{CDZ11}.} would be respectively of order $o(\|\xii\|^2_\alpha)$ and $o(\|\xii\|^{p}_\alpha)$, and would hold almost surely\footnote{using here the Landau notations considered in \SS~\ref{subsec:approx}.}; see Theorem ~\ref{App CMF}.

We are now in position to characterize,  as mentioned above, the approximating manifolds associated with \eqref{Eq_BWapp} or \eqref{Eq_Duan_alapp} as pullback limiting objects.

Regarding the approximation $h_B$ given in \eqref{Eq_BWapp}, let us introduce the following problem 
that we integrate forward with initial data taken at $-T$ for some $T>0$:
\bea  \label{eq:BW}
\frac{\d u}{\d s} = L_{\s} u + P_{\s} B(\xi,\xi), \quad \xi\in \mathcal{H}^{\c}, \; u \in \mathcal{H}^{\s}; \quad u(-T)=0.
\eea 

Then obviously 
\bea \label{eq:hB-char}
h_B(\xi)=\lim_{T\rightarrow\infty} S_{\xi}(0,-T)0,
\eea 
where $S_{\xi}(t',s')$ denotes the solution operator\footnote{Since the system is autonomous, the two-time description of the dynamics $S_{\xi}(t',s')$ reduces  obviously to that given by  a semigroup {\HL $\tilde{S}_{\xi}(t' - s')$}. We adopt however this way of writing the solution operator for the sake of unifying the different approximating manifolds discussed in this section.} associated with  $\frac{\d u}{\d s} = L_{\s} u + P_{\s} B(\xi,\xi),$ giving the state of the system at time $t'$, starting at time $s'<t'$.

The system associated with the approximation $h_N$ given in  \eqref{Eq_Duan_alapp} is reminiscent with the system \eqref{LLL}, to the exception that it involves an {\it autonomous} backward-forward system which  can be  written  as follows 
\bea \label{Eq:Duan_etal_shape}
& \frac{\mathrm{d} u_{\c}}{\mathrm{d} s} = L_\lambda^{\c} u_{\c},   && u_{\c}(0) = \xii,  && s\in [-T, 0], \\
& \frac{\mathrm{d} u_{\s}}{\mathrm{d} s} = L_\lambda^{\s} u_{\s} +  P_{\s} N_p(u_{\c} (s)), && u_{\s}(-T) = 0, && s \in [-T, 0].
\eea
Here the first equation is integrated backward from 0 to $-T$, and the second is then integrated forward from $-T$ to $0.$  If we denote by $T_{\xi}(t',s')$ the solution operator associated with the second equation using the same conventions as above, it can then be proved that
\bea \label{eq:hN-char}
h_N(\xi)=\lim_{T\rightarrow\infty} T_{\xi}(0,-T)0.
\eea
It is interesting to note that $h_N(\xi)$ corresponds to the pullback limiting object given in Lemma ~\ref{lem:pullback} for $\sigma=0$ which turns out to provide an approximating manifold  of the resulting (deterministic) PDE.
Actually, the approximation formula \eqref{Eq_Duan_alapp} may be found as a particular case of  \cite[Lemma 6.2.4]{Hen81} and \cite[Appendix A]{MW} in the context of reduction to the center manifolds for dissipative PDEs with discrete spectrum. We recall  that Theorem ~\ref{App CMF} provides error estimates for this case as well by simply setting $\sigma=0$; see Remark ~\ref{rmk finiteness of xi}.

Remarkably, the approximating manifolds associated respectively with $h_B$ given in \eqref{eq:hB-char} and $h_N$ given in \eqref{eq:hN-char} can actually be related when $N = B$. This relationship is in fact a particular case of the following lemma; see also \cite[Thm. ~3.8]{MW05}.

\bl  \label{lem:det}
 Under the conditions of Corollary ~\ref{Cor:Taylor} and assuming that
\be \label{beta-c recall}
\beta_1(\lambda) = \cdots = \beta_m(\lambda) =: \beta_\ast (\lambda), \quad \Forall \lambda \in \Lambda_{2k},
\ee
the following identity holds for the approximation formula $\widehat{h}^\mathrm{app}_\lambda$ defined in \eqref{AF_man} of the critical manifolds when $\sigma = 0$ and when $\lambda$ is sufficiently close to $\lambda_c$:
\bea
\widehat{h}^\mathrm{app}_\lambda(\xii) & := \int_{-\infty}^0  e^{-s L^{\s}_\lambda} P_{\s} F_k(
e^{s L^{\c}_\lambda}\xi) \,\mathrm{d}s \\
& = (-L^{\s}_\lambda)^{-1} P_{\s}F_k(\xi) + O(|\beta_\ast(\lambda)| \cdot \|\xi\|^k_\alpha), \quad \Forall \xii \in  \mathcal{H}^{\c}.
\eea

\el

\bp

As in the proof of Corollary ~\ref{Cor:Taylor}, we can write for $\sigma = 0$ the corresponding approximating manifold function $\widehat{h}^\mathrm{app}_\lambda(\xii)$ as:
\begin{equation} \label{eq:h expansion}
\widehat{h}^\mathrm{app}_\lambda(\xii) = \sum_{n=m+1}^\infty  \widehat{h}^{\mathrm{app},n}_\lambda(\xii) e_n,  \quad \Forall \, \xi \in \mathcal{H}^{\c},
\end{equation}
where
\beas 
& \widehat{h}^{\mathrm{app},n}_\lambda(\xii) =  \langle \widehat{h}^\mathrm{app}_\lambda(\xii, \omega), e_n \rangle, \quad \Forall n \ge m+1, \; \xii \in \mathcal{H}^{\c}.
\eeas
In this deterministic setting, by following the same arguments as in \eqref{eq:cor taylor-1}--\eqref{eq:cor taylor} {\hh and using the condition \eqref{beta-c recall}}, it can be shown that
\bea \label{eq:h-111}
\widehat{h}^\mathrm{app,n}_\lambda(\xii) &=  {\hh \sum_{(i_1, \cdots, i_k )\in \mathcal{I}^k }} \xii_{i_1} \cdots \xii_{i_k}\langle F_k(e_{i_1}, \cdots, e_{i_k}), e_n \rangle \M_n(\lambda) = \langle F_k(\xi), e_n \rangle \M_n(\lambda),
\eea
where {\hh $\mathcal{I}=\{1,\cdots,m\}$, $\xi_i =  \langle \xi, e_i \rangle$, and}
\beas
\M_n(\lambda) = \int_{-\infty}^0 e^{[k\beta_{\ast}(\lambda) -
\beta_n(\lambda)]s}\, \mathrm{d}s = \frac{1}{k\beta_{\ast}(\lambda) -
\beta_n(\lambda)}.
\eeas
Now, by taking $\lambda$ sufficiently close to $\lambda_c$, $\beta_\ast(\lambda)$ is close to zero {\HL while $\beta_n(\lambda)$ is bounded away from zero (see Remark ~\ref{rmk finiteness of xi})}. {\hh As a consequence,}
\beas
\frac{1}{k\beta_{\ast}(\lambda) -
\beta_n(\lambda)}  = - \frac{1}{\beta_n(\lambda)} + O(|\beta_\ast(\lambda)|),
\eeas
which, when used in \eqref{eq:h-111}, leads to
\bea
\widehat{h}^\mathrm{app}_\lambda(\xii) & = \sum_{n=m+1}^\infty - \frac{1}{\beta_n(\lambda)} \langle F_k(\xi), e_n \rangle e_n + O(|\beta_\ast(\lambda)| \cdot \|\xi\|^k_\alpha)\\
& = (-L^{\s}_\lambda)^{-1} P_{\s}F_k(\xi) + O(|\beta_\ast(\lambda)| \cdot \|\xi\|^k_\alpha),
\eea
where $O(|\beta_\ast(\lambda)| \cdot \|\xi\|^k_\alpha)$ results from the continuous $k$-linear property of $F_k$.

\ep

\br
 It is worth mentioning that manifolds built from $(-L^{\s}_\lambda)^{-1} P_{\s}F$ have been first  considered in \cite{FMT88} for  the theory of {\it approximate inertial manifolds} in the context of two-dimensional turbulence.\footnote{where $\lambda=0$ and $DF(0)$ is not necessarily equal to zero in the framework adopted by  \cite{FMT88}; see also \cite{TZ08} for application of such FMT manifolds (formula (1.34) therein) in the context  of one-parameter family of hyperbolic conservation laws.} 
\er

To summarize we have understood that any of the approximating manifolds derived in  \cite{BW10, CDZ11} or in the present article can be obtained as {\HL pullback limiting  objects} of {\HL auxiliary systems} which may involve a backward-forward integration procedure as  employed for system \eqref{LLL} or system \eqref{Eq:Duan_etal_shape}.

\br 
We also mention that in \cite{SDL10}, the authors studied the approximation of stochastic invariant manifold for certain type of SPDEs\footnote{with global Lipschitz nonlinearities which do not cause a loss of regularity compared to the ambient space $\mathcal{H}$.} {\it via} perturbation techniques, where the noise amplitude $\sigma$ serves as the small parameter. Formal expansion in power series of $\sigma$ in the integral equation \eqref{h eqn}  satisfied by the stochastic invariant manifold function $h$, are used to derive integral equations for the zeroth-order and the first-order  terms (in power of $\sigma$) arising in the formal expansion of $h$. In particular, the zeroth-order term matches the invariant manifold function of the corresponding deterministic PDE ({\it i.e.} for $\sigma=0$). However, no analytic formulas were derived, and the computational efforts to solve these integral equations can be quite involved; see also \cite{SKD12}.
\er

\subsection{Stochastic parameterizing manifolds} \label{ss:PM_def}
The need of considering other manifolds than the invariant (or approximating) ones 
relies in part on the fact that when the latter are not inertial, only a very specific class of initial data is concerned by such manifolds so that in particular the longtime behavior cannot be appropriately captured; see  Remark \ref{rmk: cm dynamical interpretation} (1). This need is supported by theoretical as well as practical considerations which support the idea that an inertial manifold may not exist for a broad class of dynamical systems depending on the dimension of $\mathcal{H}^{\c}$; see {\it e.g.} \cite{Chekroun_al13,CLW13, Chorin_al02,DT95, DT96,KCG13,MTV01}. 

Numerous candidates to a substitute of inertial manifolds have  been introduced in the literature in that purpose; {\mkk see {\it e.g.} \cite{FMT88,FT94,DebTem94,DT95, DT96, NTW01,Tit90}.} In all the cases, the idea was to relax the requirements of the inertial manifold theory so that  the longterm dynamics can be at  least described in some approximate sense, still by some finite-dimensional manifold. Particular efforts have been devoted to developping efficient methods to determine  in practice such manifolds. 
This has led in particular to the theory of approximate inertial manifolds (AIMs) and the so-called nonlinear Galerkin methods; see {\it e.g.} \cite{BJKT90, Chu93,DMT93, FMT88, FST89, JT94, JKT90, JRT01, MT89, NTW01, Tem97, Tit90}.  Approximate inertial manifolds have also been considered in a stochastic context \cite{Chu95}, {\mkk but contrarily to the deterministic case, very few algorithms are available to compute stochastic AIMs in practice; see {\it e.g.} \cite{Kan_al12}}. 
 
In that general spirit, we introduce here an alternative concept   to inertial manifolds formulated below as the notion of {\it parameterizing manifolds} for  stochastic dissipative evolution equations such as Eq. ~\eqref{SEE}. 
The idea underlying this concept  is motivated by the following {\mkk {\it least-square minimization problem}} associated with  the problem of parameterization of the high modes by the low modes of solutions evolving on the global random attractor, $\mathcal{A}.$\footnote{Here various notions of random attractor could be used ~\cite{CDS09}, for simplicity we can keep in mind the more standard one \cite{CF94}.} When an inertial manifold does not exist it is indeed natural to consider --- for a given $\omega \in \Omega,\; T > 0 $ and $u_0 \in \mathcal{A}(\omega)$ --- the following {\it cost functional}
\be
J(h,u_0, T;\omega):= \frac{1}{T}\int_0^T \bigl \|u (t,\omega;u_0) -  (u_{\c}(t,\omega; u_0)+h(u_{\c}(t,\omega; u_0), \theta_t \omega) \bigr \|_\alpha^2 \, \d t  
\ee
 and the corresponding (formal) multi-objective minimization problem 
\be \label{minimization pb}
\textrm{Find $h\in \mathfrak{E}$ which realizes }\underset{h \in \mathfrak{E}}\min \Big\{J(h,u_0, T;\omega) \; | \; u_0 \in \mathcal{A}(\omega), \; \omega \in \Omega,\; T > 0 \Big\}.
\ee
 where 
\bes
\mathfrak{E}:= \{h: \mathcal{H}^{\c} \times \Omega \rightarrow \mathcal{H}^{\s}_\alpha \vert   \text{ $h$ is measurable, $h(\cdot, \omega)$ is {\mkk continuous} for all $\omega \in \Omega$ and $h(0, \omega) = 0$} \},
\ees
Here $\mathcal{H}^{\c}$ denote as before the subspace spanned by the resolved modes\footnote{{\mkk Typically, the $m$ first eigenmodes of the linear part.}}, $\mathcal{H}^{\s}_\alpha$ denotes its complement in $\mathcal{H}_\alpha$, and $u_{\c}(t,\omega;u_0) =P_{\c} u(t,\omega;u_0)$, where $u$ is the solution of ~\eqref{SEE} emanating --- for a given realization $\omega$ ---  from $u_0$.

The minimization problem \eqref{minimization pb} is clearly challenging to solve in general due to its infinite dimensional nature.\footnote{Note that when a stochastic inertial manifold $\Phi$ exists  for a given $\mathcal{H}^{\c}$, then $\Phi$ provides  obviously an optimal solution to this problem {\mkk since for $u=u_{\s}+u_{\c} $ on the global attractor, we have $u_{\s} (t,\omega)=\Phi(u_{\c}(t,\omega),\theta_t \omega)$ in that case}.} 
For this reason, we {\mkk consider a weaker version of} this problem by looking instead for parameterizations $h$ that 
 improve our partial knowledge of the full solution $u$ when only the resolved modes $u_{\c}$ are known. This idea is formalized  below in Definition \ref{def:PM}. For this purpose, we assume again that a global RDS acting on an appropriate interpolated space $\mathcal{H}_\alpha$ is associated with the SPDE under consideration.  Such an RDS is furthermore assumed to possess a random absorbing ball  which absorbs deterministic bounded sets of  $\mathcal{H}_{\alpha}$. 
 In what follows we denote also by $\lambda_{c}$ the $\lambda$-value at which the trivial steady state starts losing its stability.\footnote{{\mkk The PES is here no-longer required}. }

 \bd  \label{def:PM}

A stochastic  manifold $\mathfrak{M}$ of the form
\beas
\mathfrak{M}(\omega) := \{\xii + h^{\mathrm{pm}}(\xii, \omega) \mid \xii \in \mathcal{H}^{\c}\},  \; \, \omega \in \Omega,
\eeas
with $h^{\mathrm{pm}}: \mathcal{H}^{\c} \times \Omega \rightarrow \mathcal{H}^{\s}_\alpha$ being a measurable mapping, is called a stochastic parameterizing manifold (PM) associated with the SPDE \eqref{SEE} (for some fixed $\lambda$ and $\sigma$)  if the following conditions are satisfied:

\bi
\item[(i)] For each $\omega$, $h(\cdot, \omega)$ is continuous and $h(0, \omega) = 0$. 
\item[(ii)] For any $u_0 \in \mathcal{H}_\alpha $, there exists a positive random variable $\omega \mapsto T_0(\omega; u_0)$, such that  the following inequality holds:
\bea \label{PM condition}
\int_0^T \bigl \|u_{\s}(t,\omega;u_0)  & -  h^{\mathrm{pm}}(u_{\c}(t,\omega; u_0),\theta_t \omega) \bigr \|_\alpha^2 \, \d t  \\
&   \leq \int_0^T  \|u_{\s}(t,\omega;u_0)\|_\alpha^2 \, \d t, \quad  \Forall \omega \in \Omega, \; T > T_0(\omega; u_0),
\eea
where $u_{\c}(t,\omega;u_0)$ and $u_{\s}(t,\omega;u_0)$ are the projections to respectively the subspaces $\mathcal{H}^{\c}$ and $\mathcal{H}^{\s}_\alpha$ of the solution $u(t,\omega;u_0)$ for the SPDE ~\eqref{SEE} emanating from $u_0$.

\ei

In the case where the trivial steady state is unstable {\it a strict inequality} in \eqref{PM condition} is required, and the parameterization defect of $\mathfrak{M}$ for a given realization $\omega$ and a given initial datum 
$u_0 \neq 0$ is defined as the following time-dependent ratio:
\be\label{Eq_QualPM}
Q(T, \omega; u_0) :=  \frac{\int_0^T  \bigl \|u_{\s}(t,\omega;u_0) - h^{\mathrm{pm}}(u_{\c}(t,\omega; u_0), \theta_t \omega) \bigr \|_\alpha^2 \, \d t }{ \int_0^T \|u_{\s}(t,\omega;u_0)\|_\alpha^2 \, \d t},  \quad T > T_0(\omega; u_0).
\ee

\ed 

\medskip

\br
{\mkk We note that in the case where the trivial steady state is unstable and a PM is found, $1-Q$ will refer to its parameterization quality in some forthcoming discussions.  In particular, the more $Q$ will get close to zero, the better the parameterization quality will be. }
\er

\br
We mention that the notion of stochastic parameterizing manifold introduced in Definition \ref{def:PM} should not be confused with the notion of stochastic parameterization as encountered in stochastic mode reduction techniques and other  multiscale modeling of complex nonlinear systems; see {\it e.g.} \cite{Chorin_al02, Chorin_Hald-book, GKS04,  KCG13, MTV01, Majda_etal05, Majda_al06, Majda_etal08, Stinis06} (and Remark \ref{MZ_rmk} below). The latter 
are designed for the parameterization problem  of the small scales  by the large ones for  (mainly) deterministic systems. Stochastic parameterizing manifolds as introduced above deal instead with such a problem for (exclusively) stochastic systems such as the SPDE \eqref{SEE}.  When $\sigma=0$, stochastic PMs reduce to deterministic parameterization  (in a mean square sense) of the small scales by the large ones; see \SS ~\ref{sec_PMdet} and \cite{CLW13}. 
\er

We note that condition \eqref{PM condition} in the above definition means that the parameterizing manifold function $h^{\mathrm{pm}}$  provides an approximate parameterization of the ``small'' scales by the ``large'' ones,  which lead, for any given realization $\omega$, to a  mean squared error, $\int_0^T \bigl \|u_{\s}(t,\omega) -  h^{\mathrm{pm}}(u_{\c}(t,\omega),\theta_t \omega) \bigr \|_\alpha^2 \, \d t $, smaller than the variance of $u_{\s}$,  $\int_0^T  \|u_{\s}(t,\omega)\|_\alpha^2 \; \d t$, whenever $T$ is sufficiently large.   {\mkk In statistical terms, a PM function $h^{\mathrm{pm}}$, is thus such that the fraction of variance of $u_{\s}$ unexplained by $h^{\mathrm{pm}}(u_{\c})$ is less than 1, or similarly $h^{\mathrm{pm}}$ is such that 
\be\label{Eq_PM_reduces_ERR}
\int_0^T \bigl \|u (t,\omega) -  u_{\c}(t,\omega)- h^{\mathrm{pm}}(u_{\c}(t,\omega),\theta_t \omega) \bigr \|_\alpha^2 \, \d t \leq \int_0^T \bigl \|u(t,\omega) -  u_{\c}(t,\omega) \bigr \|_\alpha^2 \, \d t,
\ee
for $T$ sufficiently large. The requirement for a manifold to be a PM consists thus of reducing, in a mean-square sense,  the  error committed on $u$ when $u_{\c}$ is known. }

A parameterizing manifold  provides thus a relaxation of the concept of stochastic inertial manifold, which allows for {\mkk (pathwise)} non-exact parameterization of the small scales by the large ones. 
{\mkk  As mentioned above, such an idea of relaxing the concept of an inertial manifold has been abundantly investigated  for PDEs (see \cite{FMT88,DebTem94,DT95, DT96,JRT01,NTW01,Tit90} for references), but has been much less explored for SPDEs \cite{Chu95,Kan_al12}.}
{\mkk  Parameterizing manifolds can be also formulated  for PDEs,  and can be seen as an alternative to the concept of an AIM, where the notion of order of an AIM \cite{DebTem94} is replaced by the notion of parameterization defect,  and the distance of a solution to an AIM is considered in a mean-square sense  leading to a control of $\int_{0}^T \mathrm{dist}(u(t;u_0), \mathfrak{M})^2  \d t$ instead of $\mathrm{dist}(u(t;u_0), \mathfrak{M})$; see  \cite{CLW13} for more details.}

Note also that the time-dependent ratio  $Q(T, \omega; u_0)$ given in \eqref{Eq_QualPM} is just one convenient way of  quantifying the quality of parameterization  achieved by a {\mkk stochastic PM  when compared with the the  error committed (in  a mean-square sense) on $u$ when $u_c$ is known}. For other choices of  ``metrics'' to measure the parameterization defect   we refer to \cite{CLW13}.
{\mkk We will see that theoretical estimates of   $Q(T, \omega; u_0)$  can be derived for certain {\it PM candidates}; see  Theorems \ref{thm:hG_PM} and \ref{thm:h1_PM} for the stochastic context, and Theorem \ref{thm:h1_PMdet} for the deterministic situation. In all the cases, $Q(T, \omega; u_0)$ can be estimated via appropriate numerical simulations; see Sections \ref{s:Burgers}  and \ref{s:Fly}. Furthermore, such notion of parameterization defect allows us to get easily access to modeling error estimates associated with reduced stochastic models based on a  given PM; see Proposition \ref{lem:PM error} below.}

\br
 It is worthwhile noting that other definitions of PMs could have been considered by requiring for instance $ \int_0^T\mathbb{E} \bigl \|u_{\s}(t,\cdot;u_0)   -  h^{\mathrm{pm}}(u_{\c}(t,\cdot; u_0),\theta_t \cdot ) \bigr \|_\alpha^2 \, \d t$ to be  dominated by $\int_0^T \mathbb{E}  \|u_{\s}(t,\cdot;u_0)\|_\alpha^2 \, \d t$ for $T$ sufficiently large, or by requiring \eqref{PM condition} to hold in probability.  These relaxed notions of PMs should be investigated elsewhere,  exploring in the present article the theory of PMs such as built on Definition \ref{def:PM}. The next proposition shows for instance that stochastic inertial manifolds provided by Corollary ~\ref{attractiveness thm SEE} are (for $\lambda>\lambda_{\c}$) stochastic PMs in the sense of such a definition.

\er

\bt\label{Prop_IMisPM}

Assume that the assumptions of Corollary ~\ref{attractiveness thm SEE} hold, leading thus to the existence of a stochastic inertial manifold $\widehat{\mathfrak{M}}_{\lambda}$ associated with the SPDE \eqref{SEE} for $\lambda$ in some interval $\Lambda$. Let us denote the corresponding inertial manifold function by
$h_\lambda: \mathcal{H}^{\c} \times \Omega \rightarrow \mathcal{H}^{\s}_\alpha$. 

Then, for any given initial datum $u_0 \in \mathcal{H}_\alpha \backslash\{0\}$ and $\lambda \in \Lambda$,  the following inequality holds:
\bea\label{Eq_PMcontrol}
\int_0^T \bigl \|u_{\s}(t,\omega;u_0) & - h_{\lambda}(u_{\c}(t,\omega; u_0),\theta_t \omega) \bigr \|_\alpha^2 \, \d t \\
& \le R(T,\omega, u_0) \int_0^T \|u_{\s}(t,\omega;u_0)\|_\alpha^2 \, \d t, \quad \Forall \omega \in \Omega, \; T > 0,
\eea
where
\bea\label{Eq_RT}
R(T,\omega, u_0) := \frac{1-e^{2(\eta + \epsilon) T}}{2|\eta + \epsilon|} \frac{\Big((1+\mathrm{Lip}(h_\lambda)) \widehat{C}_{\epsilon,\sigma}(\omega) \|P_{\s} u_0 - P_{\s} \overline{u}_0(\omega)\|_\alpha\Big)^2}{\int_0^T \|u_{\s}(t,\omega;u_0)\|_\alpha^2 \, \d t}, \quad \epsilon \in (0, |\eta|).
\eea
Recall that $\mathrm{Lip}(h_\lambda)$ is the global Lipschitz constant of $h_\lambda(\cdot, \omega)$, which admits a uniform bound independent of $\omega$ as given by \eqref{Lip h}. The positive random constant $\widehat{C}_{\epsilon,\sigma}$ is defined by:
\bea
\widehat{C}_{\epsilon,\sigma}(\omega) := \sup_{t\ge 0} \frac{K e^{-\epsilon t + \diffusion W_t(\omega)}}{1 - \Upsilon_{1}(F)},
\eea
see \eqref{eq:widehat_C}; and $\overline{u}_0(\omega)$ corresponds to the random initial datum on $\widehat{\mathfrak{M}}_{\lambda}$   such that:
\bea \label{eq:asy ineq}
& \|\overline{u}(t,\omega; \overline{u}_0(\omega)) - u(t, \omega ; u_0)\|_\alpha \le \widehat{C}_{\epsilon,\sigma}(\omega) e^{(\eta + \epsilon) t} \|P_{\s}\overline{u}_0(\omega) - P_{\s}u_0\|_\alpha, \; \Forall t \ge 0,\, \omega \in \Omega.
\eea

In the case where $\mathcal{H}^{\c}$ contains at least one critical mode, then  for $\lambda > \lambda_{\c}$  (and $\lambda \in \Lambda$)  
\bes
\lim_{T\rightarrow \infty} R(T,\omega, u_0) = 0.
\ees
This implies in particular that a stochastic inertial manifold as provided by Corollary ~\ref{attractiveness thm SEE} is always a stochastic parameterizing manifold for $\lambda > \lambda_{\c}$. 
\et

\bp
Let $u$ be the solution of   Eq. ~\eqref{SEE} emanating from $u_0 \in \mathcal{H}_\alpha$. 
The proof of Corollary ~\ref{attractiveness thm SEE} ensures the existence of  a solution $\overline{u}$ of Eq. ~\eqref{SEE} which evolves on the corresponding stochastic inertial manifold and such that \eqref{eq:asy ineq} holds.   Let us denote by $p(t,\omega)$ the projection $P_{\c} \overline{u}(t,\omega, \overline{u}_0(\omega))$.


By applying respectively $P_{\c}$ and $P_{\s}$ to $\overline{u} - u$,  since $\overline{u}(t,\omega, \overline{u}_0(\omega)) = p(t,\omega) + h_\lambda(p(t,\omega), \theta_t\omega)$, we obtain  from \eqref{eq:asy ineq} 
\be \label{eq: asy ineq proj0}
\|p(t,\omega) - u_{\c}(t, \omega, u_0)\|_\alpha \le \widehat{C}_{\epsilon,\sigma}(\omega) e^{(\eta + \epsilon) t} \|P_{\s}\overline{u}_0(\omega) - P_{\s}u_0\|_\alpha, \; \Forall t \ge 0,\, \omega \in \Omega, 
\ee
and 
\be \label{eq: asy ineq proj}
\|h_\lambda(p(t,\omega), \theta_t\omega) - u_{\s}(t, \omega, u_0)\|_\alpha \le \widehat{C}_{\epsilon,\sigma}(\omega) e^{(\eta + \epsilon) t} \|P_{\s}\overline{u}_0(\omega) - P_{\s}u_0\|_\alpha, \; \Forall t \ge 0,\, \omega \in \Omega.
\ee

By simply noting that $\| h_\lambda(p(t,\omega), \theta_t\omega) - h_\lambda(u_{\c}(t, \omega, u_0), \theta_t\omega) \|_\alpha \le \mathrm{Lip}(h_\lambda) \|p(t,\omega) -u_{\c}(t, \omega, u_0)\|_\alpha$, we deduce from \eqref{eq: asy ineq proj0} and \eqref{eq: asy ineq proj}
\bea \label{eq:error_IM}
\|u_{\s}(t, \omega, u_0) - h_\lambda(u_{\c}(t, \omega, u_0), \theta_t\omega)\|_\alpha & \le \|u_{\s}(t, \omega, u_0) - h_\lambda(p(t,\omega), \theta_t\omega) \|_\alpha \\
& + \| h_\lambda(p(t,\omega), \theta_t\omega) - h_\lambda(u_{\c}(t, \omega, u_0), \theta_t\omega) \|_\alpha \\ 
& \leq  (1+ \mathrm{Lip}(h_\lambda)) \widehat{C}_{\epsilon,\sigma}(\omega) e^{(\eta + \epsilon) t} \|P_{\s}\overline{u}_0(\omega) - P_{\s}u_0\|_\alpha.
\eea

We have thus  proved that
\bea
\int_0^T & \|u_{\s}(t, \omega, u_0) - h_\lambda(u_{\c}(t, \omega, u_0), \theta_t\omega)\|^2_\alpha \d t \\
& \le \frac{1-e^{2(\eta + \epsilon) T}}{2|\eta + \epsilon|} \Big((1+ \mathrm{Lip}(h_\lambda)) \widehat{C}_{\epsilon,\sigma}(\omega)\|P_{\s}\overline{u}_0(\omega) - P_{\s}u_0\|_\alpha \Big)^2,
\eea
and \eqref{Eq_PMcontrol} follows.

In the case where $\mathcal{H}^{\c}$ contains at least one critical mode, then  for $\lambda > \lambda_{\c}$ the trivial steady state is unstable and it can be shown that
\be\label{Eq_IM_is_PM}
\int_0^T \|u_{\s}(t,\omega;u_0)\|_\alpha^2 \, \d t \underset{T\rightarrow \infty}\longrightarrow  \infty,
\ee
holds for all $\omega$.  The proof is complete.
\ep

\br 
 In what follows we will often refer a stochastic PM to be simply a PM, the stochastic attribute of such an object  being implicit in the context of this article; only \SS ~\ref{sec_PMdet} below deals with the deterministic situation. 
\er

\subsection{Parameterizing manifolds as pullback limits} \label{ss:PM}

 We examine now the important problem of the determination of stochastic PMs in practice.   The main message is here that the approach introduced in \SS ~\ref{Sec_app_man+PB}
to provide approximating manifolds of critical manifolds, can be in fact generalized to get access to a variety of PMs in practice.

In particular, we will illustrate in applications that auxiliary backward-forward systems such as   \eqref{LLL} provide an efficient way to determine in practice PM candidates; see Sections ~\ref{s:Burgers} and  \ref{s:Fly}.  
This is supported by the rigorous results formulated below in Theorem  \ref{thm:h1_PM},  in Theorem \ref{thm:hG_PM} and in Theorem \ref{thm:h1_PMdet}.  
For instance, by taking the pullback limit of $\widehat{u}_{\s}^{(1)}[\xii]$ solving \eqref{Eq:spde-II}, it will be shown that this limit may give access to a parameterizing manifold useful for applications (Theorems  \ref{thm:h1_PM} and  \ref{thm:h1_PMdet}); see also Section ~\ref{s:Burgers}.  This pullback limit {\mkk exists and is always} given by the  Lyapunov-Perron integral $\mathfrak{I}_{\lambda}$ provided  that the condition 
\be \label{gap for LP integral-1b} 
\eta_{\s} < k \eta_{\c}< 0,
\ee
is satisfied for $\lambda \in \Lambda$. 
Such a pullback characterization of $\mathfrak{I}_{\lambda}$ follows indeed from the same arguments given in the proof of {\mkk Proposition}  ~\ref{lem:pullback}  since an inspection of the latter 
shows that only the existence of the integral $\mathfrak{I}_{\lambda}$ --- ensured by \eqref{gap for LP integral-1b}  ---  is required for the pullback characterization to hold. To summarize,  when condition \eqref{gap for LP integral-1b} is satisfied, 
\bea  \label{Eq_PBA1uc}
\widehat{h}^{(1)}_\lambda(\xii, \omega) & := \lim_{T \rightarrow +\infty} \widehat{u}_{\s}^{(1)}[\xii](T, \theta_{-T}\omega; 0) \\
&  = {\hh \int_{-\infty}^0  e^{ \diffusion (k-1)W_s(\omega)\Id}e^{-s L^{\s}_\lambda} P_{\s} F_k(
e^{s L^{\c}_\lambda}\xii) \,\mathrm{d}s, \quad \Forall \xii \in  \mathcal{H}^{\c}, \omega \in \Omega,}
\eea
provides a natural candidate to be a parameterizing manifold function, where $\widehat{u}_{\s}^{(1)}[\xii]$ solves \eqref{Eq:spde-II}.

\br\label{Rem_NR_forPM} 
Note  that in the case where $L_\lambda$ is self-adjoint,  the pullback characterization of $\mathfrak{I}_{\lambda}$ holds under the even weaker condition  \eqref{NR} as given in Section ~\ref{s:hyperbolic}.  Indeed the \eqref{NR}-condition  is necessary and sufficient for $\mathfrak{I}_{\lambda}$ to be well defined in this case.  
\er

\br
Note that $\widehat{h}^{(n)}_\lambda$ given in \eqref{h_n} as the pullback limit of $\widehat{u}^{(n)}_{\s}$ associated with the system \eqref{n_layer} is $(\mathcal{B}(\mathcal{H}^{\c})\otimes \mathcal{F}; \mathcal{B}(\mathcal{H}^{\s}_\alpha))$-measurable and is continuous in $\xi$ since $\widehat{u}^{(n)}_{\s}[\xi](T, \theta_{-T} \omega; 0)$ shares these properties.

\er

In Section \ref{s:Fly} it will be illustrated  that  other auxiliary backward-forward systems than \eqref{LLL} can be designed  to get access to  parameterizing manifolds with different parameterization qualities than those associated with ~\eqref{LLL}. In particular, we will illustrate that a parameterizing manifold with better parameterization quality than another  leads to a reduced stochastic system which is able to better capture the dynamics on the resolved low-modes.


Such an auxiliary backward-forward system turned out to be given by:
\begin{subequations} \label{LLL_new}
\begin{empheq}[box=\fbox]{align}
& \mathrm{d} \widehat{u}^{(1)}_{\c} =   L_\lambda^{\c} \widehat{u}^{(1)}_{\c}  \mathrm{d} s +\sigma \widehat{u}^{(1)}_{\c} \circ \mathrm{d} W_s,  && s \in [-T, 0],  \\ 
& \mathrm{d} \widehat{u}^{(2)}_{\c} =  \bigl( L_\lambda^{\c} \widehat{u}^{(2)}_{\c}  + P_{\c}F(  \widehat{u}^{(1)}_{\c} ) \bigr) \mathrm{d} s +\sigma \widehat{u}^{(2)}_{\c} \circ \mathrm{d} W_s,  && s \in [-T, 0],  \label{L2-b}\\
& \mathrm{d} \widehat{u}_{\s} ^{(2)}= \bigl( L_\lambda^{\s} \widehat{u}_{\s}^{(2)} + P_{\s} F(\widehat{u}^{(2)}_{\c}(s-T, \omega)) \bigr) \mathrm{d} s + \sigma \widehat{u}_{\s} ^{(2)}\circ \mathrm{d} W_{s-T} ,  &&  s \in [0, T],\\ 
&\textrm{with }  \widehat{u}^{(1)}_{\c}(s, \omega)\vert_{s=0} = \xii,  \;\; \widehat{u}^{(2)}_{\c}(s, \omega)\vert_{s=0} = \xii, \textrm{ and } \widehat{u}_{\s}^{(2)}(s, \theta_{-T}\omega)\vert_{s=0}= 0. 
\end{empheq}
\end{subequations}
where as we will see in  Section \ref{s:Fly}, $\widehat{h}^{(2)}_\lambda$ given by 
\bea \label{Eq_PBA2uc}
\widehat{h}^{(2)}_\lambda(\xi, \omega) & := \lim_{ T \rightarrow \infty} \widehat{u}^{(2)}_{\s}[\xi](T, \theta_{-T} \omega; 0) \\
& = \int_{-\infty}^{0}  e^{ - \diffusion  W_{\tau}(\omega) \Id} e^{-\tau L^{\s}_\lambda}  P_{\s} F\bigl( \widehat{u}^{(2)}_{\c}(\tau, \omega; \xi) \bigr) \,\mathrm{d}\tau, \quad \xi \in \mathcal{H}^{\c}, \, \omega \in \Omega.
\eea
provides a good candidate to be a parameterizing manifold function, when the latter limit exists\footnote{which is equivalent to the existence of the integral in \eqref{Eq_PBA2uc}.}.  In Section  ~\ref{ss:Burgers-fly},  non-resonance conditions under which this  limit exists are  exhibited  for a stochastic Burgers-type equation, in the case where the resolved modes consist of the first two eigenmodes; see \eqref{NR2}  therein.

 System \eqref{LLL_new} introduces a new type of nonlinear self-interactions  between the resolved modes (in $\mathcal{H}^{\c}$) which take into account some of nonlinear effects (through  $P_{\c}F$) that were absent  in the auxiliary system  \eqref{LLL}.
As it will be described in  Section ~\ref{Sec_Existence_PQ},   these new interactions  lead to a new type of memory terms  in the construction of $\widehat{h}^{(2)}_\lambda$ {\it via} \eqref{Eq_PBA2uc} which will help improve the parameterization quality offered by $\widehat{h}^{(1)}_\lambda$;  we refer to Section \ref{Sec_PM-based} for a description of the memory effects conveyed by   $\widehat{h}^{(1)}_\lambda$; see Lemma \ref{lem:Mn}.

 From a practical point of view,  backward-forward  systems such as \eqref{LLL} or \eqref{LLL_new}, benefit from the fact that they come with a partial-coupling between their variables, which in particular simplifies their numerical treatment  when compared with the fully-coupled backward-forward  systems considered in \cite{Kan_al12} (based on \cite{DPD96})  for the approximation problem of stochastic inertial manifolds. 
 Furthermore, we mention that our particular choice of multiplicative noise  allows us to consider --- {\it via} the cohomology approach ---  transformed versions  ($\omega$ by $\omega$) such as \eqref{Eq:bvp-rpde} so that  we do not have to deal with problems which arise  in trying to solve more general stochastic equations backward in time. 
 
As we will see in Sections  ~\ref{s:Burgers} and  ~\ref{s:Fly} in the case of a stochastic Burgers-type equation,  although  the systems \eqref{LLL} and \eqref{LLL_new} are each partially-coupled, they turn out to give access to parameterizing manifolds with quite remarkable performances when used for the modeling of the dynamics on the low modes.

Note that $\widehat{h}^{(2)}_\lambda$ can be actually written as the sum of $\widehat{h}^{(1)}_\lambda$ and some correction terms. To make it transparent, let us assume that the nonlinearity is bilinear: $F(u) = B(u,u)$. {\mkk Then, according to \eqref{Eq_PBA2uc},} we have
\bea \label{h1 correction}
 {\hspace*{-2em}} \widehat{h}^{(2)}_\lambda(\xi, \omega)  & = \int_{-\infty}^{0}  e^{ - \diffusion  W_{\tau}(\omega) \Id} e^{-\tau L^{\s}_\lambda}  P_{\s} B \bigl( \widehat{u}^{(1)}_{\c}(\tau, \omega; \xi) + C_B(\tau, \omega; \xi), \widehat{u}^{(1)}_{\c}(\tau, \omega; \xi) + C_B(\tau, \omega; \xi) \bigr) \,\mathrm{d}\tau \\
 & = \int_{-\infty}^{0}  e^{ - \diffusion  W_{\tau}(\omega) \Id} e^{-\tau L^{\s}_\lambda}  P_{\s} B \bigl( \widehat{u}^{(1)}_{\c}(\tau, \omega; \xi), \widehat{u}^{(1)}_{\c}(\tau, \omega; \xi) \bigr) \,\mathrm{d}\tau  \\
& + \int_{-\infty}^{0}  e^{ - \diffusion  W_{\tau}(\omega) \Id} e^{-\tau L^{\s}_\lambda}  P_{\s} B \bigl( \widehat{u}^{(1)}_{\c}(\tau, \omega; \xi), C_B(\tau, \omega; \xi) \bigr) \,\mathrm{d}\tau  \\
& +\int_{-\infty}^{0}  e^{ - \diffusion  W_{\tau}(\omega) \Id} e^{-\tau L^{\s}_\lambda}  P_{\s} B \bigl( C_B(\tau, \omega; \xi), \widehat{u}^{(1)}_{\c}(\tau, \omega; \xi) \bigr) \,\mathrm{d}\tau \\ 
& + \int_{-\infty}^{0}  e^{ - \diffusion  W_{\tau}(\omega) \Id} e^{-\tau L^{\s}_\lambda}  P_{\s} B \bigl( C_B(\tau, \omega; \xi),  C_B(\tau, \omega; \xi) \bigr) \,\mathrm{d}\tau,
\eea
where $C_B$ is given by
\bes
C_B(\tau, \omega; \xi) := \int_{\tau}^0 e^{\diffusion ( W_{\tau}(\omega) - W_{s'}(\omega))\Id} e^{(\tau-s') L^{\c}_\lambda} P_{\c} B \bigl( \widehat{u}^{(1)}_{\c}(s', \omega; \xi) \bigr) \,\mathrm{d}s'.
\ees

The first term $\int_{-\infty}^{0}  e^{ - \diffusion  W_{\tau}(\omega) \Id} e^{-\tau L^{\s}_\lambda}  P_{\s} B \bigl( \widehat{u}^{(1)}_{\c}(\tau, \omega; \xi), \widehat{u}^{(1)}_{\c}(\tau, \omega; \xi) \bigr) \,\mathrm{d}\tau$ in \eqref{h1 correction} above  corresponds to $\widehat{h}^{(1)}_\lambda(\xi, \omega)$, and the remaining three terms are correction terms of $\widehat{h}^{(2)}_\lambda(\xi, \omega)$ based on $\widehat{h}^{(1)}_\lambda(\xi, \omega)$ brought by the nonlinear term $P_{\c}F(  \widehat{u}^{(1)}_{\c} )$ in \eqref{L2-b}. As we will see in Section ~\ref{Sec_Existence_PQ}, these three terms bring a new type of memory terms compared with  those conveyed by $\widehat{h}^{(1)}_\lambda$.

The pullback characterization of $\widehat{h}^{(1)}_\lambda$ and $\widehat{h}^{(2)}_\lambda$ given above leads naturally to consider  the following generalization  
associated to the  $n$-layer backward-forward system ($n \ge 2$):
\begin{subequations} \label{n_layer}
\begin{empheq}[box=\fbox]{align}
& \mathrm{d} \widehat{u}^{(1)}_{\c} =   L_\lambda^{\c} \widehat{u}^{(1)}_{\c}  \mathrm{d} s +\sigma \widehat{u}^{(1)}_{\c} \circ \mathrm{d} W_s,  \hspace{17.5em} s \in [-T, 0],  \label{Ln-a}\\
& \mathrm{d} \widehat{u}^{(2)}_{\c} =  \bigl( L_\lambda^{\c} \widehat{u}^{(2)}_{\c}  + P_{\c}F(  \widehat{u}^{(1)}_{\c} ) \bigr) \mathrm{d} s +\sigma \widehat{u}^{(2)}_{\c} \circ \mathrm{d} W_s,  \hspace{11.2em} s \in [-T, 0],  \label{Ln-b}\\
& {\hspace*{4em}} \vdots \hspace{6em} \vdots \hspace{6em} \vdots  \nonumber \\
& \mathrm{d} \widehat{u}^{(n)}_{\c} =  \bigl( L_\lambda^{\c} \widehat{u}^{(n)}_{\c}  + P_{\c}F(  \widehat{u}^{(n-1)}_{\c} ) \bigr) \mathrm{d} s +\sigma \widehat{u}^{(n)}_{\c} \circ \mathrm{d} W_s,  \hspace{9.8em} s \in [-T, 0],  \label{Ln-n}\\
& \mathrm{d} \widehat{u}_{\s} ^{(n)}= \bigl( L_\lambda^{\s} \widehat{u}_{\s}^{(n)} + P_{\s} F(\widehat{u}^{(n)}_{\c}(s-T, \omega)) \bigr) \mathrm{d} s + \sigma \widehat{u}_{\s} ^{(n)}\circ \mathrm{d} W_{s-T} , \hspace{6em}  s \in [0, T],  \label{Ln-c} \\ 
&\textrm{with }  \widehat{u}^{(1)}_{\c}(s, \omega)\vert_{s=0} =  \cdots = \widehat{u}^{(n)}_{\c}(s, \omega)\vert_{s=0} = \xii, \textrm{ and } \widehat{u}_{\s}^{(n)}(s, \theta_{-T}\omega)\vert_{s=0}= 0. 
\end{empheq}
\end{subequations}

A basic  recurrence argument leads to the following representation of the solutions of this system:
\begin{subequations}  \label{soln:Ln}
\begin{align}
& \widehat{u}^{(1)}_{\c}(s, \omega; \xi)  = e^{\diffusion W_s(\omega)\Id} e^{s L^{\c}_\lambda} \xi, \\
& \widehat{u}^{(n)}_{\c}(s, \omega; \xi)  = \widehat{u}^{(1)}_{\c}(s, \omega; \xi) \\
& {\hspace*{7em}} - \int_s^0  e^{\diffusion ( W_s(\omega)  - W_{s'}(\omega))\Id} e^{(s-s') L^{\c}_\lambda} P_{\c} F\bigl( \widehat{u}^{(n-1)}_{\c}(s', \omega; \xi) \bigr) \,\mathrm{d}s', \qquad  n \ge 2,  \nonumber \\
& \widehat{u}^{(n)}_{\s}[\xi](s, \theta_{-T} \omega; 0)  = \int_{0}^{s}  e^{\diffusion (W_{s-T}(\omega) - W_{s'-T}(\omega)) \Id} e^{(s-s') L^{\s}_\lambda}  P_{\s} F\bigl( \widehat{u}^{(n)}_{\c}(s'-T, \omega; \xi) \bigr) \,\mathrm{d}s' \\
&\hspace{7.6em} = \int_{-T}^{s-T}  e^{\diffusion( W_{s-T}(\omega) - W_{\tau}(\omega)) \Id} e^{(s-T-\tau) L^{\s}_\lambda}  P_{\s} F\bigl( \widehat{u}^{(n)}_{\c}(\tau, \omega; \xi) \bigr) \,\mathrm{d}\tau. \nonumber
\end{align}
\end{subequations}
The resulting $\widehat{u}^{(n)}_{\s}$ involves thus the past history of the noise path and of $\widehat{u}^{(n)}_{\c}$.   The latter conveys a {\it ``matriochka'' of  self-interactions between the low modes},  $e_i \in \mathcal{H}^{\c}$,   through the collection  of terms
\be
\big( \langle P_{\c} F(\widehat{u}^{(p)}_{\c}), e_i \rangle \big)_{1\leq p\leq n-1; 1\leq i\leq m} \; ;
\ee
see also Section ~\ref{Sec_Existence_PQ}.
 
 The corresponding pullback limit of $\widehat{u}^{(n)}_{\s}$, assuming its existence, provides thus a new parameterization of the unresolved modes as functional of the past history of the resolved modes as modeled by $\widehat{u}^{(n)}_{\c}$:
\bea \label{h_n} 
\widehat{h}^{(n)}_\lambda(\xi, \omega) & := \lim_{ T \rightarrow \infty} \widehat{u}^{(n)}_{\s}[\xi](T, \theta_{-T} \omega; 0), \qquad \Forall \xi \in \mathcal{H}^{\c}, \; \omega \in \Omega,
\eea

As it will be illustrated in Section ~\ref{Sec_Existence_PQ} for a stochastic Burgers-type equation,   this operation leads to a {\it hierarchy of memory terms} generalizing those conveyed by $\widehat{h}^{(1)}_\lambda$. In the self-adjoint case, these memory terms  arise as random coefficients in the analytic expression of $\widehat{h}^{(n)}_\lambda$ and  come with the aforementioned  matriochka of nonlinear interactions. 

For situations where the sequence $(\widehat{u}^{(n)}_{\s})_{n\in \mathbb{N}^{*}}$ is well-defined,  the analysis of the convergence problem when $n\rightarrow \infty$ is reported for a next study; see however the numerical results of Section \ref{PM_quality_CVE}.  For the moment, we give  in the next section conditions under which  the limiting stochastic manifold (when it exists) and the one associated with $\widehat{h}^{(1)}_\lambda$, provide both stochastic  parameterizing manifolds; see Theorems  \ref{thm:hG_PM}  and  \ref{thm:h1_PM} below.

\subsection{Existence of stochastic parameterizing manifolds via backward-forward systems}\label{Sec_Exist_PMstoch}
{\mkk In this section we consider the following (formal) limit system of the backward-forward system \eqref{n_layer} as the number of layers for the resolved modes tends to infinity}:
\begin{subequations} \label{Galerkin_PB}
\begin{empheq}[box=\fbox]{align}
& \mathrm{d} \widehat{u}^{G}_{\c} =  \bigl( L_\lambda^{\c} \widehat{u}^{G}_{\c}  + P_{\c}F(  \widehat{u}^{G}_{\c} ) \bigr) \mathrm{d} s +\sigma \widehat{u}^{G}_{\c} \circ \mathrm{d} W_s,  &&  s \in [-T, 0],  \label{Galerkin_uc}\\
& \mathrm{d} \widehat{u}_{\s} ^{G}= \bigl( L_\lambda^{\s} \widehat{u}_{\s}^{G} + P_{\s} F(\widehat{u}^{G}_{\c}(s-T, \omega)) \bigr) \mathrm{d} s + \sigma \widehat{u}_{\s} ^{G}\circ \mathrm{d} W_{s-T} , &&  s \in [0, T],  \label{Galerkin_us} \\ 
&\textrm{with }  \widehat{u}^{G}_{\c}(s, \omega)\vert_{s=0} =  \xii, \textrm{ and } \widehat{u}_{\s}^{G}(s, \theta_{-T}\omega)\vert_{s=0}= 0. 
\end{empheq}
\end{subequations}
{\mkk Such  a system can be naturally obtained for $\sigma=0$, by assuming that  $\widehat{u}_{\s} ^{(n)}$ possesses a limit  (for instance in the uniform topology of $C^1((-T,0);\mathcal{H}^{\c})$) and by passing to the limit in 
\eqref{Ln-n}. In such a deterministic case the equation \eqref{Galerkin_uc}, obtained by passing to the limit,  turns out  to be the Galerkin reduced system\footnote{associated with the SPDE \eqref{SEE} and the subspace $\mathcal{H}^{\c}$.} integrated backward.}
 
{\mkk In the stochastic case, our particular choice of multiplicative noise allows us (by use of  appropriate cohomology)  to give sense of backward (pathwise) solutions of  \eqref{Galerkin_uc} when the latter exist.
By assuming that such solutions exist and obey to an appropriate growth as $T\rightarrow -\infty$,  conditions under which  the pullback limit associated with system \eqref{Galerkin_PB} 
exists and provides a PM for  SPDEs of type \eqref{SEE} (with $F(u)$ reduced to bilinear terms), can be then derived.  This is the content of the following theorem. }

\bt \label{thm:hG_PM}

{\mkk Let the SPDE \eqref{SEE} be considered within  the framework of Section \ref{ss:SEE},  where the nonlinearity is assumed to be bilinear and continuous ({\it i.e.} $F(u) = B(u,u)$), and $\lambda$ to be in some interval $\Lambda$ on which the uniform spectrum decomposition \eqref{gap 3} holds}.  

Assume {\mkk furthermore}  that {\mkk for some $\sigma>0$} the solution $\widehat{u}^{G}_{\c}(s, \omega; \xi)$ to \eqref{Galerkin_uc} emanating from $\xi$ in fiber $\omega$ exists on $[-T, 0]$ for all $\xi \in \mathcal{H}^{\c}$ and all $T > 0$. Assume also that there exists a positive random variable $\gamma_{\lambda,\sigma}(\omega)$ such that
\be  \label{uc growth control}
\lim_{T\rightarrow +\infty} e^{- \gamma_{\lambda,\sigma}(\omega) T} \|\widehat{u}^{G}_{\c}(-T, \omega; \xi)\|_\alpha  = 0, \qquad \Forall \xi \in \mathcal{H}^{\c}, \, \omega \in \Omega,\; \lambda \in \Lambda,
\ee
and that there exists $\delta > 0$, such that 
\be \label{random gap}
-2 \gamma_{\lambda,\sigma}(\omega) - \eta_2 > \delta, \qquad \omega \in \Omega,
\ee
where $\eta_2 < 0$ is chosen so that the partial-dichotomy estimates \eqref{Proj-I}--\eqref{Proj-II} hold.

Then, the following pullback limit exists:
\bea \label{hG}
 \widehat{h}^{G}_\lambda(\xi, \omega) & := \lim_{ T \rightarrow \infty} \widehat{u}^{G}_{\s}[\xi](T, \theta_{-T} \omega; 0) \\
& = \int_{-\infty}^{0}  e^{ -\diffusion W_{\tau}(\omega)) \Id} e^{-\tau L^{\s}_\lambda}  P_{\s} B \bigl( \widehat{u}^{G}_{\c}(\tau, \omega; \xi), \widehat{u}^{G}_{\c}(\tau, \omega; \xi) \bigr) \,\mathrm{d}\tau, 
\eea
for all $\xi \in \mathcal{H}^{\c}$, $\omega \in \Omega$, and $\lambda \in \Lambda$.

Assume furthermore {\mkk that a global tempered  random attractor $\mathcal{A}_\lambda$ exists for the SPDE  \eqref{SEE}}. Then, there exists a positive random variable $C_{\lambda,\sigma}(\omega)$ such that
\bea \label{error hG}
& \|u_{\s}(t,\omega; u_0) - \widehat{h}^{G}_\lambda(u_{\c}(t,\omega; u_0), \theta_t \omega) \|_\alpha \\
 & \le \frac{C_{\lambda,\sigma}(\theta_t\omega)}{(-\eta_2 - 2 \gamma_{\lambda,\sigma}(\theta_t \omega) - \frac{3\delta}{4})^{1-\alpha}}, \; \Forall u_0 \in \mathcal{A}_\lambda(\omega), \; t \ge 0, \; \lambda \in \Lambda,
\eea
where $u_{\s}(t,\omega; u_0) := P_{\s}u_\lambda(t,\omega; u_0)$ and $u_{\c}(t,\omega; u_0) := P_{\c}u_\lambda(t,\omega; u_0)$ with $u_\lambda(t,\omega; u_0)$ being the {\mkk  complete SPDE trajectory} on $\mathcal{A}_\lambda$ taking value $u_0$ in fiber $\omega$.

If additionally,  for each $\omega \in \Omega$, $\lambda \in \Lambda$, and $u_0 \in \mathcal{A}_\lambda(\omega)$, there exists $T_0(\omega, u_0) \ge 0$, such that
\bea  \label{PM condition-hG}
\int_0^T \frac{\bigl (C_{\lambda,\sigma}(\theta_t \omega) \bigr) ^2 }{(-\eta_2 - 2 \gamma_{\lambda,\sigma}(\theta_t \omega) - \frac{3 \delta}{4})^{2(1-\alpha)}} \d t \le  \int_0^T \|u_{\s}(t, \omega; u_0)\|_\alpha^2 \d t, \qquad T \ge  T_0(\omega, u_0),  
\eea
then $\widehat{h}^{G}_\lambda$ given by \eqref{hG} is a PM function for the SPDE \eqref{SEE},   when the initial data in 
condition \eqref{PM condition} of Definition ~\ref{def:PM} are taken on the global attractor $\mathcal{A}_\lambda$.

\et

\bp

We proceed in two steps.

\medskip

{\bf Step 1.} Let us first examine the existence of the pullback limit $\widehat{h}^{G}_\lambda$. Since 
the {\mkk backward solution $\widehat{u}^{G}_{\c}$} to \eqref{Galerkin_uc} exists according to our assumption, we obtain then for all {\mkk $s \in [0,T]$} that
\beas
\widehat{u}^{G}_{\s}[\xi](s, \theta_{-T} \omega; 0) & = \int_{0}^{s}  e^{\diffusion (W_{s-T}(\omega) - W_{s'-T}(\omega)) \Id} e^{(s-s') L^{\s}_\lambda}  P_{\s} B \bigl( \widehat{u}^{G}_{\c}(s'-T, \omega; \xi), \widehat{u}^{G}_{\c}(s'-T, \omega; \xi) \bigr) \,\mathrm{d}s' \\
& = \int_{-T}^{s-T}  e^{\diffusion( W_{s-T}(\omega) - W_{\tau}(\omega)) \Id} e^{(s-T-\tau) L^{\s}_\lambda}  P_{\s} B \bigl( \widehat{u}^{G}_{\c}(\tau, \omega; \xi), \widehat{u}^{G}_{\c}(\tau, \omega; \xi) \bigr) \,\mathrm{d}\tau.
\eeas
It follows {\mkk then by application of  \eqref{Proj-II} and  by making $s=T$}
\be \label{us_G control}
\|\widehat{u}^{G}_{\s}[\xi](T, \theta_{-T} \omega; 0)\|_\alpha \le K M \int_{-T}^{0} \frac{e^{- \sigma W_{\tau}(\omega) -\tau \eta_2}}{|\tau|^\alpha} \|\widehat{u}^{G}_{\c}(\tau, \omega; \xi)\|_\alpha^2 \,\mathrm{d}\tau,
\ee
where we have used the {\mkk existence of  $M>0$  satisfying}
\be \label{bilinear control}
\|B(u,v)\| \le M \|u\|_\alpha \|v\|_\alpha, \qquad \Forall u, v \in \mathcal{H}_\alpha,
\ee
{\mkk given that   $B$ is bilinear and continuous}.

{\mkk Due  to assumption \eqref{uc growth control}}, there exists a positive random variable $C_{1}(\omega,\lambda)$ such that
\be  \label{uc growth control-2}
\|\widehat{u}^{G}_{\c}(s, \omega; \xi)\|_\alpha  \le C_{1}(\omega,\lambda) e^{\gamma_{\lambda,\sigma}(\omega) |s|}, \qquad \Forall \xi \in \mathcal{H}^{\c}, \;  s \in (-\infty, 0], \; \omega \in \Omega.
\ee
We obtain then from \eqref{us_G control} that
\be
\|\widehat{u}^{G}_{\s}[\xi](T, \theta_{-T} \omega; 0)\|_\alpha \le K M \bigl( C_{1}(\omega,\lambda) \bigr)^2 \int_{-T}^{0} \frac{e^{- \sigma W_{\tau}(\omega) -\tau \eta_2 - 2 \tau \gamma_{\lambda,\sigma}(\omega)}}{|\tau|^\alpha}  \,\mathrm{d}\tau. 
\ee
Note that by the {\mkk sublinear growth property \eqref{Wt control} of the Wiener process}, there exists  a positive random variable $C_{\sigma}(\omega)$ such that
\be \label{Wt control}
| \sigma W_{t}(\omega)| \le {\mkk  \frac{\delta}{4}|t|} + C_{\sigma}(\omega), \qquad t \in \mathbb{R}, \omega \in \Omega,
\ee
where $\delta$ is given by \eqref{random gap}; see Lemma ~\ref{Lem:OU}. 

It follows then that
\bea
\|\widehat{u}^{G}_{\s}[\xi](T, \theta_{-T} \omega; 0)\|_\alpha & \le K M \bigl( C_{1}(\omega,\lambda) \bigr)^2 C_{\sigma}(\omega) \int_{-T}^{0} \frac{e^{  (-\eta_2 - 2 \gamma_{\lambda,\sigma}(\omega) - \delta/4) \tau}}{|\tau|^\alpha}  \,\mathrm{d}\tau \\
& \le \frac{K M \bigl( C_{1}(\omega,\lambda) \bigr)^2 C_{\sigma}(\omega) \Gamma(1-\alpha)}{(-\eta_2 - 2 \gamma_{\lambda,\sigma}(\omega) - \delta/4)^{1-\alpha}}, \quad \Forall \xi \in \mathcal{H}^{\c}, \;  T > 0, \; \omega \in \Omega,
\eea
where we have used $-\eta_2 - 2 \gamma_{\lambda,\sigma}(\omega) - \delta/4 > 0$ {\mkk due} to \eqref{random gap}. Consequently, $\|\widehat{u}^{G}_{\s}[\xi](T, \theta_{-T} \omega; 0)\|_\alpha$ is uniformly bounded with respect to $T$. Following the same type of estimates, it can also be shown that 
\beas
& {\hspace*{-3em}} \Bigl \|\int_{-\infty}^{-T}  e^{ -\diffusion W_{\tau}(\omega) \Id} e^{-\tau L^{\s}_\lambda}  P_{\s} B \bigl( \widehat{u}^{G}_{\c}(\tau, \omega; \xi), \widehat{u}^{G}_{\c}(\tau, \omega; \xi) \bigr) \,\mathrm{d}\tau \Bigr \|_\alpha \\
& \le \frac{K M \bigl( C_{1}(\omega,\lambda) \bigr)^2 C_{\sigma}(\omega) e^{- (-\eta_2 - 2 \gamma_{\lambda,\sigma}(\omega) - \delta/4) T}}{T^\alpha (-\eta_2 - 2 \gamma_{\lambda,\sigma}(\omega) - \delta/4)},
\eeas
where the RHS converges to zero as $T \rightarrow +\infty$. Therefore, the pullback limit of $\widehat{u}^{G}_{\s}[\xi](T, \theta_{-T} \omega; 0)$ exists and is given by \eqref{hG}.

\medskip
{\bf Step 2.} We show now that the estimate \eqref{error hG} holds. Let $u_\lambda(t,\omega; u_0)$ {\mkk ($t \in \mathbb{R}$)  be a complete trajectory on the global attractor $\mathcal{A}_\lambda$ taking value $u_0$ in fiber $\omega$}. Note that for any $t>t_0$, it holds that
\bea \label{us integral representation-1}
u_{\s}(t,\omega; u_0) & = e^{\sigma (W_t(\omega) - W_{t_0}(\omega))\Id}e^{(t-t_0)L_\lambda^{\s}}u_{\s}(t_0, \omega; u_0) \\
& {\hspace*{2em}} + \int_{t_0}^{t} e^{\sigma (W_t(\omega) - W_{s}(\omega))\Id}e^{(t-s)L_\lambda^{\s}} P_{\s}B(u_\lambda(s, \omega; u_0), u_\lambda(s, \omega; u_0)) \d s.
\eea

Let us first show that for each fixed $t \in \mathbb{R}$, 
\be \label{us linear part control}
\lim_{t_0\rightarrow -\infty}\|e^{\sigma (W_t(\omega) - W_{t_0}(\omega))\Id}e^{(t-t_0)L_\lambda^{\s}}u_{\s}(t_0, \omega; u_0)\|_\alpha = 0.
\ee

Since $\mathcal{A}_\lambda$ is tempered, there exists a positive random variable $C_3(\omega,\lambda)$ such that 
\be \label{Attractor control}
\|u_\lambda(s, \omega; u_0)\|_\alpha < C_3(\omega, \lambda) e^{\frac{\delta}{4}|s|}, \qquad \Forall s \in \mathbb{R}, \; u_0 \in \mathcal{A}_\lambda(\omega), \; \omega \in \Omega.
\ee
This together with the partial-dichotomy estimate \eqref{Proj-I} implies that
\beas
\|e^{\sigma (W_t(\omega) - W_{t_0}(\omega))\Id}e^{(t-t_0)L_\lambda^{\s}}u_{\s}(t_0, \omega; u_0)\|_\alpha  &  \le  K C_3(\omega, \lambda) e^{\sigma (W_t(\omega) - W_{t_0}(\omega)) + (t-t_0) \eta_2 + \frac{\delta}{4} |t_0|} \\
&   \le  K C_3(\omega, \lambda) e^{C_{\sigma}(\omega) + \sigma W_t(\omega) + t \eta_2} e^{\frac{\delta}{4} |t_0| - t_0 \eta_2 + \frac{\delta}{4} |t_0|},
\eeas
where we applied \eqref{Wt control} to control $\sigma W_{t_0}(\omega)$. 

Since $-\eta_2 - \frac{\delta}{2} > 0$ thanks to \eqref{random gap}, it follows that 
\bes
\lim_{t_0\rightarrow -\infty} e^{\frac{\delta}{4} |t_0| - t_0 \eta_2 + \frac{\delta}{4} |t_0|} = \lim_{t_0\rightarrow -\infty} e^{ (- \eta_2 - \frac{\delta}{2}) t_0} = 0;
\ees
and \eqref{us linear part control} follows. We obtain then from \eqref{us integral representation-1} the following representation of $u_{\s}(t,\omega; u_0)$:
\bea \label{us representation}
u_{\s}(t,\omega; u_0)  = \int_{-\infty}^{t} e^{\sigma (W_t(\omega) - W_{s}(\omega))\Id}e^{(t-s)L_\lambda^{\s}} P_{\s}B(u_\lambda(s, \omega; u_0), u_\lambda(s, \omega; u_0)) \d s.
\eea

Let us now derive \eqref{error hG} for $t = 0$. By setting $t=0$ in \eqref{us representation}, we obtain
\beas
u_{\s}(0,\omega; u_0)  = \int_{-\infty}^{0} e^{- \sigma W_{s}(\omega))\Id}e^{-s L_\lambda^{\s}} P_{\s}B(u_\lambda(s, \omega; u_0), u_\lambda(s, \omega; u_0)) \d s.
\eeas
Let us denote
\bea
v_\lambda(s,\omega;u_0)  := u_\lambda(s,\omega;u_0)  - \widehat{u}^{G}_{\c}(s, \omega; P_{\c}u_0), \qquad s \le 0.
\eea
It follows then from the bilinearity of $B$ and the representation of $\widehat{h}_\lambda^{G}$ given in \eqref{hG} that
\bea  \label{us control-1}
u_{\s}(0,\omega; u_0)  = \widehat{h}_\lambda^{G}(P_{\c}u_0, \omega) & + \int_{-\infty}^{0} e^{- \sigma W_{s}(\omega))\Id}e^{-s L_\lambda^{\s}} P_{\s}B(\widehat{u}^{G}_{\c}(s, \omega; P_{\c}u_0), v_\lambda(s, \omega; u_0)) \d s \\
& + \int_{-\infty}^{0} e^{- \sigma W_{s}(\omega))\Id}e^{-s L_\lambda^{\s}} P_{\s}B(v_\lambda(s, \omega; u_0), \widehat{u}^{G}_{\c}(s, \omega; P_{\c}u_0)) \d s \\
& + \int_{-\infty}^{0} e^{- \sigma W_{s}(\omega))\Id}e^{-s L_\lambda^{\s}} P_{\s}B(v_\lambda(s, \omega; u_0), v_\lambda(s, \omega; u_0)) \d s. 
\eea 
We proceed now with the control of the three  integral terms {\mkk  appearing in \eqref{us control-1}}. For this purpose, let us note that according to \eqref{uc growth control-2} and \eqref{Attractor control} it holds {\mkk trivially that  for all $s\le 0$,}
\be
\|v_\lambda(s,\omega;u_0)\|_\alpha \le \|u_\lambda(s,\omega;u_0)\|_\alpha + \|\widehat{u}^{G}_{\c}(s, \omega; P_{\c}u_0)\|_\alpha \le C_3(\omega,\lambda) e^{\frac{\delta}{4}|s|} +  C_{1}(\omega,\lambda) e^{\gamma_{\lambda,\sigma}(\omega) |s|}.
\ee

{\mkk By using again the partial-dichotomy estimate \eqref{Proj-II}, the growth control of $W_t(\omega)$, and the control of $B$ given in \eqref{bilinear control} with the estimates of $\|v_\lambda \|_\alpha$ and $\|\widehat{u}^{G}_{\c} \|_\alpha$; the first integral term in \eqref{us control-1} is estimated as follows:}
\bea \label{us control-2}
& \Bigl \|\int_{-\infty}^{0} e^{- \sigma W_{s}(\omega))\Id}e^{-s L_\lambda^{\s}} P_{\s}B(\widehat{u}^{G}_{\c}(s, \omega; P_{\c}u_0), v_\lambda(s, \omega; u_0)) \d s \Bigr \|_\alpha  \\
& \le K M \int_{-\infty}^{0} \frac{e^{- \sigma W_{s}(\omega)-s \eta_2}}{|s|^\alpha} \|\widehat{u}^{G}_{\c}(s, \omega; P_{\c}u_0)\|_\alpha \|v_\lambda(s, \omega; u_0) \|_\alpha \d s \\
& \le K M  \int_{-\infty}^{0} \frac{e^{C_{\sigma}(\omega) - \frac{\delta}{4} s -s \eta_2}}{|s|^\alpha} C_{1}(\omega,\lambda) e^{\gamma_{\lambda,\sigma}(\omega) |s|}  \Bigl( C_3(\omega,\lambda) e^{\frac{\delta}{4}|s|} +  C_{1}(\omega,\lambda) e^{\gamma_{\lambda,\sigma}(\omega) |s|} \Bigr)  \d s \\
& = \frac{K M  e^{C_{\sigma}(\omega)} C_1(\omega,\lambda)C_3(\omega,\lambda)\Gamma(1-\alpha)}{(-\eta_2 - \gamma_{\lambda,\sigma}(\omega) - \frac{\delta}{2})^{1-\alpha}} +  \frac{K M  e^{C_{\sigma}(\omega)} \bigl( C_1(\omega,\lambda) \bigr)^2 \Gamma(1-\alpha)}{(-\eta_2 - 2\gamma_{\lambda,\sigma}(\omega) - \frac{\delta}{4})^{1-\alpha}}.
\eea
The {\mkk second and third integrals in \eqref{us control-1} can be controlled using {\mkk similar arguments}, to lead respectively to}
\bea \label{us control-3}
& \Bigl \|\int_{-\infty}^{0} e^{- \sigma W_{s}(\omega))\Id}e^{-s L_\lambda^{\s}} P_{\s}B(\widehat{u}^{G}_{\c}(s, \omega; P_{\c}u_0), v_\lambda(s, \omega; u_0)) \d s \Bigr \|_\alpha  \\
& \le \frac{K M  e^{C_{\sigma}(\omega)} C_1(\omega,\lambda)C_3(\omega,\lambda)\Gamma(1-\alpha)}{(-\eta_2 - \gamma_{\lambda,\sigma}(\omega) - \frac{\delta}{2})^{1-\alpha}} +  \frac{K M  e^{C_{\sigma}(\omega)} \bigl( C_1(\omega,\lambda) \bigr)^2 \Gamma(1-\alpha)}{(-\eta_2 - 2\gamma_{\lambda,\sigma}(\omega) - \frac{\delta}{4})^{1-\alpha}},
\eea
{\mkk and to}
\bea  \label{us control-4}
& \Bigl \|\int_{-\infty}^{0} e^{- \sigma W_{s}(\omega))\Id}e^{-s L_\lambda^{\s}} P_{\s}B(v_\lambda(s, \omega; u_0), v_\lambda(s, \omega; u_0)) \d s \Bigr \|_\alpha  \\
& \le \frac{2 K M  e^{C_{\sigma}(\omega)} \bigl( C_1(\omega,\lambda)\bigr)^2 \Gamma(1-\alpha)}{(-\eta_2  - 2\gamma_{\lambda,\sigma}(\omega) - \frac{\delta}{4})^{1-\alpha}} +  \frac{2 K M  e^{C_{\sigma}(\omega)} \bigl( C_3(\omega,\lambda) \bigr)^2 \Gamma(1-\alpha)}{(-\eta_2 -  \frac{3 \delta}{4})^{1-\alpha}}.
\eea

Now, by \eqref{us control-1} and the estimates derived in \eqref{us control-2}--\eqref{us control-4}, we {\mkk have}
\bes
\|u_{\s}(0,\omega; u_0) - \widehat{h}^{G}_\lambda(u_{\c}(0,\omega; u_0), \omega) \|_\alpha \le \frac{C_\lambda(\omega)}{(-\eta_2 - 2 \gamma_{\lambda,\sigma}(\omega) - \frac{3\delta}{4})^{1-\alpha}}, \; \Forall u_0 \in \mathcal{A}_\lambda(\omega), 
\ees
where
\be \label{C in error hG}
C_\lambda(\omega) := K M  e^{C_{\sigma}(\omega)}  \Gamma(1-\alpha) \Bigl ( 4 \bigl( C_1(\omega,\lambda)\bigr)^2 + 2  C_1(\omega,\lambda) C_3(\omega,\lambda) + 2 \bigl( C_3(\omega,\lambda)\bigr)^2 \Bigr).
\ee
The estimate \eqref{error hG} is thus derived when $t = 0$. The result for $t>0$ follows by simply {\mkk changing the fiber $\omega$ by the fiber $\theta_t\omega$.}

\medskip

Finally, if  \eqref{PM condition-hG} is satisfied {\mkk it is clear that for each $\lambda \in \Lambda$, $\widehat{h}^{G}_\lambda$ consists of a PM function} when the initial data in 
condition \eqref{PM condition} required by  Definition ~\ref{def:PM} are taken on the global attractor $\mathcal{A}_\lambda$.
 
\ep

\br

Note that the existence of the backward solution $\widehat{u}^{G}_{\c}(\cdot, \omega; \xi)$ for \eqref{Galerkin_uc} on the negative real line and the condition \eqref{uc growth control} can be both guaranteed if for instance the bilinear term  is furthermore energy preserving, {\it i.e.}, $\langle B(u,u), u \rangle = 0$ for all $u \in \mathcal{H}_\alpha$; see \cite{FJ05} for such results  in the deterministic context.

\er

The same ideas than  used in the proof of Theorem ~\ref{thm:hG_PM} can be applied to show that $\widehat{h}^{(1)}_\lambda$  associated with the backward-forward systems \eqref{LLL} provides a PM under appropriate assumptions as stated in the following theorem.
 
\bt \label{thm:h1_PM}

Consider an SPDE \eqref{SEE} whose cocycle  possesses a global tempered  random attractor $\mathcal{A}_\lambda$.  Let us assume  that the nonlinearity $F$ is a continuous bilinear term,  and that $\lambda$  lives  in some interval $\Lambda$ for which a uniform spectral decomposition \eqref{gap 3} over $\Lambda$ holds;  and  $\eta_1$ and $\eta_2$ can be chosen such that 
\be\label{Eq_cond_golh1}
\eta_2 < 2 \eta_1, \textrm { and } \eta_{\s} < \eta_2 < \eta_1 < \eta_{\c}.
\ee

Then, there exists a positive random variable $C_\lambda(\omega)$ such that
\be \label{error hG}
\|u_{\s}(t,\omega; u_0) - \widehat{h}^{(1)}_\lambda(u_{\c}(t,\omega; u_0), \theta_t \omega) \|_\alpha \le \frac{C_\lambda(\theta_t\omega)}{((2 \eta_1 -\eta_2)/2)^{1-\alpha}}, \; \Forall u_0 \in \mathcal{A}_\lambda(\omega), \; t \ge 0, \; \lambda \in \Lambda,
\ee
where $u_{\s}(t,\omega; u_0) := P_{\s}u_\lambda(t,\omega; u_0)$ and $u_{\c}(t,\omega; u_0) := P_{\c}u_\lambda(t,\omega; u_0)$ with $u_\lambda(t,\omega; u_0)$ being the complete trajectory on $\mathcal{A}_\lambda$ taking value $u_0$ in fiber $\omega$.

If additionally, for each $\omega \in \Omega$, $\lambda \in \Lambda$, and $u_0 \in \mathcal{A}_\lambda(\omega)$, there exists $T_0(\omega, u_0) \ge 0$, such that
\bea  \label{PM condition-hG}
  \frac{1}{((2 \eta_1 -\eta_2)/2)^{1-\alpha}} \int_0^T \bigl ( C_\lambda(\theta_t \omega) \bigr) ^2 \d t  \le  \int_0^T \|u_{\s}(t, \omega; u_0)\|_\alpha^2 \d t, \qquad T \ge  T_0(\omega, u_0), 
\eea
then $\widehat{h}^{(1)}_\lambda$ is a PM function for the SPDE \eqref{SEE} for each $\lambda \in \Lambda$, when the initial data in 
condition \eqref{PM condition} of Definition ~\ref{def:PM} are taken on the global attractor $\mathcal{A}_\lambda$.

\et

\bp

The results can be derived by following the same arguments presented in the proof of Theorem ~\ref{thm:hG_PM} since analogues of the additional assumptions \eqref{uc growth control} and \eqref{random gap} required therein can be verified for the case considered here as shown below.

As mentioned in \SS~\ref{ss:PM}, $\widehat{h}^{(1)}_\lambda$ is well-defined under the condition $\eta_{\s} < 2 \eta_{\c}$. Note also that for any fixed $\epsilon > 0$, there exists a positive random variable $C_\epsilon(\omega)$ such that the backward solution $\widehat{u}_{\c}^{(1)}(s,\omega; \xi)$ for \eqref{Eq:spde-I} has the following growth control:
\be
\|\widehat{u}_{\c}^{(1)}(s,\omega; \xi)\|_\alpha \le K e^{\eta_1 s + \sigma W_s(\omega)}\|\xi\|_\alpha \le C_\epsilon(\omega) e^{\eta_1 s + \epsilon |s|}\|\xi\|_\alpha.
\ee
We derive then in analogue to the condition \eqref{uc growth control} required in Theorem ~\ref{thm:hG_PM} that:
\bes
\lim_{T\rightarrow +\infty} e^{- (|\eta_1| + \epsilon) T} \|\widehat{u}^{(1)}_{\c}(-T, \omega; \xi)\|_\alpha  = 0, \qquad \Forall \xi \in \mathcal{H}^{\c}, \, \omega \in \Omega,\; \lambda \in \Lambda, \; \epsilon > 0,
\ees
where the rate of the growth control, $\gamma_\epsilon:= |\eta_1| + \epsilon$, is deterministic in contrast with the random rate $\gamma_{\lambda,\sigma}(\omega)$ in \eqref{uc growth control}.  Since $\eta_2 < 2 \eta_1$, in analogue to \eqref{random gap}, we can choose the positive constants $\delta$ and $\epsilon$ properly so that 
\bes
-\eta_2 - 2 \gamma_\epsilon > \delta.
\ees
We can then follow the same argument as presented in Step 2 of the proof of Theorem ~\ref{thm:hG_PM} to obtain the desired results, and the proof is complete.

\ep

Clearly condition \eqref{Eq_cond_golh1} is far to be necessary  for  $\widehat{h}^{(1)}_\lambda$ to be a PM. This is for instance supported by the numerical results  of  Section 
~\ref{ss:Burgers reduction} 
 on the stochastic Burgers-type equation analyzed therein. 
In the deterministic context, actually  much weaker conditions on the spectrum of the linear part can be identified to  
ensure $\widehat{h}^{(1)}_\lambda$ to be a PM. This is the purpose of Theorem \ref{thm:h1_PMdet} proved in the next section.

\br
It is worth mentioning that the $n$-layer backward-forward system \eqref{n_layer}  can  actually give access to PMs,  by considering the random graph associated with $\widehat{u}^{(n)}_{\s}[\xi](T, \theta_{-T} \omega; 0)$ 
($\xi$ varying in $\mathcal{H}^{\c}$) for some finite fixed time $T$.  We refer to Section ~\ref{s:Fly} for illustrations on the stochastic Burgers-type equation analyzed therein. 
\er

\subsection{Existence of  parameterizing manifolds in the deterministic case}\label{sec_PMdet}


In this section we study, in a deterministic setting {\hh ($\sigma = 0$)}, conditions under which the pullback limit associated with the backward-forward system \eqref{n_layer} provides a PM.\footnote{
Note that, in the deterministic case, the corresponding backward-forward system \eqref{n_layer} becomes autonomous and the equation for the $\widehat{u}_{\s}^{(n)}$ is solved over the interval $[-T, 0]$ instead of $[0, T]$, see for instance \eqref{Eq:Duan_etal_shape}. Here, the definition of a parameterizing manifold for a deterministic PDE follows the same lines as Definition ~\ref{def:PM} with the $\omega$-dependence removed and the requirements regarding the measurability dropped.} For simplicity, we consider $\widehat{h}^{(1)}_\lambda$ only; and throughout this subsection, we assume that the nonlinearity $F$ is $k$-linear for some $k\ge 2$.

The {\mkk evolution} equation under consideration takes the following form 
\be \label{PDE}
\frac{\mathrm{d} u}{\d t} =  L_\lambda u + F(u),
\ee
which corresponds to Eq. ~\eqref{SEE} by setting $\sigma = 0$.

Note that when $\sigma = 0$ and $F = F_k$, we obtain from \eqref{Eq_PBA1uc} the following expression for $\widehat{h}^{(1)}_\lambda(\xi)$:
\be \label{det_h1}
\widehat{h}^{(1)}_\lambda(\xi) = \int_{-\infty}^0 e^{-s L^{\s}_\lambda} P_{\s} F(
e^{s L^{\c}_\lambda}\xii) \,\mathrm{d}s, \quad \Forall \xii \in  \mathcal{H}^{\c}.
\ee

{\hh Since $F$ is $k$-linear, it follows from \cite[Lem.~6.2.4]{Hen81} that $\widehat{h}^{(1)}_\lambda$ satisfies the following equation:}
\be \label{INV_h1}
D \widehat{h}^{(1)}_\lambda(\xi) L_\lambda^{\c} \xi - L_\lambda^{\s} \widehat{h}^{(1)}_\lambda(\xi) - P_{\s} F(\xi)  = 0.
\ee
In particular, the {\it invariance defect} {\mkk  of $\widehat{h}^{(1)}_\lambda$ (associated with \eqref{PDE})} is given by:\footnote{ 
{\mkk Note that the invariance defect (associated with \eqref{PDE}) of a manifold is zero if and only if  this manifold is invariant for the semigroup generated by \eqref{PDE}}.} 
\bea
M [\widehat{h}^{(1)}_\lambda](\xi) & := D \widehat{h}^{(1)}_\lambda(\xi) \Bigl(L_\lambda^{\c} \xi  + P_{\c} F \Bigl ( \xi + \widehat{h}^{(1)}_\lambda (\xi) \Bigr) \Bigr) -  \Bigl( L_\lambda^{\s} \widehat{h}^{(1)}_\lambda(\xi) +  P_{\s} F \Bigl (\xi + \widehat{h}^{(1)}_\lambda(\xi) \Bigr)\Bigr) \\
& = R[\widehat{h}^{(1)}_\lambda]  (\xi),
\eea
with
\be
R[\widehat{h}^{(1)}_\lambda]  (\xi) := D \widehat{h}^{(1)}_\lambda(\xi) P_{\c} F \Bigl( \xi + \widehat{h}^{(1)}_\lambda (\xi) \Bigr ) -  P_{\s} F\Bigl(\xi + \widehat{h}^{(1)}_\lambda(\xi) \Bigr) + P_{\s} F(\xi).
\ee
We will show that this invariance defect (when $\widehat{h}^{(1)}_\lambda$ is replaced by $u_{\s}(t)$) plays an important role to ensure $\widehat{h}^{(1)}_\lambda$ to be a PM. 

To do so, let us introduce $v(t; u_0):= u_{\s}(t;u_0) - \widehat{h}^{(1)}_\lambda(u_{\c}(t;u_0))$, where $u_{\c}(t;u_0) := P_{\c} u_\lambda(t;u_0)$ and $u_{\s}(t;u_0) := P_{\s}u_\lambda(t;u_0)$ with $u_\lambda(t;u_0)$ being the solution to Eq. ~\eqref{PDE} emanating from some $u_0 \in \mathcal{H}_{\alpha}$. Note that
\bea
\frac{\d \widehat{h}^{(1)}_\lambda(u_{\c})}{\d t} & = D \widehat{h}^{(1)}_\lambda(u_{\c}) \frac{\d u_{\c}}{\d t} = D \widehat{h}^{(1)}_\lambda(u_{\c}) \Bigl(L_\lambda^{\c} u_{\c}  + P_{\c} F( u_{\c} + u_{\s})\Bigr ) \\
& = L_\lambda^{\s} \widehat{h}^{(1)}_\lambda(u_{\c}) + P_{\s} F_k(u_{\c}) + D \widehat{h}^{(1)}_\lambda(u_{\c}) P_{\c} F (u_{\c} + u_{\s}),
\eea
where we have used \eqref{INV_h1} to derive the last equality. It follows then that 
\bea  \label{eq:v} 
\frac{\d v}{\d t} = L_\lambda^{\s} v + P_{\s} F\bigl (u_{\c} + v + \widehat{h}^{(1)}_\lambda(u_{\c}) \bigr) -  P_{\s} F_k(u_{\c}) - 
D \widehat{h}^{(1)}_\lambda(u_{\c}) P_{\c} F \bigl(u_{\c} +v + \widehat{h}^{(1)}_\lambda (u_{\c}) \bigr). 
\eea

The following theorem identifies conditions under which $\widehat{h}^{(1)}_\lambda$ provides an AIM for a broad class of PDEs. These conditions are subject to the theory of time-analyticity properties of the solutions $u_\lambda(t)$ of the underlying PDE as extended in \cite{Pro91} from the original works \cite{foias1979some,FT89}, that we adapt to our framework; see also \cite{GK03}. Similar time-analyticity arguments for the theory of AIMs were first used in \cite{FMT88}. We first recall some assumptions (adapted to  our framework) required in \cite{Pro91} in order to establish this property.

\bi

\item[{\bf (A1)}] The linear operator $L_\lambda= -A + B_\lambda$ is self-adjoint and satisfies the assumptions of Section ~\ref{s:preliminary}. 

\item[{\bf (A2)}] The nonlinearity $F:\mathcal{H}_\alpha \rightarrow \mathcal{H}$  is continuous for some $\alpha \in [0, 1/2]$; and $F$ is analytic on finite dimensional subspaces of $\mathcal{H}_\alpha$ and has analytic extensions to the complexifications of these subspaces.

\item[{\bf (A3)}]  There exist $1 \le r  < \infty$, $C > 0$, and a continuous function $c: (0, \infty) \rightarrow \mathbb{R}^{+}$ such that for all $r > 0$, the following estimate holds 
\be
| \langle \widetilde{F}(u) + \widetilde{B}_\lambda u, \widetilde{A} u  + u \rangle | \le  r \|\widetilde{A} u \|^2 + c(r)\|u\|^{2p} + C, \qquad \Forall u \in \widetilde{\mathcal{H}}_1,
\ee
where $\widetilde{A}, \widetilde{B}_\lambda, \widetilde{F}$ are the complexification of $A, B_\lambda$, and $F$ respectively; $\widetilde{\mathcal{H}}_1$ is the complexification of $\mathcal{H}_1$ defined by \eqref{eq: complex H}; $\langle \cdot, \cdot \rangle$ is the inner product on $\widetilde{\mathcal{H}}$ induced by the corresponding inner product on $\mathcal{H}$; and $\lambda > \lambda_c$  is fixed.

\item[{\bf (A4)}]  Moreover, $R:=\widetilde{F} + \widetilde{B}_\lambda$ satisfies the following weak continuity property for $\mathcal{H}_\alpha$-valued functions of the complex variables: 

If $v_m : \mathbb{C} \rightarrow \mathcal{H}_\alpha$ converges weakly to $v : \mathbb{C} \rightarrow \mathcal{H}_\alpha$ in the sense that $l(v_m(z)) \rightarrow l(v(z))$ as $m\rightarrow \infty$ for every $l \in (\mathcal{H}_\alpha)^{\ast}$ and uniformly in $z$ on every compact subset $K \subset \mathbb{C}$, then 
\bes
\int_{\mathbb{C}} \Bigl( R(v_m(z)) - R(v(z)) \Bigr) \phi(z) \d z \mathop{\longrightarrow}_{m\rightarrow \infty} 0, \qquad \Forall \phi \in C_{c}^{\infty}(\mathbb{C}).
\ees

\ei


\bt\label{thm:h1_PMdet}

Let us assume that {\bf (A1)}--{\bf (A4)} hold. Let us assume furthermore  that $A$ and $L_\lambda$ share the same eigenfunctions (in $\mathcal{H}_1$), with the corresponding eigenvalues $\gamma_n$ and $\beta_n(\lambda)$ verifying 
\be \label{eq:eigen_relation}
\exists \; C_\lambda > 0,\;\forall \; n \in \mathbb{N},\;  |\gamma_n - \beta_n(\lambda)| < C_\lambda |\beta_n(\lambda)|.
\ee 


The projector $P_{\c}$ is assumed to be chosen so that $-L_\lambda^{\s}$ defined in \eqref{Lc Ls}  is positive, and so that the cross non-resonance condition 
\eqref{NR},  between the eigenvalues of $L_\lambda^{\s}$ and $L_\lambda^{\c}$,  holds.

Finally, we assume that a semigroup can be associated with ~\eqref{PDE} in $\mathcal{H}_{\alpha}$, and that this semigroup possesses  a global attractor  $\mathcal{A}_\lambda$ whose $\epsilon$-neighborhood is given by
\bes
\mathcal{A}^\epsilon_\lambda:= \{ \psi \in \mathcal{H}_\alpha  : \mathrm{dist}(\psi, \mathcal{A}_\lambda) \le \epsilon \}.
\ees
Then for any given solution $u_\lambda(t; u_0)$ of \eqref{PDE}, there exist $t_0(u_0)> 0$ and $\gamma(\lambda,u_0)>0$ such that $v(t; u_0) := u_{\s}(t; u_0) - \widehat{h}^{(1)}_\lambda(u_{\c}(t; u_0) )$, with $\widehat{h}^{(1)}_\lambda$  given by \eqref{det_h1}, satisfies:
\bea\label{v_est_1}
\|v(t; u_0)\|_\alpha & \le \frac{( 1 + C_\lambda)^\alpha \gamma(\lambda, u_0)}{|\beta_{m+1}(\lambda)|^{1-\alpha}} +  
\frac{( 1 + C_\lambda)^\alpha \mathrm{Lip}(F)\vert_{\mathcal{A}^\epsilon_\lambda}}{|\beta_{m+1}(\lambda)|^{1-\alpha}} \| u_{\s}(t; u_0)\|_\alpha, \quad t \ge t_0(u_0),
\eea
where $\mathrm{Lip}(F)\vert_{\mathcal{A}^\epsilon_\lambda}$ is the Lipschitz constant of $F$ restricted to $\mathcal{A}^\epsilon_\lambda$ and $C_\lambda$ is as given in \eqref{eq:eigen_relation}.
  \et

\bp

Let $u_\lambda(t; u_0)$ be the solution to Eq. ~\eqref{PDE} emanating from $u_0 \in \mathcal{H}_\alpha$. First note that under the assumptions {\bf (A1)}--{\bf (A4)} given above, we can apply \cite[Thm. ~1.1]{Pro91} to the solution $u_\lambda$ to obtain that there exists a positive constant $\gamma_1:=\gamma_1(\lambda, u_0)$ such that
\be 
\Big \|\frac{\d u_\lambda(t; u_0)}{\d t}\Big\|_\alpha \le \gamma_1(\lambda, u_0), \quad t \ge t_0(u_0),
\ee
where $t_0(u_0)$ is chosen so that $u_\lambda(t; u_0) \in \mathcal{A}^\epsilon_\lambda$ for all $t \ge t_0$. It follows then that
\bea  \label{eq:v_t control}
\Big\|\frac{\d v(t; u_0)}{\d t} \Big \| & \le \Big \|\frac{\d u_{\s}(t; u_0)}{\d t} \Big\| + \Big \|\frac{\d \widehat{h}^{(1)}_\lambda(u_{\c} (t; u_0) )}{\d t} \Big \|  \\
& \le C_1 \Big \|\frac{\d u_{\s}(t; u_0)}{\d t} \Big \|_\alpha + C_1 \Big \|D \widehat{h}^{(1)}_\lambda(u_{\c} (t; u_0) ) \Big \|  \Big \|\frac{\d u_{\c}(t; u_0) }{\d t} \Big\|_\alpha  \\
& \le C_1 \gamma_1(\lambda, u_0) (1 + \sup_{\xi \in P_{\c}\mathcal{A}^\epsilon_\lambda} \|D \widehat{h}^{(1)}_\lambda (\xi)\|), \quad t \ge t_0,
\eea
where $C_1$ is the generic constant verifying $\|\psi\| \le C_1 \|\psi\|_\alpha$ for all $\psi \in \mathcal{H}_\alpha$. Note that the supremum $\sup_{\xi \in P_{\c}\mathcal{A}^\epsilon_\lambda} \|D \widehat{h}^{(1)}_\lambda (\xi)\|$ is finite since $P_{\c}\mathcal{A}^\epsilon_\lambda$ is a bounded closed set in the finite dimensional subspace $\mathcal{H}^{\c}$, and $D \widehat{h}^{(1)}_\lambda$ is continuous due to the fact that, with $L_\lambda$ being self-adjoint, $\widehat{h}^{(1)}_\lambda$ is a homogenous (deterministic) polynomial in the sense of Definition ~\ref{Def:random polynomial} (for $\sigma=0$).

Note also that according to \eqref{eq:v}, we have 
\bea \label{eq:v_identity}
v = (-L_{\lambda}^{\s})^{-1}\Bigl( \frac{\d v}{\d t} - P_{\s} F\bigl (u_{\c} + v +  \widehat{h}^{(1)}_\lambda(u_{\c}) \bigr) + P_{\s} F_k(u_{\c})   +  
D \widehat{h}^{(1)}_\lambda(u_{\c}) P_{\c} F \bigl(u_{\c} +v + \widehat{h}^{(1)}_\lambda (u_{\c}) \bigr)  \Bigr). 
\eea
The desired estimate on $v$ can be derived now based on the following two inequalities. Thanks to the condition \eqref{eq:eigen_relation}, we get
\bea
\|(-L_{\lambda}^{\s})^{-1} \psi \|_\alpha = \|A^{\alpha} (-L_{\lambda}^{\s})^{-1} \psi\|  & \le \frac{\gamma_n^\alpha}{|\beta_{m+1}(\lambda)|} \|\psi\| \\
 & \le \frac{( |\beta_{m+1}(\lambda)| + |\gamma_n -\beta_{m+1}(\lambda)|)^\alpha}{|\beta_{m+1}(\lambda)|} \|\psi\| \\
& \le \frac{( 1 + C_\lambda)^\alpha}{|\beta_{m+1}(\lambda)|^{1-\alpha}} \|\psi\|, &&  \psi \in \mathcal{H}_\alpha^{\s}.
\eea
Note also that 
\bea\label{Eq_controlv}
& \Big \| \frac{\d v (t; u_0)}{\d t} - P_{\s} F\bigl (u_{\c}(t; u_0) + v(t; u_0) +  \widehat{h}^{(1)}_\lambda(u_{\c}(t; u_0)) \bigr) + P_{\s} F_k(u_{\c}(t; u_0))   \\
& {\hspace{3em}} +  
D \widehat{h}^{(1)}_\lambda(u_{\c}(t; u_0)) P_{\c} F \bigl(u_{\c}(t; u_0) +v(t; u_0) + \widehat{h}^{(1)}_\lambda (u_{\c}(t; u_0)) \bigr) \Big \| \\
& \le \Big \| \frac{\d v(t; u_0)}{\d t} \Big\| + \mathrm{Lip}(F)\vert_{\mathcal{A}^\epsilon_\lambda} \| u_{\s}(t; u_0) \|_\alpha  \\
& {\hspace{3em}} + \Big\| D \widehat{h}^{(1)}_\lambda(u_{\c}(t; u_0)) P_{\c} F \bigl(u_{\c}(t; u_0) +v(t; u_0) + \widehat{h}^{(1)}_\lambda (u_{\c}(t; u_0) ) \bigr) \Big\| \\
& \le \gamma(\lambda, u_0) +  \mathrm{Lip}(F)\vert_{\mathcal{A}^\epsilon_\lambda} \| u_{\s}(t; u_0)\|_\alpha,   \qquad t \ge t_0(u_0),
\eea 
where 
\bea
\gamma(\lambda, u_0) & := C_1 \gamma_1(\lambda, u_0) \Bigl(1 + \sup_{\xi \in P_{\c}\mathcal{A}^\epsilon_\lambda}\|D \widehat{h}^{(1)}_\lambda (\xi)\| \Bigr) \\
& {\hspace{3em}} +  \mathrm{Lip}(F)\vert_{\mathcal{A}^\epsilon_\lambda}  \sup_{\xi \in P_{\c}\mathcal{A}^\epsilon_\lambda} \|D \widehat{h}^{(1)}_\lambda (\xi)\|  \sup_{\psi \in \mathcal{A}^\epsilon_\lambda}  \|\psi\|_\alpha,
\eea
and where the first  term in the RHS of the second inequality of \eqref{Eq_controlv} has been controlled by using  \eqref{eq:v_t control}, and the third term of the same RHS,
by using 
\beas
& \Big\| D \widehat{h}^{(1)}_\lambda(u_{\c}(t; u_0)) P_{\c} F \bigl(u_{\c}(t; u_0) +v(t; u_0) + \widehat{h}^{(1)}_\lambda (u_{\c}(t; u_0) ) \bigr) \Big\| \\
& \le \sup_{\xi \in P_{\c}\mathcal{A}^\epsilon_\lambda} \|D \widehat{h}^{(1)}_\lambda (\xi)\|   \cdot \mathrm{Lip}(F)\vert_{\mathcal{A}^\epsilon_\lambda}  \|u_{\lambda}(t; u_0) \|_\alpha \\
& \le \mathrm{Lip}(F)\vert_{\mathcal{A}^\epsilon_\lambda}  \sup_{\xi \in P_{\c}\mathcal{A}^\epsilon_\lambda} \|D \widehat{h}^{(1)}_\lambda (\xi)\|  \sup_{\psi \in \mathcal{A}^\epsilon_\lambda}  \|\psi\|_\alpha.
\eeas

We obtain then from \eqref{eq:v_identity} that
\bea  \label{v_est}
\|v(t; u_0)\|_\alpha & \le \frac{( 1 + C_\lambda)^\alpha \gamma(\lambda, u_0)}{|\beta_{m+1}(\lambda)|^{1-\alpha}} +  
\frac{( 1 + C_\lambda)^\alpha \mathrm{Lip}(F)\vert_{\mathcal{A}^\epsilon_\lambda}}{|\beta_{m+1}(\lambda)|^{1-\alpha}} \| u_{\s}(t; u_0)\|_\alpha, \qquad t \ge t_0(u_0).
\eea
The proof is now complete.

\ep

\br

The results in the above theorem shows that $\widehat{h}^{(1)}_\lambda$ consists of an approximate inertial manifold for Eq. ~\eqref{PDE} under the given conditions, since for all $t \ge t_0(u_0)$ and  $u_0 \in \mathcal{H}_\alpha$, it holds that
\beas
\mathrm{dist}(u(t; u_0), \mathfrak{M}^{(1)}_\lambda) & \le \|v(t; u_0)\|_\alpha  \\
& \le \frac{(1 + C_\lambda)^\alpha \gamma(\lambda, u_0)}{|\beta_{m+1}(\lambda)|^{1-\alpha}} +  
\frac{( 1 + C_\lambda)^\alpha \mathrm{Lip}(F)\vert_{\mathcal{A}^\epsilon_\lambda}}{|\beta_{m+1}(\lambda)|^{1-\alpha}} \| u_{\s}(t; u_0)\|_\alpha, 
\eeas
where $\mathfrak{M}^{(1)}_\lambda$ is the manifold associated with $\widehat{h}^{(1)}_\lambda$. 

\er

The following corollary provides additional conditions under which $\widehat{h}^{(1)}_\lambda$ is a PM.

\bc\label{Cor_h1_is_PM}
Let us assume  that $\lambda>\lambda_c$, and that the assumptions of Theorem \ref{thm:h1_PMdet} hold, with $t_0(u_0)$ such as given in \eqref{v_est_1} for an arbitrary $u_0 \in \mathcal{H}_\alpha$.

Then  $\widehat{h}^{(1)}_\lambda$ defined by \eqref{det_h1}  provides a  parameterizing manifold for Eq. ~\eqref{PDE}, if 
there exists $t_1 \ge t_0(u_0)$ such that, 
\be \label{eq:PM_hypothesis}
\frac{2 (T- t_0) \big((1 + C_\lambda)^\alpha \gamma(\lambda, u_0)\big)^2}{|\beta_{m+1}(\lambda)|^{2(1-\alpha)} \int_0^T \|u_{\s}(t; u_0)\|_\alpha^2 \d t} + \frac{2 \big((1 + C_\lambda)^\alpha \mathrm{Lip}(F)\vert_{\mathcal{A}^\epsilon_\lambda}\big)^2}{|\beta_{m+1}(\lambda)|^{2(1-\alpha)}} < 1, \qquad \Forall T \ge t_1. 
\ee

\ec

\bp

Let $u_0 \in \mathcal{H}_\alpha$  be fixed and $T \ge t_0 :=t_0(u_0)$.   Let us first rewrite the time-dependent ratio $Q(T; u_0)$ in Definition ~\ref{def:PM} (adapted to the deterministic case) as follows:
\bea \label{eq:Q_h1}
Q(T; u_0) :=\frac{\int_0^T \|v(t; u_0)\|^2_\alpha \d t}{\int_0^T \|u_{\s}(t; u_0)\|^2_\alpha \d t} = \frac{\int_0^{t_0} \|v(t; u_0)\|^2_\alpha \d t}{\int_0^T \|u_{\s}(t; u_0)\|^2_\alpha \d t} + \frac{\int_{t_0}^T \|v(t; u_0)\|^2_\alpha \d t}{\int_0^T \|u_{\s}(t; u_0)\|^2_\alpha \d t},
\eea
where $v(t; u_0) = u_{\s}(t; u_0) - \widehat{h}^{(1)}_\lambda(u_{\c}(t; u_0) )$.

Since $\lambda>\lambda_c$, we have $\int_0^T \|u_{\s}(t; u_0)\|^2_\alpha \d t \rightarrow \infty$ as $T \rightarrow \infty$. The first term on the RHS of \eqref{eq:Q_h1} can be thus made sufficiently close to zero by choosing $T$ sufficiently large.

The control of the second term in the RHS is based on the estimate \eqref{v_est_1} from Theorem ~\ref{thm:h1_PMdet} and the assumption \eqref{eq:PM_hypothesis}. Indeed, thanks to \eqref{v_est_1}, we get for all $T \ge t_0$ that
\beas
\int_{t_0}^T \|v(t; u_0)\|^2_\alpha \d t & \le \int_{t_0}^T \Biggl[\frac{( 1 + C_\lambda)^\alpha \gamma(\lambda, u_0)}{|\beta_{m+1}(\lambda)|^{1-\alpha}} +
\frac{( 1 + C_\lambda)^\alpha \mathrm{Lip}(F)\vert_{\mathcal{A}^\epsilon_\lambda}}{|\beta_{m+1}(\lambda)|^{1-\alpha}} \| u_{\s}(t; u_0)\|_\alpha \Biggr]^2 \d t \\
& \le 2 \int_{t_0}^T \Biggl[\frac{( 1 + C_\lambda)^\alpha \gamma(\lambda, u_0)}{|\beta_{m+1}(\lambda)|^{1-\alpha}}\Biggr]^2 +
\Biggl[ \frac{( 1 + C_\lambda)^\alpha \mathrm{Lip}(F)\vert_{\mathcal{A}^\epsilon_\lambda}}{|\beta_{m+1}(\lambda)|^{1-\alpha}} \| u_{\s}(t; u_0)\|_\alpha \Biggr]^2 \d t \\
& = \frac{ 2 (T-t_0) [( 1 + C_\lambda)^ {\alpha} \gamma(\lambda, u_0)]^2}{|\beta_{m+1}(\lambda)|^{2(1-\alpha)}} +
\frac{2 [(1 + C_\lambda)^\alpha \mathrm{Lip}(F)\vert_{\mathcal{A}^\epsilon_\lambda}]^2}{|\beta_{m+1}(\lambda)|^{2(1-\alpha)}} \int_{t_0}^T \| u_{\s}(t; u_0)\|_\alpha^2 \d t.
\eeas
This together with \eqref{eq:PM_hypothesis} implies that $\frac{\int_{t_0}^T \|v(t; u_0)\|^2_\alpha \d t}{\int_0^T \|u_{\s}(t; u_0)\|^2_\alpha \d t} < 1$ for all $T \ge t_1$. We obtain thus that $Q(T; u_0) < 1$ when $T$ is sufficiently large. The proof is complete.

\ep

\section{Non-Markovian Stochastic Reduced Equations} \label{s:reduction}

As mentioned in Introduction, the practical aspects of the reduction problem of a deterministic dynamical system to its corresponding (local) center or center-unstable manifold has been well investigated in certain  finite- and infinite-dimensional settings; see {\it e.g.} \cite{BK98, Car81, DS06, EvP04, Far01, Kuznetsov04, Hen81, HI10,MW05, MR09, Pot11, PR06}. {\mkk In  \cite{AI98, AX95,Boxler89, Boxler91,NL91,XR96} certain extensions of such reduction techniques have been considered for  finite dimensional RDSs generated by random differential equations (RDEs) or SDEs}. {\mkk Effective reduction procedures to  stochastic critical or other local invariant  manifolds}, have been nevertheless much less explored for RDSs associated with SPDEs; {\mkk see however \cite{CLR01,WD07}}.  

{\mkk A major drawback of any reduction procedure based on local invariant manifolds relies in their local nature which, as explained in Introduction,  is somewhat incompatible with large excursions  of the SPDE solutions caused by white  noise. This is particularly constraining away from the critical value ($\lambda>\lambda_{c}$) where the large excursions are typically further amplified by the nonlinear effects. As we will see in applications the theory of parameterizing manifolds introduced in Section ~\ref{ss:PM_def} in this article  becomes then an interesting substitutive concept. }

In this section, an efficient stochastic reduction procedure is presented based on stochastic parameterizing manifolds introduced {\mkk in Sections  ~\ref{ss:PM}  and \ref{Sec_Exist_PMstoch}}. 
The goal is here to derive efficient reduced models to describe the main dynamical features of the amplitudes of the critical modes contained in $\mathcal{H}^{\c}$. This stochastic reduction procedure developed below can be seen as an alternative to the nonlinear Galerkin method where the approximate inertial manifolds (AIMs) used therein are replaced here by the {\mkk parameterizing manifolds introduced in Sections ~\ref{ss:PM} and \ref{Sec_Exist_PMstoch}.} 

Among the differences with the AIM approach, the PM approach seeks for manifolds which provide modeling error of the evolution of $u_{\c}$ in a {\it  mean-square sense}; see {\mkk Proposition ~\ref{lem:PM error} below}. This modeling error is controlled by the product of three terms: the energy of the unresolved modes ({\it i.e.} the unresolved information), the nonlinear effects (associated with the size of the global random attractor), and the quality of the PM. 
{\mkk In particular, in the case of a stochastic inertial manifold, the corresponding parameterization defect $Q(T,\omega;u_0)$ and the modeling error decays to zero as $T\rightarrow \infty$, a manifestation of the fact the parameterization of the small scales by the large ones can be made ``exact'' up to a reminder that becomes negligible for large $T$; see Theorem \ref{Prop_IMisPM}. }

 {\mkk We will see in  Sections  ~\ref{s:Burgers} and  ~\ref{s:Fly}, that parameterizing} manifolds which can be built from the pullback characterizations \eqref{Eq_PBA1uc} or \eqref{Eq_PBA2uc} play a key role in the derivation of efficient stochastic reduced models in the case of large amplitudes. To simplify the presentation, only reduction  associated with parameterizing manifolds based on the pullback characterization \eqref{Eq_PBA1uc}  associated with the {\it one-layer backward-forward system}  \eqref{LLL}, are described below.  Reduced equations {\mkk based on parameterizing manifolds obtained via the pullback} characterization \eqref{Eq_PBA2uc} associated with the {\it two-layer backward-forward system}  \eqref{LLL_new} {\mkk are dealt with in Section ~\ref{s:Fly} in the case where $\mathcal{H}^{\c}$ contains both} stable and unstable modes.

The  reduced equations {\mkk based on PMs associated with the backward-forward system \eqref{LLL},} are low-dimensional SDEs arising typically with random coefficients which convey  {\it extrinsic memory effects} \cite{HO07, Hai09} expressed in terms of decay of correlations (see Lemma ~\ref{lem:Mn}),  making the stochastic reduced equations genuinely non-Markovian; {\mkk see Eqns. ~\eqref{bilinear pm} below}.  These random coefficients  involve the {\mkk the past of the noise path} and  exponentially decaying terms depending in the self-adjoint case on the gap between some linear combinations of the eigenvalues associated with the low modes and the eigenvalues  associated with high modes. These gaps  correspond exactly to those arising in the \eqref{NR}-condition; see  \eqref{general Mn recall} and Remark \ref{Rem_decaycorr_gene} below.

{\mkk Throughout} this section it is assumed that the RDS associated with the SPDE \eqref{SEE}, possesses a random global attractor   which pullback attracts  deterministic bounded sets of  $\mathcal{H}_{\alpha}$ so that in particular  the framework of Section \ref{ss:PM} can be applied.

\subsection{Low-order stochastic reduction procedure based on parameterizing manifolds}\label{ss_reduced_pm}

In this subsection, we present a stochastic reduction procedure based on parameterizing manifolds associated with ~\eqref{SEE}, which is intended to be used to model the dynamics of a given number of resolved modes.  As before, the $m$-dimensional subspace $\mathcal{H}^{\c}$ is spanned by the $m$ resolved modes, and we assume that these modes are  associated with the first $m$ eigenvalues. 


We first note that by projecting Eq. ~\eqref{SEE} onto the subspace $\mathcal{H}^{\c}$ spanned by the resolved modes,  we get:
\be\label{Eq_reduced_exact}
\mathrm{d} u_{\c}=\Big( L_\lambda u_{\c} +P_{\c} F(u_{\c} +u_{\s})\Big) \mathrm{d}t + \diffusion u_{\c} \circ \mathrm{d} W_t,
\ee
where $u_{\s}=P_{\s} u_{\lambda}$. This equation is the {\it exact} reduced equation of Eq. ~\eqref{SEE} but is not written  in a closed form since $u_{\s}$ needs to be known in order to determine $u_{\c}$.

Various parameterization strategies of the unresolved variable $u_{\s}$ can be imagined to derive a closed version of \eqref{Eq_reduced_exact} from which an approximation of the dynamics of $u_{\c}$  is sought. For the reduction procedure considered here, the unresolved variable $u_{\s}$ is parameterized in terms of the resolved variable $u_{\c}$ through a PM function $\widehat{h}^{\mathrm{pm}}_\lambda: \mathcal{H}^{\c} \times \Omega \rightarrow \mathcal{H}^{\s}_\alpha$ to be determined. 

When such a PM function is given,  the PM-based  (abstract) reduced equation for the resolved modes is then given  by:\footnote{ When $\widehat{h}^{\mathrm{pm}}_\lambda$ corresponds to a stochastic inertial manifold then the asymptotic behavior of   $P_{\c} u_\lambda(t, \omega; u_0)$ can be derived from  the asymptotic behavior of the solutions of Eq. ~\eqref{SDE abstract-pm}. In this case, Eq.  ~\eqref{SDE abstract-pm} corresponds  to the {\it inertial form} of Eq. ~\eqref{SEE}  in the language  of inertial manifold theory; see {\it e.g.} \cite{CFNT89,Tem97}  for $\sigma=0$.}
\begin{equation} \label{SDE abstract-pm}
\mathrm{d}\xii = \Big( L_\lambda^{\c} \xii +P_{\c} F\big(\xii + \widehat{h}^{\mathrm{pm}}_\lambda(\xii, \theta_t\omega)\big)\Big) \mathrm{d}t + \diffusion \xii \circ \mathrm{d} W_t.
\end{equation}

Before proceeding to the derivation of operational form of the above stochastic reduced system, we provide a way to assess the modeling error associated with this reduced system. For this purpose, we rely  on the cohomology relation introduced in \SS ~\ref{sec. conjugacy}. The modeling error is estimated  for the RPDE \eqref{REE 1}, and its corresponding reduced system. 

When a modeling error is small for the RPDE case, it means that the projected RPDE dynamics  on $\mathcal{H}^{\c}$  is well modeled  by its reduced system. The implications are then that according to the cohomology transfer principle mentioned in \SS ~\ref{sec. conjugacy}, any qualitative features of the resolved modes captured by the reduced system for the  RPDE correspond then also to features captured by the reduced system for the SPDE. 

Recall that Eq. ~\eqref{REE 1} is associated with the SPDE \eqref{SEE} via the random transformation \eqref{random transform}. By projecting Eq. ~\eqref{REE 1} onto the subspace $\mathcal{H}^{\c}$, we derive the following exact reduced system associated with Eq. ~\eqref{REE 1}:
\begin{equation} \label{REE reduced}
\frac{\mathrm{d} v_{\c}}{\mathrm{d} t} =L_\lambda^{\c} v_{\c} + z_\diffusion(\theta_t\omega)v_{\c} + P_{\c}G(\theta_t\omega,\, v_{c} + v_{\s}), 
\end{equation}
where $v_{\c} := P_{\c} v_\lambda$, $v_{\s}:= P_{\s}v_\lambda$. Recall also that $G(\omega,\, v) := e^{-z_\diffusion(\omega)} F(e^{z_\diffusion(\omega)}v)$. 

For any given PM function $\widehat{h}^{\mathrm{pm}}_\lambda: \mathcal{H}^{\c} \times \Omega \rightarrow \mathcal{H}^{\s}_\alpha$ associated with the SPDE \eqref{SEE}, let us define $h_\lambda^{\mathrm{pm}}:\mathcal{H^{\c}} \times \Omega \rightarrow \mathcal{H}^{\s}_\alpha$ by:
\bea
h_\lambda^{\mathrm{pm}}(\xi, \omega) := e^{-z_\sigma(\omega)} \widehat{h}_\lambda^{\mathrm{pm}}(e^{z_\sigma(\omega)}\xi, \omega), \quad \xi \in \mathcal{H}^{\c}, \; \omega \in \Omega. 
\eea
Similar to the reduced equations \eqref{SDE abstract-pm} associated with $\widehat{h}^{\mathrm{pm}}_\lambda$, we have the following reduced equation (for the RPDE) based on  the parameterization  function $h_\lambda^{\mathrm{pm}}$:
\begin{equation} \label{RDE abstract-pm}
\frac{\mathrm{d}\xii}{\d t} = L_\lambda \xii +  z_\diffusion(\theta_t\omega) \xi + P_{\c} G\big(\theta_t \omega, \xii + h^{\mathrm{pm}}_\lambda(\xii, \theta_t\omega)\big).
\end{equation}

The following result provides {\mkk an estimate of} the modeling error ---  in a {\mkk mean-square sense}  --- associated with the reduced system \eqref{RDE abstract-pm} {\mkk for the modeling  of the RPDE  dynamics projected onto the resolved modes, when $\lambda>\lambda_{\c}$}.

\bprop  \label{lem:PM error}

Let $v_\lambda(t,\omega)$ be a solution of the RPDE \eqref{REE 1} on its global random attractor {\mkk $\widetilde{\mathcal{A}}_{\lambda}$ for $\lambda>\lambda_{\c}$}. Then, the following  modeling error estimate holds:
\bea\label{Eq_model_err}
\frac{1}{T}\int_0^T \Bigl \| \frac{\mathrm{d} v_{\c}(t,\omega)}{\d t} -&\Bigl( L_{\c} v_{\c}(t,\omega)   +  z_\sigma(\theta_t\omega) v_{\c}(t,\omega)  +P_{\c} G \bigl( \theta_t \omega, v_{\c}(t,\omega) + h^{\mathrm{pm}}_\lambda(v_{\c}(t,\omega), \theta_t\omega) \bigr) \Bigr) \Bigr \|^2 \d t\\
&\leq \frac{1}{T} \underset{t\in[0,T]}\max\Big( l(G,t,\omega )\Big) \widetilde{Q}(T,\omega; v_0) |v_{\s}(\cdot,\omega)|^2_{L^2(0,T;\mathcal{H}_{\alpha})},
\eea
where $v_{\c} = P_{\c}v_\lambda$, $v_{\s} = P_{\s}v_\lambda$, $v_0 = v_{\lambda}(0,\omega)$, $l(G,t,\omega ):=\mathrm{Lip}^2(G|_{\widetilde{\mathcal{A}}(\theta_t \omega)})$, and
\be \label{eq:REE_Q}
\widetilde{Q}(T, \omega; v_0) := \frac{\int_0^T  \bigl \|v_{\s}(t,\omega;v_0) - h_\lambda^{\mathrm{pm}}(v_{\c}(t,\omega; v_0), \theta_t \omega) \bigr \|_\alpha^2 \, \d t }{ \int_0^T \|v_{\s}(t,\omega;v_0)\|_\alpha^2 \, \d t},  \quad T > 0.
\ee

\eprop

\bp 
First note that  $|v_{\s}(\cdot,\omega)|_{L^2(0,T;\mathcal{H}_{\alpha})}\neq 0$ when $\lambda >\lambda_{\c}$ and $T>0$, since the trivial steady state is unstable.

The proof is then based on the following basic inequality:
\bea \label{Eq_PM error}
\frac{1}{T}\int_0^T  & \bigl \| P_{\c} G\big(\theta_t\omega, v_{\c}(t,\omega) + v_{\s}(t,\omega)\big) - P_{\c} G\big(\theta_t \omega, v_{\c}(t,\omega) + h^{\mathrm{pm}}_\lambda(v_{\c}(t,\omega), \theta_t\omega)\big) \bigr \|^2 \d t  \\
& \le \frac{1}{T} \underset{t\in[0,T]}\max\Big( l(G,t,\omega )\Big) \widetilde{Q}(T,\omega; u_0) |v_{\s}(\cdot,\omega)|^2_{L^2(0,T;\mathcal{H}_{\alpha})},
\eea
which results  from a straightforward application of a Lipschitz estimate of the nonlinearity assessed on the global attractor $\widetilde{\mathcal{A}}_{\lambda}$, and by using the definition of the ratio $\widetilde{Q}(T, \omega; u_0)$ given by \eqref{eq:REE_Q}. 

Note also that from the exact reduced equation \eqref{REE reduced} for the RPDE \eqref{REE 1}, we have 
\beas
P_{\c} G\big(v_{\c}(t,\omega) + v_{\s}(t,\omega)\big)  = \frac{\mathrm{d} v_{\c}(t,\omega)}{\d t} - L_{\c} v_{\c}(t,\omega) - z_\sigma(\theta_t\omega) v_{\c}(t,\omega), \quad t \ge 0, \; \omega \in \Omega.
\eeas
This identity used  in \eqref{Eq_PM error} leads then trivially to \eqref{Eq_model_err}. 
\ep

\br
It is important to note that the upper bound  derived in \eqref{Eq_model_err} splits the modeling error estimate {\mkk(over any finite time interval)} into  the product of three terms whose each of them takes its source in different aspects of the reduction problem: the L$^2$-average of the $\|\cdot\|_{\alpha}$-energy contained in the unresolved modes ({\it i.e.} the {\mkk unknown information}\footnote{related to the dimension of $\mathcal{H}^{\c}$.}), the nonlinear effects (associated with the size of the global random attractor), and the quality of the PM. 
\er
\medskip

We explain now how to derive an operational form of Eq. ~\eqref{SDE abstract-pm} in the case when an explicit analytic formula of $\widehat{h}^{\mathrm{pm}}_\lambda$ is available.  To simplify the presentation, we consider the case  where $L_\lambda$ is self-adjoint, the  \eqref{NR}-condition is satisfied and where $\widehat{h}^{(1)}_\lambda$ given by \eqref{Eq_PBA1uc} provides a PM.
We will see  in Section \ref{s:Burgers} that  $\widehat{h}^{(1)}_\lambda$ can indeed provides a PM for a broad class of regimes in the case of  a stochastic Burgers-type equation.

The derivation of an effective reduced equations from Eq. ~\eqref{SDE abstract-pm} based on $\widehat{h}^{(1)}_\lambda$,  can be then articulated according to the three steps outlined  below for the {\it self-adjoint case} when the \eqref{NR}-condition is satisfied.

\bi

\item[{\bf Step 1}.] {\bf Expansion of  $\widehat{h}^{(1)}_\lambda(\xii, \theta_t\omega)$.} 
From the discussion of Section \ref{s:hyperbolic}, since $\widehat{h}^{(1)}_\lambda$ takes the same form as $\mathfrak{I}_\lambda$ given by \eqref{LP integral}, we deduce that
$\widehat{h}^{(1)}_\lambda$  can be decomposed as  in  \eqref{diagonal case SEE} given in Corollary ~\ref{Cor:Taylor}.
Namely, $\widehat{h}^{(1)}_\lambda$ can be written as
\bea \label{h1-expansion-a}
\widehat{h}^{(1)}_\lambda(\xii, \omega)= \sum_{n=m+1}^\infty \widehat{h}^{\mathrm{(1)},n}_\lambda(\xii,\,
\omega) e_n,  \quad \Forall \, \xi \in \mathcal{H}^{\c}, \, \omega \in \Omega,
\eea
where $m$ is the dimension of the subspace $\mathcal{H}^{\c}$ spanned by the resolved modes, and $\widehat{h}^{\mathrm{(1)},n}_\lambda(\xii,\, \omega)$ is given for each $n \ge m+1$ by
\bea  \label{h1-expansion}
\widehat{h}^{\mathrm{(1)},n}_\lambda(\xii, \omega) &:= \sum_{(i_1, \cdots, i_k )\in \mathcal{I}^k }y_{i_1} \cdots y_{i_k} \langle F_k(e_{i_1}, \cdots, e_{i_k}), e_n \rangle \M_n^{i_1,\cdots,i_k}(\omega, \lambda).
\eea
Here, $\mathcal{I}=\{1,\cdots,m\}$,   $y_i =  \langle \xii, e_i \rangle$ for each $i \in \mathcal{I}$ with  $\langle \cdot, \cdot \rangle$ denoting  the inner-product in the ambient Hilbert space $\mathcal{H}$, and for each $(i_1, \cdots, i_k ) \in \mathcal{I}^k$, the $\M_n^{i_1,\cdots,i_k}(\omega, \lambda)$-term is given by:
\be  \label{general Mn recall}
 \boxed{\M_n^{i_1,\cdots,i_k}(\omega, \lambda) := \int_{-\infty}^0 e^{ \bigl( \sum_{j=1}^{k}\beta_{i_j}(\lambda) -
\beta_n(\lambda) \bigr)s + \diffusion (k-1)W_s(\omega)}\, \mathrm{d}s,}
\ee
where  each $\beta_{i_j}(\lambda)$ denotes the eigenvalue associated with the corresponding mode $e_{i_j}$ in $\mathcal{H}^{\c}$, and $\beta_n(\lambda)$  denotes the eigenvalue associated with the unresolved mode $e_n$ in $\mathcal{H}^{\s}$. {\mkk Note that each such terms is well-defined due to the  \eqref{NR}-condition and \eqref{Wt control}.}

\item[{\bf Step 2}.] {\bf Effective form of the reduced equation  \eqref{SDE abstract-pm}.}
We seek for $\xii(t,\omega)$, solution of \eqref{SDE abstract-pm},  into its  Fourier expansion:
\bea \label{expansion}
\xii(t, \omega) = \sum_{i=1}^m y_i(t, \omega) e_i, 
\eea
where again $m$ is the dimension of the subspace $\mathcal{H}^{\c}$. 
The effective reduced equation which rules the evolution of the amplitude variables  $y_i(t, \omega)$ is sought  
by projecting Eq. ~\eqref{SDE abstract-pm} onto each of the $m$ resolved modes\footnote{The projection is done by taking the inner product in  the ambient Hilbert space $\mathcal{H}$ on both sides of Eq. ~\eqref{SDE abstract-pm} with each of the $m$ resolved modes $e_1, \cdots, e_m$.}. In that respect,  the expression of $\widehat{h}^{(1)}_\lambda(\xii, \theta_t\omega)$ 
given by Eq. ~\eqref{h1-expansion-a}
 plays a key role in the derivation of the resulting $m$-dimensional reduced system of SDEs. Its use in the nonlinear terms introduce indeed nonlinear interactions between the unresolved and resolved modes as well as self-interactions between the unresolved modes which are {\mkk both} absent when $\widehat{h}^{\mathrm{pm}}=0$ in  \eqref{SDE abstract-pm}. 
 \ei

The nonlinear interactions just mentioned are clearly the crucial point to resolve in order to derive efficient reduced models. {\mkk This explains in part why  a PM is required to perform better than the ``linear manifold''\footnote{{\mkk corresponding to $h^{\mathrm{pm}}=0$, {\it i.e.} $\mathfrak{M}=\mathcal{H}^{\c}$}.} from  its definition;
see \eqref{Eq_PM_reduces_ERR}.} We will see in applications that the parameterization defect \eqref{Eq_QualPM} is a useful indicator to measure how good these nonlinear interactions are resolved by a PM candidate.

The above reduction procedure is illustrated below in \SS ~\ref{Sec_abstract-example} on an abstract example when the nonlinearity consists   of a bilinear term. A stochastic Burgers-type equation will serve as a  concrete example in Section ~\ref{s:Burgers}.

In the case where no analytic expression for $\widehat{h}_\lambda^{\mathrm{pm}}$ is available, the vector field in 
\eqref{SDE abstract-pm} is built simultaneously with $\xi(t,\omega)$, as the time $t$ evolves. The parameterizing manifold function is thus computed ``on the fly."  This procedure is described  in details in Section ~\ref{s:Fly} on an example, where the parameterizing manifold is taken to be $\widehat{h}_\lambda^{(2)}$ given by \eqref{Eq_PBA2uc}.

\br  \label{rmk:det case}
It is important to note that the low-order reduction procedure described above also applies when Eq. ~\eqref{SEE} is a deterministic PDE, {\it i.e.}, $\diffusion  =0$. In this case, all the expressions become deterministic. We refer to \cite{CLW13} for the description of the deterministic  theory of parameterizing manifolds.

\er

\subsection{An abstract example of PM-based reduced system} \label{Sec_abstract-example}
To illustrate how to apply the procedure described above, we consider the case where the nonlinearity consists of a bilinear term, {\it i.e.} when $F = F_2$. 
Here again to fix the ideas,  we will use $\widehat{h}^{(1)}_\lambda$ associated with the {\it one-layer backward-forward system}   ~\eqref{LLL} to be the PM candidate, and we assume the same working assumptions as above.

 First, we expand $\xi(t, \omega)$ and $\widehat{h}^{(1)}_\lambda(\xi, \theta_t\omega)$  according to \eqref{expansion} and \eqref{h1-expansion-a}, respectively. Then by taking the $\mathcal{H}$-inner product on both sides of Eq. ~\eqref{SDE abstract-pm} with the resolved mode $e_j$  for $1\le j \le m$, we obtain for $\xi = \xi(t,\omega)$:
\bea \label{bilinear y}
\mathrm{d}y_j = \beta_{j}(\lambda) y_j \mathrm{d}t  & +   \Big \langle P_{\c} F_2 \Big(\sum_{i=1}^m y_i e_i + \sum_{n=m+1}^\infty \widehat{h}^{(1),n}_\lambda({\HL \xi},\, \theta_t\omega) e_n\Big), e_j  \Big \rangle \mathrm{d}t   +  \diffusion y_j \circ \mathrm{d}W_t.
\eea

Let us now introduce the following {\it nonlinear interaction coefficients}:
\bea \label{coef B}
B^l_{pq} := \langle F_2(e_{p}, e_{q}), e_l \rangle, \Forall p, \ q, \  l \in \mathbb{N}^\ast.
\eea
These coefficients  describe three types of  nonlinear interactions when projected against a given $e_l$-mode:  the  self-interactions of the resolved modes (when $p,q \in \{1,...,m\}$), the cross-interactions between the resolved and unresolved modes (when $p \in \{1,...,m\}$  and $q\geq m+1$, or {\it vice versa}), and  the self-interactions of the unresolved modes (when $p,q\geq m+1$).

We obtain thus for any $j\in\{1,...,m\}$, that:
\bea \label{F2 expansion}
&\Big \langle P_{\c} F_2 \Big(\sum_{i=1}^m y_i e_i + \sum_{n=m+1}^\infty \widehat{h}^{(1),n}_\lambda({\HL \xi},\, 
\theta_t\omega) e_n\Big), e_j  \Big \rangle\\
 & = \sum_{i_1, i_2 = 1}^m y_{i_1}y_{i_2}\langle F_2(e_{i_1}, e_{i_2}), e_j \rangle \\
& \hspace{1em} + \sum_{i = 1}^m \sum_{n = m+1}^\infty y_{i} \widehat{h}^{(1),n}_\lambda(\xi, \theta_t\omega) \Big( \langle F_2(e_{i}, e_{n}), e_j \rangle  + \langle F_2(e_{n}, e_{i}), e_j \rangle \Big) \\
& \hspace{1em} + \sum_{n_1, n_2 = m+1}^\infty \widehat{h}_\lambda^{(1), n_1}(\xi, \theta_t\omega) \widehat{h}_\lambda^{(1),n_2}({\HL \xi}, \theta_t\omega) \langle F_2(e_{n_1}, e_{n_2}), e_j \rangle \\
 & = \sum_{i_1, i_2 = 1}^m  B^j_{i_1 i_2} y_{i_1}y_{i_2} + \sum_{i = 1}^m \sum_{n = m+1}^\infty ( B^j_{i n}   + B^j_{n i}) y_{i} \widehat{h}_\lambda^{(1),n}({\HL \xi}, \theta_t\omega)  \\
& \hspace{1em} + \sum_{n_1, n_2 = m+1}^\infty B^j_{n_1 n_2} \widehat{h}_\lambda^{(1), n_1}({\HL \xi}, \theta_t\omega) \widehat{h}_\lambda^{(1),n_2}({\HL \xi}, \theta_t\omega).
\eea

Using the  expression of the nonlinear interactions  coefficients given in \eqref{coef B}, we  get from \eqref{h1-expansion} that
\bea \label{coef recast}
\widehat{h}_\lambda^{(1),n}({\HL \xi}, \theta_t\omega) = \sum_{i_1, i_2 =1}^m B^n_{i_1i_2} \M^{i_1i_2}_n(\theta_t\omega, \lambda) y_{i_1}y_{i_2}.
\eea
Now, using  \eqref{coef recast} in \eqref{F2 expansion}, we find {\HL by rearranging the terms:}
\bea \label{F2 expansion-2}
&\Big \langle P_{\c} F_2 \Big(\sum_{i=1}^m y_i e_i + \sum_{n=m+1}^\infty \widehat{h}^{(1),n}_\lambda({\HL \xi},\, 
\theta_t\omega) e_n\Big), e_j  \Big \rangle\\
 & = \sum_{i_1, i_2 = 1}^m  B^j_{i_1 i_2} y_{i_1}y_{i_2} + \sum_{i, i_1, i_2=1}^m \sum_{n = m+1}^\infty ( B^j_{i n}   + B^j_{n i})B^n_{i_1 i_2} \M_n^{i_1i_2}(\theta_t\omega, \lambda) y_{i_1}y_{i_2} y_{i} \\
& \hspace{1em} + \sum_{n_1, n_2 = m+1}^\infty \sum_{i_1,i_2 = 1}^m \sum_{i_3, i_4 = 1}^m B^j_{n_1 n_2}B_{i_1 i_2}^{n_1} B_{i_3 i_4}^{n_2} \M_{n_1}^{i_1i_2}(\theta_t\omega,\lambda) \M_{n_2}^{i_3i_4}(\theta_t\omega,\lambda)y_{i_1}y_{i_2}y_{i_3}y_{i_4}.
\eea

By reporting \eqref{F2 expansion-2} in \eqref{bilinear y}, the following PM-based system of \emph{effective reduced equations} is now  obtained {\mkk as}:
\begin{equation} \label{bilinear pm} 
\boxed{
\begin{aligned}
 \mathrm{d}y_j  & = \bigg (\beta_{j}(\lambda) y_j   + \overbrace{\sum_{i_1, i_2 = 1}^m  B^j_{i_1 i_2} y_{i_1}y_{i_2}}^{(a)} \\
& + \overbrace{\sum_{i, i_1, i_2=1}^m \sum_{n = m+1}^\infty ( B^j_{i n}   + B^j_{n i})B^n_{i_1 i_2} \M^{i_1 i_2}_n(\theta_t\omega, \lambda) y_{i_1}y_{i_2} y_{i}}^{(b)} \\ 
&+ \underbrace{\sum_{n_1, n_2 = m+1}^\infty \sum_{i_1,i_2 = 1}^m \sum_{i_3, i_4 = 1}^m B^j_{n_1 n_2}B_{i_1 i_2}^{n_1} B_{i_3 i_4}^{n_2} \M^{i_1i_2}_{n_1}(\theta_t\omega,\lambda) \M^{i_3i_4}_{n_2}(\theta_t\omega,\lambda)y_{i_1}y_{i_2}y_{i_3}y_{i_4}
}_{(c)} \bigg) \, \mathrm{d} t  \\
& + \diffusion y_j \circ \mathrm{d}W_t, \quad \ 1 \le j \le m.
\end{aligned}
}
\end{equation}
where $B^l_{pq}$ are defined in \eqref{coef B}, and $\M_n^{i_k i_l}(\theta_t\omega, \lambda)$ are given in \eqref{general Mn recall}.

Note that in many situations, only finitely many terms $B^r_{pq}$ appearing in $(b)$ and $(c)$ of the above equations \eqref{bilinear pm} are nonzero; see Section ~\ref{s:Burgers} 
and \cite{CLW13,CGLW}.  In the case  where infinitely many  of such terms would be nonzero, the summation should be obviously truncated {\mkk for} sufficiently large $n$ in practice.

\subsection{PM-based reduced systems  as non-Markovian SDEs}\label{Sec_PM-based}
In the effective reduced equations \eqref{bilinear pm} described above, the terms in $(a)$ reflect nonlinear self-interactions between the resolved modes in  $\mathcal{H}^{\c}$ via the coefficients $B^j_{i_1 i_2}$. On the other hand, the terms in $(b)$ {\mkk (coming from expansion of $F_2(\xii(t,\omega), \widehat{h}_\lambda^{(1)}(\xii(t,\omega), \theta_t\omega))$, $\xii \in \mathcal{H}^{\c}$)} reflect nonlinear  cross-interactions between the resolved  modes in $\mathcal{H}^{\c}$ and the unresolved ones\footnote{{\mkk as parameterized by $\widehat{h}_\lambda^{(1)}$.}} in $H^{\s}_\alpha$; and the terms in $(c)$ {\mkk (coming from expansion of $F_2(\widehat{h}_\lambda^{(1)}(\xii(t,\omega), \theta_t\omega), \widehat{h}_\lambda^{(1)}(\xii(t,\omega), \theta_t\omega))$)} reflect self-interactions between the unresolved modes. Both $(b)$ and $(c)$ involve the random coefficients $\M^{i_k i_l}_n(\omega, \lambda)$.

 These random coefficients  involve the history of the random forcing and  exponentially decaying terms, depending in the self-adjoint case, on the gap between some linear combinations of the eigenvalues associated with the low modes and the eigenvalues  associated with high modes. These gaps  correspond exactly to those arising in the \eqref{NR}-condition; see  \eqref{general Mn recall} and Remark \ref{Rem_decaycorr_gene} below. 
 
 The resulting reduced equations are  thus low-dimensional SDEs arising typically with random coefficients which convey  {\it extrinsic memory effects} \cite{HO07, Hai09} expressed in terms of decay of correlations,  making the stochastic reduced equations genuinely non-Markovian \cite{Hai09} {\mkk provided that $\sigma$ lives in some admissible range; see Lemma ~\ref{lem:Mn} below}. {\mkk It is worthwhile noting that in such cases, the ``noise bath" is essential for these coefficients  to exhibit decay of correlations, the $\M^{i_k i_l}_n(\omega, \lambda)$-terms being reduced to a constant when $\sigma=0$;  see Remark \ref{rmk finiteness of xi} (4).}

The precise results about the decay of correlations  of these terms are presented below in the case where 
 the resolved modes in $\mathcal{H}^{\c}$ correspond to a same eigenvalue with multiplicity $m$, {\it i.e.,}    when 
\be\label{Eq_mult=m}
 \beta_1(\lambda) = \cdots \beta_m(\lambda)=:\beta_{\ast}(\lambda),
\ee 
in \eqref{general Mn recall}. This is just for convenience, to simplify the notations and the proof provided in Appendix ~\ref{sec:Mn proof}. The changes corresponding to the general case are indicated in Remark ~\ref{Rem_decaycorr_gene}.


\bl  \label{lem:Mn}
{\mkk  
Let $g(\lambda) := k\beta_\ast(\lambda) - \beta_{n}(\lambda),$ and  
\bea \label{eq:g and sigma}
 \qquad  \sigma_*(\lambda): = \frac{\sqrt{2 g(\lambda)}}{k-1}, \qquad \sigma_{\#}(\lambda) := \frac{\sqrt{g(\lambda)}}{k-1}.
\eea


Then for  $\lambda$ such that $g(\lambda)>0$, the $\M_n(\omega, \lambda)$-term defined by
\be \label{Mn recall2}
\M_n(\omega, \lambda) := \int_{-\infty}^0 e^{g(\lambda) s + \diffusion (k-1)W_s(\omega)}\, \mathrm{d}s,
\ee
satisfies the following properties:}

\bi

\item[(i)] The expectation of $\M_n(\cdot, \lambda)$ exists if and only if $\sigma < \sigma_{*}$, and is given {\hh in that case} by
\bea \label{expectation Mn}
\mathbb{E}(\M_{n}(\cdot, \lambda)) = \frac{1}{g(\lambda) - (k-1)^2\diffusion^2/2}, \quad  \sigma < \sigma_{*}.
\eea 

\item[(ii)] The variance of $\M_n(\cdot, \lambda)$ exists if and only if $\sigma < \sigma_{\#}$, and is given {\hh in that case} by
\bea \label{variance Mn}
\mathrm{Var}(\M_{n}(\cdot, \lambda))  = \frac{(k-1)^2\sigma^2}{2(g(\lambda) - (k-1)^2\sigma^2/2 )^2(g(\lambda) - (k-1)^2 \sigma^2)}, \quad  \sigma < \sigma_{\#}.
\eea 

\item[(iii)] The autocorrelation of $\M_n(\theta_t \cdot, \lambda)$ exists if and only if  $\sigma < \sigma_{\#}$, and is given  in that case by
\bea \label{eq:autocorrelation}
R(t)  & :=  \frac{\mathrm{Cov}(\M_{n}(\theta_{s+t} \cdot, \lambda), \M_{n}(\theta_{s} \cdot, \lambda))}{\mathrm{Var}(\M_{n}(\cdot, \lambda))} \\
&  = \exp\Bigl( -\Bigl ( g(\lambda) - (k-1)^2 \sigma^2/2 \Bigr ) |t|\Bigr), \quad t \in \mathbb{R}, \,  0 < \sigma < \sigma_{\#},
\eea
where 
\beas
\mathrm{Cov}(\M_{n}(\theta_{s+t} \cdot, \lambda), \M_{n}(\theta_{s} \cdot, \lambda)) &   =  \mathbb{E}(\M_{n}(\theta_{s+t} \cdot, \lambda)  \M_{n}(\theta_{s} \cdot, \lambda)) \\
& \hspace{3em} -  \mathbb{E}(\M_{n}(\theta_{s+t} \cdot, \lambda)) \mathbb{E}(\M_{n}(\theta_{s} \cdot, \lambda)).
\eeas

\ei

\el

This lemma results from direct application of the  Fubini Theorem, the independent increment property of the Wiener process, and the fact that $\mathbb{E}(e^{\diffusion W_t(\cdot)}) = e^{\diffusion^2|t|/2}$ for any $t\in \mathbb{R}$, as expectation of the geometric Brownian motion generated by $\d S_t = \frac{\sigma^2}{2} S_t \d t + \sigma S_t \d W_t$; see {\it e.g.} ~\cite[Sect.~5.1]{Oksendal98}. For the reader's convenience, a proof is provided in Appendix ~\ref{sec:Mn proof}.

\br\label{Rem_decaycorr_gene}

The results presented in Lemma ~\ref{lem:Mn} also hold for the more general $\M_n^{i_1,\cdots,i_k}$-terms given in \eqref{general Mn recall} with the suitable  changes $g(\lambda)=\sum_{j=1}^{k}\beta_{i_j}(\lambda) -
\beta_n(\lambda) >0$ from \eqref{NR}.  We will see in Section \ref{ss:extreme}, that the size of   $\sum_{j=1}^{k}\beta_{i_j}(\lambda) -
\beta_n(\lambda)$ as  $(i_1, \cdots, i_k )$  varies in  $\mathcal{I}^k$ --- called hereafter the {\mkk NR-gaps} ---  play a dominant role in the contribution of the $\M_n^{i_1,\cdots,i_k}$-terms  
to achieve good modeling performances of $u_{\c}$ by  reduced systems  such as \eqref{bilinear pm}.  
\er

As we will see, the  memory terms \eqref{general Mn recall} can turn out to play an important role in applications in order to derive efficient closed models for  the  dynamics of the low modes;  see Section ~\ref{s:Burgers} and \cite{CGLW}.   We provide for the moment another representation of such terms, which will be useful regarding their computations. It is indeed worth noting that, for each $n \ge m+1$, the $\M_n^{i_1,\cdots,i_k}$-term given by \eqref{general Mn recall}   for a fixed $k$-tuple $(i_1,\cdots,i_k)$,  
corresponds  in fact to the unique stationary solution of the following auxiliary scalar SDE: 
\bea \label{xi eqn} 
\mathrm{d} M = \left(1 - \Big(\sum_{j=1}^{k}\beta_{i_j}(\lambda) -
\beta_n(\lambda)\Big)M \right )\mathrm{d} t -  (k-1) \diffusion M \circ \mathrm{d} W_t.
\eea
Since $\sum_{j=1}^{k}\beta_{i_j}(\lambda) -
\beta_n(\lambda) > 0$, the stationary solution $\M_n^{i_1,\cdots,i_k}(\omega, \lambda)$ of Eq. ~\eqref{xi eqn}  is measurable with respect to the past $\sigma$-algebra $\mathcal{F}_{-}$ generated by the mappings $\{\omega  \mapsto \varphi(\tau, \theta_{-t}\omega) \; \vert \; 0 \le \tau \le t\}$, where $\varphi$ is the cocycle associated with such an SDE; see \cite[Example 2.4.4]{Chueshov02}. The fact that each memory term  as given by \eqref{general Mn recall}  can be represented as the stationary solution to its corresponding Eq. ~\eqref{xi eqn}  simplifies clearly  its computation in practice\footnote{when compared to a direct evaluation using its integral representation \eqref{general Mn recall}.}; cf. Section ~\ref{s:Burgers}.

\br 
In the case where $M_n^{lk}=0$, for $n\geq N>m+1$ $(l,k\in \{1,..,m\})$, by supplementing the equations \eqref{bilinear pm} (satisfied by the  $y_i$-variables) {\attn with} the equations  \eqref{xi eqn} given below for the corresponding coefficients $M_n^{lk}$, 
it can be shown that such an augmented system  in the $(y_i;M_n^{lk})$-variables, leads to well-defined pathwise (local) solutions using the standard theory of existence of solutions to SDEs \cite{Oksendal98}, provided that the initial data for the $M_n^{lk}$-equations are taken to be their corresponding stationary solutions given in \eqref{general Mn recall}. Existence of pathwise solutions to the non-Markovian SDE \eqref{bilinear pm}, can be then deduced from the existence of such solutions to the corresponding augmented (Markovian) system.
\er

\needspace{1\baselineskip}

\br \label{rmk:Mn}
Although the memory terms given by \eqref{general Mn recall}  have exponential decay of autocorrelation (Lemma ~\ref{lem:Mn} (iii)),  they exhibit typically non-Gaussian statistics; see Fig. ~\ref{fig:Mn}. Note that in the reduction of stochastic systems, other extrinsic memory terms with similar properties have been found in the literature \cite{BM13, FS09,  LR12,  Rob08}, although we are not aware of an equivalent of Lemma ~\ref{lem:Mn} for those terms.

\er

\begin{figure}[!hbtp]
   \centering
   \includegraphics[height=8cm, width=13cm]{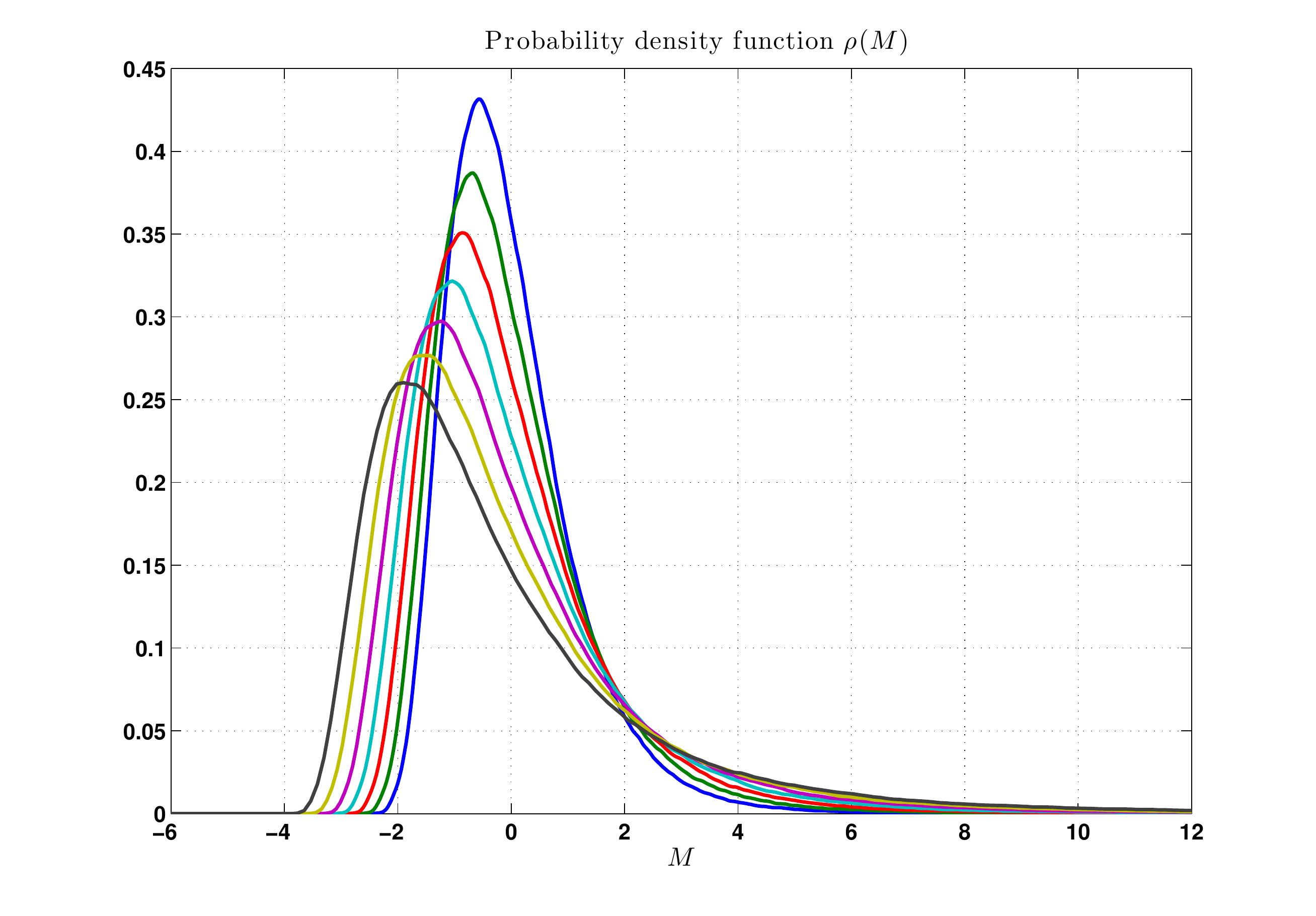}
  \caption{{\footnotesize {\mkk Numerically estimated probability density functions (PDFs) of the $\M_n$-term, with parameters $\beta_\ast(\lambda)=10$ in \eqref{Eq_mult=m}, $k=2$, and  $\beta_n(\lambda)=19.7$, and for various values of $\sigma$.} 
  The range of $\sigma$-values used for this plot is $[0.5 \sigma^{+}, \sigma^{+}]$ where 
  $\sigma^{+} =  \sigma_*(\lambda)/2$.
  The PDF corresponding to the smallest {\mkk  $\sigma$ used for this plot} has the biggest maximum {\mkk value}. The more $\sigma$ increases, the more this maximum decreases and the skewness of the PDF becomes pronounced.}}   \label{fig:Mn}
\end{figure}

As mentioned in the above remark,   extrinsic memory terms of different forms have been already encountered in reduction strategies of finite-dimensional SDEs to random center manifolds; see {\it e.g.} ~\cite{Boxler89}. Extrinsic memory effects also arise in procedures  which consist of deriving  simultaneously both normal forms and center manifold reductions of  SDEs; see for instance \cite{AI98, AX95, NL91} and \cite[Sect. ~8.4.5]{Arnold98}. In {\mkk such a two-in-one strategy},  anticipating terms may arise  --- as integrals involving the future of the noise path  ---  in  both  the  corresponding  random change of coordinates and the resulting normal form.

In \cite{FS09, LR12, Rob08},  pursuing the works of \cite{AI98,AX95}, reduced stochastic equations involving also extrinsic memory terms have been derived mainly in the context of the stochastic slow manifold; see also \cite{BM13}. By seeking for a random change of variables, which typically involves repeated stochastic convolutions,  stochastic reduced equations (different from those  derived in \eqref{bilinear pm}) are obtained to model the dynamics of the slow variables. These reduced equations are also non-Markovian but require a special care in their derivation  to push the anticipative terms (arising in such an approach) to higher order albeit not eliminating them \cite{FS09,Rob08}.

As a comparison, our reduction strategy  is naturally associated with the theory of (stochastic) parameterizing manifolds introduced in this article, and in particular it does not require the existence of  a stochastic slow (or inertial)  manifold. Our approach prevents furthermore the emergence of anticipative terms to any order in the corresponding reduced SDEs.   Memory terms of more elaborated structures than described in Lemma  ~\ref{lem:Mn} (see e.g. \eqref{N_2} or \eqref{N_3} in Section \ref{Sec_Existence_PQ}) can  also arise in  our stochastic reduced equations    built from stochastic PMs  defined as  pullback limits associated with the {\it multilayer backward-forward systems} introduced  in Section  \ref{ss:PM}.

 As illustrated for the stochastic Burgers-type equation analyzed in Section \ref{s:Fly}, stochastic reduced equations  \eqref{SDE abstract-pm}  built from $\widehat{h}^{(2)}_\lambda$ defined by the pullback limit \eqref{Eq_PBA2uc}, convey such memory terms.  As we will see in Section \ref{Sec_Existence_PQ}, these terms arise as repeated compositions of functions involving integrals depending on the past of the noise path driving the SPDE (see e.g. \eqref{N_2} or \eqref{N_3} in Section \ref{Sec_Existence_PQ}), and come with a {\it ``matriochka'' of nonlinear self-interactions} between the low eigenmodes.  We will see  in Sections ~\ref{s:Burgers} and ~\ref{s:Fly} that the memory terms
 arising in  the  PM-based reduced systems --- either from  $\widehat{h}^{(1)}_\lambda$ or $\widehat{h}^{(2)}_\lambda$ ---  play an important  role in the modeling performance achieved by such systems, regarding the SPDE dynamics projected onto the resolved modes.

 \br\label{MZ_rmk}
  The extrinsic memory terms that arise in our approach (or the aforementioned ones) should  not be confused with  {\it intrinsic memory terms} built from the past history of the low modes such as appearing in other reduction strategies; see {\it e.g.}    \cite{Chorin_al02, CKG11, Chorin_al06, Chorin_Hald-book, EMS01, KCG13, kondrashov2013low, harlim, MTV01, McWilliams12, Stinis06, WL13}.   
  Extrinsic memory terms such as \eqref{general Mn recall} take their sources in the ``noise bath'' and in the nonlinear, leading-order self-interactions of the low modes, as projected against the high modes.
In contrast, the intrinsic memory terms such as arising in the Mori-Zwanwig approach \cite{Chorin_al02,Chorin_Hald-book} may be present even when $\sigma=0$ and result from a decomposition of the reduced vector field in an averaged part plus fluctuating components; see \cite{KCG13,Chorin_Hald-book}.
 \er

\section{Application to a stochastic Burgers-type equation: Numerical results} \label{s:Burgers}
We illustrate below, how  the PM-based stochastic reduction procedure introduced above performs in the case of a Burgers-type equation perturbed by a linear multiplicative white noise. 
The performance of the reduction procedure will be assessed mainly on a quantitative level hereafter.  For a dynamically oriented study based on this approach  we refer to  \cite{CGLW}.

The purpose is here  to illustrate that the one-layer backward-forward  system \eqref{LLL} introduced in Section \ref{s:pullback} provides,  for a broad class of regimes,  candidates of parameterizing manifolds, which  lead to reduced systems that model the SPDE dynamics projected onto  the resolved modes with good performances.  The use of the two-layer backward-forward  system \eqref{LLL_new} is analyzed in Section \ref{s:Fly}. In particular, this system will be shown to give access to  PMs of good quality for parameter regimes where the one-layer backward-forward  system \eqref{LLL} fails to provide a PM.

The problem analyzed hereafter consists of the following stochastic initial-boundary value problem on the interval $(0, l)$:
\bea \label{eq:Burgers}
& \mathrm{d} u = \big( \nu u_{xx}  + \lambda u  - \gamma u  u_x\big) \mathrm{d} t + \sigma u \circ \mathrm{d}W_t, \\
& u(0,t) = u(l,t) = 0, \;\;\; t\geq 0,\\
& u(x, 0) = u_0(x), \qquad x\in (0,l),
\eea
where $\nu, \lambda, \gamma$, and $\sigma$ are positive parameters, and $u_0$ is some appropriate initial datum. This problem can be cast into the abstract form \eqref{SEE} which fits the framework of Section \ref{s:preliminary} with $\mathcal{H}=\{u \in L^2(0,l)\;:\; u(0)=u(l)=0 \}$ and $\mathcal{H}_1=H^2(0,l)\cap H_0^1(0,l)$, and we refer to \cite{HW06} for more details when $\sigma=0$. 
In particular, using the cohomology  approach recalled in \SS ~\ref{sec. conjugacy} (see also \cite{FL05}) it can be shown that a random version of the estimates performed in \cite{HW06} ensures the existence of a global dissipative RDS in $H_0^1(0,l)$ associated with \eqref{eq:Burgers} so that the reduction techniques introduced in Section \ref{s:reduction} can be applied. From a physical point of view, we mention that \eqref{eq:Burgers} is a stochastic version of problems arising  in the modeling of  flame front propagation; see {\it e.g.} \cite{berestycki2001meta}. It can also be seen  as inscribed in the long tradition for the study of Burgers turbulence subject to random forces; see {\it e.g.} ~\cite{BFK00,DDT94,E_al00,EV00,Frisch95,Mitra_al05}.

Note that the eigenvalues of the associated deterministic linear problem are given by 
\be\label{Eq_egein_burg}
\beta_n(\lambda) := \lambda - \frac{\nu n^2\pi^2}{l^2},  \;\; n \in \mathbb{N}^\ast,
\ee
 and the corresponding eigenvectors are 
 \be  \label{Eq_eigenmodes}
 e_n(x) := \sqrt{\frac{2}{l}}\sin\Big(\frac{n\pi}{l} x\Big), \;\; x\in(0,l).
 \ee
 The critical value of the control parameter $\lambda$ at which the trivial steady state loses its stability is given by $\lambda_c = \frac{\nu\pi^2}{l^2}$, and the corresponding critical mode is $e_1$. 

We consider below the case where the subspace $\mathcal{H}^{\c}$ is not only spanned by the critical mode $e_1$, but also by the stable mode $e_2$, namely:
\bea  \label{Burgers Hc}
\mathcal{H}^{\c} := \mathrm{span}\{e_1, e_2\}.
\eea
In particular, this choice of $\mathcal{H}^{\c}$ leads to examples of  hyperbolic stochastic {\mkk parameterizing} manifolds mentioned in Section ~\ref{s:hyperbolic}.  

\subsection{Parameterization defect of $\widehat{h}_\lambda^{(1)}$: Numerical estimates} \label{ss:Burgers reduction} 
In this subsection, we first determine analytically the PM candidate {\mkk $\widehat{h}^{(1)}_\lambda$  as given by \eqref{h1-expansion-a}--\eqref{h1-expansion}
for the SPDE \eqref{eq:Burgers}} when $\mathcal{H}^{\c}$ is given by \eqref{Burgers Hc} and then analyze numerically  its parameterization quality. For the regimes {\mkk for which} $\widehat{h}^{(1)}_\lambda$ is inferred to constitute a PM from this numerical analysis,  the corresponding stochastic reduced equations for the resolved modes $e_1$ and $e_2$ are then derived; {\mkk see Eqns. ~\eqref{eq:Burgers reduced} in \SS ~\ref{ss: Reduced equations based on h1}}. The modeling performances of the evolution of $u_{\c}=P_{\c}u$ from  these reduced equations are assessed in \SS ~\ref{ss:extreme}. 

Since the linear part of the SPDE considered here is self-adjoint, the pullback limit $\widehat{h}^{(1)}_\lambda$ associated with the auxiliary system \eqref{LLL} exists if and only if the \eqref{NR}-condition is satisfied;  see \SS ~\ref{ss:PM}.  In that respect, note first that the nonlinear self-interactions of the resolved modes take the following form when projected against the unresolved modes:
\bea  \label{nonlinear interaction}
& \langle e_1(e_2)_x,e_3\rangle = \frac{\sqrt{2} \pi }{l^{3/2}}, \quad \langle e_2(e_1)_x,e_3\rangle = \frac{\pi}{\sqrt{2}l^{3/2}},\quad \langle e_2(e_2)_x,e_4\rangle =  \frac{\sqrt{2} \pi }{l^{3/2}}, \\
& \langle e_i(e_j)_x,e_n\rangle = 0, \qquad i, j \in \{1,\ 2 \}, \; n \ge 5,
\eea
so that with $\mathcal{H}^{\c}$ given in \eqref{Burgers Hc}, the \eqref{NR}-condition reduces here to:
\bea
\beta_1(\lambda) + \beta_2(\lambda) - \beta_3(\lambda) > 0, \qquad 2\beta_2(\lambda) - \beta_4(\lambda) > 0.
\eea
Now, by using the explicit expression of the eigenvalues provided in \eqref{Eq_egein_burg}, it follows that the \eqref{NR}-condition is satisfied for any $\lambda > - \frac{4\nu \pi^2}{l^2}$. For each such $\lambda$, the pullback limit $\widehat{h}_\lambda^{(1)}$  can be obtained by using the formula \eqref{h1-expansion-a}; and thanks to \eqref{nonlinear interaction}, $\widehat{h}_\lambda^{(1)}$ reduces to the following:
\bea \label{eq:h Burgers}
\widehat{h}_\lambda^{(1)}(y_1e_1 + y_2e_2, \omega) = -\frac{3\gamma \pi}{\sqrt{2}l^{3/2}} \M^{12}_3(\omega, \lambda) y_1y_2e_3 -\frac{\sqrt{2} \gamma \pi }{l^{3/2}} \M^{22}_4(\omega, \lambda)  y_2^2 e_4,  \quad y_1, y_2 \in \mathbb{R},
\eea
where the {\it extrinsic memory terms} $\M^{12}_3(\omega, \lambda)$ and $\M^{22}_4(\omega, \lambda)$ are given by
\bea \label{M3M4}
\M^{12}_3(\omega, \lambda) & = \int_{-\infty}^0 e^{[\beta_{1}(\lambda) + \beta_{2}(\lambda) -
\beta_3(\lambda)]s + \diffusion W_s(\omega)}\, \mathrm{d}s, \\
\M^{22}_4(\omega, \lambda) & = \int_{-\infty}^0 e^{[2 \beta_{2}(\lambda) -
\beta_4(\lambda)]s + \diffusion W_s(\omega)}\, \mathrm{d}s.
\eea

{\mkk These terms come here with the nonlinear, self-interactions between the low modes $e_1$ and $e_2$, projected against  respectively the high modes $e_3$ and $e_4$; see \eqref{h1-expansion}.}

\medskip

The  parameterization defect of $\widehat{h}_\lambda^{(1)}$, $Q(T, \omega; u_0)$ as defined by \eqref{Eq_QualPM},  is  {\mkk numerically investigated}  for various values of $\lambda$ and $\sigma$.   In that purpose,  we describe below the numerical schemes used to integrate  Eq.~\eqref{eq:Burgers}  and to  compute the memory terms $\M^{12}_3(\omega, \lambda)$ and $\M^{22}_4(\omega, \lambda)$ present  in the analytic expression of  $\widehat{h}_\lambda^{(1)}$ given by \eqref{eq:h Burgers} above.

 Eq.~\eqref{eq:Burgers} is solved by a semi-implicit Euler scheme where for each time step the nonlinearity $uu_x = (u^2)_x/2$ and the noise term $\sigma u \circ \d W_t$ are treated explicitly, and the other terms are treated implicitly. The Laplacian operator is discretized using the standard second-order central difference approximation. The resulting semi-implicit scheme now reads as follows:
\bea \label{Burgers discrete-1}
u_{j}^{n+1} - u_{j}^{n} = 
 \Big( \nu \Delta_d u_{j}^{n+1} + \lambda  u_{j}^{n+1} + \frac{\sigma^2}{2}u_{j}^{n} - \frac{\gamma}{2} \nabla_d\big((u_{j}^{n})^2 \bigr)\Big)\delta t  + \sigma \zeta_n   u_{j}^{n} \sqrt{\delta t},
\eea
where $u_{j}^{n}$ is the discrete approximation of $u(j\delta x, n\delta t)$, $\delta x$ denotes  the mesh size of the spatial discretization; $\delta t$, the time step;  the discretized Laplacian $\Delta_d$  and the discretized spatial derivative $\nabla_d$   are given by 
\bes
\Delta_d u_{j}^{n}= \frac{u_{j-1}^{n} - 2u_{j}^{n} + u_{j+1}^{n}}{(\delta x)^2}; \qquad \nabla_d \big( (u_j^{n})^2\big) = \frac{ (u_{j+1}^{n})^2 - (u_j^{n})^2 }{\delta x}, \quad j \in \{1, \cdots, N_x - 2\},
\ees  
supplemented by the boundary conditions
$$
u_0^n=u_{N_x}^n=0, 
$$
where $N_x$ is the total number of grid points used for the discretization of the spatial domain $[0, l]$.
The $\zeta_n$  are random variables drawn independently from a normal distribution $\mathcal{N}(0,1)$. 
Note that the additional drift term $\diffusion^2 u_{j}^{n}/2$ in the RHS of \eqref{Burgers discrete-1} is due to the conversion of the Stratonovich noise term $\sigma u\circ \mathrm{d}W_t$ into its  It\^o form. 
 In this section, the simulations are performed for $\delta t=0.01$ and $N_x = 132$ with $l = 2.5\pi$ so that $\delta x \approx 0.06$. The parameters $\nu$ and  $\gamma$ will remain fixed to be respectively $2$ and $0.5$.  The values of the parameters $\lambda$ and $\sigma$ will be specified when necessary.

Let us denote the solution to \eqref{Burgers discrete-1} at time $t = n \delta t$ to be $\mathbf{U}^n$, then after rearranging the terms, equation \eqref{Burgers discrete-1} can be rewritten into the following algebraic system:
\bea \label{Burgers discrete}
\bigl( (1- \lambda \delta t) \mathbf{I} - \nu  \delta t \mathbf{A}  \bigr)\mathbf{U}^{n+1} = ( 1+\frac{\sigma^2}{2}\delta t +\sigma \zeta_n  \sqrt{\delta t}) \mathbf{U}^{n} -  \frac{\gamma}{2} \delta t \mathbf{B}[\mathbf{S}(\mathbf{U}^{n})],
\eea
where $\mathbf{I}$ is the $(N_x-2) \times (N_x-2)$ identity matrix, $\mathbf{A}$ is the tridiagonal matrix associated with the discrete Laplacian $\Delta_d$, $\mathbf{B}$ is the matrix associated with the discrete spatial derivative $\nabla_d$, and $\mathbf{S}(\mathbf{U}^{n})$ denotes the vector whose entries are the square of the corresponding entries of $\mathbf{U}^{n}$. 

Note that the eigenvalues of $\mathbf{A}$ are given by $\frac{2}{(\delta x)^2} \Big( \cos( \frac{j \pi \delta x}{l}) - 1 \Big)$, $j = 1, \cdots, N_x - 2$, and the corresponding eigenvectors are the discretized version of the first $N_x-2$ sine modes $e_1, \cdots, e_{N_x-2}$ given in \eqref{Eq_eigenmodes}. Consequently, the eigenvalues of the matrix $\mathbf{M}:=(1- \lambda \delta t) \mathbf{I} - \nu  \delta t \mathbf{A}$ on the LHS of \eqref{Burgers discrete} can be obtained easily, and the corresponding eigenvectors are still the discretized sine functions. At each time step and for a fixed realization $\omega$ of the $\zeta_n$'s, this algebraic system \eqref{Burgers discrete} can thus be solved using the {\it discrete sine transform}. More precisely, we first compute the discrete sine transform of the RHS; we then divide the elements of the transformed vector by the eigenvalues of $\mathbf{M}$; the inverse discrete sine transform is then performed to find $\mathbf{U}^{n+1}$; see {\it e.g.} ~\cite[{\att Sect. ~3.2}]{Eyre98} for more details. All the numerical experiments performed in this article have been carried out by using the Matlab version \texttt{7.13.0.564} (R2011b), where the discrete sine transform has been handled by using the built-in function \texttt{dst.m}.

\br
Cross checking has been carried out regarding the semi-implicit scheme \eqref{Burgers discrete-1} used to solve the SPDE \eqref{eq:Burgers}.  Simulations of \eqref{eq:Burgers} based on this scheme have been compared with those based on {\mkk a finite difference discretization of the conservative form, which has been adapted to the stochastic context.} For the parameters used, the relative error under the $L^2$-norm between the numerical solutions produced by these two schemes has been found to be  about $1\%$. The scheme \eqref{Burgers discrete-1} has been chosen here since it allows for relatively larger time steps  than required by an energy-preserving scheme to achieve the same accuracy. {\mkk We refer the reader to \cite{AG06, BJ13,burkardt2007reduced,Hou_al06,LR04} for other numerical approximation schemes  of nonlinear SPDEs.}
\er

Given a realization $\omega$ and an initial datum $u_0$, to compute $Q(T, \omega; u_0)$ as defined by \eqref{Eq_QualPM},  the quantity $\widehat{h}^{(1)}_{\lambda}(u_{\c}(t,\omega; u_0), \theta_t \omega)$ has to be evaluated  for  $t \in [0,T]$.  For that purpose and based on \eqref{eq:h Burgers}, the memory terms $\M^{12}_{3}(\theta_t\omega, \lambda)$ and $\M^{22}_{4}(\theta_t\omega, \lambda)$ need to be simulated.   As mentioned in Section \ref{Sec_PM-based},  this problem is reduced to the approximation problem of stationary solutions of the corresponding auxiliary SDEs of type \eqref{xi eqn}. For instance in the case of $\M^{12}_{3}(\theta_t\omega, \lambda)$, an approximation of the stationary solution of 
\bea \label{eq: M3 eqn} 
\mathrm{d} M = \left(1 - (\beta_{1}(\lambda) + \beta_{2}(\lambda) - \beta_{3}(\lambda))M \right )\mathrm{d} t - \diffusion M \circ \mathrm{d} W_t,
\eea
is obtained as follows.  A semi-implicit scheme is adopted to simulate \eqref{eq: M3 eqn}, where  the term $(\beta_{1}(\lambda) + \beta_{2}(\lambda) - \beta_{3}(\lambda))M \mathrm{d} t$ is treated implicitly, and the other terms are treated explicitly; see \cite[Sect. ~12.2]{KP92}. The simulation of Eq. ~\eqref{eq: M3 eqn} is  then performed on an interval $[-s, t]$ for a sufficiently large $s$ to ensure that the stationary regime is  reached on $[0,t]$  so that a good approximation of $\M^{12}_{3}(\theta_t\omega, \lambda)$ is guaranteed. An approximation of $\M^{22}_{4}$ is obtained in a similar way.

\medskip

Having clarified the simulations of $u(t,\omega; u_0)$, 
and of $\widehat{h}^{(1)}_{\lambda}(u_{\c}(t,\omega; u_0), \theta_t \omega)$,  the parameterization defect $Q_{\lambda,\sigma}(T, \omega; u_0)$ defined by \eqref{Eq_QualPM} is first  averaged over $T \in I=[T_1,T_2]=[400, 1000]$,   for $u_0 = 0.1e_1 + 0.2 e_2 + 0.1 e_5$, an arbitrary realization $\omega$ of the Wiener process and for various values of $\lambda$ and $\sigma$ as specified in the caption and legends of Fig. ~\ref{fig:ave_Q}.    The resulting time average $\overline{Q}_{\lambda,\sigma}(\omega;u_0,I)$ is intended to give a first idea about the parameterization defect for different regimes as $\lambda $ and $\sigma$ are varied.

The results are reported in Fig. ~\ref{fig:ave_Q} where the $\lambda$-dependence (for a fixed value of $\sigma$) is observed to be  linear and the $\sigma$-dependence (for a fixed value of $\lambda$)   is observed to be of quadratic type  for small values of  $\sigma$, while  linear for larger values. As one can notice,  the (time average) parameterization quality gets deteriorated as $\sigma$ increases  and $\lambda$ moves away from (above) its critical value.   The time average of the parameterization defect remains however below $1$ for a broad class of parameter regimes, and in particular for a relatively large amount of noise. 

For the corresponding ranges of $\sigma$- and $\lambda$-values where $\overline{Q}_{\lambda,\sigma} <1$, it has been numerically observed that the SPDE \eqref{eq:Burgers} exhibits  two stable stationary solutions $u_\lambda^{\pm}$ which attract --- in a pullback sense --- all other solutions except the trivial steady state, where $u_\lambda^{+}$ has positive amplitude and $u_\lambda^{-}$ has negative amplitude. The initial datum $u_0$ used for Fig. ~\ref{fig:ave_Q}  is in the basin of attraction of  $u_\lambda^{+}$, and other choice of initial data in this basin lead to similar results.  The global shapes of the curves reported in  Fig. ~\ref{fig:ave_Q} have been also observed to be not sensitive to the realization $\omega$ used in the simulations\footnote{Only slight changes of the slopes at the intercept with $y=1$ have been observed when $\omega$ is changed.}. The intercept of these curves with the line $y=1$ has been also observed to do not vary  sensitively as $T_2$ increases.

\begin{figure}[h]
    \subfigure[For different values of $\lambda$]{
    \begin{minipage}[b]{0.4\textwidth}
    \centering
   \includegraphics [height=.8\textwidth, width=1\textwidth, trim = 0 0 0 0, clip]{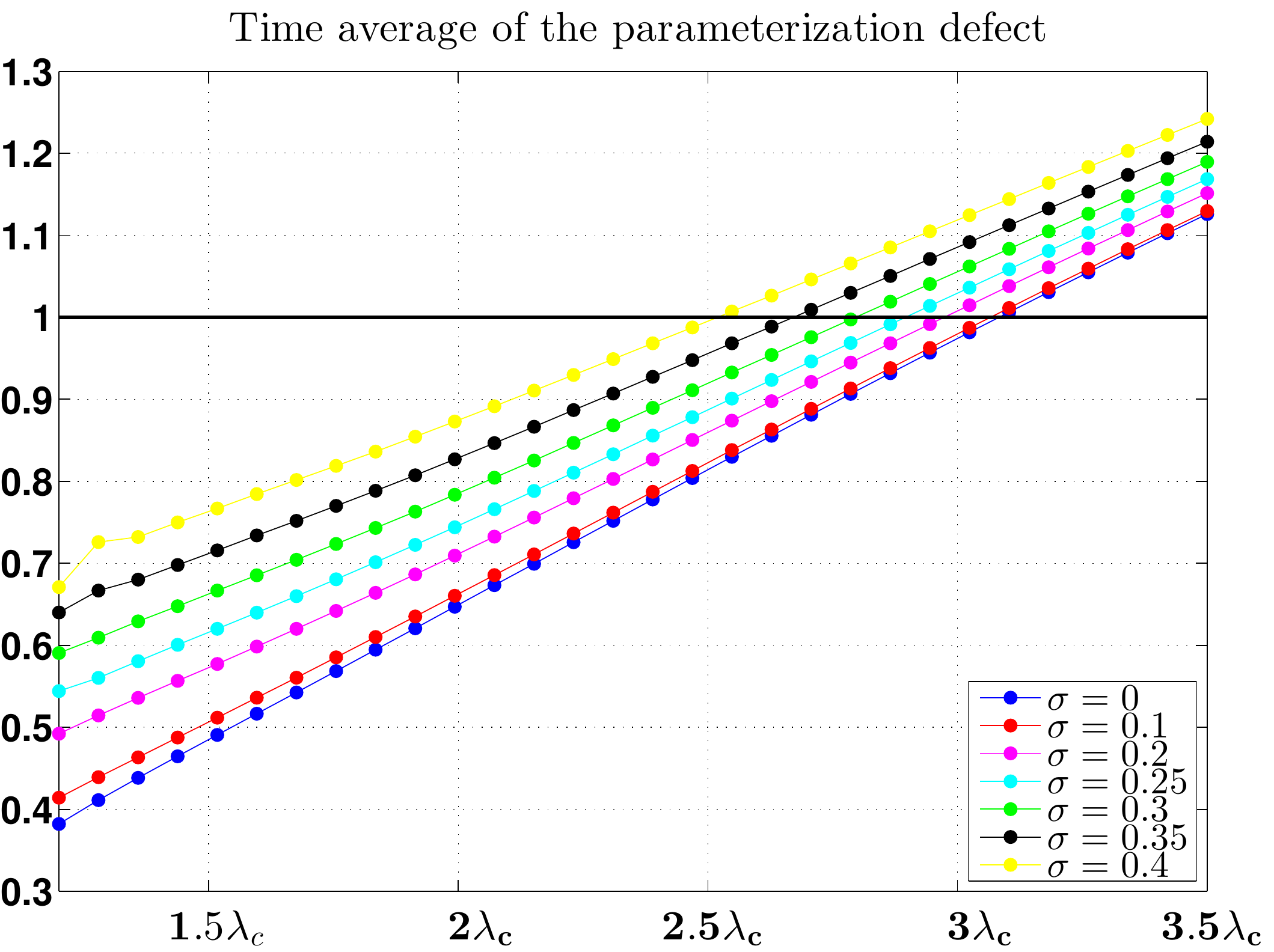}
    \end{minipage}}
    \subfigure[For different values of $\sigma$]{
    \begin{minipage}[b]{0.4\textwidth}
    \centering
   \includegraphics [height=.8\textwidth, width=1\textwidth, trim = 0 0 0 0, clip]{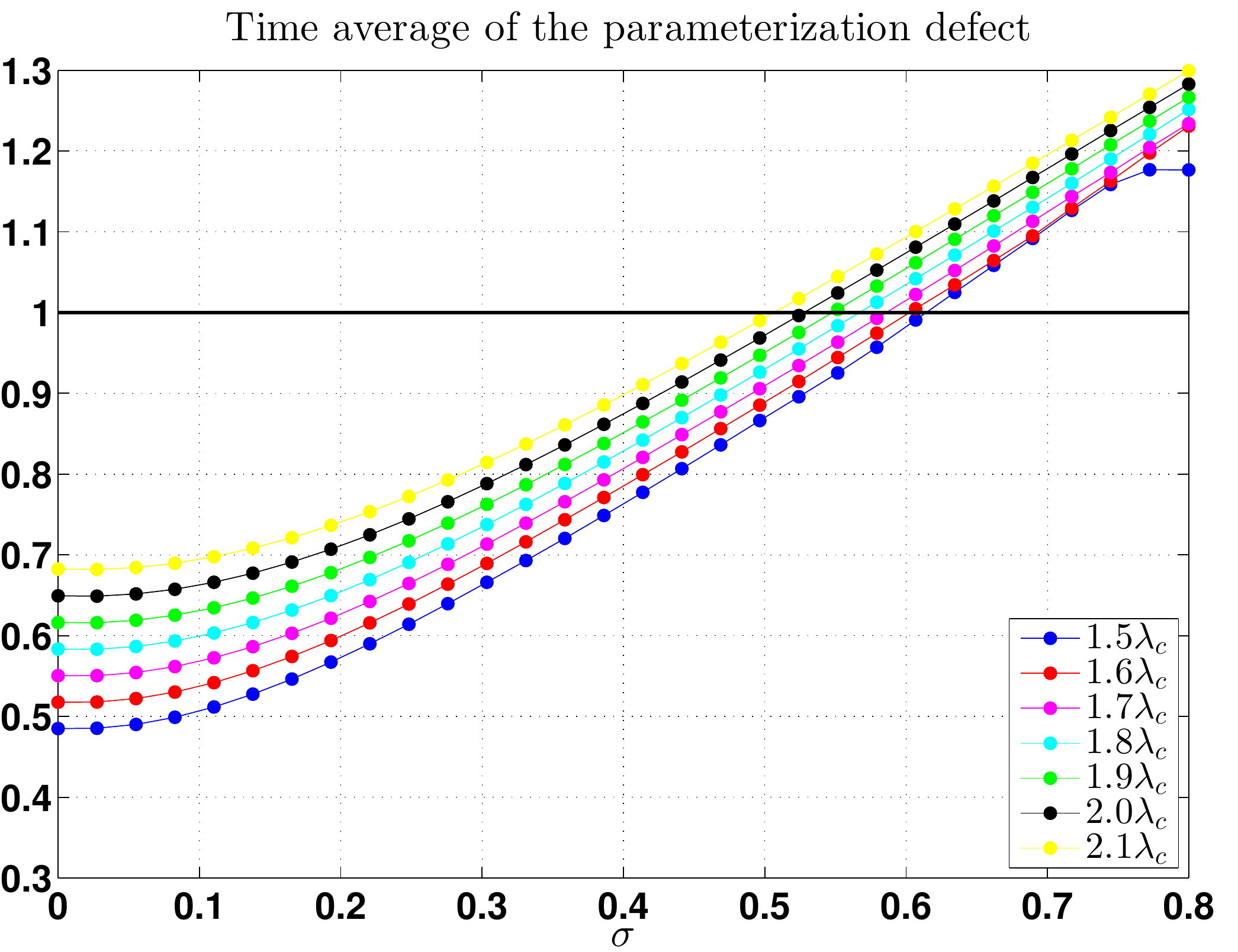}
    \end{minipage}}
\caption{{\footnotesize {\bf Time-average of the parameterization defect over $[T_1, T_2]$.}   Here, $\frac{1}{T_2 - T_1}\int_{T_1}^{T_2} Q(T, \omega; u_0)\, \d T$ is computed  for  $u_0 = 0.1e_1 + 0.2 e_2 + 0.1 e_5$, $T_1= 400$, $ T_2 = 1000$ and a fixed realization: (a) as  $\lambda$ varies in $[1.2\lambda_c, 3.5 \lambda_c]$ for $\sigma$ as indicated in the legend; (b) as $\sigma$ varies in $[0, 0.8]$ for $\lambda$ as indicated in the legend. Other parameters in Eq. ~\eqref{eq:Burgers} are chosen to be  $\gamma = 0.5$, $l = 2.5\pi$ and $\nu = 2$, leading to a critical value of  $\lambda_c = 0.32$. As explained in the text, these results are robust w.r.t. the choices of the realization, the initial datum and $T_2>1000$.}}
\label{fig:ave_Q}
\end{figure}

\begin{figure}[!h]
    \subfigure[For $u_0 = u_0^{+}$]{
    \begin{minipage}[b]{0.4\textwidth}
    \centering
   \includegraphics [height=1\textwidth, width=1.1\textwidth, trim = 0 0 0 0, clip]{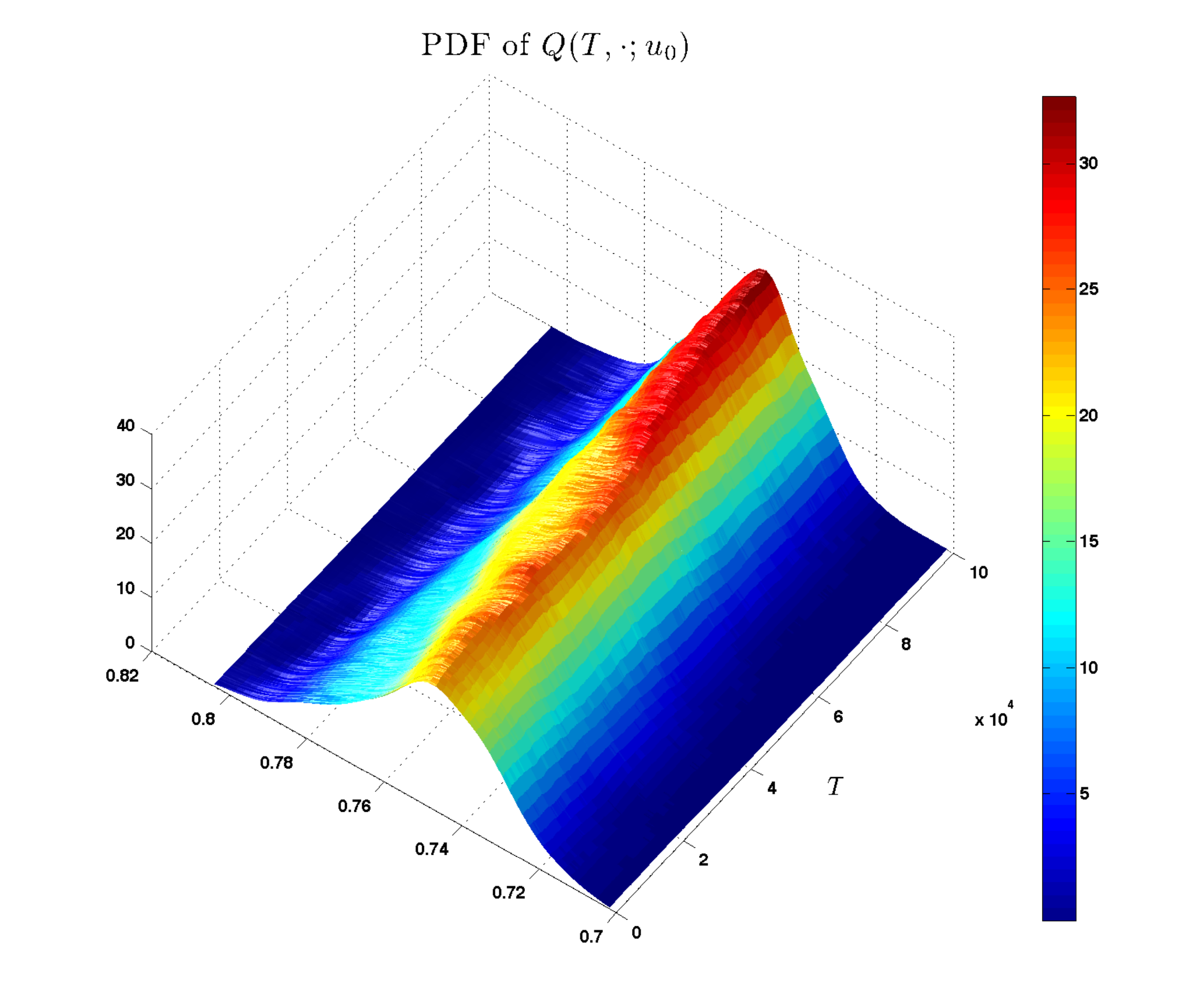}
    \end{minipage}}  
    \subfigure[For $u_0 = u_0^{-}$]{
    \begin{minipage}[b]{0.4\textwidth}
    \centering
   \includegraphics [height=1\textwidth, width=1.1\textwidth, trim = 0 0 0 0, clip]{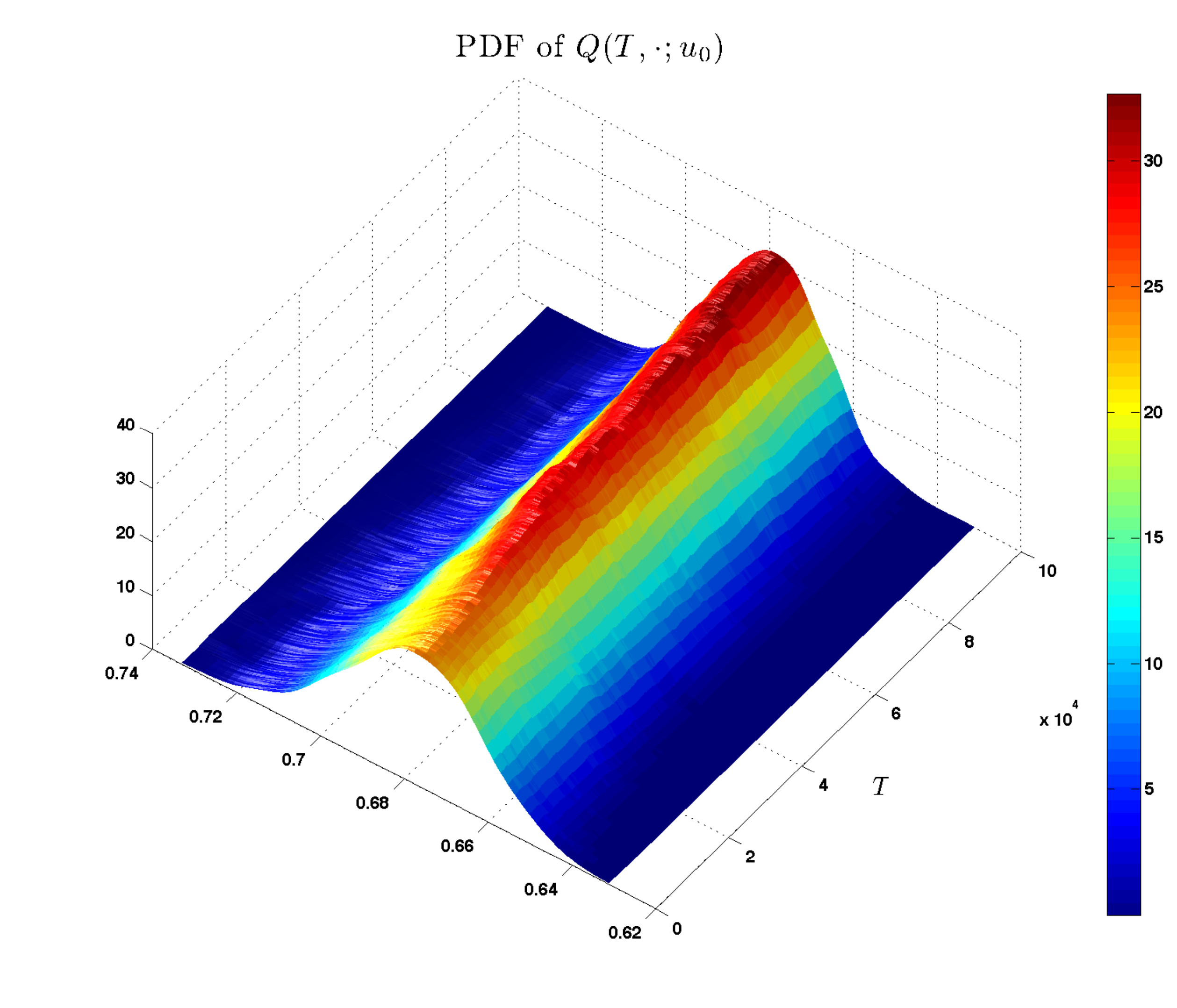}
    \end{minipage}} 
\caption{{\footnotesize {\bf Empirical probability density functions (PDFs) of the parameterization defect}:  Here for each fixed $T>0$,  the PDFs have been estimated using a standard kernel density estimator from $10^4$ realizations in (a) and (b). The values of  $\gamma$, $l$,  and $\nu$ are those used for Fig. ~\ref{fig:ave_Q}. Here $\lambda = 1.7 \lambda_c$, and $\sigma = 0.35$.
  The initial datum for the SPDE is taken to be $u_0^{+} = 0.1e_1 + 0.2 e_2 + 0.1 e_5$ in (a); and is taken to be $u_0^{-} = -0.1e_1 + 0.2 e_2 + 0.1 e_5$ in (b).}} 
\label{fig:PM_ratio_PDFs}
\end{figure}

To complete the parameterization defect analysis we computed --- for  each $T\in [0,10^4]$ --- an empirical estimation of the probability density functions (PDFs) of 
$Q_{\lambda,\sigma}(T, \omega; u_0)$  from $\omega$ living in  an empirical set $\Omega_n$ with card($\Omega_n$)$=n=10^4$ and for two distinct initial data, one living in the basin of attraction of  $u_\lambda^{+}$ and one living  in the basin of attraction of  $u_\lambda^{-}$. The results are reported in Fig. ~\ref{fig:PM_ratio_PDFs} for $\lambda = 1.7 \lambda_c$, and $\sigma = 0.35$, parameter values which will be used  later on to analyze the modeling performances achieved by the reduction strategy based on $\widehat{h}_\lambda^{(1)}$; see  Fig. ~\ref{fig:extreme} below. 
As one can observe, the support of these PDFs remains contained in a narrow subinterval  of $(0,1)$ as $T$ evolves.

These empirical facts  strongly indicate that  the pullback limit $\widehat{h}_\lambda^{(1)}$ (as given by its analytic form \eqref{eq:h Burgers}) provides a PM  for the SPDE \eqref{eq:Burgers} for a large set  of values of the control parameter $\lambda$ and the noise amplitude $\sigma$.  It is worthwhile to note that as $\lambda$ increases, the amplitude of the resolved modes (as well as of the full solution) grows, and this growth is further subject to fluctuations driven by the noise effects as $\sigma$ increases. We will show later in this section that even for parameter regimes corresponding to  large amplitudes of the solutions, the corresponding reduced equations based on $\widehat{h}_\lambda^{(1)}$ are still able to achieve good modeling performances  of   the SPDE dynamics projected onto  the resolved modes.

\medskip

\subsection{Stochastic reduced equations based on $\widehat{h}_\lambda^{(1)}$}\label{ss: Reduced equations based on h1}

The stochastic reduced equations governing the evolution of the amplitudes $y_1(t, \omega)$ and $y_2(t,\omega)$ --- associated with respectively the resolved  modes $e_1$ and $e_2$ --- are obtained by applying the results of Section  ~\ref{Sec_abstract-example} to the SPDE \eqref{eq:Burgers}. The reduced equations \eqref{bilinear pm} become then in this case:
\bea \label{eq:Burgers reduced}
& \d y_1 = \biggl( \beta_1(\lambda) y_1 + \frac{\gamma \pi}{\sqrt{2}l^{3/2}} y_1y_2 - \frac{3\gamma^2 \pi^2 \M^{12}_3(\theta_t\omega, \lambda)}{2l^3}y_1y_2^2  \\
& \hspace{15em} + \frac{3\gamma^3\pi^3 \M^{12}_3(\theta_t\omega, \lambda) \M^{22}_4(\theta_t\omega, \lambda)}{ \sqrt{2}l^{9/2}}y_1 y_2^3 \biggr) \d t + \sigma y_1 \circ \d W_t, \\
& \d y_2 = \Bigl( \beta_2(\lambda) y_2 - \frac{\gamma \pi}{\sqrt{2}l^{3/2}} y_1^2  - \frac{3\gamma^2 \pi^2 \M^{12}_3(\theta_t\omega, \lambda)}{l^3}y_1^2y_2  -  \frac{2\gamma^2\pi^2 \M^{22}_4(\theta_t\omega, \lambda)}{l^3} y_2^3 \Bigr) \d t + \sigma y_2 \circ \d W_t.
\eea
Note that based on the parameterization defect analysis  that precedes, this system should be considered in practice (for instance) for $\lambda \in [1.2 \lambda_c, 2.1 \lambda_c]$, and $\sigma \in [0,0.4]$ for which $\widehat{h}_\lambda^{(1)}$ has been empirically shown to constitute a good PM candidate.

In order to assess the role of the memory terms,  the modeling performances of the dynamics of $u_{\c}$ achieved by  ~\eqref{eq:Burgers reduced}  will be compared with those achieved by  the following averaged version of ~\eqref{eq:Burgers reduced}:
\bea \label{eq:Burgers averaged}
& \d v_1 = \biggl( \beta_1(\lambda) v_1 + \frac{\gamma \pi}{\sqrt{2}l^{3/2}} v_1v_2 - \frac{3\gamma^2 \pi^2 \mathbb{E}(\M^{12}_3(\cdot, \lambda))}{2l^3}v_1v_2^2  \\
& \hspace{15em} + \frac{3\gamma^3\pi^3 \mathbb{E}(\M^{12}_3(\cdot, \lambda))\mathbb{E}(\M^{22}_4(\cdot, \lambda))}{\sqrt{2}l^{9/2}}v_1 v_2^3 \biggr) \d t + \sigma v_1 \circ \d W_t, \\
& \d v_2 = \Bigl( \beta_2(\lambda) v_2 - \frac{\gamma \pi}{\sqrt{2}l^{3/2}} v_1^2  - \frac{3\gamma^2 \pi^2 \mathbb{E}(\M^{12}_3(\cdot, \lambda))}{l^3}v_1^2v_2  -  \frac{2\gamma^2\pi^2 \mathbb{E}(\M^{22}_4(\cdot, \lambda))}{l^3} v_2^3 \Bigr) \mathrm{d} t + \sigma v_2 \circ \mathrm{d} W_t.
\eea
This system consists in replacing the memory terms $\M^{12}_3(\cdot, \lambda)$ and $\M^{22}_4(\cdot, \lambda)$ by their corresponding expected  values $\mathbb{E}(\M^{12}_3(\cdot, \lambda))$ and $\mathbb{E}(\M^{22}_4(\cdot, \lambda))$.  For practical purposes, we recall that the latter are given by 
\beas
\mathbb{E}(\M^{12}_3(\cdot, \lambda)) = \frac{1}{\beta_{1}(\lambda) + \beta_{2}(\lambda) -
\beta_3(\lambda) - \sigma^2/2}, \qquad \mathbb{E}(\M^{22}_4(\cdot, \lambda))  = \frac{1}{2\beta_{2}(\lambda) -
\beta_4(\lambda) - \sigma^2/2};
\eeas
provided that $\sigma < \sigma_{\#}(\lambda) = \sqrt{\min(\beta_1(\lambda) + \beta_2(\lambda) - \beta_3(\lambda),2\beta_{2}(\lambda) -
\beta_4(\lambda))}$ according to Lemma \ref{lem:Mn} (here $k=2$). 

\br \label{Rmk_initialdata}
{\mkk
In practice such systems need to be integrated from an initial datum $\phi$ well-chosen.
In applications considered below, it is enough to take $\phi$ in either  \eqref{eq:Burgers reduced} or \eqref{eq:Burgers averaged}, to be $P_{\c} u_0$, where $u_0$ is the initial datum used for the SPDE simulations. In a modeling perspective of the SPDE dynamics as projected onto the resolved modes, this choice is not always appropriate in particular in presence of multiple local attractors and other choices have to be made for $\phi$.\footnote{{\mkk In that regards techniques from {\it e.g.} ~\cite{CJ12, Roberts89} can be adapted to serve our purpose and will be discussed elsewhere; see also \cite{CGLW}.}}}
\er



\subsection{Modeling performance achieved by the stochastic reduced equations based on $\widehat{h}_\lambda^{(1)}$}  \label{ss:extreme}
We  report here on the modeling performances of the SPDE dynamics on the $\mathcal{H}^{\c}$-modes achieved  by  numerical simulations of the reduced system \eqref{eq:Burgers reduced}.  A particular attention is paid on the contribution  of the memory terms to these modeling performances.  In that respect,  an appropriate parameter regime for the SPDE  \eqref{eq:Burgers} is carefully selected. 


The numerical integration of the reduced system \eqref{eq:Burgers reduced} is performed using a standard  Euler-Maruyama scheme; see {\it e.g.} ~\cite[p.~305]{KP92}. At each time step, the extrinsic memory terms $\M^{12}_{3}(\theta_t\omega, \lambda)$ and $\M^{22}_{4}(\theta_t\omega, \lambda)$ are simulated using the corresponding auxiliary SDE of type \eqref{xi eqn}. 
 Obviously, the realization $\omega$ of the Wiener process employed in this operation has to be the one that forces the effective reduced equations ~\eqref{eq:Burgers reduced} as well as the SPDE given in \eqref{eq:Burgers}.

\medskip

\vspace{1ex} 
{\bf Choice of the parameter regime.} For $\gamma$ small and $\lambda$ away from the critical value $\lambda_c$ ($\lambda>\lambda_{c}$), the amplitude of the solutions to the SPDE \eqref{eq:Burgers} get typically large so that for such choices of $\gamma$ and $\lambda$, an interesting class of regimes is selected in order to  assess  the modeling performances achieved by the reduced system \eqref{eq:Burgers reduced} based on   $\widehat{h}_\lambda^{(1)}$.


The size of the domain and the value of $\lambda$ have also been selected so that the NR-gaps $\beta_{1}(\lambda) + \beta_{2}(\lambda) - \beta_3(\lambda)$ and $2\beta_{2}(\lambda) - \beta_4(\lambda)$ are small enough in order that the memory terms $\M^{12}_3$ and $\M^{22}_4$ do not exhibit too fast decay of correlations; see Lemma \eqref{lem:Mn} (iii). The amount of the noise has been also calibrated so that 
its amount has a significant impact on the dynamics when compared with the deterministic situation where a steady state is typically reached asymptotically. Note that however $\diffusion$ cannot be chosen to be too large in order to  ensure the first two moments of $\M^{12}_3$ and $\M^{22}_4$ to exist; see again Lemma \ref{lem:Mn}. 


Such regimes --- small $\gamma$, small NR-gaps, and large enough $\sigma$ ---  constitute  thus an appropriate laboratory to assess the role of the memory effects in the stochastic modeling by  \eqref{eq:Burgers reduced}, of the SPDE dynamics  projected onto  the first two modes.  We use hereafter the averaged reduced system \eqref{eq:Burgers averaged} as a basis for comparison.  This system comes with the same nonlinear terms to model the cross-interactions between the resolved and unresolved modes where however the random fluctuations of these terms\footnote{brought by  $\M^{12}_3$ and $\M^{22}_4$ in the reduced system \eqref{eq:Burgers reduced}.} have been averaged out. This averaged reduced system is integrated also by using the Euler-Maruyama scheme with the same time step.

\vspace{1ex} 
{\bf  Large excursions from reduced equations: the role of memory effects.} 
Based on the previous discussion, we choose the SPDE parameters to be $\gamma = 0.5$, $\sigma = 0.35$, $\lambda = 1.7 \lambda_c$, $l = 2.5\pi$ and $\nu = 2$; {\mkk referred hereafter as \bf regime A}}. For this set of parameters, the values of the two NR-gaps involved in the memory terms $\M^{12}_3(\cdot, \lambda)$ and $\M^{22}_4(\cdot, \lambda)$ are respectively $\beta_{1}(\lambda) + \beta_{2}(\lambda) - \beta_3(\lambda) = 1.8$ and $2\beta_{2}(\lambda) - \beta_4(\lambda)=3.1$.

\begin{figure} [!hbtp]
   \centering
  \vspace{-2ex} \includegraphics[height=10cm, width=16cm]{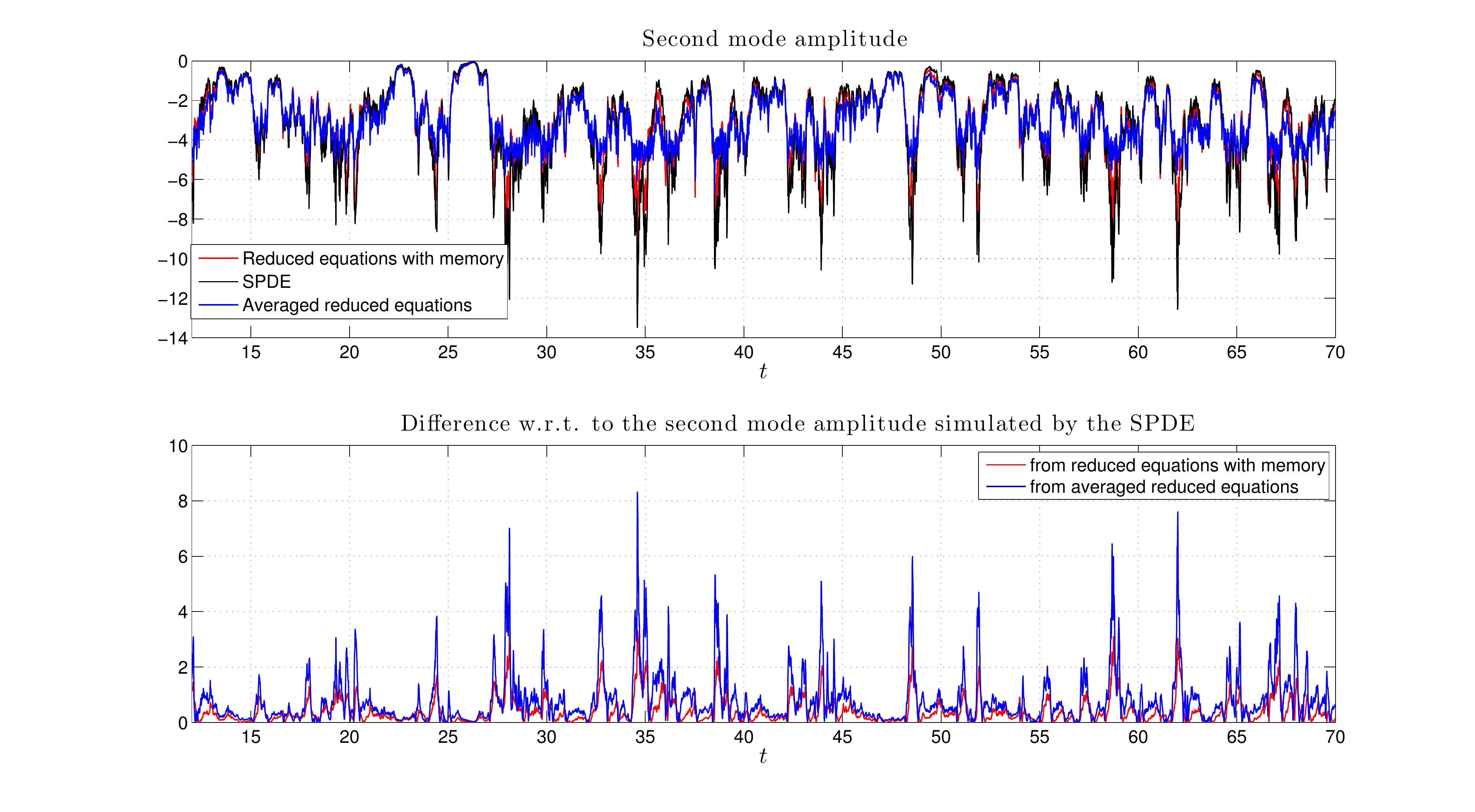}
  \vspace{-6ex}
  \caption{{\footnotesize {\bf Upper panel}: Second mode amplitude  associated with respectively the non-Markovian reduced system \eqref{eq:Burgers reduced} (red curve), the averaged reduced system \eqref{eq:Burgers averaged} (blue curve), and the SPDE \eqref{eq:Burgers} (black curve). {\bf Lower panel}: In red (resp. blue), the absolute value of the difference between the second mode amplitude $y_2$ (resp. $v_2$) simulated from the non-Markovian reduced system \eqref{eq:Burgers reduced} (resp. from the averaged system  \eqref{eq:Burgers averaged}) and $u_{\c,2}$ simulated from the SPDE. The parameters are chosen to be $\gamma = 0.5$, $\sigma = 0.35$, $\lambda = 1.7 \lambda_c$, $l=2.5\pi$, and $\nu =2$.  
The initial datum for the SPDE \eqref{eq:Burgers} is $u_0 = 0.1e_1 + 0.2e_2 + e_5$ and the initial datum for each the reduced system is taken to be $(y_{10}, y_{20}) =(0.1, 0.2)$. The values of the two NR-gaps  involved in the memory terms $\M^{12}_3(\cdot, \lambda)$ and $\M^{22}_4(\cdot, \lambda)$ are respectively $\beta_{1}(\lambda) + \beta_{2}(\lambda) - \beta_3(\lambda) = 1.8$ and $2\beta_{2}(\lambda) - \beta_4(\lambda)=3.1$.}}   \label{fig:extreme}
    \vspace{-.1ex}
\end{figure}

In what follows we focus on the modeling of  $u_{\c,2}$ --- the SPDE solution  projected onto the second mode --- from respectively  the non-Markovian  system \eqref{eq:Burgers reduced} and the  system \eqref{eq:Burgers averaged}.  The results are reported in Fig. ~\ref{fig:extreme} for a typical realization of the noise.  We can observe that while the reduced dynamics as simulated by $v_2$ from the average system \eqref{eq:Burgers averaged}  gives a reasonable  average behavior of $u_{\c,2}$,  this system fails to simulate the large deviations exhibited by $u_{\c,2}$.  To the contrary the reduced dynamics as simulated by $y_2$ from the  non-Markovian  system \eqref{eq:Burgers reduced} shows much better performances in modeling these large excursions. The lower panel of Fig. ~\ref{fig:extreme} shows an overall better modeling performance achieved by $y_2$ rather than $v_2$. Regarding the first mode amplitude, $u_{\c,1}$, less striking modeling improvement  has been observed from  $y_1$ when compared with $v_1$, although an overall better modeling performance is still achieved by the former (not shown).

These numerical results demonstrate  the importance of the the memory terms $\M^{12}_3(\cdot, \lambda)$ and $\M^{22}_4(\cdot, \lambda)$ in the modeling of the large deviations exhibited by $u_{\c,2}$.
Nevertheless it can be observed that such a modeling performance  achieved by the non-Markovian  system \eqref{eq:Burgers reduced}  (based on $\widehat{h}_\lambda^{(1)}$),  deteriorates as the noise amplitude $\sigma$  or the size $l$ of the domain increases.

This deterioration comes with two complementary observations. First, it is not difficult to deduce analytically from \eqref{Eq_egein_burg} that  the NR-gaps associated with $\M^{12}_3(\cdot, \lambda)$ and $\M^{22}_4(\cdot, \lambda)$ are decreasing functions of $l$.  Second,  it can be observed numerically that a deterioration of  the parameterization quality of $\widehat{h}_\lambda^{(1)}$ (for fixed values of $\sigma$ and $\lambda$) takes place as $l$ increases; compare for instance Fig. ~\ref{fig:ave_Q} with Fig. ~\ref{fig:ave_Q_sec10} below. 
 
 The decrease in the NR-gaps combined with the deterioration of the  parameterization quality of $\widehat{h}_\lambda^{(1)}$ underpin the need of designing parameterizing manifolds which exhibit better  
  parameterization quality  while conveying memory effects of more elaborated structure than {\mkk conveyed} by $\M^{12}_3(\cdot, \lambda)$ and $\M^{22}_4(\cdot, \lambda)$.  As we will see in the next section, this can be reached by parameterizing manifolds based on  the pullback limit $\widehat{h}^{(2)}_\lambda$ given by \eqref{Eq_PBA2uc}, when the latter exist. However,  explicit expressions of  $\widehat{h}^{(2)}_\lambda$ may be out of reach and  the determination of the reduced dynamics becomes more involved.  The next section presents how to still achieve efficient reduced systems in such  cases.

\section{Non-Markovian Stochastic Reduced Equations on the Fly}  \label{s:Fly}

{\mkk As pointed out earlier,  the pullback limits introduced in Section \ref{ss:PM} offer various} possibilities to construct manifolds that can be used to improve the parameterization quality offered by $\widehat{h}^{(1)}_\lambda$ associated with the auxiliary system \eqref{LLL}. In this section,  we consider one such manifold obtained by the pullback limit \eqref{Eq_PBA2uc} associated with the {\mkk two-layer backward-forward} system \eqref{LLL_new}; and we illustrate that the corresponding reduced system, when applied to the stochastic Burgers-type equation introduced in Section ~\ref{s:Burgers},  provides indeed better performances in modeling the dynamics of the resolved  modes when compared {\mkk with} those achieved by the reduced system based on $\widehat{h}^{(1)}_\lambda$.  {\mkk As it will be emphasized later, this feature is particularly noticeable when  the memory effects and the nonlinear cross-interactions, as well as the self-interactions between the unresolved modes become important for the modeling of the SPDE dynamics projected onto the resolved modes.  }

Methodological aspects are presented in \SS ~\ref{ss:fly-method} {\mkk and \SS ~\ref{ss:Burgers-fly}} below, where a numerical procedure is described to determine  ``on the fly''  the  reduced random vector field (based on $\widehat{h}^{(2)}_\lambda$) along a trajectory $\xii(t,\omega)$ generated by the latter as the time is advanced. This method is particularly useful when no analytic formulas of $\widehat{h}^{(2)}_\lambda$ are available. The substitutive cornerstone in this case, is the pullback characterization of  $\widehat{h}^{(2)}_\lambda$ given by \eqref{Eq_PBA2uc}, which allows to update the reduced vector field once $\xii(t,\omega)$  is known at a particular time instance $t$. 

The method is  illustrated in \SS ~\ref{ss:Burgers-fly} on the stochastic Burgers-type equation introduced in Section ~\ref{s:Burgers}. As {\mkk main numerical results, it is shown in \SS ~\ref{ss:num results} that} the statistics of the large excursions present in the SPDE dynamics projected onto the resolved modes 
are reproduced with a very good accuracy from simulations of the reduced dynamics based  on $\widehat{h}^{(2)}_\lambda$.

\subsection{Reduced system on the fly based on $\widehat{h}^{(2)}_\lambda$} \label{ss:fly-method}

{\mkk The existence of $\widehat{h}^{(2)}_\lambda$ as given by  \eqref{Eq_PBA2uc} as well as its parameterization defect are examined in  \SS ~\ref{Sec_Existence_PQ},   in the case of the stochastic Burgers-type equation.}
By assuming {\mkk for the moment} that  $\widehat{h}^{(2)}_\lambda$ is well-defined {\mkk and constitutes a PM}, the corresponding reduced system given by \eqref{SDE abstract-pm} reads as follows:
\begin{equation} \label{SDE abstract layered}
\mathrm{d} \xii = \Big( L_\lambda^{\c} \xii +P_{\c} F\big(\xii + \widehat{h}^{(2)}_\lambda(\xii, \theta_t\omega)\big)\Big) \mathrm{d}t + \diffusion \xii \circ \mathrm{d} W_t. 
\end{equation}


{\mkk As mentioned above, it is in general more involved to derive explicit analytic expressions of $\widehat{h}_\lambda^{(2)}$.  The key idea is to use  \eqref{LLL_new} in order to provide an approximation of $\widehat{h}^{(2)}_\lambda(\xii, \theta_t\omega)$ via the pullback characterization \eqref{Eq_PBA2uc}, but this has to be executed with care. 

It  is indeed too cumbersome to use \eqref{Eq_PBA2uc} directly to approximate  the vector field $P_{\c} F\big(\xii + \widehat{h}^{(2)}_\lambda(\xii, \theta_t\omega)\big)$ as $\xi$ varies in $\mathcal{H}^{\c}$.  We adopt instead 
a ``Lagragian approach'' which consists of approximating ``on the fly'' this vector field along a trajectory   $ \xi(t,\omega)$ of interest, as the time $t$ flows.\footnote{This trajectory is determined in practice by the initial datum $\phi$ used for the reduced system,  which itself depends on the initial datum $u_0$  used for  the SPDE simulation; see also Remark \ref{Rmk_initialdata}.}

As we will see, this  is much more manageable and leads naturally to consider, instead of \eqref{SDE abstract layered},  
the following reduced equation} 
 \begin{equation} \label{Reduced_fly_approx}
{\mkk \mathrm{d} \xii_t = \Big( L_\lambda^{\c} \xii_t +P_{\c} F\big(\xii_t + \widehat{u}^{(2)}_{\s}[\xi(t,\omega)](t+T, \theta_{-T}\omega; 0)\big) \Big) \mathrm{d}t + \diffusion \xii_t \circ \mathrm{d} W_t, }
\end{equation}
{\mkk where we have stressed the $t$-dependence of the variable $\xii$ to better characterize an important feature of \eqref{Reduced_fly_approx}.  This feature hinges on the fact that \eqref{Reduced_fly_approx} is actually coupled with the backward-forward system recalled in \eqref{layered_uc1}--\eqref{layered_ini} below, {\it via} the term $\widehat{u}^{(2)}_{\s}[\xi(t,\omega)](t+T, \theta_{-T}\omega;0)$. Such a coupling is however of a particular type since it requires to update, at each time $t$, the vector field in \eqref{Reduced_fly_approx} along the trajectory $ \xi(t,\omega)$. 
This is achieved at each time $t$ by first  integrating
Eqns.  ~\eqref{layered_uc1}--\eqref{layered_uc2} backward from $\xi(t,\omega)$ (in fiber $\theta_t \omega$) up to fiber $\theta_{t-T} \omega$, and then by coming back to fiber $\theta_t \omega$ by forward integration Eq. ~\eqref{layered_us}. This operation provides $\widehat{u}^{(2)}_{\s}[\xi(t,\omega)](t+T, \theta_{-T}\omega;0)$ necessary to properly determine $\d \xi_t$. 
A schematic of this procedure is depicted in Figure \ref{fig:fly_integ} below.

 }

\begin{figure}[!hbtp]
   \centering
   \includegraphics[height=6.5cm, width=13cm]{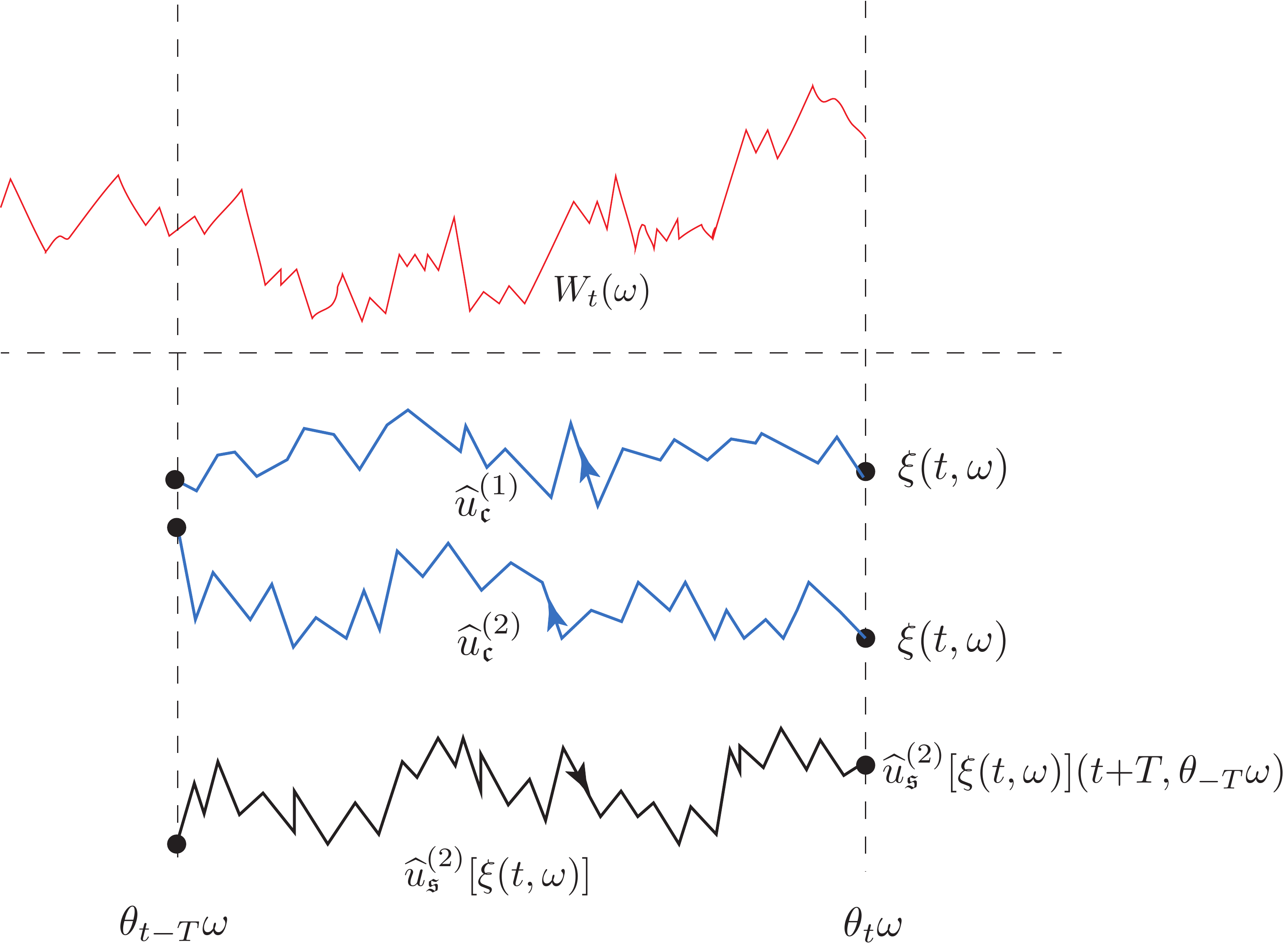}
  \caption{{\footnotesize 
{\mkk  {\bf Schematic of the on-the-fly reduction  procedure}. For a given path of the noise $W_t(\omega)$ (red curve),  the system \eqref{layered_uc1}--\eqref{layered_ini}   is first integrated backward from $\xi(t,\omega)$ (blue curves) --- the current state of the reduced equation \eqref{eq:fly_xi} (also driven by $W_t(\omega)$)  --- and then forward (black curve). The resulting $\widehat{u}^{(2)}_{\s}$-component in fiber $\theta_t \omega$, $\widehat{u}^{(2)}_{\s}[\xi(t,\omega)](t+T, \theta_{-T}\omega; 0)$,  is then used to update the vector field in \eqref{eq:fly_xi} in the same fiber.}}}  \label{fig:fly_integ}
\end{figure}

%

{\mkk The resulting coupled system can be then written  as below and will be called the {\it reduced system on the fly}  associated with  $\widehat{u}^{(2)}_{\s}[ \cdot]$ (and $T$). This system is aimed to provide an approximation of 
the PM-based reduced system \eqref{SDE abstract layered} associated with $\widehat{h}^{(2)}_\lambda$, for $T$ sufficiently large. }

\begin{subequations} \label{Reduced_fly_sys}
\begin{empheq}[box=\fbox]{align}
& {\mkk \mathrm{d} \xii_t = \Big( L_\lambda^{\c} \xii_t +P_{\c} F\big(\xii_t + \widehat{u}^{(2)}_{\s}[\xi(t,\omega)](t+T, \theta_{-T}\omega; 0)\big) \Big) \mathrm{d}t + \diffusion \xii_t \circ \mathrm{d} W_t},   && t > 0,  \label{eq:fly_xi} \\
&\xi(0, \omega) = \phi,         \label{eq:fly_xi_ini} \\ 
& \text{where $\widehat{u}^{(2)}_{\s}[\xi](t+T, \theta_{-T}\omega; 0)$ is obtained via}  \nonumber \\ 
& \mathrm{d} \widehat{u}^{(1)}_{\c} =   L_\lambda^{\c} \widehat{u}^{(1)}_{\c}  \mathrm{d} s +\sigma \widehat{u}^{(1)}_{\c} \circ \mathrm{d} W_s,  &&  s\in [t-T, t],  \label{layered_uc1}\\
& \mathrm{d} \widehat{u}^{(2)}_{\c} =  \bigl( L_\lambda^{\c} \widehat{u}^{(2)}_{\c}  + P_{\c}F(  \widehat{u}^{(1)}_{\c} ) \bigr) \mathrm{d} s +\sigma \widehat{u}^{(2)}_{\c} \circ \mathrm{d} W_s,  && s\in [t-T, t],  \label{layered_uc2}\\
& \mathrm{d} \widehat{u}^{(2)}_{\s} = \bigl( L_\lambda^{\s} \widehat{u}^{(2)}_{\s} + P_{\s} F(\widehat{u}^{(2)}_{\c}(s-T, \omega)) \bigr) \mathrm{d} s  + \sigma \widehat{u}^{(2)}_{\s} \circ \mathrm{d} W_{s-T}, && s\in [t, T+t],  \label{layered_us} \\
&  \text{with}\quad \widehat{u}^{(1)}_{\c}(s, \omega)\vert_{s=t} = \xi(t, \omega), \quad \widehat{u}^{(2)}_{\c}(s, \omega)\vert_{s=t} = \xi(t, \omega),  && \hspace{-3.75em} \widehat{u}^{(2)}_{\s}(s, \theta_{-T}\omega)\vert_{s=t}= 0. \label{layered_ini}
\end{empheq}
\end{subequations}
{\mkk
In applications considered below, it is enough to take $\phi$ in \eqref{eq:fly_xi_ini} to be $P_{\c} u_0$, where $u_0$ is the initial datum used for the SPDE simulations; see also Remark \ref{Rmk_initialdata}.}

\br
As we will see in \SS ~\ref{Sec_Existence_PQ}, such reduced systems arise typically with random coefficients  that convey memory effects of more elaborated structure 
than present in reduced systems  based on $\widehat{h}^{(1)}_\lambda$  such as \eqref{bilinear pm}. Such memory terms are exhibited analytically in \SS ~\ref{Sec_Existence_PQ}, and arise from a ``matriochka'' of nonlinear self-interactions between the low modes, which lead to a description of these memory terms  as repeated compositions of functions involving integrals depending on the past of the noise path driving the SPDE; see {\it e.g.} ~\eqref{N_2} and  \eqref{N_3} below.   For such reasons, an on-the-fly reduced system  \eqref{Reduced_fly_sys} will be qualified to be non-Markovian in what follows.

\er

For the remaining of this section, we show how the above reduced system can be efficiently  {\mkk determined, by projecting each equation of the system \eqref{Reduced_fly_sys} onto the appropriate eigenmodes in the spirit of Section \ref{Sec_abstract-example}.  Its numerical integration is described in \SS ~\ref{sss:Numerical integration}. We will see also in \SS ~\ref{ss:num results} that good modeling performance can be still achieved from \eqref{Reduced_fly_sys}, when $T$ is not necessarily large.

}

\subsection{Reduced system on the fly for the stochastic Burgers-type equation} \label{ss:Burgers-fly} 


We apply the above reduction strategy based on $\widehat{h}^{(2)}_\lambda$ to the problem \eqref{eq:Burgers} introduced in Section ~\ref{s:Burgers}. {\mkk 
 In particular,  the linear part $L_\lambda u$ in \eqref{Reduced_fly_sys} is given by $\nu \partial_{xx}^2 u + \lambda u$ and the nonlinear term by $F(u) = - \gamma u \partial_x u$. 
}


As in Section ~\ref{s:Burgers}, we take $\mathcal{H}^{\c}$ to be the subspace spanned by the first two eigenmodes:
\bes
{\mkk \mathcal{H}^{\c}:=\mathrm{span} \{ e_1,e_2\}.}
\ees
 


\subsubsection{\bf Reduced equations on the fly  in coordinate form: Derivation}

{\mkk We derive here  in the basis $\{e_1,e_2\}$, the reduced system on the fly  \eqref{Reduced_fly_sys}  associated with the problem \eqref{eq:Burgers}.


Let us first consider the system of SDEs associated with \eqref{layered_uc1} -- \eqref{layered_us}: }
\begin{subequations} \label{eq:layered-t}
\begin{align} 
& \mathrm{d} \widehat{u}^{(1)}_{\c} =   L_\lambda^{\c} \widehat{u}^{(1)}_{\c}  \mathrm{d} s +\sigma \widehat{u}^{(1)}_{\c} \circ \mathrm{d} W_s,  &&  s\in [t-T, t],  \label{layered a-t}\\
& \mathrm{d} \widehat{u}^{(2)}_{\c} =  \bigl( L_\lambda^{\c} \widehat{u}^{(2)}_{\c}  + P_{\c}F(  \widehat{u}^{(1)}_{\c} ) \bigr) \mathrm{d} s +\sigma \widehat{u}^{(2)}_{\c} \circ \mathrm{d} W_s,  && s\in [t-T, t],  \label{layered b-t}\\
& \mathrm{d} \widehat{u}^{(2)}_{\s} = \bigl( L_\lambda^{\s} \widehat{u}^{(2)}_{\s} + P_{\s} F(\widehat{u}^{(2)}_{\c}(s-T, \omega)) \bigr) \mathrm{d} s  + \sigma \widehat{u}^{(2)}_{\s} \circ \mathrm{d} W_{s-T}, && s\in [t, T+t],  \label{layered c-t} \\
& \text{with } \widehat{u}^{(1)}_{\c}(s, \omega)\vert_{s=t} = \xi(t, \omega),   \qquad \widehat{u}^{(2)}_{\c}(s, \omega)\vert_{s=t} = \xi(t, \omega), &&   \widehat{u}^{(2)}_{\s}(s, \theta_{-T}\omega)\vert_{s=t}= 0.  \label{layered ini}
\end{align}
\end{subequations}

Since $\widehat{u}^{(1)}_{\c}$ takes value in $\mathcal{H}^{\c}$, by projecting \eqref{layered a-t} against $e_1$ and $e_2$, we obtain 
\bea \label{layered a-2}
& \d y^{(1)}_1 = \beta_1(\lambda) y^{(1)}_1 \d s + \sigma y^{(1)}_1 \circ \d W_s, &&  s \in [t-T, t],  \\ 
& \d y_2^{(1)} = \beta_2(\lambda) y_2^{(1)}  \d s + \sigma y_2^{(1)} \circ \d W_s, && s \in [t-T, t], \\
\eea
with
\bea  \label{layered a-2-ini}
y^{(1)}_1(s, \omega)\vert_{s=t} = \xi_{1}(t,\omega), \qquad y_2^{(1)}(s, \omega)\vert_{s=t} = \xi_2(t,\omega),
\eea
where 
\bea \label{layered y0 1st layer}
y_i^{(1)}(s, \omega) := \langle \widehat{u}^{(1)}_{\c}(s, \omega), e_i \rangle,  \quad \xi_{i}(t,\omega) := \langle \xi(t,\omega), e_i \rangle, \quad i=1, 2,
\eea
and the initial condition \eqref{layered a-2-ini} {\mkk results} from $\widehat{u}^{(1)}_{\c}(s, \omega)\vert_{s=t} = \xi(t, \omega)$ according to \eqref{layered ini}. Recall that $\xi(t,\omega)$ is the current state of the solution associated with \eqref{eq:fly_xi}.

For \eqref{layered b-t}, first note that 
\bea\label{Eq_interactions for uc1}
\langle P_{\c} F(\widehat{u}^{(1)}_{\c}), e_1 \rangle & = - \gamma \bigl \langle P_{\c} \bigl(y_1^{(1)} e_1 + y_2^{(1)} e_2 \bigr) \bigl( y_1^{(1)}e_1 + y_2^{(1)}e_2\bigr)_x, e_1 \bigr \rangle  = \frac{\gamma\pi}{\sqrt{2}l^{3/2}}y_1^{(1)}y_2^{(1)}, \\
\langle P_{\c} F(\widehat{u}^{(1)}_{\c}), e_2 \rangle  &= - \gamma \bigl \langle P_{\c} \bigl(y_1^{(1)}e_1 + y_2^{(1)}e_2
\bigr) \bigl(y_1^{(1)}e_1 + y_2^{(1)}e_2\bigr)_x, e_2 \bigr \rangle  = -\frac{\gamma\pi}{\sqrt{2}l^{3/2}}\bigl(y_1^{(1)}\bigr)^2.
\eea
It follows then by projecting  \eqref{layered b-t} against $e_1$ and $e_2$ that
\bea \label{layered b-2}
& \d y_1^{(2)} = \Bigl( \beta_1(\lambda) y_1^{(2)} + \frac{\gamma \pi}{\sqrt{2}l^{3/2}} y_1^{(1)}y_2^{(1)} \Bigr) \d s + \sigma y_1^{(2)} \circ \d W_s, && s \in [t-T, t],  \\ 
& \d y_2^{(2)} = \Bigl( \beta_2(\lambda) y_2^{(2)} - \frac{\gamma \pi}{\sqrt{2}l^{3/2}} \bigl(y_1^{(1)}\bigr)^2  \Bigr) \d s + \sigma y_2^{(2)} \circ \d W_s, && s \in [t-T, t], \\
\eea 
with  
\bea
y_1^{(2)}(s, \omega)\vert_{s=t} = \xi_1(t,\omega), \qquad y_2^{(2)}(s, \omega)\vert_{s=t} = \xi_2(t,\omega),
\eea
where $y_i^{(2)}(s,\omega) := \langle \widehat{u}^{(2)}_{\c}(s, \omega), e_i \rangle$, $i=1, 2$; $y_1^{(1)}$ and $y_2^{(1)}$ are solutions to \eqref{layered a-2}; and $\xi_1(t,\omega)$ and $\xi_2(t,\omega)$ are given in \eqref{layered y0 1st layer}.

For \eqref{layered c-t}, note that 
\begin{subequations}\label{Eq_interactions for us2}
\begin{align}
& \langle P_{\s} F(\widehat{u}^{(2)}_{\c}), e_3 \rangle = - \gamma \bigl \langle P_{\s} \bigl(y_1^{(2)} e_1 + y_2^{(2)} e_2 \bigr) \bigl( y_1^{(2)} e_1 + y_2^{(2)}e_2 \bigr)_x, e_3 \bigr \rangle  = -\frac{3 \gamma\pi}{\sqrt{2}l^{3/2}}y_1^{(2)}y_2^{(2)},  \label{Burgers_F1}\\
& \langle P_{\s} F(\widehat{u}^{(2)}_{\c}), e_4 \rangle = - \gamma \bigl \langle P_{\s} \bigl( y_1^{(2)} e_1 + y_2^{(2)} e_2 \bigr) \bigl( y_1^{(2)}e_1 + y_2^{(2)} e_2 \bigr)_x, e_4 \bigr \rangle  = -\frac{\sqrt{2}\gamma\pi}{l^{3/2}}\bigl(y_2^{(2)}\bigr)^2, \label{Burgers_F2} \\
& \langle P_{\s} F(\widehat{u}^{(2)}_{\c}), e_n \rangle = 0, \qquad \Forall \, n \ge 5,  \label{Burgers_F3}
\end{align}
\end{subequations}
where we have used again \eqref{nonlinear interaction}.


We obtain then by projecting  \eqref{layered c-t} against $e_3$ and $e_4$ that
\bea \label{layered c-2}
& \d y_3^{(2)}= \Bigl( \beta_3(\lambda) y_3^{(2)} - \frac{3\gamma \pi}{\sqrt{2}l^{3/2}} y_1^{(2)}(s-T, \omega) y_2^{(2)}(s-T, \omega) \Bigr) \d s + \sigma y_3^{(2)} \circ \d W_{s-T}, &&  s \in [t, T+t],  \\ 
& \d y_4^{(2)} = \Bigl( \beta_4(\lambda) y_4^{(2)} - \frac{\sqrt{2}\gamma \pi}{l^{3/2}} [y_2^{(2)}(s-T, \omega)]^2  \Bigr) \d s + \sigma y_4^{(2)} \circ \d W_{s-T}, &&   s \in [t, T+t],
\eea 
with
\bea \label{layered us0}
y_3^{(2)}(s, \theta_{-T}\omega)\vert_{s=t} = 0,  \quad y_4^{(2)}(s, \theta_{-T} \omega)\vert_{s=t} = 0,
\eea
where $y_i^{(2)}(s, \theta_{-T} \omega) := \langle \widehat{u}^{(2)}_{\s}(s, \theta_{-T} \omega), e_i \rangle$, $i = 3, 4$, and $y_1^{(2)}$ and $y_2^{(2)}$ are solutions to \eqref{layered b-2}. Note that the time shift $s-T$ for $y_1^{(2)}$, $y_2^{(2)}$, and $\d W$ appeared in \eqref{layered c-2} is due again to the fact that $y_3^{(2)}$ and $y_4^{(2)}$ are initialized in fiber $\theta_{t-T} \omega$ whereas $y_1^{(2)}$ and $y_2^{(2)}$ are initialized in fiber $\theta_{t}\omega$.  

Similarly, by projecting  \eqref{layered c-t} against $e_n$ for $n\ge 5$ and using \eqref{Burgers_F3}, we obtain that
\bea \label{layered c-3}
\d y_n^{(2)} = \beta_n(\lambda) y_n^{(2)} \d s + \sigma y_n^{(2)} \circ \d W_{s-T}, && y_n^{(2)}(s, \theta_{-T} \omega)\vert_{s=t} = 0,  \quad \; s \in [t, T+t],  \quad  n \ge 5.
\eea
It is clear that the solution to \eqref{layered c-3} is identically zero. Thus, for the example considered here,  the solution to \eqref{layered c-t} takes the following form:
\be  \label{us two layer}
\boxed{
\widehat{u}^{(2)}_{\s}[\xi_t](s, \theta_{-T}\omega; 0) = y_3^{(2)}[\xi_{t}](s, \theta_{-T}\omega) e_3 + y_4^{(2)}[\xi_t](s, \theta_{-T}\omega) e_4, \quad s \in [t, T+t],
}
\ee
where $(y_3^{(2)}[\xi_t](s, \theta_{-T}\omega), y_4^{(2)}[\xi_t](s, \theta_{-T}\omega))$ is the solution to \eqref{layered c-2}--\eqref{layered us0}. As before, the dependence on $\xi(t,\omega)$ brought by $y_1^{(2)}$ and $y_2^{(2)}$ are emphasized here, and we have used $\xi_t$ to denote this dependence for brevity. In the following, we will recall this dependence when necessary, and will otherwise suppress it in most places.

{\mkk By using the expression of  $\widehat{u}^{(2)}_{\s}$ given in \eqref{us two layer},  we are now in position to derive  an operational version of  \eqref{eq:fly_xi}  for the  modeling of SPDE dynamics projected onto the resolved modes}.   {\mkk In that respect}, let us denote 
\be \label{eq:uc-components}
\xi(t,\omega) := \xi_1(t,\omega) e_1 +\xi_2(t,\omega)e_2,  \quad \text{ with } \quad   \xi_i(t,\omega) :=  \langle \xi(t,\omega), e_i \rangle,  \; i = 1, 2.
\ee 

Note that by using \eqref{us two layer}, we obtain
\beas
 \Bigl \langle P_{\c} F(\xi + \widehat{u}^{(2)}_{\s}), e_1 \Bigr \rangle & = - \gamma \Bigl \langle P_{\c} \bigl(\xi_1e_1 + \xi_2e_2 + y_3^{(2)} e_3 + y_4^{(2)} e_4\bigr)\bigl(\xi_1 e_1 + \xi_2e_2 + y_3^{(2)} e_3 + y_4^{(2)} e_4\bigr)_x, e_1 \Bigr \rangle \\
&   = \frac{\gamma\pi}{\sqrt{2}l^{3/2}} \Bigl (\xi_1\xi_2 + \xi_2 y_3^{(2)} + y_3^{(2)} y_4^{(2)} \Bigr). 
\eeas
Similarly, we have
\beas
& \Bigl \langle P_{\c} F(\xi + \widehat{u}^{(2)}_{\s}), e_2 \Bigr \rangle  = \frac{\sqrt{2} \gamma \pi}{l^{3/2}} \Bigl( - \frac{1}{2} (\xi_1)^2 + \xi_1 y_3^{(2)} + \xi_2 y_4^{(2)} \Bigr).
\eeas

Now, by projecting the system \eqref{eq:fly_xi}--\eqref{eq:fly_xi_ini} against $e_1$ and $e_2$ respectively, and by using the above two equalities, we obtain:
\bea \label{eq:uc-Burgers}
& \d \xi_1 =\Bigl( \beta_1(\lambda) \xi_1 + \frac{\gamma\pi}{\sqrt{2}l^{3/2}} (\xi_1\xi_2 + \xi_2 y_3^{(2)} + y_3^{(2)} y_4^{(2)}) \Bigr )\d t  + \sigma \xi_1 \circ \d W_t,  \\ 
& \d \xi_2 =\Bigl( \beta_2(\lambda) \xi_2 + \frac{\sqrt{2} \gamma \pi}{l^{3/2}} \bigl( - \frac{1}{2} (\xi_1)^2 + \xi_1 y_3^{(2)} + \xi_2 y_4^{(2)} \bigr) \Bigr) \d t  + \sigma \xi_2 \circ \d W_t, 
\eea
and
\be
\xi_1(0, \omega) = \xi_{1,0}, \qquad  \xi_2 (0, \omega) = \xi_{2,0},
\ee
where $\xi_1$ and $\xi_2$  are the projection of $\xi(t,\omega)$ against $e_1$ and $e_2$ respectively as given by \eqref{eq:uc-components}, and 
\bea \label{layered y0}
 \xi_{i,0} := \langle \phi, e_i \rangle, \quad i=1, 2,
\eea
with $\phi$ being the initial datum for \eqref{eq:fly_xi}.

To summarize, the on-the-fly reduced system  \eqref{Reduced_fly_sys} associated with  the SPDE \eqref{eq:Burgers}, leads, in coordinate form (for $\mathcal{H}^{\c}$ given by \eqref{Burgers Hc}),  to the  corresponding set of equations:
\begin{subequations} \label{Reduced_AE}
\begin{empheq}[box=\fbox]{align}
& \d \xi_1 =\Bigl( \beta_1(\lambda) \xi_1 + \frac{\gamma\pi}{\sqrt{2}l^{3/2}} \bigl (\xi_1\xi_2 + \xi_2 y_3^{(2)} + y_3^{(2)} y_4^{(2)} \bigr ) \Bigr)\d t  + \sigma \xi_1 \circ \d W_t,  \hspace{2em} t > 0,  \label{SDE:xi_1} \\ 
& \d \xi_2 =\Bigl( \beta_2(\lambda) \xi_2 + \frac{\sqrt{2} \gamma \pi}{l^{3/2}} \bigl( - \frac{1}{2} (\xi_1)^2 + \xi_1 y_3^{(2)} + \xi_2 y_4^{(2)}\bigr) \Bigr) \d t  + \sigma \xi_2 \circ \d W_t, \hspace{1.3em} t > 0, \label{SDE:xi_2}\\
& \text{with} \hspace{4em} \xi_1(0, \omega) = \xi_{1,0}, \hspace{7em} \xi_2 (0, \omega) = \xi_{2,0},   \label{SDE:xi_ini}\\
& \text{where $\xi_t = \xi_1(t,\omega) e_1 + \xi_2(t,\omega)e_2$ and $y_3^{(2)}$ and $y_4^{(2)}$ are obtained via} \nonumber \\
& \d y^{(1)}_1 = \beta_1(\lambda) y^{(1)}_1 \d s + \sigma y^{(1)}_1 \circ \d W_s,   \hspace{15em}  s \in [t-T, t],  \label{SDE_uc1_y1}\\ 
& \d y_2^{(1)} = \beta_2(\lambda) y_2^{(1)}  \d s + \sigma y_2^{(1)} \circ \d W_s,  \hspace{15em} s \in [t-T, t], \label{SDE_uc1_y2} \\
& \d y_1^{(2)} = \Bigl( \beta_1(\lambda) y_1^{(2)} + \frac{\gamma \pi}{\sqrt{2}l^{3/2}} y_1^{(1)}y_2^{(1)} \Bigr) \d s + \sigma y_1^{(2)} \circ \d W_s,  \hspace{6.5em} s \in [t-T, t],  \label{SDE_uc2_y1} \\ 
& \d y_2^{(2)} = \Bigl( \beta_2(\lambda) y_2^{(2)} - \frac{\gamma \pi}{\sqrt{2}l^{3/2}} \bigl(y_1^{(1)}\bigr)^2  \Bigr) \d s + \sigma y_2^{(2)} \circ \d W_s,  \hspace{6.8em} s \in [t-T, t],  \label{SDE_uc2_y2} \\
& \d y_3^{(2)} =  \Bigl[  \beta_3(\lambda) y_3^{(2)} - \frac{3\gamma \pi}{\sqrt{2}l^{3/2}} y_1^{(2)}(s - T, \omega) y_2^{(2)}(s-T, \omega) \Bigr] \d s  \label{SDE_us_y1} \\
& \hspace{8em} + \sigma y_3^{(2)} \circ \d W_{s-T},  \hspace{14.3em} s \in [t, T+t], \nonumber \\ 
& \d y_4^{(2)} = \Bigl( \beta_4(\lambda) y_4^{(2)} - \frac{\sqrt{2}\gamma \pi}{l^{3/2}} [y_2^{(2)}(s-T, \omega)]^2  \Bigr) \d s + \sigma y_4^{(2)} \circ \d W_{s-T}, \hspace{2em} s \in [t, T+t], \label{SDE_us_y2}\\ 
& \text{with} \hspace{4em} y^{(1)}_1(s, \omega)\vert_{s=t} = \xi_{1}(t,\omega), \hspace{3em} y_2^{(1)}(s, \omega)\vert_{s=t} = \xi_2(t,\omega),  \label{SDE_ini}\\
& \hspace{6em} y_1^{(2)}(s, \omega)\vert_{s=t} = \xi_1(t,\omega),  \hspace{3em}  y_2^{(2)}(s, \omega)\vert_{s=t} = \xi_2(t,\omega). \nonumber \\
& \hspace{6em} y_3^{(2)}(s, \theta_{-T}\omega)\vert_{s=t} = 0,  \hspace{4em}  y_4^{(2)}(s, \theta_{-T}\omega)\vert_{s=t} = 0. \nonumber
\end{empheq}
\end{subequations}

\subsubsection{\bf Reduced equations on the fly in coordinate form: Numerical integration}\label{sss:Numerical integration}

We describe now the numerical integration of the reduced system  \eqref{Reduced_AE}. For the sake of clarity, we make explicit the dependence of $y_3^{(2)}$ and $y_4^{(2)}$ on $\xi_t := \xi(t,\omega)$.

 {\mkk The crucial point is the determination, as the time $t$ flows,  of approximations  of $y_3^{(2)}[\xi_t] (T+t, \theta_{-T}\omega)$ and $y_4^{(2)}[\xi_t](T+t, \theta_{-T}\omega).$\footnote{{\mkk The latter give indeed access to approximations of $\widehat{u}^{(2)}_{\s}[\xi_t](T+t, \theta_{-T}\omega; 0)$, which in turn is aimed to approximate (when $T$ is sufficiently large) the parameterizing manifold function $\widehat{h}^{(2)}_\lambda$ evaluated at $( \xi(t,\omega), \theta_t \omega)$.}} This  is accomplished numerically by the use of two discrete times, one for the  numerical integration of \eqref{SDE:xi_1}--\eqref{SDE:xi_2}, and a second one (conditioned on the first)  
for the numerical integration of  \eqref{SDE_uc1_y1}--\eqref{SDE_us_y2}. The resulting double loop algorithm is described below. 
 }


{\mkk In that respect, we first provide the Euler-Maruyama  discretization of \eqref{SDE:xi_1}--\eqref{SDE:xi_2}:} 
\bea \label{eq:xi_disc}
& \xi_1^{n+1} = \xi_1^n + \Bigl[ \beta_1(\lambda) \xi_1^{n} + \frac{\sigma^2}{2} \xi_1^n + \frac{\gamma\pi}{\sqrt{2}l^{3/2}} (\xi_1^{n} \xi_2^{n} + \xi_2^{n} y_{3}^{(2),N}[\xi^{n}] + y_{3}^{(2),N}[\xi^{n}] y_{4}^{(2),N}[\xi^{n}]) \Bigr ] \delta t  + \sigma \zeta_n \sqrt{\delta t}\xi_1^{n},  \\ 
& \xi_2^{n+1} =\xi_2^n + \Bigl( \beta_2(\lambda) \xi_2^{n} + \frac{\sigma^2}{2} \xi_2^n + \frac{\sqrt{2} \gamma \pi}{l^{3/2}} \bigl( - \frac{1}{2} (\xi_1^n)^2 + \xi_1^n y_{3}^{(2),N}[\xi^{n}] + \xi_2^n y_{4}^{(2),N}[\xi^{n}]\bigr) \Bigr) \delta t  + \sigma \zeta_n \sqrt{\delta t} \xi_2^n, 
\eea
where $\xi_i^n$ is the approximation of $\xi_i(n \delta t, \omega)$, $i=1,2$,  associated with the time step $\delta t$; and $\xi^n$ denotes the corresponding approximation of $\xi(n\delta t, \omega) = \xi_1(n \delta t, \omega) e_1 + \xi_2(n \delta t, \omega)e_2$.   Similarly,  $y_{j}^{(2),N}[\xi^{n}]$  are the approximation of $y_j^{(2)}[\xi_{n\delta t}]((n+N)\delta t, \theta_{-N\delta t} \omega)$, for $j = 3,4$, where the latter are obtained {\mkk for each $n$} by {\mkk the backward-forward numerical integration of \eqref{SDE_uc1_y1}--\eqref{SDE_ini} as described below}.  The terms $\frac{\sigma^2}{2} \xi_i^n \delta t$ are due again to the conversion of the Stratonovich noise term into its It\^o form and the $\zeta_n$  are random variables drawn independently from a normal distribution $\mathcal{N}(0,1)$. The above system \eqref{eq:xi_disc} is initialized using the initial data \eqref{SDE:xi_ini}:
\bes
\xi_1^0 = \xi_{1,0}, \qquad \xi_2^0 = \xi_{2,0},
\ees
 corresponding to the components of $P_{\c} u_0$, with $u_0$ the  initial datum used for the SPDE simulation.

At each iteration $n$, once $\xi_1^n$ and  $\xi_2^n$ have been determined,  the values $y_{3}^{(2),N}[\xi^{n}]$ and $y_{4}^{(2),N}[\xi^{n}]$ are computed by backward-forward numerical integration of  \eqref{SDE_uc1_y1}--\eqref{SDE_ini} in order to determine $\xi_1^{n+1}$ and  $\xi_2^{n+1}$ via  \eqref{eq:xi_disc}.  This is organized in two steps as described below. 

First, \eqref{SDE_uc1_y1}--\eqref{SDE_uc2_y2} are integrated backward from $n\delta t$ up to $(n-N) \delta t$. Second,   \eqref{SDE_us_y1}--\eqref{SDE_us_y2} are integrated forward from $(n-N) \delta t$ to $n\delta t$. 
Both integrations are performed using an Euler-Maruyama scheme adapted to the ``arrow of time.''
 More precisely, Eqns. ~\eqref{SDE_uc1_y1}--\eqref{SDE_uc2_y2} are first integrated backward using the following scheme:
\bea \label{uc_disc}
& y_{1}^{(1),p}  - y_{1}^{(1),p+1} =  \bigl( \beta_1(\lambda)+ \frac{\sigma^2}{2}\bigr) y_{1}^{(1),p+1}\delta t   +  \sigma y_{1}^{(1),p+1} {\sqrt{\delta t}} \zeta_{n-(p+1)},  \\ 
& y_{2}^{(1),p}  - y_{2}^{(1),p+1}  =  \bigl( \beta_2(\lambda) + \frac{\sigma^2}{2}\bigr) y_{2}^{(1),p+1} \delta t + \sigma y_{2}^{(1),p+1} {\sqrt{\delta t}} \zeta_{n-(p+1)}, \\
& y_{1}^{(2),p}  - y_1^{(2),p+1}  =  \Bigl( \bigl( \beta_1(\lambda)+ \frac{\sigma^2}{2}\bigr) y_1^{(2), p+1} 
+ \frac{\gamma \pi}{\sqrt{2}l^{3/2}} y_1^{(1), p+1}y_2^{(1), p+1} \Bigr) \delta t  + \sigma y_1^{(2),p+1} {\sqrt{\delta t}} \zeta_{n-(p+1)},  \\ 
& y_{2}^{(2),p}  - y_2^{(2),p+1}  =   \Bigl( \bigl( \beta_2(\lambda) + \frac{\sigma^2}{2}\bigr) y_2^{(2),p+1} 
- \frac{\gamma \pi}{\sqrt{2}l^{3/2}} \Bigl(y_1^{(1),p+1}\Bigr)^2 \Bigr)\delta t + \sigma y_2^{(2),p+1} {\sqrt{\delta t}} \zeta_{n-(p+1)}, 
\eea
where $p = 0, \cdots, N-1$; $y_{i}^{(1),p}$ and $y_i^{(2),p}$ are respectively the discrete approximations of $y_i^{(1)}((n-p)\delta t, \omega)$ and $y_i^{(2)}((n-p)\delta t, \omega)$, $i = 1, 2$ (emanating both backward from $\xi_{n \delta t}$).

  The term $\sqrt{\delta t} \zeta_{n-(p+1)}$ is aimed to provide an approximation of the increment $W_{(n-p) \delta t}(\omega) - W_{(n-(p+1)) \delta t}(\omega)$.  The system is initialized at $p=0$ with 
\bes
y_{1}^{(1),0} = y_{1}^{(2),0}= \xi_1^n, \qquad  y_{2}^{(1),0} = y_{2}^{(2),0} = \xi_2^n,
\ees 
{\mkk with $\xi^n=(\xi_1^n,\xi_2^n)$ representing the state of \eqref{eq:xi_disc} at the $n^{\mathrm{th}}$ iteration.}

Once $\bigl \{ \bigl( y_{1}^{(2),p}, y_{2}^{(2),p}\bigr) : p = 1, \cdots, N \bigr\}$ {\mkk is obtained for a given $n^{\mathrm{th}}$ iteration of \eqref{eq:xi_disc},}
the {\mkk forward} solution of \eqref{SDE_us_y1}--\eqref{SDE_us_y2} {\mkk is then approximated according to the following scheme:}
\bea  \label{us_disc} 
y_{3}^{(2), p+1}[\xi^n]  - y_3^{(2),p}[\xi^n] & =   \Bigl( \bigl( \beta_3(\lambda)+ \frac{\sigma^2}{2}\big) y_3^{(2),p}[\xi^n] 
- \frac{3 \gamma \pi}{\sqrt{2}l^{3/2}} y_1^{(2),N-p}y_2^{(2),N-p} \Bigr) \delta t  \\
& {\hspace*{12em}}+ \sigma y_3^{(2),p}[\xi^n] {\sqrt{\delta t}} \zeta_{n-(N-p)}, \\ 
y_{4}^{(2),p+1}[\xi^n]  - y_4^{(2),p}[\xi^n] & =  \Bigl( \bigl( \beta_4(\lambda) + \frac{\sigma^2}{2}\big) y_4^{(2),p}[\xi^n] 
- \frac{\sqrt{2} \gamma \pi}{l^{3/2}} \Bigl ( y_2^{(2), N-p}\Bigr )^2 \Bigr)  \delta t \\
& {\hspace*{12em}}  + \sigma y_4^{(2),p}[\xi^n] {\sqrt{\delta t}} \zeta_{n-(N-p)}, 
\eea
where $p = 0, \cdots, N-1$; {\mkk $y_{j}^{(2), p}$ are the discrete approximations of $y_j^{(2)}[\xii_{n\delta t}]((p+n) \delta t, \theta_{-N \delta t}\omega)$, $j = 3, 4$}; and {\mkk \eqref{us_disc} is} initialized at $p=0$ with 
\bes
y_{3}^{(2),0}[\xi^n] = y_{4}^{(2),0}[\xi^n] = 0. 
\ees


Note that in practice, $N$ has to be chosen  sufficiently large so that  $y_{3}^{(2),N}[\xi^n] e_3 + y_{4}^{(2),N}[\xi^n] e_4 $ provides a good approximation of $\widehat{h}^{(2)}_\lambda(\xi(n \delta t, \omega), \theta_{n \delta t}\omega)$. 


The numerical  integration of the reduced system \eqref{Reduced_AE},  can be thus summarized in the pseudocode given below.

\medskip 
\needspace{2\baselineskip}
\noindent{\bf Pseudocode for the reduced system on the fly \eqref{Reduced_AE}:}

\begin{framed}
\begin{enumerate}[leftmargin=0.5cm]

\item [] $\zeta_{j} := \mathrm{randn(1,1)}$,   $j = -N, \cdots, 0, 1, \cdots, n_{\mathrm{max}}$;  \% generating the noise path.

\vspace*{1ex}

\item[] $(\xi_1^0, \xi_2^0) = (\xi_{1,0}, \xi_{2,0})$;  \hspace{5em}  \% initialization of  the system \eqref{eq:xi_disc} using \eqref{SDE:xi_ini}.

\vspace*{1ex}

\item[] {\bf For} $ n  = 1:n_{\mathrm{max}}$  \hspace{6.75em} \% main loop to solve the {\mkk difference equations} \eqref{eq:xi_disc}. 
\begin{equation*}
\begin{array}{r}
\left.\begin{aligned}
 \hspace*{4.5em}	y_1^{(1), 0} = \xi_1^{n-1}; \quad	y_2^{(1), 0} = \xi_2^{n-1};  \\
 \hspace*{4.5em}	y_1^{(2), 0} = \xi_1^{n-1}; \quad	y_2^{(2), 0} = \xi_2^{n-1};  
	\end{aligned} \right \}  \hspace{0.5em} \text{\% initialization of the {\mkk difference equations} \eqref{uc_disc}.}
	\end{array}
	\end{equation*} 

\bi 
\item[] {\bf For} $p=1:N$ \hspace{10em}  \% solving  the {\mkk difference equations} \eqref{uc_disc}.
\begin{equation*}
\begin{aligned}
y_1^{(1),p}= \frac{ y_1^{(1),p-1}}{1 + \bigl( \beta_1(\lambda)+ \frac{\sigma^2}{2}\bigr) \delta t + \sigma {\sqrt{\delta t}} \zeta_{n-p}}; \\
y_2^{(1),p}= \frac{ y_2^{(1),p-1}}{1 + \bigl( \beta_2(\lambda)+ \frac{\sigma^2}{2}\bigr) \delta t + \sigma {\sqrt{\delta t}} \zeta_{n-p}}; \\
y_1^{(2),p}= \frac{ y_1^{(2),p-1} - \frac{\gamma \pi}{\sqrt{2}l^{3/2}} y_1^{(1),p} y_2^{(1),p} \delta t }{1 + \bigl( \beta_1(\lambda)+ \frac{\sigma^2}{2}\bigr) \delta t + \sigma {\sqrt{\delta t}} \zeta_{n-p}}; \\
y_2^{(2),p}= \frac{ y_2^{(2),p-1} - \frac{\gamma \pi}{\sqrt{2}l^{3/2}} (y_1^{(1),p})^2 \delta t}{1 + \bigl( \beta_2(\lambda)+ \frac{\sigma^2}{2}\bigr) \delta t + \sigma {\sqrt{\delta t}} \zeta_{n-p}};      
\end{aligned} 
\end{equation*}
{\bf End For}

\vspace{2ex}

$y_3^{(2),0}= 0; \; \;  y_4^{(2),0}= 0;$   \hspace{4em} \% initialization of the {\mkk difference equations} \eqref{us_disc}.

\vspace*{1ex}

{\bf For} $p=1:N$   \hspace{9em} \% solving the {\mkk difference equations} \eqref{us_disc}.
\begin{equation*}
\begin{aligned}
\hspace*{6.5em} y_3^{(2),p} & =  y_3^{(2),p-1} + \Bigl( \bigl( \beta_3(\lambda)+ \frac{\sigma^2}{2}\bigr) y_3^{(2),p-1}
- \frac{3 \gamma \pi}{\sqrt{2}l^{3/2}} y_1^{(2),N-(p-1)}y_2^{(2),N-(p-1)} \Bigr) \delta t  \\
& {\hspace*{12em}}+ \sigma y_3^{(2),p-1} {\sqrt{\delta t}} \zeta_{n-(N-(p-1))};  \\
\hspace*{6.5em} y_{4}^{(2),p} &=  y_4^{(2),p-1} + \Bigl( \bigl( \beta_4(\lambda) + \frac{\sigma^2}{2}\bigr) y_4^{(2),p-1}
- \frac{\sqrt{2} \gamma \pi}{l^{3/2}} \Bigl ( y_2^{(2), N-(p-1)}\Bigr )^2 \Bigr)  \delta t \\
& {\hspace*{12em}}  + \sigma y_4^{(2),p-1} {\sqrt{\delta t}} \zeta_{n-(N-(p-1))};
\end{aligned} 
\end{equation*}
{\bf End For}


\vspace*{1ex}

\item[]  \% update $(\xi_1^{n}, \xi_2^{n})$:
\beas
& {\hspace*{3em}}  \xi_1^{n} = \xi_1^{n-1} + \Bigl( \bigl ( \beta_1(\lambda)  + \frac{\sigma^2}{2} \bigr) \xi_1^{n-1} + \frac{\gamma\pi}{\sqrt{2}l^{3/2}} (\xi_1^{n-1} \xi_2^{n-1} + \xi_2^{n-1} y_{3}^{(2),N} + y_{3}^{(2),N} y_{4}^{(2),N}) \Bigr) \delta t  \\
& {\hspace*{18em}} + \sigma \zeta_n {\sqrt{\delta t}}\xi_1^{n-1}; \\ 
& {\hspace*{3em}} \xi_2^{n} =\xi_2^{n-1} + \Bigl( \bigl( \beta_2(\lambda) + \frac{\sigma^2}{2} \bigr) \xi_2^{n-1} + \frac{\sqrt{2} \gamma \pi}{l^{3/2}} \bigl( - \frac{1}{2} (\xi_1^{n-1})^2 + \xi_1^{n-1} y_{3}^{(2),N} + \xi_2^{n-1} y_{4}^{(2),N}\bigr) \Bigr) \delta t  \\
& {\hspace*{18em}} + \sigma \zeta_n {\sqrt{\delta t}} \xi_2^{n-1};
\eeas

\ei

\item[] {\bf End For}

\end{enumerate} 

\end{framed}

\subsection{Existence of $\widehat{h}^{(2)}_\lambda$ as pullback limit, new memory terms, and non-resonance conditions}\label{Sec_Existence_PQ}  

We turn now to identify the condition under which the pullback limit as given in \eqref{Eq_PBA2uc} exists {\mkk for the SPDE \eqref{eq:Burgers} and $\mathcal{H}^{\c}$ as given in \eqref{Burgers Hc}}. This is made possible by noting that the auxiliary systems given in \eqref{SDE_uc1_y1}--\eqref{SDE_ini} can actually be solved analytically. Indeed, the solution of \eqref{SDE_uc1_y1}--\eqref{SDE_uc1_y2} together with the corresponding inital datum given in \eqref{SDE_ini} is given as follows:
\bea \label{layered a-sol} 
& y_1^{(1)}(s,\omega) = \xi_1(t,\omega) e^{\beta_1(\lambda)  (s-t) + \sigma (W_s(\omega) - W_t(\omega))}, \quad  s \in [t-T, t], \\
& y_2^{(1)}(s,\omega) = \xi_2(t,\omega) e^{\beta_2(\lambda) (s-t) + \sigma ( W_s(\omega) - W_t(\omega)) },  \quad  s \in [t-T, t].
\eea

For system \eqref{SDE_uc2_y1}--\eqref{SDE_uc2_y2}, we have
\bea \label{layered b-sol}
y_1^{(2)}(s,\omega) &= y_1^{(1)}(s,\omega) - \frac{\gamma \pi}{\sqrt{2}l^{3/2}} \int_s^t e^{\beta_1(\lambda)  (s-\tau) + \sigma (W_s(\omega) - W_\tau(\omega))}y_1^{(1)}(\tau,\omega) y_2^{(1)}(\tau,\omega) \,\d \tau, \\
y_2^{(2)}(s,\omega) &= y_2^{(1)}(s,\omega) + \frac{\gamma \pi}{\sqrt{2}l^{3/2}} \int_s^t e^{\beta_2(\lambda)  (s-\tau) + \sigma (W_s(\omega) - W_\tau(\omega))}[y_1^{(1)}(\tau,\omega)]^2 \,\d \tau,
\eea
where $y_1^{(1)}$ and $y_2^{(1)}$ are given in \eqref{layered a-sol}, and $s \in [t-T, t]$.

The solution to system \eqref{SDE_us_y1}--\eqref{SDE_us_y2} is then given by:
\bea \label{layered c-sol}
y_3^{(2)}(s,\theta_{-T}\omega) &=  -\frac{3\gamma \pi}{\sqrt{2}l^{3/2}} \int_t^s e^{\beta_3(\lambda)  (s-\tau) + \sigma (W_{s-T}(\omega) - W_{\tau-T}(\omega))}y_1^{(2)}(\tau-T,\omega) y_2^{(2)}(\tau-T,\omega) \,\d \tau, \\
y_4^{(2)}(s,\theta_{-T}\omega) &=  -\frac{\sqrt{2}\gamma \pi}{l^{3/2}} \int_t^s e^{\beta_4(\lambda)  (s-\tau) + \sigma (W_{s-T}(\omega) - W_{\tau-T}(\omega))}[y_2^{(2)}(\tau-T,\omega)]^2 \,\d \tau,
\eea
where $y_1^{(2)}$ and $y_2^{(2)}$ are given in \eqref{layered b-sol}, and $s \in [t, t+T]$. 

By performing the change of variables $\tau' = \tau - T$ setting $s$ to $t+T$ in the integrals involved in \eqref{layered c-sol}, 
we obtain that $y_3^{(2)}(t+T,\theta_{-T}\omega)$ and $y_4^{(2)}(t+T,\theta_{-T}\omega)$ are given by
\bea\label{Eq_golint}
y_3^{(2)}(t+T,\theta_{-T}\omega) &=  - \frac{3\gamma \pi}{\sqrt{2}l^{3/2}}  \int_{t-T}^{t} e^{\beta_3(\lambda) (t - \tau') + \sigma ( W_{t}(\omega) - W_{\tau'}(\omega))}y_1^{(2)}(\tau',\omega) y_2^{(2)}(\tau',\omega) \,\d \tau', \\
y_4^{(2)}(t+T,\theta_{-T}\omega) &=  - \frac{\sqrt{2}\gamma \pi}{l^{3/2}} \int_{t-T}^t e^{\beta_4(\lambda) (t- \tau') - \sigma( W_{t}(\omega) - W_{\tau'}(\omega))}[y_2^{(2)}(\tau',\omega)]^2 \,\d \tau'.
\eea

 By using the expressions of $y_1^{(2)}(\tau',\omega)$ and $y_2^{(2)}(\tau',\omega)$ given by \eqref{layered b-sol} (using also \eqref{layered a-sol}), the integrands 
\be
a(t,\tau';\omega)=e^{\beta_3(\lambda) (t - \tau') + \sigma ( W_{t}(\omega) - W_{\tau'}(\omega))}y_1^{(2)}(\tau',\omega) y_2^{(2)}(\tau',\omega),
\ee
and 
\be\label{Eq_bterm}
 b(t,\tau';\omega)= e^{\beta_4(\lambda) (t- \tau') - \sigma( W_{t}(\omega) - W_{\tau'}(\omega))}[y_2^{(2)}(\tau',\omega)]^2,
\ee
can be expanded as the sum of seven terms: four coming by expansion of $a$,  and three by expansion of $b$.   

 For instances, the four terms coming from $a$ are given by
\be \label{eq:exponential terms1}
a_1(t,\tau';\omega)=\xi_1(t,\omega) \xi_2(t,\omega) e^{ -(\beta_1(\lambda) + \beta_2(\lambda) - \beta_3(\lambda)) (t-\tau') - \sigma (  W_{t}(\omega) - W_{\tau'}(\omega))}, \\
\ee
\bea\label{eq:exponential terms2}
a_2(t,\tau';\omega)=&- \frac{\gamma \pi}{\sqrt{2}l^{3/2}} \xi_1(t,\omega) [\xi_2(t,\omega)]^2 e^{ -(\beta_1(\lambda) + 2 \beta_2(\lambda) - \beta_3(\lambda))t - 2 \sigma W_{t}(\omega)} \\
&{\hspace*{5em}} \times e^{(\beta_1(\lambda) + \beta_2(\lambda) - \beta_3(\lambda)) \tau' + \sigma  W_{\tau'}(\omega)} \int_{\tau'}^t e^{\beta_2(\lambda) \tau + \sigma  W_{\tau}(\omega)} \, \d \tau,
\eea 
\bea\label{eq:exponential terms3}
 a_3(t,\tau';\omega)=&\frac{\gamma \pi}{\sqrt{2}l^{3/2}} [\xi_1(t,\omega)]^3 e^{ -(3\beta_1(\lambda) - \beta_3(\lambda))t - 2 \sigma W_{t}(\omega)} \\
& {\hspace*{5em}} \times e^{ (\beta_1(\lambda) + \beta_2(\lambda) - \beta_3(\lambda)) \tau' + \sigma W_{\tau'}(\omega)} \int_{\tau'}^t e^{( 2\beta_1(\lambda) - \beta_2(\lambda) ) \tau + \sigma W_{\tau}(\omega)} \, \d \tau,
\eea
and 
\bea\label{eq:exponential terms4}
a_4(t,\tau';\omega)=&- \frac{\gamma^2 \pi^2}{2 l^{3}} [\xi_1(t,\omega)]^3 \xi_2(t,\omega) e^{ -(3\beta_1(\lambda) + \beta_2(\lambda) - \beta_3(\lambda))t - 3 \sigma W_{t}(\omega)} \\
& {\hspace*{1em}} \times e^{ (\beta_1(\lambda) + \beta_2(\lambda) - \beta_3(\lambda)) \tau' + \sigma W_{\tau'}(\omega)} \int_{\tau'}^t e^{\beta_2(\lambda) \tau + \sigma W_{\tau}(\omega)} \, \d \tau \int_{\tau'}^t e^{( 2\beta_1(\lambda) - \beta_2(\lambda)) \tau + \sigma W_{\tau}(\omega)} \, \d \tau,
\eea
with similar expressions for the three terms coming from $b$, also containing terms of degree 2, 3 and 4 in the variables $\xi_1$ and $\xi_2$.

 From \eqref{Eq_golint}, the existence of the pullback limits $\lim_{T\rightarrow +\infty}y_3^{(2)}(t+T,\theta_{-T}\omega)$ and $\lim_{T\rightarrow +\infty}y_4^{(2)}(t+T,\theta_{-T}\omega)$  reduces then to the existence of  the integral  over $(-\infty, t]$ (in the $\tau'$-variable) of each of the aforementioned seven terms.  

The condition under which this holds can be obtained by using the following growth control of the Wiener process
\bea\label{W-s control}
|W_s(\omega)| \le C_{\epsilon}(\omega) + \epsilon |s|, \quad \Forall s\in \mathbb{R}, \; \omega \in \Omega,
\eea
where $\epsilon$ is an arbitrarily fixed positive constant, and $C_{\epsilon}$ is a positive random variable depending on $\epsilon$; see Lemma ~\ref{Lem:OU}. Indeed, by controlling $W_s$ in all the terms of the type \eqref{eq:exponential terms1}-\eqref{eq:exponential terms4} by $C_{\epsilon}(\omega) + \epsilon |s|$,  it can be checked that the set of conditions
\begin{equation}  \label{NR2} \tag{NR2}
\begin{aligned} 
& \beta_1(\lambda) + \beta_2(\lambda) - \beta_3(\lambda) > 0, \qquad && \beta_1(\lambda) + 2 \beta_2(\lambda) - \beta_3(\lambda)  > 0, \\
&3 \beta_1(\lambda)  - \beta_3(\lambda)  > 0 , \qquad && 3 \beta_1(\lambda) + \beta_2(\lambda) - \beta_3(\lambda) > 0, \\
& 2 \beta_1(\lambda)  + \beta_2(\lambda) - \beta_4(\lambda) > 0,  && 4 \beta_1(\lambda) - \beta_4(\lambda) > 0,\\
&  2 \beta_2(\lambda) - \beta_4(\lambda) > 0,
\end{aligned}
\end{equation}
are necessary and sufficient for the aforementioned pullback limits to exist.

Now, in virtue of \eqref{Eq_PBA2uc} and \eqref{us two layer}, and under the above \eqref{NR2}-condition, $\widehat{h}^{(2)}_\lambda(\xi(t, \omega), \theta_{t}\omega)$ is obtained as the following pullback limit:
\bea \label{eq:happ2 PB}
\widehat{h}^{(2)}_\lambda(\xi(t, \omega), \theta_{t}\omega)  = \lim_{T \rightarrow +\infty} \Big(y_3^{(2)}(T+t, \theta_{-T}\omega) e_3 + y_4^{(2)}(T+t, \theta_{-T}\omega) e_4\Big).
\eea

From \eqref{eq:exponential terms1}-\eqref{eq:exponential terms4}, it can be inferred that the resulting manifold function $\widehat{h}^{(2)}_\lambda$ is a random polynomial function of degree 4 in the variables $\xi_1$ and $\xi_2$. The  analytic expression of $\widehat{h}^{(2)}_\lambda(\xi, \theta_t \omega)$ contains actually the degree-two monomials  constituting the expression of $\widehat{h}^{(1)}_\lambda(\xi, \theta_t \omega)$  
given in \eqref{eq:h Burgers}: $\int_{-\infty}^t a_1(t,\tau';\omega) \d \tau'$ from \eqref{eq:exponential terms1} corresponding to $-\frac{3\gamma \pi}{\sqrt{2}l^{3/2}} \M^{12}_3(\theta_t\omega, \lambda) \xi_1\xi_2$,  and $\frac{\sqrt{2} \gamma \pi }{l^{3/2}} \M^{22}_4(\theta_t\omega, \lambda)  \xi_2^2$ from $b$, with $\M^{12}_3$ and $\M^{22}_4$ as given in \eqref{M3M4}. We recall that the fact that the expression of $\widehat{h}^{(2)}_\lambda$ contains the one of $\widehat{h}^{(1)}_\lambda$  is  not limited to this stochastic Burgers-type equation, and can be deduced in the general case from the abstract definition of $\widehat{h}^{(2)}_\lambda$ given by \eqref{Eq_PBA2uc} as it was noted in \eqref{h1 correction}.


\vspace{1ex}
{\bf Non-resonances, matriochka of nonlinear interactions and hierarchy of memory effects.} The five extra terms contained in  $\widehat{h}^{(2)}_\lambda(\xi, \theta_t \omega)$ compared with  $\widehat{h}^{(1)}_\lambda(\xi, \theta_t \omega)$, come with 
new type of nonlinear self-interactions  between the low modes $e_1$ and $e_2$. This is visible from the  self-interactions between the low modes such as described by \eqref{Eq_interactions for us2}  used in the construction of the system  \eqref{layered c-2} that provides $y_3^{(2)}$  and  $ y_4^{(2)}$. As can be seen on \eqref{Eq_interactions for us2}, these interactions involve the variables $y_1^{(2)}$ and $y_2^{(2)}$ that result from (backward) integration of  \eqref{layered b-2}, and thus convey themselves self-interactions  between the low modes as described in \eqref{Eq_interactions for uc1}.

Interestingly such a {\it ``matriochka'' of nonlinear self-interactions} comes with a new type of memory terms of more elaborated structures than $\M^{12}_3$ and $\M^{22}_4$ conveyed by $\widehat{h}^{(1)}_\lambda$. To understand these structures, a closer look at the coefficient associated with the monomial $- \frac{\gamma \pi}{\sqrt{2}l^{3/2}} \xi_1 (\xi_2)^2$ arising in the analytic expression of 
$\widehat{h}^{(2)}_\lambda$, is illuminating.  The latter is obtained by integration of \eqref{eq:exponential terms2} over $(-\infty, t]$ (in the $\tau'$-variable) as:
\be \label{N_2}
N_2(t,\omega)=e^{ -(\beta_1(\lambda) + 2 \beta_2(\lambda) - \beta_3(\lambda))t - 2 \sigma W_{t}(\omega)}\int_{-\infty}^t \alpha_2(t,\tau';\omega) \d \tau',
\ee
where
\be \label{alpha_2}
\alpha_2(t,\tau';\omega):=e^{ (\beta_1(\lambda) + \beta_2(\lambda) - \beta_3(\lambda)) \tau' + \sigma W_{\tau'}(\omega)} \int_{\tau'}^t e^{\beta_2(\lambda) \tau + \sigma W_{\tau}(\omega)} \, \d \tau, \qquad \omega \in \Omega.
\ee

Note that $N_2(t,\omega)$ is well defined thanks to the conditions $\beta_1(\lambda) + \beta_2(\lambda) - \beta_3(\lambda) > 0$ and $\beta_1(\lambda) + 2 \beta_2(\lambda) - \beta_3(\lambda) > 0$ from the non-resonance condition \eqref{NR2}, and again the growth control of the Wiener process given by \eqref{W-s control}. 

Interestingly, $N_2(t,\omega)$ can be obtained as the pullback limit of the $N$-component of the following backward-forward system:
\begin{subequations}\label{eq:aux memory two-layer}
\begin{empheq}[box=\fbox]{align}
& \d M= \left(1 - \beta_{2}(\lambda) M \right )\mathrm{d} s - \diffusion M \circ \mathrm{d} W_s, && s \in[ t-T, t], \label{eq:aux memory two-layer_1}\\
& \d N = \bigl( - (\beta_{1}(\lambda) + 2 \beta_{2}(\lambda) - \beta_{3}(\lambda))N - M(s-T, \omega) \bigr) \d s -2\sigma N \circ \d W_{s-T}, && s \in [t, t+T], \label{eq:aux memory two-layer_2} \\
& \textrm{with } M(s, \omega)\vert_{s=t} = 0, \textrm{ and } N(s, \theta_{-T}\omega)\vert_{s=t}= 0,
\end{empheq}
\end{subequations}
where Eq. ~\eqref{eq:aux memory two-layer_1} is integrated backward from fiber $\theta_t \omega$ up to fiber $\theta_{t-T}\omega$, and Eq. ~\eqref{eq:aux memory two-layer_2} is integrated forward from fiber $\theta_{t-T} \omega$ up to fiber $\theta_{t}\omega$.


Indeed, by solving \eqref{eq:aux memory two-layer_1} backward from fiber $\theta_t \omega$ up to fiber $\theta_{t-T}\omega$,  with initial datum $M(s, \omega)\vert_{s=t} = 0$, we obtain that
\be \label{eq:M}
M(s, \omega) = -\int_s^t e^{-\beta_2(\lambda)(s - \tau) - \sigma (W_s(\omega) - W_\tau(\omega))} \d \tau, \qquad s \in [t - T, t].
\ee

By integrating then \eqref{eq:aux memory two-layer_2} forward from fiber $\theta_{t-T}\omega$   up to fiber  $\theta_{t}\omega$, we obtain
\beas
N(s, \theta_{-T}\omega) & = - \int_t^s e^{-(\beta_{1}(\lambda) + 2 \beta_{2}(\lambda) - \beta_{3}(\lambda))(s - \tau) - 2 \sigma (W_{s-T}(\omega) - W_{\tau - T}(\omega))} M(\tau - T, \omega) \d \tau, \qquad s \in [t, t + T],
\eeas
where the initial datum is taken to be zero  in fiber $\theta_{t-T} \omega$. 

A change of variables $\tau - T = \tau'$, leads then to 
\beas
N(t+T, \theta_{-T}\omega) & = - \int_{t}^{t+T} e^{-(\beta_{1}(\lambda) + 2 \beta_{2}(\lambda) - \beta_{3}(\lambda))(t+T - \tau) - 2 \sigma (W_{t}(\omega) - W_{\tau - T}(\omega))} M(\tau - T, \omega) \d \tau \\
& = - \int_{t-T}^{t} e^{-(\beta_{1}(\lambda) + 2 \beta_{2}(\lambda) - \beta_{3}(\lambda))(t - \tau') - 2 \sigma (W_{t}(\omega) - W_{\tau'}(\omega))} M(\tau', \omega) \d \tau'.
\eeas

Now, by reporting the expression of $M(\cdot, \omega)$ given by \eqref{eq:M} in the above equality and using the definition of $\alpha_2(t,\tau';\omega)$ given by \eqref{alpha_2}, we obtain after rearranging the terms,
\beas
N(t+T, \theta_{-T}\omega) = e^{-(\beta_{1}(\lambda) + 2 \beta_{2}(\lambda) - \beta_{3}(\lambda))t - 2 \sigma W_{t}(\omega) } \int_{t-T}^{t} \alpha_2(t,\tau';\omega) \d \tau',
\eeas
which leads to the desired result for $N_2(t,\omega)$ defined in \eqref{N_2}, {\it i.e.}
\bes
N_2(t,\omega) = \lim_{T \rightarrow +\infty} N(t+T, \theta_{-T}\omega).
\ees

Similarly, the memory term associated with \eqref{eq:exponential terms3} given by,
\bea\label{N_3}
N_3(t,\omega) & := e^{ -(3\beta_1(\lambda) - \beta_3(\lambda))t - 2 \sigma W_{t}(\omega)} \\
& \hspace{5em} \times \int_{-\infty}^{t} \Big( e^{ (\beta_1(\lambda) + \beta_2(\lambda) - \beta_3(\lambda)) \tau' + \sigma W_{\tau'}(\omega)} \int_{\tau'}^t e^{( 2\beta_1(\lambda) - \beta_2(\lambda) ) \tau + \sigma W_{\tau}(\omega)} \d \tau \Big) \; \d \tau',
\eea
can be obtained as the pullback limit of the $N$-component of the following backward-forward system:
\begin{subequations}\label{eq:aux memory two-layer_3}
\begin{empheq}[box=\fbox]{align}
& \d M= \left(1 - (2\beta_{1}(\lambda) - \beta_{2}(\lambda)) M \right )\mathrm{d} s - \diffusion M \circ \mathrm{d} W_s, && s \in[ t-T, t], \label{eq:aux memory two-layer_3-1}\\
& \d N = \bigl( - (\beta_{1}(\lambda) + \beta_{2}(\lambda) - \beta_{3}(\lambda))N - M(s-T, \omega) \bigr) \d s -2\sigma N \circ \d W_{s-T}, && s \in [t, t+T], \label{eq:aux memory two-layer_3-2} \\
& \textrm{with } M(s, \omega)\vert_{s=t} = 0, \textrm{ and } N(s, \theta_{-T}\omega)\vert_{s=t}= 0,
\end{empheq}
\end{subequations}
and the backward-forward system for the memory term associated with \eqref{eq:exponential terms4},
\beas
N_4(t,\omega) & := e^{ -(3\beta_1(\lambda) + \beta_2(\lambda) - \beta_3(\lambda))t - 3 \sigma W_{t}(\omega)} \\
& \times \int_{-\infty}^{t}\Big(e^{ (\beta_1(\lambda) + \beta_2(\lambda) - \beta_3(\lambda)) \tau' + \sigma W_{\tau'}(\omega)} \int_{\tau'}^t e^{\beta_2(\lambda) \tau + \sigma W_{\tau}(\omega)} \, \d \tau \int_{\tau'}^t e^{( 2\beta_1(\lambda) - \beta_2(\lambda)) \tau + \sigma W_{\tau}(\omega)} \d \tau\Big) \; \d \tau',
\eeas
is given by
\begin{subequations}\label{eq:aux memory two-layer_4}
\begin{empheq}[box=\fbox]{align}
& \d M_1= \left(1 - \beta_{2}(\lambda) M_1 \right )\mathrm{d} s - \diffusion M_1 \circ \mathrm{d} W_s, && s \in[ t-T, t], \label{eq:aux memory two-layer_4-1}\\
& \d M_2 = \left(1 - (2\beta_{1}(\lambda) - \beta_{2}(\lambda)) M_2 \right )\mathrm{d} s - \diffusion M_2 \circ \mathrm{d} W_s, && s \in[ t-T, t], \label{eq:aux memory two-layer_4-2}\\
& \d N = \bigl( - (\beta_{1}(\lambda) + \beta_{2}(\lambda) - \beta_{3}(\lambda))N + M_1(s-T, \omega)M_2(s-T, \omega) \bigr) \d s \label{eq:aux memory two-layer_4-3} \\
& \hspace{16em} -2\sigma N \circ \d W_{s-T}, && s \in [t, t+T], \nonumber \\
& \textrm{with } M_1(s, \omega)\vert_{s=t} = M_2(s, \omega)\vert_{s=t} = 0, \textrm{ and } N(s, \theta_{-T}\omega)\vert_{s=t}= 0.
\end{empheq}
\end{subequations}

Similar statements hold for two of the three extra terms  arising  from expansion of $b$ given by \eqref{Eq_bterm}.
In the case where $\widehat{h}_\lambda^{(3)}$ obtained as pullback limit associated with the {\it three-layer system} \eqref{n_layer} (with $n=3$) is considered, new memory terms (and new non-resonance conditions) also arise.  As in the case of $\widehat{h}_\lambda^{(2)}$, these terms are defined by means of pullback limits associated with auxiliary backward-forward systems of type just introduced above.
Such auxiliary backward-forward systems are also associated with their corresponding non-resonances. The resulting memory terms built in this way   
arise also  from {\it nonlocal information} coming from the past of the noise path. As for $N_2$, $N_3$ and $N_4$ described above, these terms are constituted by {\it integrals of compositions of functions}, where the latter are  themselves defined by means of integrals involving  
the past of the noise, except that the length of such compositions increase by one compared to those involved in $N_2$, $N_3$ and $N_4$.

An increase in  the number of layers in  \eqref{n_layer} --- as long as the related pullback limits are well-defined --- leads thus to a {\it hierarchy of memory terms} {\it hierarchy of memory terms} obtained via repeated compositions of functions involving integrals depending on the past of the noise path driving the SPDE.  Such a hierarchy arise  with a matriochka of nonlinear self-interactions between the low modes, as well as with  a sequence of non-resonance conditions,  both of increasing complexity.

Such features are actually related to the problem of convergence of the sequence of manifold functions $\{\widehat{q}_\lambda^{(n)}\}_{q\in \mathbb{N}^{*}}$. This problem will be rigorously analyzed elsewhere for the general case. In the case of our stochastic Burgers-type equation,  we report below on the numerical  results that strongly indicate that such a convergence takes place here for  a broad class of regimes as $\lambda$ varies.

For the sake of completeness, we conclude this subsection with the pseudocode corresponding to a reduction on the fly based on $\widehat{h}_\lambda^{(q)}$ obtained as the pullback limit  \eqref{h_n}  (for $q\geq 2$). As for $\widehat{h}_\lambda^{(2)}$, 
we write actually this pseudocode for the corresponding version in coordinate form  (for $\mathcal{H}^{\c}=\mathrm{span}(e_1,e_2)$)  of the reduced system  \eqref{Reduced_fly_sys},  with $u_{\s}^{(2)}[\xi(t,\omega)](t+T, \theta_{-T}\omega; 0) $ substituted by $u_{\s}^{(q)}[\xi(t,\omega)](t+T, \theta_{-T}\omega; 0) $ therein.

\vspace{2ex}

\noindent{\bf Pseudocode for the reduced system on the fly based on a $q$-layer system  \eqref{n_layer}}

\begin{framed}
\begin{enumerate}[leftmargin=0.5cm]

\item [] $\zeta_{j} := \mathrm{randn(1,1)}$,   $j = -N, \cdots, 0, 1, \cdots, n_{\mathrm{max}}$;  \% generating the noise path.

\vspace*{1ex}

\item[] $(\xi_1^0, \xi_2^0) = (\xi_{1,0}, \xi_{2,0})$;  \hspace{5em}  \% initialization of  the system \eqref{Eq_reduced_gen} (given below).

\vspace*{1ex}
\item[]{\bf For} $ n  = 1:n_{\mathrm{max}}$  \hspace{2.5em} \% main loop to solve the {\mkk difference equations} \eqref{Eq_reduced_gen}. 

\item[] \% solving an analogue of the difference equations \eqref{uc_disc} with $q$ layers.
\begin{equation*}
\begin{array}{r}
\left.\begin{aligned}
\hspace*{4.5em}	y_1^{(1), 0} = \cdots = y_1^{(q), 0} = \xi_1^{n-1}; \\
\hspace*{4.5em}	y_2^{(1), 0} = \cdots = y_2^{(q), 0} = \xi_2^{n-1};
\end{aligned} \right \} \hspace{0.5em} \text{\% initialization.}
\end{array}
\end{equation*}

\bi 
\item[] {\bf For} $p = 1:N$ \hspace{1em}
\begin{equation*}
\begin{aligned}
& y_1^{(1),p}= \frac{ y_1^{(1),p-1}}{1 + \bigl( \beta_1(\lambda)+ \frac{\sigma^2}{2}\bigr) \delta t + \sigma {\sqrt{\delta t}} \zeta_{n-p}}; \qquad \text{ \% the first layer} \\
& y_2^{(1),p}= \frac{ y_2^{(1),p-1}}{1 + \bigl( \beta_2(\lambda)+ \frac{\sigma^2}{2}\bigr) \delta t + \sigma {\sqrt{\delta t}} \zeta_{n-p}}; \\
\end{aligned}
\end{equation*}

\bi

\item[] {\bf For} $l = 2: q$ \hspace{16em} \% the $l^{\mathrm{th}}$ layer
\begin{equation*}
\begin{aligned}
y_1^{(l),p}= \frac{ y_1^{(l),p-1} - \frac{\gamma \pi}{\sqrt{2}l^{3/2}} y_1^{(l-1),p} y_2^{(l-1),p} \delta t }{1 + \bigl( \beta_1(\lambda)+ \frac{\sigma^2}{2}\bigr) \delta t + \sigma {\sqrt{\delta t}} \zeta_{n-p}}; \\
y_2^{(l),p}= \frac{ y_2^{(l),p-1} - \frac{\gamma \pi}{\sqrt{2}l^{3/2}} (y_1^{(l-1),p})^2 \delta t}{1 + \bigl( \beta_2(\lambda)+ \frac{\sigma^2}{2}\bigr) \delta t + \sigma {\sqrt{\delta t}} \zeta_{n-p}};
\end{aligned}
\end{equation*}

\item[] {\bf End For}

$y_3^{(2),0}= 0; \; \;  y_4^{(2),0}= 0;$   \hspace{4em} \% initialization of the {\mkk difference equations} \eqref{us_disc} associated with . 

\vspace*{1ex}

{\bf For} $p=1:N$   \hspace{.4em} \% solving the forward difference equations \eqref{us_disc} with\\
  \hspace*{7em} \% the following changes

\beas
y_{1}^{(2),p}\leftarrow y_{1}^{(q),p},\\
y_{2}^{(2),p}\leftarrow y_{2}^{(q),p},\\
y_{3}^{(2),p}\leftarrow y_{3}^{(q),p},\\
y_{4}^{(2),p}\leftarrow y_{4}^{(q),p}\\
\eeas

\vspace{1ex}
{\bf End For}

\ei

\item[]  \% update $(\xi_1^{n}, \xi_2^{n})$ via the following analogue of the difference equations \eqref{eq:xi_disc}:
\bea\label{Eq_reduced_gen}
& {\hspace*{3em}}  \xi_1^{n} = \xi_1^{n-1} + \Bigl( \bigl ( \beta_1(\lambda)  + \frac{\sigma^2}{2} \bigr) \xi_1^{n-1} + \frac{\gamma\pi}{\sqrt{2}l^{3/2}} (\xi_1^{n-1} \xi_2^{n-1} + \xi_2^{n-1} y_{3}^{(q),N} + y_{3}^{(q),N} y_{4}^{(q),N}) \Bigr) \delta t  \\
& {\hspace*{18em}} + \sigma \zeta_n {\sqrt{\delta t}}\xi_1^{n-1}; \\ 
& {\hspace*{3em}} \xi_2^{n} =\xi_2^{n-1} + \Bigl( \bigl( \beta_2(\lambda) + \frac{\sigma^2}{2} \bigr) \xi_2^{n-1} + \frac{\sqrt{2} \gamma \pi}{l^{3/2}} \bigl( - \frac{1}{2} (\xi_1^{n-1})^2 + \xi_1^{n-1} y_{3}^{(q),N} + \xi_2^{n-1} y_{4}^{(q),N}\bigr) \Bigr) \delta t  \\
& {\hspace*{18em}} + \sigma \zeta_n {\sqrt{\delta t}} \xi_2^{n-1};
\eea

\vspace{1ex}
\item[] {\bf End For}

\ei

\item[] {\bf End For}

\end{enumerate}

\end{framed}

\subsection{Parameterization defect of $\widehat{h}_\lambda^{(q)}$, for $2\leq q\leq 10$: Numerical estimates}\label{PM_quality_CVE}

As mentioned in the previous subsection, the pseudocode such as introduced above is built on $u_{\s}^{(q)}[\xi(t,\omega)](t+\tau, \theta_{-\tau}\omega; 0)$.\footnote{For the sake of the discussion here, we have changed $T$ in \eqref{h_n} by $\tau$ here.}   In principle, the latter  has to be used for $\tau$ sufficiently large in order to approximate the pullback limit  \eqref{h_n}.  Interestingly,  in many cases a reasonably small $\tau$  is sufficient  to provide a good approximation in  terms of {\it time average} of the parameterization defect  of  
$\widehat{h}_\lambda^{(q)}$.  This is observed for instance for the case at hand,  as reported in Table \ref{Table_T}.

To understand the building blocks of this table,  let us introduce $\widehat{g}^{(q)}_{\lambda,\tau}(\xi,\omega):=u_{\s}^{(q)}[\xi](t+\tau, \theta_{-\tau}\omega; 0)$, and let us denote by $\mathcal{Q}_{\lambda,\tau}^{(q)}(T, \omega; u_0)$ the corresponding parameterization defect  associated with the manifold function $\widehat{g}^{(q)}_{\lambda,\tau}$.   By replacing in the above pseudocode, the current state $(\xi_1^{n}, \xi_2^{n})$ by $u_{\c}(n\delta t,\omega; u_0)$, we can numerically estimate  $\overline{\mathcal{Q}}_{\lambda,\tau}^{(q)}(T_1,T_2,\omega;u_0)=\frac{1}{T_2 - T_1}\int_{T_1}^{T_2} \mathcal{Q}_{\lambda,\tau}^{(q)}(T, \omega; u_0)\, \d T$, for $T_2>T_1$. The SPDE parameters are those specified in the caption of Fig. ~\ref{fig:ave_Q_sec10}, except the parameter $\lambda$ that is fixed to be $8\lambda_c$, and $\delta t=0.01$. As reported in Table \ref{Table_T}, a convergence up to the first four digits is observed to be reached starting from $\tau=6$, independently of the number  $q$ of layers used in \eqref{n_layer}.   These results have been observed to be robust with respect to  the choices of the realization, the initial datum $u_0$\footnote{which falls within the same basin of attraction of a stationary solution; see Section \ref{ss:Burgers reduction}.} and $T_2(>1000)$. 

 These results indicate strongly that in all the cases (except for $q=1$), the manifold function  $\widehat{g}^{(q)}_{\lambda,\tau}(\xi,\omega)$ (and thus $\widehat{h}_\lambda^{(q)}$)  for $\lambda = 8 \lambda_c$, provides a PM (at least in a time average sense).   We can thus reasonably state  that, for $\tau\geq2$,  $\widehat{g}^{(q)}_{\lambda,\tau}(\xi,\omega)$ provides a PM even far away from the critical value $\lambda_c$,  where in particular the amplitude of the solutions of the SPDE \eqref{eq:Burgers} gets large. Such a property has been checked to hold for $\lambda \leq 8 \lambda_{\c}$; see also Table   \ref{tab:h^n}.

\begin{table}  

\begin{tabular}{| c | c | c|c|c|c|c|c|c|c|c|}  
  \hline
  \diaghead(5,-3){\theadfont Diag Cl}%
{$\tau$}{$\overline{\mathcal{Q}}_{\lambda, \tau}^{(q)}$} &  $\overline{\mathcal{Q}}_{\lambda, \tau}^{(1)}$ & $\overline{\mathcal{Q}}_{\lambda, \tau}^{(2)}$& $\overline{\mathcal{Q}}_{\lambda, \tau}^{(3)}$& $\overline{\mathcal{Q}}_{\lambda, \tau}^{(4)}$ &$\overline{\mathcal{Q}}_{\lambda, \tau}^{(5)}$ & $\overline{\mathcal{Q}}_{\lambda, \tau}^{(6)}$& $\overline{\mathcal{Q}}_{\lambda, \tau}^{(7)}$ & $\overline{\mathcal{Q}}_{\lambda, \tau}^{(8)}$& $\overline{\mathcal{Q}}_{\lambda, \tau}^{(9)}$ & $\overline{\mathcal{Q}}_{\lambda, \tau}^{(10)}$ \\
\hline
$\tau=2$ & 2.0364 & 0.5583 & 0.9254 & 0.6472 & 0.6371 & 0.6485 & 0.6487 & 0.6482 & 0.6483 & 0.6483 \\
\hline
$\tau=4$ &    2.1333 & 0.5768 & 0.9921 & 0.6610 & 0.6540 & 0.6680 & 0.6679 & 0.6673 & 0.6674 & 0.6674 \\
\hline
$\tau=6$ &     2.1357 & 0.5774 & 0.9934 & 0.6612 & 0.6544 & 0.6685 & 0.6683 & 0.6678 & 0.6678 & 0.6678 \\
\hline
$\tau=8$ &     2.1357 & 0.5774 & 0.9935 & 0.6612 & 0.6544 & 0.6685 & 0.6683 & 0.6678 & 0.6678 & 0.6678 \\
\hline
$\tau=10$ &     2.1357 & 0.5774 & 0.9935 & 0.6612 & 0.6544 &  0.6685 & 0.6683 & 0.6678 & 0.6678 & 0.6678 \\
\hline
\end{tabular}
\vspace{2ex}    
\caption {{\footnotesize {\bf Time average of  parameterization defect for $\widehat{g}^{(q)}_{\lambda,\tau}$,  with $1\leq \tau \leq 10$,  and $1\leq q\leq 10$.}   Here,  $\overline{\mathcal{Q}}_{\lambda,\tau}^{(q)}(T_1,T_2,\omega;u_0)=\frac{1}{T_2 - T_1}\int_{T_1}^{T_2} \mathcal{Q}_{\lambda,\tau}^{(q)}(T, \omega; u_0)\, \d T$ is computed  for  $u_0 = 0.1e_1 + 0.2 e_2 + 0.1 e_5$, $T_1= 400$, $ T_2 = 1000$, a fixed realization $\omega$, for different pullback-time $\tau$ and for different number of layers $q$ used in \eqref{n_layer}.  The SPDE parameters are those specified  in the caption of Fig. ~\ref{fig:ave_Q_sec10}, except the parameter $\lambda$ that is fixed to be $8\lambda_c$, and $\delta t=0.01$. A convergence up to the first four digits is observed to be reached from $\tau=6$ independently of the number  $q$ of layers used in \eqref{n_layer}.   These results are robust w.r.t. the choices of the realization, the initial datum and $T_2(>1000)$, indicating that in all the cases, the manifold function  $\widehat{g}^{(q)}_{\lambda,\tau}(\xi,\omega)$ provides a PM (in an average sense).}} \label{Table_T} 
\end{table}

This analysis of the parameterization defect of PMs as obtained via a $q$-layer system  \eqref{n_layer},  is completed by a numerical  study of the convergence of $\overline{\mathcal{Q}}_{\lambda,\tau}^{(q)}$ for various values of $\lambda$, as the number of layers $q$ increases.  The results are reported in Table \ref{tab:h^n} for $\tau =6$. As one can observe, for a broad range of $\lambda$, a convergence of  $\overline{\mathcal{Q}}_{\lambda,\tau}^{(q)}$ takes place (up to four digit) after few iterations on $q$, while the dependence on $\lambda$ of $\overline{\mathcal{Q}}_{\lambda,\tau}^{(q)}$ is linear.  
From what has been observed in Table \ref{Table_T}, the (time average) parameterization defect of $\widehat{g}^{(q)}_{\lambda,6}$  approximating the  parameterization defect of  $\widehat{h}^{(q)}_\lambda$, we can state that  for a broad range of values of $\lambda$ and $q$,  the manifold function $\widehat{h}^{(q)}_\lambda$ gives access to a PM (in a time average sense).

Reminding that, given a fixed number of resolved modes,  the smaller the parameterization defect is, the smaller the modeling error is expected to be (cf. Proposition \ref{lem:PM error}), 
we have chosen for the derivation of effective reduced models for the case at hand,  the manifold function $\widehat{h}^{(2)}_\lambda$ that  exhibits here the smallest parameterization quality.

Second,  it can be observed numerically that a deterioration of  the parameterization defect of $\widehat{h}_\lambda^{(1)}$ (for fixed values of $\sigma$ and $\lambda$) takes place as $l$ increases; compare for instance Fig. ~\ref{fig:ave_Q} with Fig. ~\ref{fig:ave_Q_sec10} below.


\begin{table}  

\begin{tabular}{| c |c | c|c|c|c|c|c|c|c| }  
  \hline
  \diaghead{\theadfont Diag Column}%
{$\overline{\mathcal{Q}}_{\lambda, \tau}^{(q)}$}{{\large $\lambda$}}  &   $1.96 \lambda_c $  &   $2.71 \lambda_c $ &   $3.46 \lambda_c $ &  $ 4.22\lambda_c $ &    $4.98 \lambda_c $ &   $ 5.73\lambda_c $ &   $ 6.49\lambda_c $ &    $7.24\lambda_c $ &  $ 8.00\lambda_c $ \\
 \hline     
  $\overline{\mathcal{Q}}_{\lambda,\tau}^{(2)}$ & 0.3412 & 0.3734 & 0.4052 & 0.4332 & 0.4587 & 0.4831 & 0.5075 & 0.5325 & 0.5583 \\
\hline
$\overline{\mathcal{Q}}_{\lambda,\tau}^{(3)}$ & 0.3690 & 0.4208 & 0.4802 & 0.5444 & 0.6133 & 0.6869 & 0.7642 & 0.8441 & 0.9254 \\
\hline
$\overline{\mathcal{Q}}_{\lambda,\tau}^{(4)}$ & 0.3505 & 0.3893 & 0.4299 & 0.4688 & 0.5062 & 0.5425 & 0.5781 & 0.6130 & 0.6472 \\
\hline
$\overline{\mathcal{Q}}_{\lambda,\tau}^{(5)}$ & 0.3489 & 0.3869 & 0.4265 & 0.4641 & 0.5002 & 0.5353 & 0.5697 & 0.6037 & 0.6371 \\
\hline
$\overline{\mathcal{Q}}_{\lambda,\tau}^{(6)}$ & 0.3494 & 0.3878 & 0.4281 & 0.4667 & 0.5041 & 0.5407 & 0.5770 & 0.6129 & 0.6485 \\
\hline
$\overline{\mathcal{Q}}_{\lambda,\tau}^{(7)}$ & 0.3494 & 0.3879 & 0.4282 & 0.4669 & 0.5042 & 0.5409 & 0.5772 & 0.6131 & 0.6487 \\
\hline
$\overline{\mathcal{Q}}_{\lambda,\tau}^{(8)}$ & 0.3494 & 0.3878 & 0.4281 & 0.4668 & 0.5041 & 0.5407 & 0.5769 & 0.6128 & 0.6482 \\
\hline
$\overline{\mathcal{Q}}_{\lambda,\tau}^{(9)}$ & 0.3494 & 0.3878 & 0.4281 & 0.4668 & 0.5041 & 0.5407 & 0.5769 & 0.6128 & 0.6483 \\
\hline
$\overline{\mathcal{Q}}_{\lambda,\tau}^{(10)}$ & 0.3494 & 0.3878 & 0.4281 & 0.4668 & 0.5041 & 0.5407 & 0.5769 & 0.6128 & 0.6483 \\
\hline
\end{tabular}
{\vspace*{0.8em}}
\caption{{\footnotesize  {\bf Time average of the parameterization defect of $\widehat{g}^{(q)}_{\lambda,\tau}$, for $\tau=6$ and $2\leq q \leq 10$}. For each such a $q$, the parameterization defect is averaged over the same interval than used for Table \ref{Table_T}, as  $\lambda$  varies in $[1.96 \lambda_{c}, 8.00\lambda_c] $. The resulting $\overline{\mathcal{Q}}_{\lambda,\tau}^{(q)}$ are computed here for $\tau=6$. The realization $\omega$ used for the simulation is also the same than used for Table \ref{Table_T}.   Finally, the  SPDE parameters are those used for Fig. ~\ref{fig:ave_Q_sec10}. As one can see for  the values of $\lambda$ reported here,  a convergence of  $\overline{\mathcal{Q}}_{\lambda,\tau}^{(q)}$ takes place (up to four digit) after few iterations on $q$. This table shows also that PMs  are reached via $\widehat{g}^{(q)}_{\lambda,\tau}$ (and thus $\widehat{h}^{(q)}_\lambda$), even when $\lambda$ is far from its critical value $\lambda_c$,  where in particular the amplitude of the solutions of the SPDE \eqref{eq:Burgers} gets large. }}   \label{tab:h^n}
\end{table}

%
%
%
%

The parameterization defect of $\widehat{h}^{(2)}_\lambda$ is shown below in {\mkk Fig. ~\ref{fig:ave_Q_sec10}. 
As a comparison, the corresponding results for $\widehat{h}^{(1)}_\lambda$ derived in Section ~\ref{s:Burgers} associated with the one-layer auxiliary system \eqref{LLL} are also given. 
\begin{figure}[h]   
    \centering 
   \includegraphics [height=5cm, width=7cm]{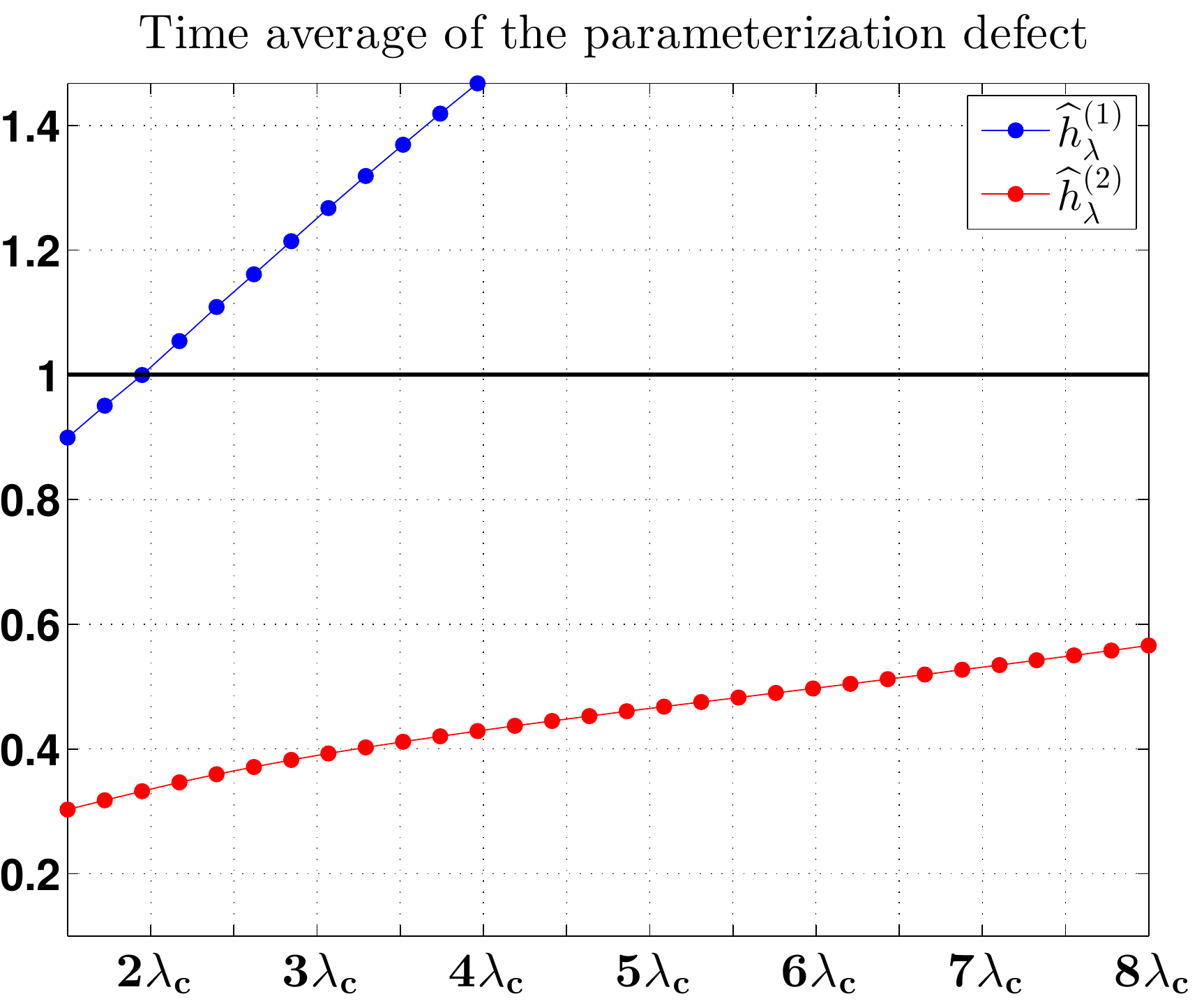}
\caption{{\footnotesize  {\bf Time average $\overline{\mathcal{Q}}_{1,\tau}$ and $\overline{\mathcal{Q}}_{2,\tau}$, associated with the  parameterization defect of $\widehat{h}^{(1)}_\lambda$ and $\widehat{h}^{(2)}_\lambda$.}
Same as Fig. ~\ref{fig:ave_Q}(a), the time average is computed for one fixed realization and over time interval $[400, 1000]$, but in a different parameter regime. The parameters for the SPDE \eqref{eq:Burgers} are chosen to be  $\gamma = 0.5$, $l = 3.5\pi$, $\nu = 2$, $\sigma = 0.4$, and various $\lambda$ in $[1.5\lambda_c, 8 \lambda_c]$ where $\lambda_c \approx 0.16$. The initial datum for all the simulations here is fixed to be $u_0 = 0.1e_1 + 0.2 e_2 + 0.1 e_5$. } }
\label{fig:ave_Q_sec10}
\end{figure}

Fig. ~\ref{fig:ave_Q_sec10} shows the time average of the parameterization defect $Q(t,\omega; u_0)$ associated with the PM function candidates $\widehat{h}^{(1)}_\lambda$ and $\widehat{h}^{(2)}_\lambda$ for an arbitrarily fixed realization of the Wiener process and for various values of $\lambda$ in $[1.5 \lambda_c, 8 \lambda_c]$. Recall that the time-dependent ratio $Q(t,\omega; u_0)$ is defined by Definition ~\ref{def:PM}.

As we can see for this parameter regime,  the parameterization defect of $\widehat{h}^{(2)}_\lambda$ associated with the two-layer backward-forward system \eqref{LLL_new} has improved by a factor  of two, compared with the one $\widehat{h}^{(1)}_\lambda$ associated with the one-layer system \eqref{LLL}. }

\subsection{Numerical results: Probability density and autocorrelation functions}\label{ss:num results} 
As mentioned at the end of Section ~\ref{ss:extreme}, the parameterization defect of $\widehat{h}^{(1)}_\lambda$ deteriorates as $\sigma$ increases (see Fig. ~\ref{fig:ave_Q}(b)), and as $l$ increases. Figure  ~\ref{fig:ave_Q_sec10} shows that the overall improvement (as $\lambda$ varies) of (time average) parameterization quality achieved by $\widehat{h}^{(2)}_\lambda$ is particularly significant.  

It is interesting to emphasize that such an improvement in the parameterization quality come with the combination of the new memory effects, as well as the matriochka of nonlinear interactions conveyed by $\widehat{h}^{(2)}_\lambda$ discussed in \SS ~\ref{Sec_Existence_PQ}. Such features are embodied with the new second, third and fourth order terms appearing in the analytic expression of $\widehat{h}^{(2)}_\lambda$ when compared with the expression of $\widehat{h}^{(1)}_\lambda$.

For the parameter regime analyzed in this section ({\bf regime B}: $\gamma = 0.5$, $\sigma = 0.4$, $\lambda = 1.7 \lambda_c$, $l = 3.5\pi$ and $\nu = 2$),  the manifold function $\widehat{h}^{(1)}_\lambda$ provides still a PM (in a time average sense),  but  of much poorer quality than $\widehat{h}^{(2)}_\lambda$ does.  It is interesting to note that for this regime, the values of the two NR-gaps  involved in the memory terms $\M^{12}_3(\cdot, \lambda)$ and $\M^{22}_4(\cdot, \lambda)$ are respectively $\beta_{1}(\lambda) + \beta_{2}(\lambda) - \beta_3(\lambda) = 0.9$ and $2\beta_{2}(\lambda) - \beta_4(\lambda)=1.6$. These gaps have been reduced essentially by half compared to regime A analyzed in  Section ~\ref{ss:extreme}.  This reduction of the NR-gaps results, according to Lemma ~\ref{lem:Mn}, corresponds to slower decay of  autocorrelations for the memory terms $\M^{12}_3$ and $\M^{22}_4$, when compared with  regime A.   The importance of memory effects such as conveyed by  the terms $\M^{12}_3$ and $\M^{22}_4$ has been illustrated in Section  ~\ref{ss:extreme} (see Fig. ~\ref{fig:extreme} again) for the  achievement  of good modeling performance by reduced systems based on  $\widehat{h}^{(1)}_\lambda$. Since here for regime B, the parameterization quality of the latter is poor\footnote{See Fig. ~\ref{fig:ave_Q_sec10} for $\lambda=1.7 \lambda_{c}$.} while the need of memory effects is still important, we turned  naturally to the reduced systems based on $\widehat{h}^{(2)}_\lambda$; the latter conveying more elaborated memory effects recalled above and discussed in  ~\ref{Sec_Existence_PQ}. Here analytic expressions of $\widehat{h}^{(2)}_\lambda$ could have been used based on the results of 
\SS ~\ref{Sec_Existence_PQ}, but we adopted the on-the-fly reduction procedure described in anterior subsections to show its performance. 

As Figure ~\ref{fig:fly} illustrates for the second mode amplitude, these new non-Markovian features conveyed by $\widehat{h}^{(2)}_\lambda$ allows us to achieve remarkable modeling performances. 
Clearly to achieve such results, the parameterizing manifold $\widehat{h}^{(2)}_\lambda$ gives access to a very good parameterization of the unresolved dynamics, from the resolved one. 

High-resolution numerical simulations were made to analyze in more details these modeling performance achieved by the reduced system \eqref{Reduced_AE} on the fly (based on $\widehat{h}^{(2)}_\lambda$). In that respect, the reproduction of statistical quantities such as probability density functions (PDFs) and the autocorrelation functions (ACFs) serves of evaluation criteria of the modeling performances.

\begin{figure}[!hbtp]
   \centering
   \includegraphics[height=8cm, width=14cm]{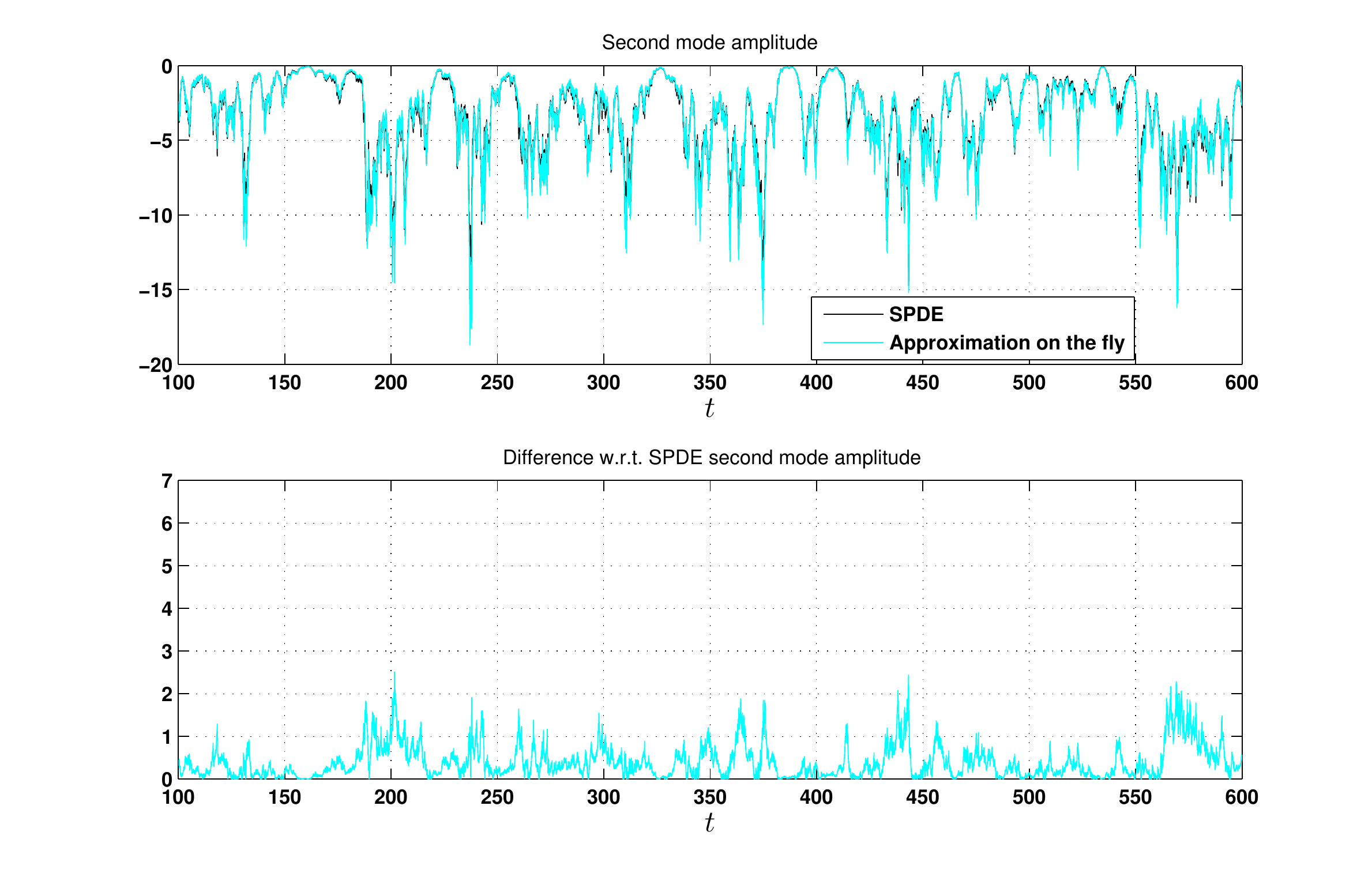}
  \caption{{\footnotesize {\bf Upper panel}: second mode amplitudes associated with respectively the two-layer reduced system \eqref{SDE abstract layered} (cyan curve) and the SPDE \eqref{eq:Burgers} (black curve). {\bf Lower panel}: The absolute value of the difference between the two curves shown in the upper panel. The system parameters are $\gamma = 0.5$, $\sigma = 0.4$, $\lambda = 1.7 \lambda_c$, $l = 3.5\pi$ and $\nu = 2$.}}   \label{fig:fly}
\end{figure}

\begin{figure}[!hbtp]
   \centering
  \hspace{-7ex}  \includegraphics[height=8cm, width=13cm]{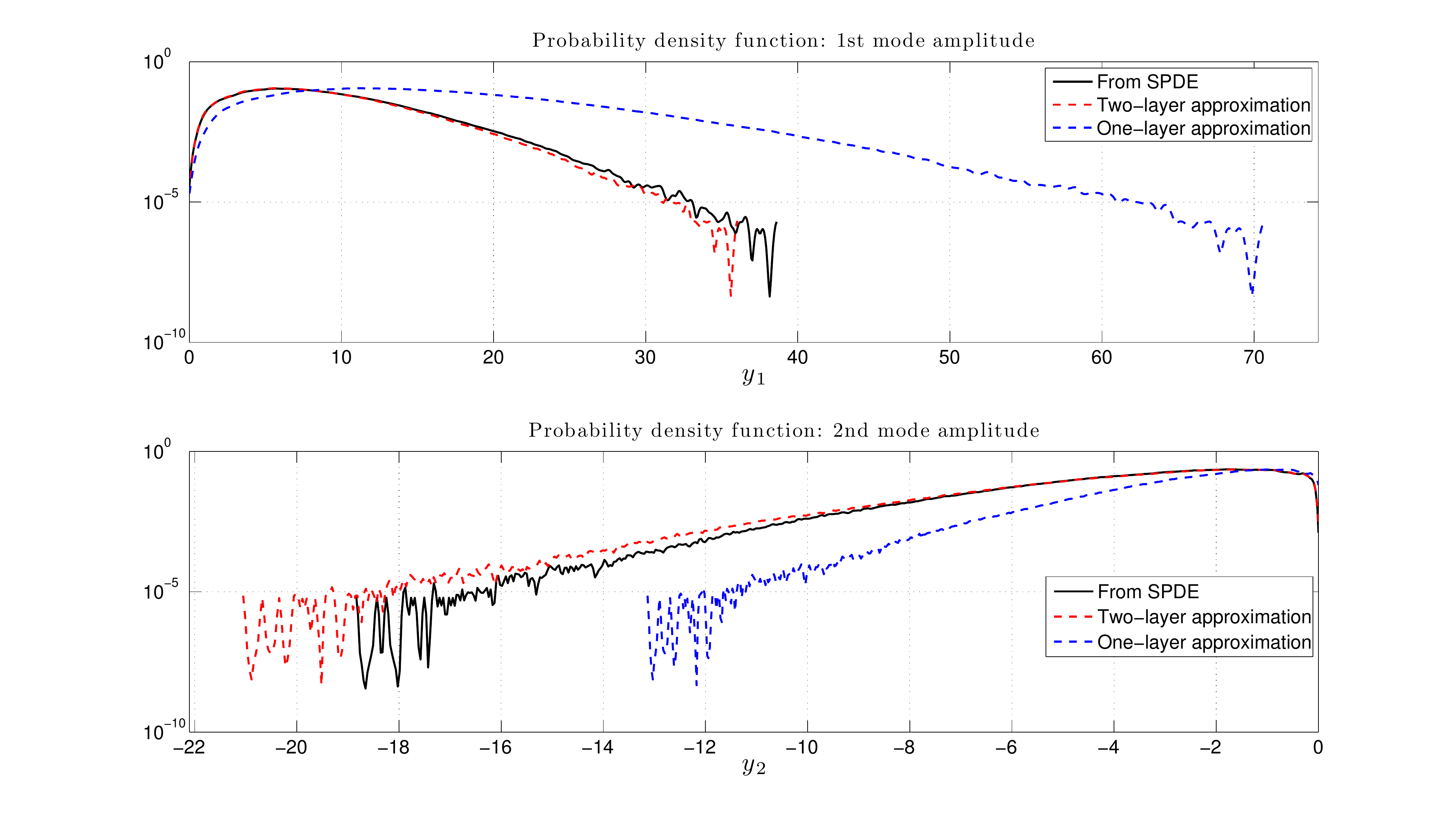}
  \caption{{\footnotesize Probability density functions of the first two modes modeled by the reduced systems \eqref{eq:Burgers reduced} (one-layer approximation) and \eqref{SDE abstract layered} (two-layer approximation) as well as the SPDE, where the system parameters are the same as used for Fig. ~\ref{fig:fly}.}}   \label{fig:PDF_approx3}
\end{figure}

\begin{figure}[!hbtp]
   \centering
  \hspace{-7ex}  \includegraphics[height=8cm, width=13cm]{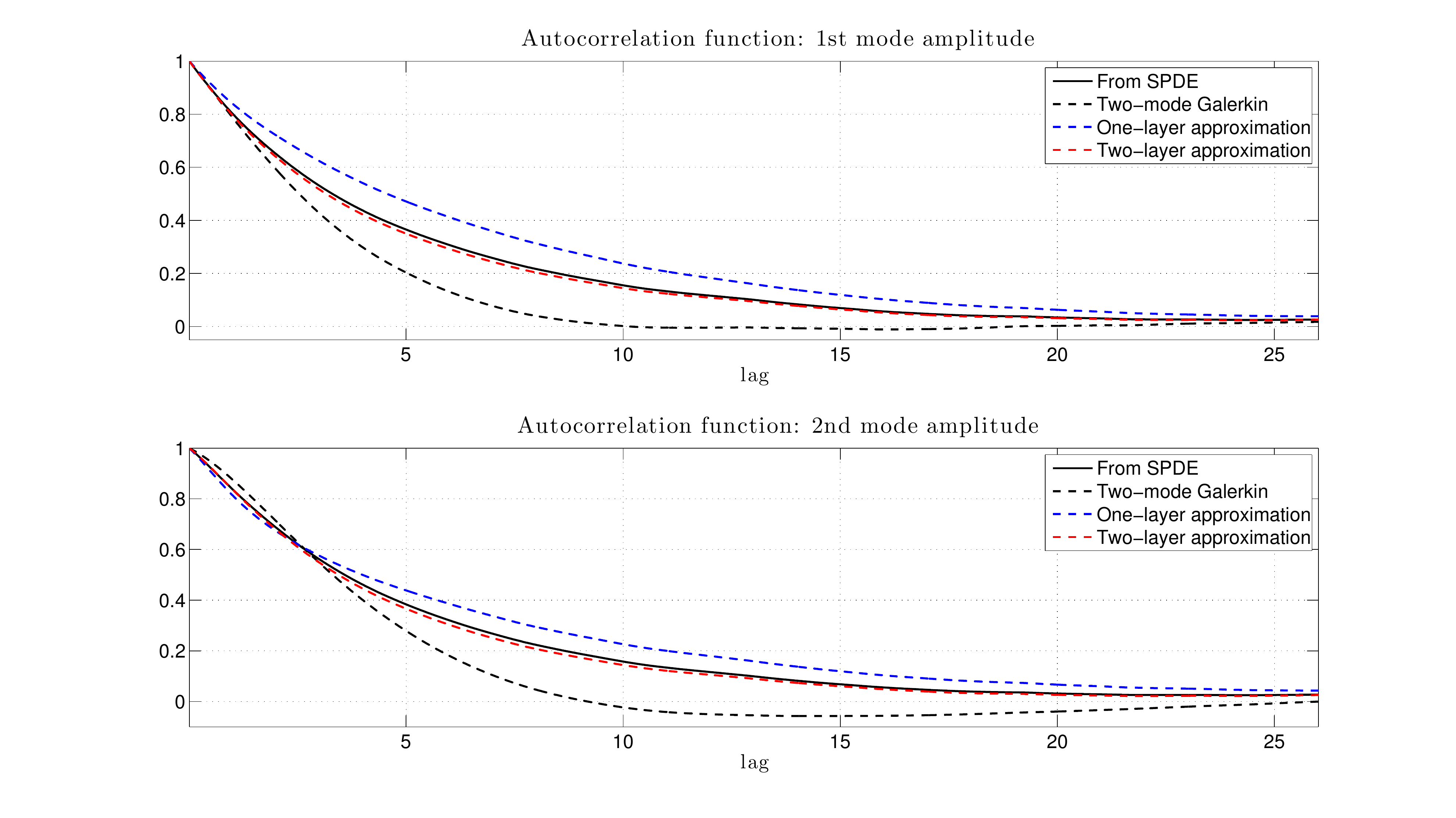}
  \caption{{\footnotesize Autocorrelation of the first two modes modeled by the two-mode Galerkin approximation, the reduced systems \eqref{eq:Burgers reduced} (one-layer approximation) and \eqref{SDE abstract layered} (two-layer approximation), and the SPDE, where the system parameters are the same as used for Fig. ~\ref{fig:fly}.}}   \label{fig:ACF}
\end{figure}

The results are reported in Fig. ~\ref{fig:PDF_approx3} and Fig. ~\ref{fig:ACF} for the first and second modes amplitudes. The corresponding PDFs and ACFs, as simulated from the reduced system \eqref{eq:Burgers reduced} based on   $\widehat{h}^{(1)}_\lambda$ and from the on-the-fly reduced system \eqref{Reduced_AE} based   on $\widehat{h}^{(2)}_\lambda$, are compared with the ones simulated by direct  integration of the SPDE according to the scheme described in Section ~\ref{ss:Burgers reduction}.  The results obtained via a basic  two-mode Galerkin reduced system (integrated via an Euler-Maruyama scheme), have been also included for reference.\footnote{The PDFs as simulated from a two-mode Galerkin reduced system are not shown in Fig. ~\ref{fig:PDF_approx3}, since, for instance, the estimated PDF of the second mode SPDE dynamics is overestimated by a factor of three, from such a reduced system.}

Three millions of  iterations (with $\delta t=0.01$) have been used, for each system,  to estimate the PDFs and ACFs, after removing the transient. 
 As it can be observed, the two-layer reduced system \eqref{SDE abstract layered}, with the more elaborated memory terms and nonlinear cross-interactions brought by $\widehat{h}^{(2)}_\lambda$,  achieves  very good performance in reproducing the statistical signature (such as PDF and ACF) of the SPDE dynamics  projected onto the resolved modes. One can observe furthermore that the modeling performance  based on $\widehat{h}^{(2)}_\lambda$  outperforms those based on $\widehat{h}^{(1)}_\lambda$, or the  two-mode Galerkin approximation.

In particular the modeling of the large excursions present in the SPDE dynamics projected onto the $1^{\textrm{st}}$ and $2^{\textrm{nd}}$ modes,  are   
 reproduced with high-accuracy from the on-the-fly reduced system   \eqref{Reduced_AE} based on $\widehat{h}^{(2)}_\lambda$. These large excursions come with an interesting spatial manifestation that is worthwhile to mention. Figure \ref{fig:snapshot} shows the SPDE solution profile at a particular time instance corresponding to a large excursion episode  as observed on the  $1^{\textrm{st}}$ and $2^{\textrm{nd}}$ modes; see Fig. ~\ref{fig:fly}. 
 As one can  observe, a concentration of the highest spatial velocity in a small region of the domain goes with such large excursions. Starting from a smooth initial condition, we have thus a steepening of the gradients localized in the field which fluctuates as the time flows.  This phenomenon is the manifestation of the  excitation of the small scales  by the noise through  the nonlinear term.  
The very good modeling performance achieved by the on-the-fly reduced system   \eqref{Reduced_AE} based on $\widehat{h}^{(2)}_\lambda$, relies thus on the ability of $\widehat{h}^{(2)}_\lambda$ to capture 
this noise-driven transfer of energy to the small scales as the time flows.  The attribution of this success to  the matriochka of self-interactions between the low modes, as well as the memory effects  such as conveyed by 
$\widehat{h}^{(2)}_\lambda$, is clear for the case at hand.  Fascinating problems are left in front of us for the generalization of such successes  regarding the reduction problem of SPDEs dynamics driven by more general noise.

\begin{figure}[!hbtp]
   \centering
   \includegraphics[height=5cm, width=9cm]{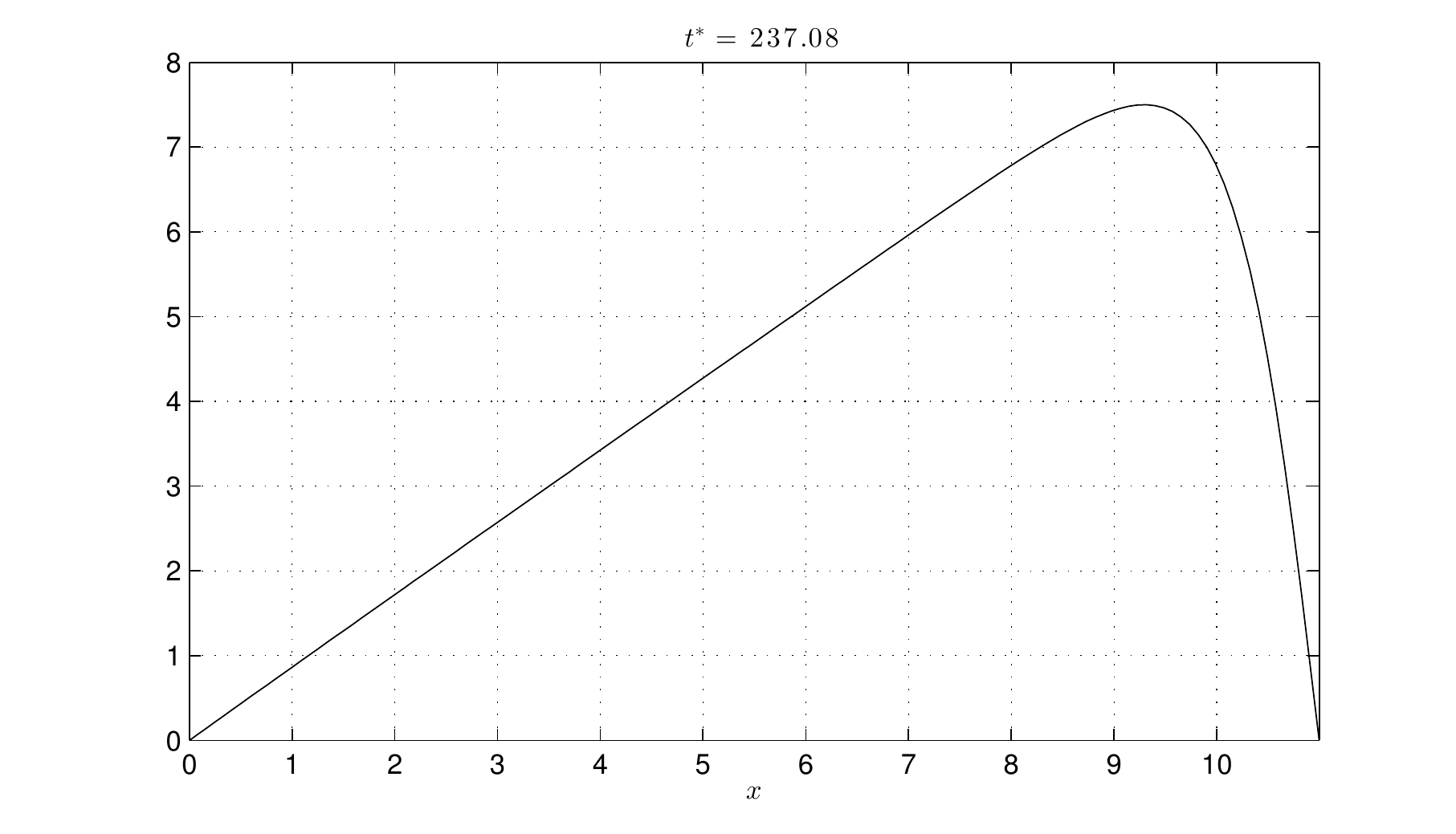}
  \caption{A snapshot of the SPDE solution profile.}   \label{fig:snapshot}
\end{figure}

\section*{Acknowledgments}

MDC and HL are supported by the National Science Foundation
 grant DMS-1049253 and Office of Naval Research grant
N00014-12-1-0911. SW is supported in part by National Science Foundation grants DMS-1211218, and DMS-1049114, and by Office of Naval Research grant N00014-11-1-0404. We thank Michael Ghil for his support and interest in this work. 

MDC is deeply grateful to James McWilliams and David Neelin for stimulating discussions on the closure problem of turbulence and stochastic parameterizations; and to Michael Ghil and Roger Temam for stimulating discussions on the slow and the ``fuzzy'' manifold, over the years.   MDC thanks also Jinqiao Duan for the reference  \cite{Kan_al12}, and for discussions on stochastic invariant manifolds.  

Preliminary versions of this work were presented by MDC and HL at the  AIMS Special Session on ``Advances in Classical and Geophysical Fluid Dynamics,'' held at the 9th AIMS Conference on Dynamical Systems, Differential Equations and Applications  in Orlando, July 2012;  at the ``workshop on Random  Dynamical Systems,'' held at the Institute of Mathematics and Applications in October 2012;  at  the ``Lunch Seminars'' held at the Center for Computational and Applied Mathematics, Purdue University, in November 2012; at the ``Applied Mathematics Seminar,''  held at the University of Illinois at Chicago in March 2013;  at AMS Special Session on ``Partial Differential Equations from Fluid Mechanics,'' held at University of Louisville in October 2013; and at the ``PDE seminar,'' held at Indiana University in October 2013.   We thank Jerry Bona, Jerome Darbon, Changbing Hu, Mike Jolly, Madelina Petcu, James Robinson, Taylan Sengul,  Jie Shen, Eric Simonnet, Yohann Tendero, Roger Temam, Florentina Tone and Kevin Zumbrun for their interest in this work and their respective invitations.

\appendix

\section{Classical and Mild Solutions of the Transformed RPDE} \label{Sect_mild}

In this appendix, we provide {\HLL for the sake of completeness}, a proof of Proposition ~\ref{prop:exist} regarding the existence and uniqueness of a measurable global classical solution to the transformed  Eq. ~\eqref{REE 1} for any given $\mathcal{H}_\alpha$-valued (random) initial datum.  We also introduce the definition of a mild solution to Eq. ~\eqref{REE 1}, and consider the existence and uniqueness of such mild solutions.

\bp[{\bf Proof of Proposition ~\ref{prop:exist}}]

We proceed in three steps.  Since Eq. ~\eqref{REE 1} becomes a non-autonomous PDE for each fixed $\omega$, the existence of a unique solution for each given deterministic initial datum in $\mathcal{H}_\alpha$ with the announced regularity given in ~\eqref{v regularity} can be proved by relying on the classical existence theory of solutions over finite time intervals and the half-line; see {\it e.g.} \cite[Thm.~3.3.3, Cor.~3.3.5]{Hen81}. The measurability property of such solutions requires more attention and details are provided in Step ~2 {\HL for the sake of completeness}. The measurability property of solutions with random initial data as claimed right after Proposition ~\ref{prop:exist} follows then from a basic composition argument {\HLL as explained} in Step ~3.

\medskip
{\bf Step 1}. For each fixed $\omega$, let us introduce
\beas
f_\omega(t, v) := z_\diffusion(\theta_t\omega)v + G(\theta_t\omega, v) = z_\diffusion(\theta_t\omega)v+ e^{-z_\diffusion(\theta_t\omega)} F(e^{z_\diffusion(\theta_t\omega)}v), \quad \Forall t \ge 0,\, v \in \mathcal{H}_\alpha.
\eeas
In order to apply the aforementioned non-autonomous theory, we first note that the following conditions hold naturally for $f_\omega(t,v)$:
\bea \label{v condition I}
\text{$f_\omega$ is locally H\"older continuous in $t$ and locally Lipschitz in $v$, $\Forall t\in \mathbb{R}^{+}, \, v \in \mathcal{H}_\alpha$}; 
\eea
and there exists a positive continuous function $Q_\omega \colon [0, \infty) \rightarrow \mathbb{R}^{+}$ depending on $\omega$ such that 
\bea \label{v condition II}
\|f_\omega(t, v)\| \le Q_\omega(t)(1+ \|v\|_\alpha), \quad \Forall  t \ge 0, \, v \in \mathcal{H}_\alpha.
\eea
These are indeed a direct consequence of the fact that the nonlinearity $F \colon \mathcal{H}_\alpha \rightarrow \mathcal{H}$ is assumed to be globally Lipschitz, and that the function $t \mapsto z_\diffusion(\theta_t\omega)$ is locally $\gamma$-H\"older continuous for all $\omega$ with any $\gamma \in (0, 1/2)$ according to Lemma ~\ref{Lem:OU}.

The continuous dependence of the solutions with respect to the initial data and the parameter $\lambda$  can be derived by using {\it e.g.} \cite[Thm.~3.4.4]{Hen81}. 

\medskip
{\bf Step 2}.   Now, we analyze the measurability of the solution for each fixed $\lambda$. First, we show that $v_\lambda(t, \cdot; v_0)$ is $(\mathcal{F}; \mathcal{B}(\mathcal{H}_\alpha))$-measurable for each fixed $t>0$ and $v_0\in \mathcal{H}_\alpha$. The treatment here is inspired by \cite{Caraballo_al09}. To this end, let us introduce for such given $t$ and $v_0$ the space 
\bea
X:= C_{v_0}([0, t]; \mathcal{H}_\alpha) :=\{ v \colon [0, t] \rightarrow \mathcal{H}_\alpha \; \vert \; v \text{ is continuous and } v(0) = v_0\}. 
\eea
Note that this space is a separable complete metric space when endowed with the metric
\bea
d_\gamma(v_1, v_2) := \sup_{t' \in [0,t]} e^{-\gamma t'}\|v_1(t') - v_2(t')\|_\alpha,  
\eea
where $\gamma > 0$.

Now, for each fixed $\omega$, let us introduce the following mapping $\mathcal{T}^{\omega, \lambda}_{t, v_0}$ defined on $X$:
\bea \label{eq:T operator}
\mathcal{T}^{\omega, \lambda}_{t, v_0} [v](s) := e^{s L_\lambda}v_0 + \int_0^{s} e^{(s-s')L_\lambda}f_\omega(s', v(s'))\d s',   \quad  \Forall \, s\in [0, t], \; v \in X. 
\eea
This mapping is well-defined due to \eqref{v condition II} and the fact that there exists $a_\lambda > 0$ and $C_1>0$ such that 
\bea \label{eq:c-b}
& \|e^{tL_\lambda}\|_{L(\mathcal{H}, \mathcal{H}_\alpha)}  \le \frac{C_1}{t^{\alpha}}
e^{a_\lambda t}, \qquad \Forall \, t > 0;
\eea
see \cite[Thm. ~44.5]{SY02} for a derivation of \eqref{eq:c-b}. The same reasons show that $\mathcal{T}^{\omega, \lambda}_{t, v_0}$ maps $X$ into itself.


We show now that $\mathcal{T}^{\omega, \lambda}_{t, v_0}$ is a contraction mapping on the space $\big(X, d_\gamma\big)$ for a sufficiently large $\gamma$ depending on $\omega$. First note that there exists $C_2>0$ such that the following inequality holds:
\bea  \label{eq:c}
& \|e^{tL_\lambda}\|_{L(\mathcal{H}_\alpha, \mathcal{H}_\alpha)} \le C_2 e^{a_\lambda t}, \qquad \Forall \, t \ge 0,
\eea
where $a_\lambda$ is the same as given in \eqref{eq:c-b}; see again \cite[Thm. ~44.5]{SY02}.

By using \eqref{eq:c-b}, \eqref{eq:c}, and the fact that $G$ and $F$ have the same Lipschitz constant $\LF$, we obtain
\bea \label{eq:contract-T}
\hspace{-1em}  d_\gamma(\mathcal{T}^{\omega, \lambda}_{t, v_0}[v_1],  \mathcal{T}^{\omega, \lambda}_{t, v_0}[v_2])  & \le \sup_{s\in [0,t]}e^{-\gamma s} \biggl \{ \int_0^{s} \|e^{(s-s')L_\lambda} z_\sigma(\theta_{s'}\omega) (v_1(s') - v_2(s'))\|_\alpha \d s' \\
& \hspace{4em} + \int_0^{s} \|e^{(s-s')L_\lambda}\bigl(
G(\theta_{s'}\omega, v_1(s')) - G(\theta_{s'}\omega, v_2(s'))\bigr)\|_\alpha \d s' \biggr \}  \\
&\le \sup_{s\in [0,t]}e^{-\gamma s}
\biggl ( C_2 \int_0^{s} e^{a_\lambda (s-s')} |z_\sigma(\theta_{s'}\omega)| \|v_1(s') - v_2(s')\|_\alpha \d s' \\
&\hspace{4em} + C_1 \LF \int_0^{s} \frac{e^{a_\lambda (s-s')}}{(s - s')^{\alpha}} \|v_1(s') - v_2(s')\|_\alpha \d s' \biggr) \\
& \le \sup_{s\in [0,t]} C d_\gamma(v_1, v_2) \int_0^{s} e^{(a_\lambda - \gamma) (s-s')} \biggl ( |z_\sigma(\theta_{s'}\omega)| + \frac{\LF}{(s - s')^{\alpha}} \biggr) \d s',
\eea
where $C = \max \{C_1, C_2\}$.

{\hh We conclude about the contractive property of  $\mathcal{T}^{\omega, \lambda}_{t, v_0}$  when $\gamma$ is sufficiently large. 
First note that from the  Lebesgue dominated convergence theorem the following property holds:
\beas
\forall s \; \in [0,t],\; g_{\gamma}(s):=\int_0^{s} e^{(a_\lambda - \gamma) (s-s')} \biggl ( |z_\sigma(\theta_{s'}\omega)| + \frac{\LF}{(s - s')^{\alpha}} \biggr) \d s' \underset{\gamma \rightarrow +\infty}\longrightarrow 0
\eeas
 Note also that  for  $\gamma'>\gamma >a_{\lambda}$, we have that $g_{\gamma'}(s)< g_{\gamma}(s)$ for all $s\in [0,t]$. Since  $[0,t]$ is furthermore  compact  it follows from the Dini Theorem \cite[Thm.~2.4.10]{Dud02} that
  $g_{\gamma}$ converges uniformly to zero on $[0, t]$ as $\gamma$ tends to $+\infty$. As a consequence, there exists $\gamma(\omega)$ sufficiently large, such that} 
\bea \label{eq:contract-T2}
{\hh \sup_{s\in [0,t]} C g_{\gamma(\omega)}(s)} = \sup_{s\in [0,t]} C \int_0^{s} e^{(a_\lambda - \gamma(\omega)) (s-s')} \biggl ( |z_\sigma(\theta_{s'}\omega)| + \frac{\LF}{(s-s')^{\alpha}} \biggr) \d s' < \frac{1}{2}.
\eea
For such a chosen $\gamma(\omega)$, we conclude then from  \eqref{eq:contract-T}--\eqref{eq:contract-T2} and the Banach fixed point theorem that there exists a unique fixed point $v^{\omega, \lambda}_{t, v_0} \in X$ of the operator $\mathcal{T}^{\omega, \lambda}_{t, v_0}$. Note also that the solution $v_\lambda(\cdot,\omega; v_0)$ obtained in Step ~1, when restricted to the time interval $[0, t]$, is in $X$ and is clearly a fixed point of $\mathcal{T}^{\omega, \lambda}_{t, v_0}$. {\HLL By uniqueness of the fixed point, we have thus:}
\be \label{eq:limit pt}
v_\lambda(s,\omega; v_0) = v^{\omega, \lambda}_{t, v_0}(s), \quad \Forall s \in [0, t].
\ee

Now we show that $v_\lambda(t,\cdot; v_0)$ is measurable. This results from the fact that $v_\lambda(t,\cdot; v_0)$ is obtained as the limit of a standard Picard scheme
 associated with $\mathcal{T}^{\,\cdot, \lambda}_{t, v_0}$ which  maps $(\mathcal{B}([0,t]) \otimes \mathcal{F}; \mathcal{B}(\mathcal{H}_\alpha))$-measurable mappings to $(\mathcal{B}([0,t]) \otimes \mathcal{F}; \mathcal{B}(\mathcal{H}_\alpha))$-measurable mappings, where $\mathcal{B}([0,t])$ denotes the trace $\sigma$-algebra of $\mathcal{B}(\mathbb{R})$ restricted to $[0,t]$. For the sake of clarity, we provide the details of such a standard argument.
 
Let us define a constant mapping $v_{t, v_0}^0 \colon  [0,t] \times \Omega \rightarrow \mathcal{H}_\alpha$ by:
\bea
v_{t, v_0}^0(s, \omega) \equiv v_0, \quad \Forall \, s\in [0, t], \, \omega \in \Omega.
\eea
Obviously, $v_{t, v_0}^0$ is $(\mathcal{B}([0,t]) \otimes \mathcal{F}; \mathcal{B}(\mathcal{H}_\alpha))$-measurable and $v_{t, v_0}^0(\cdot, \omega) \in X$ for each $\omega$.
For each $n \ge 1$, let us define recursively $v^n_{t, v_0} \colon [0,t] \times \Omega \rightarrow \mathcal{H}_\alpha$ by the following Picard scheme:
\bea \label{eq:picard}
v^n_{t, v_0}(s,\omega) := \mathcal{T}^{\omega, \lambda}_{t, v_0}[v_{t, v_0}^{n-1}(\cdot, \omega)](s), \quad \Forall \, s\in [0, t], \, \omega \in \Omega.
\eea
As noted above,  $v^n_{t, v_0}$ is $(\mathcal{B}([0,t]) \otimes \mathcal{F}; \mathcal{B}(\mathcal{H}_\alpha))$-measurable for each $n$. It then follows that  $v^n_{t, v_0}(t, \cdot)$ is $(\mathcal{F}; \mathcal{B}(\mathcal{H}_\alpha))$-measurable for each $n\ge 0$. 


Let $\gamma(\omega)>0$ be chosen such that $\mathcal{T}^{\omega, \lambda}_{t, v_0}$ is a contraction on $(X, d_{\gamma(\omega)})$. Then, we get from \eqref{eq:limit pt} that
\beas 
d_{\gamma(\omega)}\bigl( v_\lambda(\cdot,\omega; v_0), v^n_{t, v_0}(\cdot, \omega)\bigr) = d_{\gamma(\omega)}\bigl( v^{\omega, \lambda}_{t, v_0}(\cdot), v^n_{t, v_0}(\cdot, \omega)\bigr) \rightarrow 0 \quad \text{  as } \quad  n \rightarrow \infty,  \, \quad \, \Forall \omega \in \Omega.
\eeas 
This together with the definition of the metric $d_{\gamma(\omega)}$ implies that the following pointwise limit exists 
\bea \label{eq:conv RPDE}
v_\lambda(t,\omega; v_0) = \lim_{n\rightarrow \infty} v^n_{t, v_0}(t, \omega), \quad \, \Forall \omega \in \Omega, \, t > 0,
\eea 
where the limit is taken in $\mathcal{H}_\alpha$. Since $\mathcal{H}_\alpha$ is separable and the sequence $v^n_{t, v_0}(t, \cdot)$ is $(\mathcal{F}; \mathcal{B}(\mathcal{H}_\alpha))$-measurable, we conclude from \eqref{eq:conv RPDE} that the function $v_\lambda(t,\cdot; v_0) \colon \Omega \rightarrow \mathcal{H}_\alpha$ is also measurable.
 
Since $v_\lambda(\cdot, \omega; v_0)$ is continuous for each $\omega$; and $v_\lambda(t,\cdot; v_0)$ is measurable for each $t>0$, then $v_\lambda(\cdot,\cdot; v_0)$ is $\big(\mathcal{B}(\mathbb{R}^+) \otimes \mathcal{F}; \mathcal{B}(\mathcal{H}_\alpha)\big)$-measurable by using for instance \cite[Lemma III.14]{CV77}. 

Still based on \cite[Lemma III.14]{CV77}, since $v_\lambda(t, \omega; \cdot)$ is continuous for each $\omega$ and $t>0$, and $v_\lambda(\cdot,\cdot; v_0)$ is measurable for each $v_0$, then $v_\lambda$ is $\big(\mathcal{B}(\mathbb{R}^+) \otimes \mathcal{F} \otimes \mathcal{B}(\mathcal{H}_\alpha); \mathcal{B}(\mathcal{H}_\alpha)\big)$-measurable. 

\medskip
{\bf Step 3}. For any given measurable random initial datum $v_0(\omega)$, we can define $\omega$-wisely the solution $v_{\lambda, {\hh v_0(\omega)}}(t, \omega):= v_\lambda(t, \omega; v_0(\omega))$ as done in Step ~1. Clearly, the solution has the regularity given in \eqref{v regularity}. Let 
\bea
g \colon \mathbb{R}^{+} \times \Omega \rightarrow \mathbb{R}^{+} \times \Omega \times \mathcal{H}_\alpha, \quad (t, \omega) \mapsto (t, \omega, v_0(\omega)).
\eea
Then, we have $v_{\lambda, v_0(\omega)}(t,\omega) = v_\lambda \circ g(t,\omega)$. Since both $g$ and $v_\lambda$ are measurable, then $v_{\lambda,v_0(\omega)}$ is also measurable. The proof is complete.

\ep

Now, we introduce the following definition of mild solutions to Eq. ~\eqref{REE 1}.

\bd  \label{Def:mild solution}

Let $\mathcal{H}_\alpha$ be the interpolation space associated with the fractional power $A^\alpha$ for some $\alpha \in (0,1)$, where $A$ is given in \eqref{L}. Let $J:=[t_1, t_2]$ be a closed interval in $\mathbb{R}$, and $\lambda$ be a fixed value in $\mathbb{R}$. A mapping $v_\lambda \colon J \times \Omega \rightarrow \mathcal{H}_\alpha$, which is jointly measurable and continuous in $t$, is said to be a mild solution of Eq. ~\eqref{REE 1} in the space $\mathcal{H}_\alpha$ on $J$ with initial datum $v_\lambda(t_1, \omega) = v_0$, if it satisfies the following integral equation:
\bea \label{Eq:mild_soln}
v_\lambda(t,\omega) = \mathfrak{T}_{\lambda,\sigma}({t}, {t_1};\omega)  v_0 + \int_{t_1}^t \mathfrak{T}_{\lambda,\sigma}({t}, {s};\omega)   G(\theta_s\omega, v_\lambda(s, \omega))\,\mathrm{d} s, \quad \Forall \, t \in J, \, \omega \in \Omega,
\eea
where $\mathfrak{T}_{\lambda,\sigma}$ is the solution operator associated with the linearized equation of Eq. ~\eqref{REE 1} given in \SS~\ref{ss:T}.

We call $v_\lambda$ a mild solution of Eq. ~\eqref{REE 1} on $[t_1, \infty)$ if it is a mild solution on $[t_1, t_2]$ for any $t_2 > t_1$. Similarly, $v_\lambda$ is called a mild solution of Eq. ~\eqref{REE 1} on $(-\infty, t_2]$ if it is a mild solution on $[t_1, t_2]$ for any $t_1 < t_2$.

\ed

In Appendix ~\ref{appendix section 1}, we will make use of the following elementary lemma.

\bl  \label{lem:mild}

For any given $t_1 < t_2$, a measurable mapping $v_\lambda \colon [t_1, t_2] \times \Omega \rightarrow \mathcal{H}_\alpha$ is a mild solution to Eq. ~\eqref{REE 1} if and only if it satisfies the following integral equation:
\bea \label{Eq:mild_soln-3}
v_\lambda(t,\omega) &= \mathfrak{T}_{\lambda,\sigma}({t}, {t_2};\omega) P_{\c} v_\lambda(t_2, \omega)  - \int_t^{t_2} \mathfrak{T}_{\lambda,\sigma}({t}, {s};\omega) P_{\c} G(\theta_s\omega, v_\lambda(s, \omega))\,\mathrm{d} s \\
& + \mathfrak{T}_{\lambda,\sigma}({t}, {t_1};\omega) P_{\s} v_\lambda(t_1, \omega)  + \int_{t_1}^{t} \mathfrak{T}_{\lambda,\sigma}({t}, {s};\omega) P_{\s} G(\theta_s\omega, v_\lambda(s, \omega))\,\mathrm{d} s, \quad t \in [t_1, t_2].
\eea

\el

\bp


Assume that $v_\lambda$ is a mild solution of Eq. ~\eqref{REE 1} on $[t_1, t_2]$, we check that it satisfies \eqref{Eq:mild_soln-3}. First note that by the definition of mild solutions, we have
\bea
v_\lambda(t_2,\omega) = \mathfrak{T}_{\lambda,\sigma}({t_2}, {t_1};\omega)  v_\lambda(t_1, \omega) + \int_{t_1}^{t_2} \mathfrak{T}_{\lambda,\sigma}({t_2}, {s};\omega)   G(\theta_s\omega, v_\lambda(s, \omega))\,\mathrm{d} s,
\eea
which leads to
\bea
P_{\c}v_\lambda(t_2,\omega) = \mathfrak{T}_{\lambda,\sigma}({t_2}, {t_1};\omega)  P_{\c}v_\lambda(t_1, \omega) + \int_{t_1}^{t_2} \mathfrak{T}_{\lambda,\sigma}({t_2}, {s};\omega)   P_{\c} G(\theta_s\omega, v_\lambda(s, \omega))\,\mathrm{d} s,
\eea
where we used the fact that the solution operator $\mathfrak{T}_{\lambda, \sigma}$ leaves invariant the subspaces $\mathcal{H}^{\c}$ and $\mathcal{H}^{\s}$ as pointed out in \SS~\ref{ss:T}.

We then obtain from the identity above that
\bea  \label{eq:mild-111}
& \mathfrak{T}_{\lambda,\sigma}({t}, {t_2};\omega) P_{\c}v_\lambda(t_2,\omega) \\
& = \mathfrak{T}_{\lambda,\sigma}({t}, {t_2};\omega) \mathfrak{T}_{\lambda,\sigma}({t_2}, {t_1};\omega)  P_{\c}v_\lambda(t_1, \omega) \\
& \hspace{3em} + \mathfrak{T}_{\lambda,\sigma}({t}, {t_2};\omega) \int_{t_1}^{t_2} \mathfrak{T}_{\lambda,\sigma}({t_2}, {s};\omega)   P_{\c} G(\theta_s\omega, v_\lambda(s, \omega))\,\mathrm{d} s \\
&= \mathfrak{T}_{\lambda,\sigma}({t}, {t_1};\omega)  P_{\c}v_\lambda(t_1, \omega) + \int_{t_1}^{t_2} \mathfrak{T}_{\lambda,\sigma}({t}, {s};\omega)   P_{\c} G(\theta_s\omega, v_\lambda(s, \omega))\,\mathrm{d} s, \quad t \in [t_1, t_2].
\eea
By definition, we also have 
\beas
v_\lambda(t,\omega) &= \mathfrak{T}_{\lambda,\sigma}({t}, {t_1};\omega)  v_\lambda(t_1, \omega) + \int_{t_1}^{t} \mathfrak{T}_{\lambda,\sigma}({t}, {s};\omega)   G(\theta_s\omega, v_\lambda(s, \omega))\,\mathrm{d} s, \quad t\in [t_1, t_2],
\eeas
which can be rewritten as follows {\HLL by simply using the fact that $P_{\c} + P_{\s} = \Id$}:
\bea \label{eq:mild-112}
v_\lambda(t,\omega) &= \mathfrak{T}_{\lambda,\sigma}({t}, {t_1};\omega) P_{\c} v_\lambda(t_1, \omega) + \int_{t_1}^{t} \mathfrak{T}_{\lambda,\sigma}({t}, {s};\omega)  P_{\c} G(\theta_s\omega, v_\lambda(s, \omega))\,\mathrm{d} s \\ 
& \hspace{1em} + \mathfrak{T}_{\lambda,\sigma}({t}, {t_1};\omega) P_{\s}  v_\lambda(t_1, \omega) + \int_{t_1}^{t} \mathfrak{T}_{\lambda,\sigma}({t}, {s};\omega) P_{\s}  G(\theta_s\omega, v_\lambda(s, \omega))\,\mathrm{d} s, \; t\in [t_1, t_2].
\eea
{\HLL The expression of $\mathfrak{T}_{\lambda,\sigma}({t}, {t_1};\omega) P_{\c} v_\lambda(t_1, \omega)$ provided implicitly by \eqref{eq:mild-111} leads then to}:
\beas 
\mathfrak{T}_{\lambda,\sigma}({t}, {t_1};\omega) & P_{\c} v_\lambda(t_1, \omega) + \int_{t_1}^{t} \mathfrak{T}_{\lambda,\sigma}({t}, {s};\omega)  P_{\c} G(\theta_s\omega, v_\lambda(s, \omega))\,\mathrm{d} s \\
& = \mathfrak{T}_{\lambda,\sigma}({t}, {t_2};\omega) P_{\c}v_\lambda(t_2,\omega)  - \int_{t_1}^{t_2} \mathfrak{T}_{\lambda,\sigma}({t}, {s};\omega)   P_{\c} G(\theta_s\omega, v_\lambda(s, \omega))\,\mathrm{d} s \\
& \quad + \int_{t_1}^{t} \mathfrak{T}_{\lambda,\sigma}({t}, {s};\omega)  P_{\c} G(\theta_s\omega, v_\lambda(s, \omega))\,\mathrm{d} s \\
& =  \mathfrak{T}_{\lambda,\sigma}({t}, {t_2};\omega) P_{\c}v_\lambda(t_2,\omega)  - \int_{t}^{t_2} \mathfrak{T}_{\lambda,\sigma}({t}, {s};\omega)   P_{\c} G(\theta_s\omega, v_\lambda(s, \omega))\,\mathrm{d} s.
\eeas
Using this last identity in \eqref{eq:mild-112}, we obtain then \eqref{Eq:mild_soln-3}. The converse can be obtained in the same fashion, and we omit the details here.

\ep

The existence and uniqueness problem of mild solutions to Eq.~\eqref{REE 1} is easily solved by relying on 
 Proposition ~\ref{prop:exist} and the following elementary result.

\begin{prop} \label{prop_mild}

Consider the RPDE \eqref{REE 1}. The assumptions on $L_\lambda$ and $F$ are those of \SS~\ref{ss:SEE} where $F$ is assumed to be globally Lipschitz here; see \eqref{Lip F}.

Then, for each $\lambda \in \mathbb{R}$ and $t_1 < t_2$, a mapping $v_\lambda \colon [t_1, t_2] \times \Omega \rightarrow \mathcal{H}_\alpha$ is a mild solution of Eq. ~\eqref{REE 1} on the time interval $[t_1, t_2]$ if and only if $v_\lambda$ is a classical solution of Eq. ~\eqref{REE 1} in the sense given in Proposition ~\ref{prop:exist}.


\end{prop}

\bp

This is a direct generalization of \cite[Lemma ~3.3.2]{Hen81}. Indeed, by using the basic random change of variables $\widetilde{v}_\lambda(t, \omega) := e^{-\int_{t_1}^{t}  z_\sigma(\theta_\tau \omega) \, \mathrm{d}\tau \Id} v_\lambda(t, \omega)$, it can be checked that $v_\lambda$ is a classical solution of  Eq. ~\eqref{REE 1} on $(t_1, t_2)$ if and only if $\widetilde{v}_\lambda$ is a classical solution of the following equation
\bea \label{eq:REE transformed}
\frac{\mathrm{d} \widetilde{v}}{\mathrm{d} t} =  L_\lambda \widetilde{v} + \widetilde{G}(\theta_t\omega, \widetilde{v}),
\eea
where $\widetilde{G}(\theta_t\omega, \widetilde{v}) = e^{-\int_{t_1}^{t}  z_\sigma(\theta_\tau \omega) \, \mathrm{d}\tau \Id} G(\theta_t\omega, e^{\int_{t_1}^{t}  z_\sigma(\theta_\tau \omega) \, \mathrm{d}\tau \Id} \widetilde{v})$. Note also that for each $\omega$, the nonlinearity $\widetilde{G}(\theta_t\omega, \widetilde{v})$ clearly satisfies the conditions \eqref{v condition I}--\eqref{v condition II} with $\widetilde{G}$ in place of $f_\omega$ thanks to the Lipschitz property of $G$ with respect to $v$ and the H\"older continuity of the OU process $t\mapsto z_\sigma(\theta_t\omega)$. In particular, the conditions required in \cite[Lemma ~3.3.2]{Hen81} are met for Eq. ~\eqref{eq:REE transformed}, leading to the conclusion that $\widetilde{v}_\lambda$ is a classical solution to Eq. ~\eqref{eq:REE transformed} on $[t_1, t_2]$ if and only if it is a mild solution, namely if and only if it satisfies the following integral equation:
\bea \label{mild v tilde}
\widetilde{v}_\lambda(t,\omega) = e^{(t-t_1) L_\lambda} v_0 + \int_{t_1}^t e^{(t-s) L_\lambda}  \widetilde{G}(\theta_s\omega, \widetilde{v}_\lambda(s, \omega))\,\mathrm{d} s, \quad \Forall \, t \in [t_1, t_2], \, \omega \in \Omega,
\eea
where $v_0 = \widetilde{v}_\lambda(t_1,\omega)$. By the construction of $\widetilde{v}_\lambda$, one readily sees that $\widetilde{v}_\lambda$ satisfies \eqref{mild v tilde} if and only if $v_\lambda$ satisfies \eqref{Eq:mild_soln}. The proof is complete.

\ep

\section{Proof of Theorem \ref{Lip manfd}} \label{appendix section 1}

The proof is based on the Lyapunov-Perron method, and we adapt mainly the presentation {\HL of} \cite{DLS04} to our functional setting. We split the proof into four steps. In Step ~1, the {\HLL sought} random invariant manifold is characterized as the random set consisting of all {\HLL elements in $\mathcal{H}_\alpha$ such that there exists a complete trajectory of Eq. ~\eqref{REE} passing through each such element at $t=0$, which has controlled growth as $t\rightarrow -\infty$}. This characterization is shown to be equivalent to an integral equation to be satisfied by each such mild solution. The latter characterization via an integral equation is more suitable for analysis based on a fixed point argument relying on the partial-dichotomy estimates \eqref{eq:dichotomy}. A related {\HL fixed point} problem associated with this integral equation and parameterized by the critical variable $\xi \in \mathcal{H}^{\c}$ is analyzed in Step ~2 by usage of the uniform contraction mapping principle \cite[Theorems 2.1--2.2]{CH82}. In Step ~3, by taking the projection onto the non-critical space $\mathcal{H}_\alpha^{\s}$ of the solution (at $t=0$) obtained in Step ~2 for each given $\xi$, we build a random  $\mathcal{H}_\alpha^{\s}$-valued function defined on $\mathcal{H}^{\c}$, whose graph gives the sought random invariant manifold; the invariance property of the manifold being examined in Step ~4.

\bp[{\bf Proof of Theorem \ref{Lip manfd}}] We proceed in four steps as outlined above.

\medskip
{\bf Step 1. Characterization of $\mathfrak{M}_\lambda(\omega)$ via an integral equation.} For each  $\lambda \in \Lambda$ and $\omega \in \Omega$, let $\mathfrak{M}_\lambda(\omega)$ be the subset of $\mathcal{H}_\alpha$ defined by
\bea  \label{CM dynamical character}
& \mathfrak{M}_\lambda(\omega) := \bigl \{u_0 \in \mathcal{H}_\alpha \; \vert \;  \exists \; u_\lambda(\cdot, \omega; u_0) \in C_\eta^{-} \text{ with $u_\lambda(0, \omega; u_0) = u_0$}, \text{ which  is furthermore} \\
&  \hspace{6em} \text{  a mild solution of Eq. ~\eqref{REE} on $(-\infty, 0]$ in the sense of Definition ~\ref{Def:mild solution}.} \bigr \},
\eea 
where $C_\eta^{-}$ is the space defined in \eqref{space C-} with $\eta$ chosen according to condition \eqref{thm 3.1 condition}.

Note that $\mathfrak{M}_\lambda(\omega)$ thus defined is {\HL not empty}, since the origin is clearly contained in it. We will show in later steps that $\mathfrak{M}_\lambda$ so obtained is the global random invariant manifold that we look for. In the current step, we provide an equivalent characterization of this set. {\HLL Given any $u_0 \in \mathfrak{M}_\lambda(\omega)$, let $u_\lambda(\cdot, \omega; u_0)\in C_\eta^-$ be a corresponding mild solution of Eq. ~ \eqref{REE} on $(-\infty, 0]$ with  $u_\lambda(0, \omega; u_0) = u_0$. Namely, $u_\lambda$ satisfies the following:
\bea \label{Eq:bkwd-soln-2}
u_\lambda(t,\omega; u_0) &= \mathfrak{T}_{\lambda, \sigma}(t, \tau; \omega) u_\lambda(\tau, \omega; u_0) \\
& \hspace{2em} + \int_{\tau}^{t} \mathfrak{T}_{\lambda, \sigma}(t, \tau; \omega)G(\theta_s\omega, u_\lambda(s, \omega; u_0))\,\mathrm{d} s, \, \Forall \, \tau \le t \le 0,
\eea
where $\mathfrak{T}_{\lambda, \sigma}$ is the solution operator associated with the linearized stochastic flow as described in \SS~\ref{ss:T}.

According to Lemma ~\ref{lem:mild}, Eq. ~\eqref{Eq:bkwd-soln-2} can be rewritten as
\bea \label{Eq:bkwd-soln-2b}
& u_\lambda(t,\omega; u_0) = \mathfrak{T}_{\lambda,\sigma}({t}, {0};\omega) P_{\c} u_0  - \int_t^{0} \mathfrak{T}_{\lambda,\sigma}({t}, {s};\omega) P_{\c} G(\theta_s\omega, u_\lambda(s, \omega; u_0))\,\mathrm{d} s \\
& + \mathfrak{T}_{\lambda,\sigma}({t}, {\tau};\omega) P_{\s} u_\lambda(\tau, \omega; u_0)  + \int_{\tau}^{t} \mathfrak{T}_{\lambda,\sigma}({t}, {s};\omega) P_{\s} G(\theta_s\omega, u_\lambda(s, \omega; u_0))\,\mathrm{d} s, \, \Forall \, \tau \le t \le 0.
\eea

Projecting Eq. ~\eqref{Eq:bkwd-soln-2b} onto the subspace $\mathcal{H}^{\c}$, we obtain that}
\begin{equation} \label{Pu}
\begin{aligned}
 P_{\c} u_\lambda(t, \omega; u_0) & = \mathfrak{T}_{\lambda, \sigma}(t, 0; \omega) P_{\c} u_0 - \int_t^0 \mathfrak{T}_{\lambda, \sigma}(t, s; \omega) P_{\c} G(\theta_s\omega, u_\lambda(s, \omega; u_0))\,\mathrm{d} s,
\end{aligned}
\end{equation}
where we used the fact that $P_{\c}$ commutes with $L_\lambda$ and hence with $\mathfrak{T}_{\lambda, \sigma}$.

Projecting Eq. ~\eqref{Eq:bkwd-soln-2b} onto the subspace $\mathcal{H}_\alpha^{\s}$, we obtain that
\begin{equation} \label{Ps I}
\begin{aligned}
P_{\s} u_\lambda(t, \omega; u_0) & = \mathfrak{T}_{\lambda, \sigma}(t, \tau; \omega) P_{\s} u_\lambda(\tau, \omega; u_0) + \int_\tau^t \mathfrak{T}_{\lambda, \sigma}(t, s; \omega) P_{\s} G(\theta_s\omega, u_\lambda(s, \omega; u_0))\,\mathrm{d} s.
\end{aligned}
\end{equation}

Since $u_\lambda(\cdot, \omega; u_0) \in C^-_\eta$, using \eqref{d-I} we obtain for any $\tau \le t\le 0$ that
\begin{equation}
\begin{aligned}
\|\mathfrak{T}_{\lambda, \sigma}(t, \tau; \omega) P_{\s} u_\lambda(\tau, \omega; u_0)\|_\alpha & \le K e^{{\eta_2} (t -\tau) + \int_\tau^t z_\diffusion(\theta_s \omega) \,\mathrm{d}s} \|u_\lambda(\tau, \omega; u_0)\|_{\alpha} \\
& \le K e^{{\eta_2} t - \int_t ^0 z_\diffusion(\theta_s \omega)
\,\mathrm{d}s}e^{(\eta - {\eta_2})\tau} \|u_\lambda(\cdot, \omega; u_0)\|_{C_{\eta}^-}.
\end{aligned}
\end{equation}
{\HLL By recalling that $\eta_2 < \eta$ from condition \eqref{thm 3.1 condition}, we get that}, for each  $t \le 0$, the RHS above converges to zero as $\tau$ goes to $-\infty$. Now, by taking the limit $\tau \rightarrow -\infty$ in \eqref{Ps I}, we obtain that
\begin{equation} \label{Ps II}
\begin{aligned}
P_{\s} u_\lambda(t, \omega; u_0) = \int_{-\infty}^t \mathfrak{T}_{\lambda, \sigma}(t, s; \omega) P_{\s} G\big(\theta_s\omega,
u_\lambda(s, \omega; u_0)\big)\,\mathrm{d} s, \quad \Forall \, t \le 0.
\end{aligned}
\end{equation}

For each  $\xii \in \mathcal{H}^{\c}$, $\lambda \in \Lambda$, and $\omega \in \Omega$, we define an operator $ \mathcal{N}_{\xii}^{\omega, \lambda}$ as follows:
\begin{equation} \label{VCF} 
\begin{aligned}
 \mathcal{N}_{\xii}^{\omega, \lambda}[u](t) := \mathfrak{T}_{\lambda, \sigma}(t, 0; \omega) \xii & - \int_t^0 \mathfrak{T}_{\lambda, \sigma}(t, s; \omega) P_{\c} G(\theta_s\omega, u(s))\,\mathrm{d} s
\\
& + \int_{-\infty}^t \mathfrak{T}_{\lambda, \sigma}(t, s; \omega) P_{\s} G(\theta_s\omega, u(s))\,\mathrm{d} s, \quad t \le 0, \, u \in C_\eta^{-}.
\end{aligned}
\end{equation}

Combining \eqref{Pu} and \eqref{Ps II}, we obtain that $u_\lambda$ is a fixed point of $\mathcal{N}_{\xii}^{\omega, \lambda}$ with $\xii = P_{\c} u_0$. Conversely, if $u_\lambda(\cdot, \omega) \in C_\eta^-$ is a fixed point of \eqref{VCF} with $u_\lambda(0, \omega) = u_0$, we show in the following that it is also a mild solution to Eq. ~ \eqref{REE} on $(-\infty, 0]$. 

{\HLL According to Lemma ~\ref{lem:mild}, we only need to show that for any $\tau <0$, $u_\lambda(t, \omega)$ can be written into the form given in \eqref{Eq:bkwd-soln-2b} for any $t \in [\tau, 0]$.

First note that by applying $P_{\s}$ to both sides of \eqref{VCF} and setting $t = \tau$, we obtain that
\bea \label{Eq:u0-B}
P_{\s}u_\lambda(\tau, \omega) = P_{\s}\mathcal{N}_{\xii}^{\omega, \lambda}[u_\lambda](\tau) =  \int_{-\infty}^{\tau} \mathfrak{T}_{\lambda, \sigma}(\tau, s; \omega) P_{\s} G(\theta_s\omega, u_\lambda(s, \omega))\,\mathrm{d} s.
\eea
Note also that for any $t \in [\tau, 0]$, we have
\beas 
\int_{-\infty}^{t} \mathfrak{T}_{\lambda, \sigma}(t, s; \omega) P_{\s} G(\theta_s\omega, u_\lambda(s, \omega))\,\mathrm{d} s  & = \int_{\tau}^{t} \mathfrak{T}_{\lambda, \sigma}(t, s; \omega) P_{\s} G(\theta_s\omega, u_\lambda(s, \omega))\,\mathrm{d} s \\
& \hspace{1em} +  \int_{-\infty}^{\tau} \mathfrak{T}_{\lambda, \sigma}(t, s; \omega) P_{\s} G(\theta_s\omega, u_\lambda(s, \omega))\,\mathrm{d} s,
\eeas
and 
\beas
\int_{-\infty}^{\tau} \mathfrak{T}_{\lambda, \sigma}(t, s; \omega) P_{\s} G(\theta_s\omega, u_\lambda(s, \omega))\,\mathrm{d} s & = \mathfrak{T}_{\lambda, \sigma}(t, \tau; \omega) \int_{-\infty}^{\tau}  \mathfrak{T}_{\lambda, \sigma}(\tau, s; \omega) P_{\s} G(\theta_s\omega, u_\lambda(s, \omega))\,\mathrm{d} s \\
& = \mathfrak{T}_{\lambda, \sigma}(t, \tau; \omega)  P_{\s}u_\lambda(\tau, \omega),
\eeas
where we used \eqref{Eq:u0-B} to derive the last equality above. 

By combining the above two identities, we obtain
\beas
\int_{-\infty}^{t} \mathfrak{T}_{\lambda, \sigma}(t, s; \omega) P_{\s} G(\theta_s\omega, u_\lambda(s, \omega))\,\mathrm{d} s  & = \mathfrak{T}_{\lambda, \sigma}(t, \tau; \omega)  P_{\s}u_\lambda(\tau, \omega)  \\
& \hspace{2em} + \int_{\tau}^{t} \mathfrak{T}_{\lambda, \sigma}(t, s; \omega) P_{\s} G(\theta_s\omega, u_\lambda(s, \omega))\,\mathrm{d} s.
\eeas

Noting that $\xi = P_{\c} u_0$, \eqref{Eq:bkwd-soln-2b} follows then by using the above identity in the following fixed point equation satisfied by $u_\lambda$:
\bea  \label{eq:N 111}
u_\lambda(t, \omega) = \mathcal{N}_{\xii}^{\omega, \lambda}[u_\lambda](t) = \mathfrak{T}_{\lambda, \sigma}(t, s; \omega) \xii & -  \int_{t}^0 \mathfrak{T}_{\lambda, \sigma}(t, s; \omega) P_{\c} G(\theta_s\omega, u_\lambda(s,\omega))\,\mathrm{d} s \\
& + \int_{-\infty}^{t} \mathfrak{T}_{\lambda, \sigma}(t, s; \omega) P_{\s} G(\theta_s\omega, u_\lambda(s, \omega))\,\mathrm{d} s.
\eea
Hence, $u_\lambda$ is indeed a mild solution to Eq. ~\eqref{REE} on $(-\infty, 0]$ thanks to Lemma ~\ref{lem:mild}.}

Consequently, we obtain the following equivalent characterization of the set $\mathfrak{M}_\lambda(\omega)$ defined in \eqref{CM dynamical character}, that is for each  $\omega \in \Omega$ and $\lambda \in \Lambda$,
\bea \label{CM equiv character}
\bigl ( u_0 \in \mathfrak{M}_\lambda(\omega) \bigr) \Longleftrightarrow \Bigl( \exists\; u_\lambda(\cdot, \omega) \in C^-_\eta \textrm{ s.t. }  u_\lambda(0, \omega) = u_0, \;  u_\lambda = \mathcal{N}_{P_{\c} u_0}^{\omega, \lambda}[u_\lambda] \; \Bigr),
\eea
where $\mathcal{N}_{P_{\c} u_0}^{\omega, \lambda}$ is defined in \eqref{VCF} {\HL with $\xi = P_{\c}u_0$}.

\medskip
{\bf Step 2. Unique fixed point of $\mathcal{N}_{\xii}^{\omega, \lambda}$.} We show in this step that for each $\xii \in \mathcal{H}^{\c}$, $\lambda \in \Lambda$, and $\omega \in \Omega$ the operator $\mathcal{N}_{\xii}^{\omega, \lambda}$ defined in \eqref{VCF} has a unique fixed point in the space $C_\eta^{-}$.

We first check that 
\bea \label{eq:range_N}
\mathcal{N}_{\xii}^{\omega, \lambda} C_\eta^- \subset C_\eta^-.  
\eea
Note that by the partial-dichotomy estimate in \eqref{d-III} and the assumption that $\eta \in (\eta_2, \eta_1)$, we obtain that  
\beas
\sup_{t\in (-\infty, 0]} e^{-\eta t + \int_{t}^0 z_\diffusion(\theta_s\omega)\d s} \big \| \mathfrak{T}_{\lambda, \sigma}(t, 0; \omega) \xii \big \|_\alpha  & \le K \sup_{t\in (-\infty, 0]} e^{-\eta t } e^{\eta_1 t} \|\xii\|_\alpha  \le K \|\xii\|_\alpha, \quad \Forall \, \xi \in \mathcal{H}^{\c};
\eeas
and that 
\beas
& \sup_{t\in (-\infty, 0]} e^{-\eta t + \int_{t}^0 z_\diffusion(\theta_s\omega)\d s} \Bigl \| \int_t^0 \mathfrak{T}_{\lambda, \sigma}(t, s; \omega) P_{\c} G(\theta_s\omega, u(s))\,\mathrm{d} s \Bigr \|_\alpha \\
& \le \sup_{t\in (-\infty, 0]}  K \LF \|u\|_{C_\eta^-} \int_t^0 e^{({\eta_1} - \eta)(t-s)} \, \mathrm{d}s \\
& =\frac{K \LF}{{\eta_1} - \eta} \|u\|_{C_\eta^-}, \quad \qquad \Forall \, u \in C_\eta^{-}.
\eeas
Similarly, by \eqref{d-II}, we obtain that
\beas
 \sup_{t\in (-\infty, 0]} e^{-\eta t + \int_{t}^0 z_\diffusion(\theta_s\omega)\d s} & \Bigl \|\int_{-\infty}^t \mathfrak{T}_{\lambda, \sigma}(t, s; \omega) P_{\s} G(\theta_s\omega, u(s))\,\mathrm{d} s \Bigr \|_\alpha \\
&  \le K \LF \Gamma(1-\alpha) (\eta - {\eta_2})^{\alpha - 1}\|u\|_{C_\eta^-}.
\eeas
Using the above three estimates and the definition of $\mathcal{N}_{\xii}^{\omega, \lambda}$, we get
\bea
\|\mathcal{N}_{\xii}^{\omega, \lambda}[u]\|_{C_\eta^{-}} & = \sup_{t\in (-\infty, 0]} e^{-\eta t + \int_{t}^0 z_\diffusion(\theta_s\omega)\d s} \| \mathcal{N}_{\xii}^{\omega, \lambda}[u](t)\|_\alpha \\
& \le K\|\xii\|_\alpha + \Upsilon_{1}(F) \|u\|_{C_\eta^-},
\eea
where $\Upsilon_{1}(F)$ is defined in \eqref{thm 3.1 condition}. Thus, \eqref{eq:range_N} follows. 

Note that for any $u_1, u_2 \in C_\eta^-$ we have that
\begin{equation} \label{Upsilon est}
\begin{aligned}
& \hspace{-1em} \| \mathcal{N}_{\xii}^{\omega, \lambda} [u_1] - \mathcal{N}_{\xii}^{\omega, \lambda} [u_2] \|_{C_\eta^-} \\
 & \le \sup_{t\in (-\infty, 0]}\Biggl \{ e^{-\eta t + \int_t^0
z_\diffusion(\theta_s \omega) \,\mathrm{d} s} \biggl ( \Bigl \| \int_t^0
\mathfrak{T}_{\lambda, \sigma}(t, s; \omega) P_{\c}
\bigl(G(\theta_s\omega, u_1) - G(\theta_s\omega, u_2) \bigr)\,\mathrm{d} s \Bigr \|_\alpha \\
&\hspace{7em}+ \Bigl \| \int_{-\infty}^t \mathfrak{T}_{\lambda, \sigma}(t, s; \omega) P_{\s} \bigl(G(\theta_s\omega, u_1)
- G(\theta_s\omega,
u_2) \bigr)\,\mathrm{d} s \Bigr \|_\alpha \biggr ) \Biggr \} \\
& \le \sup_{t\in (-\infty, 0]} \Biggl \{ K \LF
\|u_1-u_2\|_{C_\eta^-}\biggl( \int_t^0 e^{({\eta_1} -
\eta)(t-s)} \, \mathrm{d}s + \int_{-\infty}^t \frac{e^{({\eta_2} -
\eta)(t-s)}}{(t-s)^\alpha} \,
\mathrm{d}s \biggr) \Biggr\} \\
& \le K \LF \bigl(({\eta_1} - \eta)^{-1} +  \Gamma(1-\alpha) (\eta - {\eta_2})^{\alpha - 1}
\bigr)\|u_1-u_2\|_{C_\eta^-} \\
& = \Upsilon_{1}(F)\|u_1-u_2\|_{C_\eta^-}.
\end{aligned}
\end{equation}
This together with assumption \eqref{thm 3.1 condition} implies that $\mathcal{N}_{\xii}^{\omega, \lambda}$ is a contraction mapping on $C_\eta^-$ with constant of contraction that is independent of $\omega$, $\lambda\in \Lambda$ and $\xii \in \mathcal{H}^{\c}$. Note also that by definition, $\mathcal{N}_{\xii}^{\omega, \lambda}$ is clearly Lipschitz in $\xii$, and for each fixed $u\in C_\eta^-$ the following estimate holds:
\be \label{Lip in xi}
\|\mathcal{N}_{\xii_1}^{\omega, \lambda}[u] - \mathcal{N}_{\xii_2}^{\omega, \lambda} [u] \|_{C_\eta^-} \le K \|\xii_1 - \xii_2\|_\alpha, \quad \Forall \xii_1, \xii_2 \in \mathcal{H}^{\c}.
\ee
We can apply now the uniform contraction mapping principle (see {\it e.g.} ~\cite[Theorems 2.1--2.2]{CH82}), which ensures that for each $\xii \in \mathcal{H}^{\c}$, the mapping $\mathcal{N}_{\xii}^{\omega, \lambda}$ has a unique fixed point
$u_{\lambda,\xii}(\cdot,  \omega) \in C_\eta^-$. Moreover, it follows from \eqref{Upsilon est} and \eqref{Lip in xi} that the mapping $\xii \rightarrow u_{\lambda,\xii}(\cdot,  \omega)$ is Lipschitz from $\mathcal{H}^{\c}$ to $C_\eta^-$:
\begin{equation} \label{Lip in x}
\begin{aligned}
\|u_{\lambda,\xii_1}(\cdot,  \omega) - u_{\lambda,\xii_2}(\cdot,  \omega)\|_{C_\eta^-} &= \|\mathcal{N}_{\xii_1}^{\omega, \lambda}[u_{\lambda,\xii_1}] - \mathcal{N}_{\xii_2}^{\omega, \lambda} [u_{\lambda,\xii_2}]\|_{C_\eta^-} \le \frac{K\|\xii_1-\xii_2\|_\alpha}{1- \Upsilon_{1}(F)}.
\end{aligned}
\end{equation}

Now we show that the mapping $(t, \omega, \xii) \mapsto u_{\lambda,\xii}(t,\omega)$ is jointly measurable for each fixed $\lambda \in \Lambda$. The argument is similar {\HL to the one} used in Step ~2 of the proof of Proposition ~\ref{prop:exist} {\HL and relies here on the Picard scheme associated with $u = \mathcal{N}_{P_{\c}u_0}^{\omega, \lambda}[u]$ where $u_0 = u(0, \omega)$}. Let us denote the zero mapping from $(-\infty, 0]\times \Omega$ to $\mathcal{H}_\alpha$ by $Z$; and define $u_{n} \colon (-\infty, 0]\times \Omega \rightarrow \mathcal{H}_\alpha$ for each $n\ge 1$ to be:
\bea
u_{n}(t,\omega) := (\mathcal{N}_{\xii}^{\omega, \lambda})^n [Z(\cdot, \omega)](t), \quad t \le 0, \, \omega \in \Omega.
\eea

Since $u_{\lambda,\xii}(\cdot, \omega)$ is obtained as the fixed point of $\mathcal{N}_{\xii}^{\omega, \lambda}$ in $C_\eta^{-}$ for each fixed $\omega$, and  $Z(\cdot, \omega)$ is obviously in $C_\eta^{-}$, we infer that
\bea \label{conv series}
u_{\lambda,\xii}(\cdot, \omega) = \lim_{n \rightarrow \infty} u_{n}(\cdot,\omega), \quad \Forall \, \omega \in \Omega,
\eea 
where the limit is taken in $C_\eta^{-}$. 

Note furthermore that the operator $\mathcal{N}_{\xii}^{\, \cdot, \lambda}$ maps $(\mathcal{B}((-\infty,0]) \otimes \mathcal{F}; \mathcal{B}(\mathcal{H}_\alpha))$-measurable mappings to $(\mathcal{B}((-\infty,0]) \otimes \mathcal{F}; \mathcal{B}(\mathcal{H}_\alpha))$-measurable mappings; where $\mathcal{B}((-\infty,0])$ denotes the trace $\sigma$-algebra of $\mathcal{B}(\mathbb{R})$ with respect to $(-\infty, 0]$. So, each $u_n$ has this measurability. As a consequence, $u_{\lambda,\xii}$ is also $(\mathcal{B}((-\infty,0]) \otimes \mathcal{F}; \mathcal{B}(\mathcal{H}_\alpha))$-measurable according to \eqref{conv series}. Moreover, since $u_{\lambda,\xii}(t, \omega)$ is Lipschitz in $\xii$ (hence continuous in $\xii$), by using for instance \cite[ Lemma III.14]{CV77}, we conclude that $u_{\lambda,\xii}(t, \omega)$ is jointly measurable in $t$, $\xii$ and $\omega$ for each fixed $\lambda \in \Lambda$.

\medskip
{\bf Step 3. Construction of the random invariant manifold $\mathfrak{M}_\lambda$.} Now, we show that the set  $\mathfrak{M}_\lambda$ defined in \eqref{CM dynamical character} is indeed a random invariant manifold as the graph of a random $\mathcal{H}_\alpha^{\s}$-valued function.

Let 
\bea \label{eq:h_def}
h_\lambda(\xii,\omega) := P_{\s} u_{\lambda,\xii}(0, \omega), \quad  \Forall \, \xi \in \mathcal{H}^{\c},\, \omega \in \Omega.
\eea 
Note that $h_\lambda(\xii, \omega)$ lives in $\mathcal{H}_\alpha^{\s}$ because $u_{\lambda,\xii}(0, \omega) \in \mathcal{H}_\alpha$ by construction.

By the measurability property of $u_{\lambda, \xi}$ derived in the previous step, $h_\lambda$ is measurable. Since $u_{\lambda, \xi}(\cdot, \omega)$ is the fixed point of $\mathcal{N}_{\xii}^{\omega, \lambda}$, we have in particular that
\bea \label{eq:fixed u0}
u_{\lambda, \xi}(0, \omega) = \mathcal{N}_{\xii}^{\omega, \lambda}[u_{\lambda, \xi}(\cdot, \omega)](0), \quad \Forall \, \omega \in \Omega.
\eea
{\HL Therefore, by applying $P_{\s}$ to \eqref{VCF} and setting $t$ to zero, we get from \eqref{eq:h_def}:}
\begin{equation} \label{h}
\begin{aligned}
 h_\lambda(\xii, \omega) = \int_{-\infty}^0 {\HL \mathfrak{T}_{\lambda, \sigma}(0, s; \omega)}  P_{\s} G(\theta_s\omega, u_{\lambda,\xii}(s, \omega))\,\mathrm{d}s,
\end{aligned}
\end{equation}
{\HL where we used the fact that the solution operator $\mathfrak{T}_{\lambda, \sigma}$ leaves invariant the subspaces $\mathcal{H}^{\c}$ and $\mathcal{H}^{\s}$ as pointed out in \SS~\ref{ss:T}.} {\HL By applying $P_{\c}$ to \eqref{VCF} we get also:}
\bea  \label{eq:center_part}
\xi  =  P_{\c}u_{\lambda, \xi}(0, \omega).
\eea
Now, \eqref{eq:h_def} and \eqref{eq:center_part} lead to
\bea \label{eq:u graph}
u_{\lambda,\xii}(0,\omega) = \xi + h_\lambda(\xi, \omega).
\eea
{\HL Recalling that $u_{\lambda, \xi}(\cdot, \omega) \in C_\eta^-$, we conclude} from \eqref{eq:u graph} and the equivalent characterization of $\mathfrak{M}_\lambda(\omega)$ given in \eqref{CM equiv character} that
\begin{equation} \label{eq:graph_rep}
\begin{aligned}
 \mathfrak{M}_\lambda(\omega) = \{ \xii + h_\lambda(\xii, \omega) \mid \xii \in \mathcal{H}^{\c} \}.
\end{aligned}
\end{equation}
Namely, $ \mathfrak{M}_\lambda(\omega)$ is the graph of $h_\lambda(\cdot, \omega)$ for all $\omega$.

We derive now {\HL the} properties that $h_\lambda$ satisfies as stated in the theorem. {\HL First} note that $u_{\lambda,\xi}(t,  \omega) \equiv 0$ if $\xi = 0$. This together with \eqref{h} and the fact that $G(\omega, 0) = 0$ for all $\omega$  implies that $h_\lambda(0, \omega) = 0$. {\HL By construction}, $h_\lambda$ depends continuously on $\lambda$ for $\lambda\in \Lambda$ since $u_{\lambda, \xi}$ does.
From \eqref{h}, we get for any $\xi_1, \xi_2 \in \mathcal{H}^{\c}$ that
\beas
h_\lambda(\xii_1, \omega) -  h_\lambda(\xii_2, \omega)= \int_{-\infty}^0 \mathfrak{T}_{\lambda, \sigma}(0, s; \omega)  P_{\s} \Bigl( G(\theta_s\omega, u_{\lambda,\xi_1}(s, \omega)) - G(\theta_s\omega, u_{\lambda,\xi_2}(s, \omega)) \Bigr) \mathrm{d}s.
\eeas
It then follows from \eqref{Lip in x} and \eqref{d-II} that $h_\lambda(\cdot, \omega)$ is Lipschitz in $\xi$ for each $\omega$: 
\begin{equation} \label{Lip h in proof}
\begin{aligned}
 \|h_\lambda(\xii_1, & \omega) -  h_\lambda(\xii_2, \omega)\|_\alpha \le \frac{K^2\LF(\eta-{\eta_2})^{\alpha-1}\Gamma(1-\alpha)}{1-\Upsilon_{1}(F)}\|\xii_1-\xii_2\|_\alpha, \; \; \Forall \, \xi_1, \xi_2 \in \mathcal{H}^{\c}.
\end{aligned}
\end{equation}
{\hh The estimate of the Lipschitz constant in \eqref{Lip h} follows now from \eqref{Lip h in proof}. 

By recalling from Step 1 that the fixed point $u_{\lambda,\xii}(\cdot, \omega)$ of $\mathcal{N}_{\xi}^{\omega, \lambda}$ is also the mild solution in $C_{\eta}^{-}$ to Eq. ~\eqref{REE} on $(-\infty, 0]$ with $u_{\lambda,\xii}(0, \omega) = \xi + h_\lambda(\xi, \omega)$, the integral representation \eqref{h eqn} follows then from \eqref{h} and the definition of $\mathfrak{T}_{\lambda,\sigma}$ provided in \SS ~\ref{ss:T}. }

Now we show that $\mathfrak{M}_\lambda$ is a random closed set. Since $h_\lambda$ is continuous in $\xii$, according to \eqref{eq:graph_rep}, $\mathfrak{M}_\lambda(\omega)$ is a closed subset of $\mathcal{H}_\alpha$ for each $\omega$. Then, by using the Selection Theorem \cite[Thm.~III.9]{CV77} or \cite[Thm.~2.6]{Crauel02}, in order to show that $\mathfrak{M}_\lambda$ is a random closed set, we only need to show that there exists a sequence $\{\gamma_n\}_{n \in \mathbb{N}}$ of measurable mappings
$\gamma_n \colon \Omega \rightarrow \mathcal{H}_\alpha$, such that
\beas
\mathfrak{M}_\lambda(\omega)=\overline{\{\gamma_n(\omega)\Vert n \in \mathbb{N}\}}^{\mathcal{H}_\alpha}, \quad \Forall \omega \in \Omega.
\eeas
Since $\mathcal{H}_\alpha$ is separable, there exists a sequence $\{u_n\} \in (\mathcal{H}_\alpha)^\mathbb{N}$ 
which is dense in $\mathcal{H}_\alpha$. Now, for each $u_n$, let us define $\gamma_n \colon \Omega \rightarrow \mathcal{H}_\alpha$ as follows
$$\omega \mapsto \gamma_n(\omega):=P_{\c} u_n + h_\lambda(P_{\c} u_n, \omega).$$
Clearly, $\gamma_n$ is measurable because $h_\lambda(P_{\c} u, \cdot)$ is
measurable for any fixed $u \in \mathcal{H}_\alpha$. 

Since $\mathfrak{M}_\lambda(\omega)$ is closed for each $\omega \in \Omega$, we have by the construction of $\gamma_n$ that
\beas
\overline{\{\gamma_n(\omega)\Vert n \in \mathbb{N}\}}^{\mathcal{H}_\alpha} \subset \mathfrak{M}_\lambda(\omega), \quad \Forall \omega \in \Omega.
\eeas

Now, we show that the reverse inclusion holds. Let $\omega$ be fixed in $\Omega$. For each $u \in \mathfrak{M}_\lambda(\omega)$, it can be written as $u = P_{\c} u + h_\lambda(P_{\c}u, \omega)$. By the definition of $\{u_n\}$, there exists a subsequence, stilled denoted by $\{u_n\}$, which converges to $u$, leading to the convergence of $\{P_{\c}u_n\}$  to $P_{\c}u$. Then, by the continuity of $h_\lambda(\cdot, \omega)$ and the definition of $\gamma_n$, we have that $\{\gamma_n(\omega)\}$ converges to $u$. It follows that 
\beas
\mathfrak{M}_\lambda(\omega) \subset \overline{\{\gamma_n{\omega} \mid n\in
\mathbb{N}^\ast\}}^{\mathcal{H}_\alpha}, \quad \Forall \omega \in \Omega.
\eeas
We have thus proved that $\mathfrak{M}_\lambda$ is a random closed set. 

\medskip
{\bf Step 4. Invariance property of $\mathfrak{M}_\lambda$.} We show {\HL in this last step} that $\mathfrak{M}_\lambda$ is invariant, {\it i.e.},
\bea  \label{eq:inv-goal}
S_\lambda(t, \omega) \mathfrak{M}_\lambda(\omega) \subset \mathfrak{M}_\lambda(\theta_t\omega), \quad \Forall \, t > 0, \, \omega \in \Omega,
\eea 
where $S_\lambda$ is the RDS associated with Eq. ~\eqref{REE}.

For each fixed $\omega$, $u_0 \in \mathfrak{M}_\lambda(\omega)$, and $t > 0$, let $u_\lambda(s, \omega; u_0)$, $s \in [0, t]$, be the mild solution to Eq. ~ \eqref{REE} with initial datum $u_0$ (in the fiber $\omega$). By Proposition ~\ref{prop_mild}, we know that this mild solution exists, and is also a classical solution. Thus, according to the definition of $S_\lambda$, we have
\bea
S_\lambda(t, \omega) u_0 = u_\lambda(t, \omega; u_0).
\eea
Our goal is then to show that 
\bea \label{eq:inv-goal-1}
S_\lambda(t, \omega) u_0 \in \mathfrak{M}_\lambda(\theta_t\omega).
\eea

Since $u_0\in \mathfrak{M}_\lambda(\omega)$, from the characterization of $\mathfrak{M}_\lambda(\omega)$ provided in \eqref{CM dynamical character}, we can extend the above mild solution to the interval $(-\infty, t]$ and the following property holds:
\bea \label{eq:inv-1}
u_\lambda(s, \omega; u_0)\vert_{s\in (-\infty, 0]} \in C_\eta^-(\omega).
\eea

Now, let
\bea
v_\lambda(s, \theta_t\omega; u_\lambda(t, \omega; u_0)) := u_\lambda(s+t, \omega; u_0), \quad \Forall\, s \le 0.
\eea
Note that $v_\lambda(s, \theta_t\omega; u_\lambda(t, \omega; u_0))$ is a mild solution to Eq. ~\eqref{REE} on $(-\infty, 0]$ {\HLL which takes value $u_\lambda(t, \omega; u_0)$ in the fiber $\theta_t\omega$ when $s=0$}. Note also that
\bea
v_\lambda(0, \theta_t\omega; u_\lambda(t, \omega; u_0))  = u_\lambda(t, \omega; u_0) = S_\lambda(t, \omega) u_0.
\eea
Then, according to the characterization provided in \eqref{CM dynamical character} adapted to $\mathfrak{M}_\lambda(\theta_t\omega)$,
in order to check \eqref{eq:inv-goal-1}, we only need to show that 
\bea  \label{eq:inv-goal-2}
v_\lambda(\cdot, \theta_t\omega; u_\lambda(t, \omega; u_0)) \in C_\eta^-(\theta_{t}\omega).
\eea

Note that 
\bea \label{eq:inv-2}
\|v_\lambda(\cdot, \theta_t\omega; u_\lambda(t, \omega; u_0))\|_{C_\eta^-(\theta_{t}\omega)}  &=  \sup_{s \le 0} e^{-\eta s + \int_s^0 z_\sigma(\theta_{\tau + t}\omega)\d \tau}\|v_\lambda(s, \theta_t\omega; u_\lambda(t, \omega; u_0))\|_\alpha \\
& = \sup_{s \le 0} e^{-\eta s + \int_{s+t}^t z_\sigma(\theta_{\tau}\omega)\d \tau}\|u_\lambda(s+t, \omega; u_0)\|_\alpha \\
& = \sup_{s' \le t} e^{-\eta (s'-t) + \int_{s'}^t z_\sigma(\theta_{\tau}\omega)\d \tau}\|u_\lambda(s', \omega; u_0)\|_\alpha \\
& = e^{ \eta t + \int_{0}^t z_\sigma(\theta_{\tau}\omega)\d \tau} \sup_{s' \le t} e^{-\eta s' + \int_{s'}^0 z_\sigma(\theta_{\tau}\omega)\d \tau}\|u_\lambda(s', \omega; u_0)\|_\alpha.
\eea
Since $u_\lambda(s, \omega; u_0)$ is continuous on $[0, t]$, it clearly holds that
\bea \label{eq:inv-3}
\sup_{s' \in [0, t]} e^{-\eta s' + \int_{s'}^0 z_\sigma(\theta_{\tau}\omega)\d \tau}\|u_\lambda(s', \omega; u_0)\|_\alpha < \infty.
\eea
{\HLL The control on $(-\infty, 0]$ by some finite constant is achieved thanks to \eqref{eq:inv-1}. From \eqref{eq:inv-2}, we obtain then that $\|v_\lambda(\cdot, \theta_t\omega; u_\lambda(t, \omega; u_0))\|_{C_\eta^-(\theta_{t}\omega)} < \infty$, namely \eqref{eq:inv-goal-2} is satisfied, which leads in turn to the invariance property \eqref{eq:inv-goal-1}.}


We have thus checked that the set $\mathfrak{M}_\lambda$ defined in \eqref{CM dynamical character} satisfies all the conditions required in Definition ~\ref{def_GM} in order to be a global random invariant Lipschitz manifold. The proof is now complete.

 \ep

\section{Proof of Lemma ~\ref{lem:Mn}} \label{sec:Mn proof}

\bp[{\bf Proof of Lemma ~\ref{lem:Mn}.}]


As mentioned before, the main ingredients to derive the results stated in Lemma ~\ref{lem:Mn} are Fubini's Theorem, the independent increment property of the Wiener process, and the fact that $\mathbb{E}(e^{\diffusion W_t(\cdot)}) = e^{\diffusion^2|t|/2}$ for any $t\in \mathbb{R}$ as expectation of the geometric Brownian motion generated by $\d S_t = \frac{\sigma^2}{2} S_t + \sigma S_t \d W_t$; see {\it e.g.} ~\cite[Sect.~5.1]{Oksendal98}.

Let us begin with the derivation of (i). By using the expression of $\M_{n}(\omega, \lambda)$ recalled in \eqref{Mn recall2}, we obtain that
\bea
\mathbb{E}(\M_{n}(\cdot, \lambda)) & = \mathbb{E} \biggl( \int_{-\infty}^0 e^{g(\lambda)s + \diffusion (k-1)W_s(\cdot)}\, \mathrm{d}s \biggr) \\
& =  \int_{-\infty}^0 e^{g(\lambda)s} \mathbb{E}(e^{\diffusion (k-1)W_s(\cdot)}) \mathrm{d}s \\
& =  \int_{-\infty}^0 e^{g(\lambda)s} e^{\diffusion^2 (k-1)^2 |s|/2} \mathrm{d}s.
\eea
Since the above integral is finite if and only if $g(\lambda) - \diffusion^2 (k-1)^2/2 > 0$, the result stated in (i) follows.

To prove (ii), we first compute $\mathbb{E}(\M_{n}(\cdot, \lambda)\M_{n}(\cdot, \lambda))$. In order to simplify the notations, let use introduce 
\bea \label{eq:gamma}
\gamma := \diffusion (k-1).
\eea
Using again the expression of $\M_{n}(\omega, \lambda)$ in \eqref{Mn recall2} and Fubini's Theorem, we obtain
\bea  \label{Mn 2nd moment}
\mathbb{E}(\M_{n}(\cdot, \lambda)\M_{n}(\cdot, \lambda)) & = \mathbb{E} \biggl( \int_{-\infty}^0 e^{g(\lambda)s + \diffusion (k-1)W_s(\cdot)}\, \mathrm{d}s \int_{-\infty}^0 e^{g(\lambda)s' + \diffusion (k-1)W_{s'}(\cdot)}\, \mathrm{d}s' \biggr) \\
& = \mathbb{E} \biggl( \int_{-\infty}^0 \int_{-\infty}^0  e^{ g(\lambda)s + g(\lambda)s' + \gamma W_s(\cdot) + \gamma W_{s'}(\cdot)}\, \mathrm{d}s \mathrm{d}s' \biggr) \\
& = \int_{-\infty}^0 \int_{-\infty}^0  e^{ g(\lambda)s + g(\lambda)s' } \mathbb{E} \bigl (e^{\gamma (W_s(\cdot) + W_{s'}(\cdot))}\bigr ) \mathrm{d}s \mathrm{d}s'.
\eea
Note that $W_{s'}(\cdot) - W_s(\cdot)$ is independent of $2 W_s(\cdot)$ when $s' < s\le 0$. Note also that $W_{s'}(\omega) - W_s(\omega) = W_{s' - s}(\theta_s(\omega))$, which leads to 
$\mathbb{E} \bigl (e^{\gamma (W_{s'}(\cdot) - W_s(\cdot))}\bigr) = \mathbb{E} \bigl (e^{\gamma (W_{s'-s}(\cdot))}\bigr) = e^{\frac{\gamma^2}{2}|s'-s|}$. We obtain then
\beas
\mathbb{E} \bigl (e^{\gamma (W_s(\cdot) + W_{s'}(\cdot))}\bigr) & = \mathbb{E} \bigl (e^{\gamma (W_{s'}(\cdot) - W_s(\cdot))} e^{2\gamma W_s(\cdot)}\bigr)  = \mathbb{E} \bigl (e^{\gamma (W_{s'}(\cdot) - W_s(\cdot))}\bigr) \mathbb{E} \bigl(e^{2\gamma W_s(\cdot)}\bigr) \\
& = e^{\frac{\gamma^2}{2} (s-s')} e^{- 2 \gamma^2 s}, \qquad \Forall s' < s \le 0. 
\eeas
Similarly, we have
\beas
\mathbb{E} \bigl (e^{\gamma (W_s(\cdot) + W_{s'}(\cdot))}\bigr) & = e^{\frac{\gamma^2}{2} (s'-s)} e^{- 2 \gamma^2 s'}, \qquad \Forall s < s' \le 0. 
\eeas
The expression of $\mathbb{E}(\M_{n}(\cdot, \lambda)\M_{n}(\cdot, \lambda))$ given in \eqref{Mn 2nd moment} can then be rewritten as 
\bea  \label{eq:2nd moments}
\mathbb{E}(\M_{n}(\cdot, \lambda)\M_{n}(\cdot, \lambda)) & = \int_{-\infty}^0 \int_{s'}^0  e^{ g(\lambda)s + g(\lambda)s' } \mathbb{E} \bigl (e^{\gamma (W_s(\cdot) + W_{s'}(\cdot))}\bigr ) \mathrm{d}s \mathrm{d}s' \\
& \hspace{3em} +  \int_{-\infty}^0 \int_{-\infty}^{s'} e^{ g(\lambda)s + g(\lambda)s' } \mathbb{E} \bigl (e^{\gamma (W_s(\cdot) + W_{s'}(\cdot))}\bigr ) \mathrm{d}s \mathrm{d}s' \\
& = \int_{-\infty}^0 \int_{s'}^0  e^{ g(\lambda) s + g(\lambda)s'  + \frac{\gamma^2}{2} (s-s')  - 2 \gamma^2 s} \mathrm{d}s \mathrm{d}s' \\
& \hspace{3em} +  \int_{-\infty}^0 \int_{-\infty}^{s'} e^{ g(\lambda)s + g(\lambda)s' + \frac{\gamma^2}{2} (s'-s) - 2 \gamma^2 s'} \mathrm{d}s \mathrm{d}s' \\
& =: I_1 + I_2.
\eea 
Now, $\mathbb{E}(\M_{n}(\cdot, \lambda)\M_{n}(\cdot, \lambda))$ can be evaluated by direct computations. Actually, the condition $\sigma < \sigma_{\#}$ serves as the integrability condition for $I_1$ and $I_2$ terms, and we have under this condition that
\bea \label{eq:I1}
I_1  = I_2 =  \frac{1}{2(g(\lambda) - \gamma^2/2)(g(\lambda)-\gamma^2)}.
\eea
In the following, we provide details about the calculation of $I_1$ term. The calculation for $I_2$ term is simpler and is omitted here.

Since 
\beas
I_1  = \int_{-\infty}^0 \int_{s'}^0  e^{ (g(\lambda) - 3\gamma^2/2) s + ( g(\lambda) - \gamma^2/2) s'} \mathrm{d}s \mathrm{d}s',
\eeas
there are two cases to consider. If $g(\lambda) - 3\gamma^2/2 \neq 0$, then 
\bea
I_1  & = \frac{1}{g(\lambda) - 3\gamma^2/2} \int_{-\infty}^0  (1 - e^{(g(\lambda) - 3\gamma^2/2)s'})  e^{(g(\lambda) - \gamma^2/2)s'} \mathrm{d}s' \\
& = \frac{1}{g(\lambda) - 3\gamma^2/2} \int_{-\infty}^0   \bigl( e^{(g(\lambda) - \gamma^2/2)s'} -  e^{2 (g(\lambda) - \gamma^2)s'} \bigr) \mathrm{d}s'.
\eea
It is clear that the last integral above is finite if and only if $g(\lambda) > \gamma^2$, which is the same as $\sigma < \sigma_{\#}$ according to our definition of $\sigma_{\#}$ and $\gamma$ given respectively in \eqref{eq:g and sigma} and \eqref{eq:gamma}. Under this condition, we obtain
\beas
I_1 = \frac{1}{g(\lambda) - 3\gamma^2/2} \biggl ( \frac{1}{g(\lambda) - \gamma^2/2} - \frac{1}{2(g(\lambda)-\gamma^2)} \biggr) = \frac{1}{2(g(\lambda) - \gamma^2/2)(g(\lambda)-\gamma^2)},
\eeas
which leads to \eqref{eq:I1} in this case.

If $g(\lambda) - 3\gamma^2/2 = 0$, we get still under the condition $g(\lambda) > \gamma^2$ that
\beas
I_1  = \int_{-\infty}^0 \int_{s'}^0  e^{ (g(\lambda) - 3\gamma^2/2) s + ( g(\lambda) - \gamma^2/2) s'} \mathrm{d}s \mathrm{d}s'  =  - \int_{-\infty}^0  s' e^{(g(\lambda) - \gamma^2/2)s'} \mathrm{d}s'  = \frac{1}{(g(\lambda) - \gamma^2/2)^2}  = \frac{1}{\gamma^4},
\eeas
which agrees with \eqref{eq:I1} since $\frac{1}{2(g(\lambda) - \gamma^2/2)(g(\lambda)-\gamma^2)} = \frac{1}{\gamma^4}$ when $g(\lambda) - 3\gamma^2/2 = 0$.

It follows from \eqref{eq:2nd moments} and \eqref{eq:I1} that 
\beas
\mathbb{E}(\M_{n}(\cdot, \lambda)\M_{n}(\cdot, \lambda)) = \frac{1}{(g(\lambda) - \gamma^2/2)(g(\lambda)-\gamma^2)}, \qquad \sigma < \sigma_{\#}. 
\eeas 
This together with \eqref{expectation Mn} leads to
\bea  \label{eq:variance Mn-2}
\mathrm{Var}(\M_{n}(\cdot, \lambda)) &  :=  \mathbb{E}(\M_{n}(\cdot, \lambda)  \M_{n}(\cdot, \lambda)) -  \mathbb{E}(\M_{n}(\cdot, \lambda))^2 \\
& = \frac{1}{(g(\lambda) - \gamma^2/2)(g(\lambda)-\gamma^2)} -  \frac{1}{(g(\lambda) - \gamma^2/2)^2} \\
& = \frac{\gamma^2}{2(g(\lambda) - \gamma^2/2)^2(g(\lambda)-\gamma^2)}, \qquad \sigma < \sigma_{\#},
\eea
where we used $\gamma = (k-1)\sigma$ to rewrite $\mathbb{E}(\M_{n}(\cdot, \lambda))$; and \eqref{variance Mn} follows.

To derive (iii), we only need to compute the covariance $\mathrm{Cov}(\M_{n}(\theta_{s+t} \cdot, \lambda), \M_{n}(\theta_{s} \cdot, \lambda))$ for any given $s$ and $t$. The calculations are essentially the same as those used to compute the variance of $\M_{n}(\cdot, \lambda)$, and the main idea is still to use the independent increment property of the Wiener process. We sketch the main steps below and omit the detailed calculations.

First note that by introducing $\omega' = \theta_{s}\omega$, we only need to compute
$\mathrm{Cov}(\M_{n}(\theta_{t} \cdot, \lambda), \M_{n}(\cdot, \lambda))$, and without loss of generality, we consider the case $t< 0$. As before, by using Fubini's Theorem, we obtain
\beas
\mathbb{E}(\M_{n}(\theta_{t} \cdot, \lambda)  \M_{n}(\cdot, \lambda)) & = \mathbb{E} \biggl( \int_{-\infty}^0 \int_{-\infty}^0  e^{ g(\lambda)s + g(\lambda)s' + \gamma W_s(\cdot) + \gamma W_{s'}(\theta_{t} \cdot)}\, \mathrm{d}s \mathrm{d}s' \biggr) \\
& = \int_{-\infty}^0 \int_{-\infty}^0  e^{ g(\lambda)s + g(\lambda)s' } \mathbb{E} \bigl (e^{\gamma (W_s(\cdot) + W_{s'}(\theta_{t} \cdot))}\bigr ) \mathrm{d}s \mathrm{d}s'.
\eeas
Note that 
\beas
W_s(\omega) + W_{s'}(\theta_{t} \omega) = W_{s}(\omega) + W_{s'+t}(\omega) - W_{t}(\omega). 
\eeas
We aim to write this expression as the sum of independent random variables in order to evaluate $\mathbb{E} \bigl (e^{\gamma (W_s(\cdot) + W_{s'}(\theta_{t} \cdot))}\bigr )$. Since $s'+t < t < 0$,  there are three cases to be considered depending on the relative position of $s$ with respect to $s'+t$ and $t$. 

If $s< s' + t$, we choose the independent random variables to be $W_{s}(\omega) - W_{s'+t}(\omega)$, $2(W_{s'+t}(\omega) -  W_{t}(\omega))$, and $W_{t}(\omega)$. In this case, we have
\beas
\mathbb{E} \bigl (e^{\gamma (W_s(\cdot) + W_{s'}(\theta_{t} \cdot))}\bigr) &= \mathbb{E} \bigl ( e^{\gamma( W_{s}(\cdot) - W_{s'+t}(\cdot))} \bigr) \mathbb{E} \bigl (e^{ 2\gamma(W_{s'+t}(\cdot) -  W_{t}(\cdot))} \bigr)\mathbb{E} \bigl ( e^{\gamma W_{t}(\cdot)} \bigr) \\
& = e^{\frac{\gamma^2}{2}(s'+t - s)}  e^{-2\gamma^2s'}  e^{-\frac{\gamma^2}{2}t}.
\eeas


If $s' + t < s < t$, then $W_{s'+t}(\omega) - W_{s}(\omega)$, $2(W_{s}(\omega) -  W_{t}(\omega))$, and $W_{t}(\omega)$ are  the desired independent random variables; and we have
\beas
\mathbb{E} \bigl (e^{\gamma (W_s(\cdot) + W_{s'}(\theta_{t} \cdot))}\bigr) &= \mathbb{E} \bigl ( e^{\gamma( W_{s'+t}(\cdot) - W_{s}(\cdot))} \bigr) \mathbb{E} \bigl (e^{ 2\gamma(W_{s}(\cdot) -  W_{t}(\cdot))} \bigr)\mathbb{E} \bigl ( e^{\gamma W_{t}(\cdot)} \bigr) \\
& = e^{\frac{\gamma^2}{2}(s - s' - t)}  e^{2\gamma^2 (t-s)}  e^{-\frac{\gamma^2}{2}t}.
\eeas

Finally, if $t < s \le 0$, the independent random variables are chosen to be $W_{s'+t}(\omega) - W_{t}(\omega)$ and $W_{s}(\omega)$; and 
\beas
\mathbb{E} \bigl (e^{\gamma (W_s(\cdot) + W_{s'}(\theta_{t} \cdot))}\bigr) = \mathbb{E} \bigl ( e^{\gamma( W_{s'+t}(\cdot) - W_{t}(\cdot))} \bigr) \mathbb{E} \bigl ( e^{\gamma W_{s}(\cdot)} \bigr)  = e^{-\frac{\gamma^2}{2}s'}  e^{-\frac{\gamma^2}{2}s}.
\eeas

Now, similar to \eqref{eq:2nd moments}, we can rewrite $\mathbb{E}(\M_{n}(\theta_{t} \cdot, \lambda)  \M_{n}(\cdot, \lambda))$ as the sum of three integrals:
\beas
\mathbb{E}(\M_{n}(\theta_{t} \cdot, \lambda)  \M_{n}(\cdot, \lambda))  = J_1 + J_2 + J_3,
\eeas
where
\beas
 J_1 & := \int_{-\infty}^0 \int_{-\infty}^{s'+t}  e^{ g(\lambda)s + g(\lambda)s' }e^{\frac{\gamma^2}{2}(s'+t - s)}  e^{-2\gamma^2s'}  e^{-\frac{\gamma^2}{2}t}  \mathrm{d}s \mathrm{d}s', \\
 J_2 & := \int_{-\infty}^0 \int_{s'+t}^{t}  e^{ g(\lambda)s + g(\lambda)s' } e^{\frac{\gamma^2}{2}(s - s' - t)}  e^{2\gamma^2 (t-s)}  e^{-\frac{\gamma^2}{2}t}  \mathrm{d}s \mathrm{d}s', \\
 J_3 & := \int_{-\infty}^0 \int_{t}^{0}  e^{ g(\lambda)s + g(\lambda)s' }e^{-\frac{\gamma^2}{2}s'}  e^{-\frac{\gamma^2}{2}s}  \mathrm{d}s \mathrm{d}s'.
\eeas
Following the same type of calculations as those for $I_1$ and $I_2$ terms, we obtain that all these three terms $J_1$, $J_2$, and $J_3$ are finite if and only if $\sigma < \sigma_{\#}$; and under this condition, we have
\beas
\mathbb{E}(\M_{n}(\theta_{t} \cdot, \lambda)  \M_{n}(\cdot, \lambda)) & = J_1 + J_2 + J_3 \\
& = \frac{\gamma^2 e^{(g(\lambda) - \gamma^2/2)t}}{2(g(\lambda) - \gamma^2/2)^2(g(\lambda)-\gamma^2)} + \frac{1}{(g(\lambda) - \gamma^2/2)^2} \\
& = \mathrm{Var}(\M_{n}(\cdot, \lambda)) e^{(g(\lambda) - \gamma^2/2)t} + \frac{1}{(g(\lambda) - \gamma^2/2)^2}, \qquad  \Forall t < 0,
\eeas  
where the last equality follows from the expression of $\mathrm{Var}(\M_{n}(\cdot, \lambda))$ derived in \eqref{eq:variance Mn-2}.

Now, by using $\mathbb{E}(\M_{n}(\theta_{t} \cdot, \lambda)) = \mathbb{E}(\M_{n}(\cdot, \lambda)) = \frac{1}{(g(\lambda) - \gamma^2/2)^2}$, we obtain that 
\beas
\mathrm{Cov}(\M_{n}(\theta_{t} \cdot, \lambda), \M_{n}(\cdot, \lambda)) &   =  \mathbb{E}(\M_{n}(\theta_{t} \cdot, \lambda)  \M_{n}(\cdot, \lambda))  -  \mathbb{E}(\M_{n}(\theta_{t} \cdot, \lambda)) \mathbb{E}(\M_{n}(\cdot, \lambda)) \\
& = \mathrm{Var}(\M_{n}(\cdot, \lambda)) e^{(g(\lambda) - \gamma^2/2)t},  \qquad  \Forall t < 0,
\eeas
and the expression for the autocorrelation $R(t)$ given by  \eqref{eq:autocorrelation} follows for the case $t<0$ by recalling that $\gamma = (k-1)\sigma$.

\ep

\bibliographystyle{amsalpha}
\newcommand{\etalchar}[1]{$^{#1}$}
\providecommand{\bysame}{\leavevmode\hbox to3em{\hrulefill}\thinspace}
\providecommand{\MR}{\relax\ifhmode\unskip\space\fi MR }
\providecommand{\MRhref}[2]{%
  \href{http://www.ams.org/mathscinet-getitem?mr=#1}{#2}
}
\providecommand{\href}[2]{#2}

\end{document}